\documentclass[12pt, leqno]{article}

\usepackage{wrapfig}
\usepackage{euscript}

\usepackage{tikz, tikz-3dplot, pgfplots}
\usetikzlibrary{decorations.markings, patterns, decorations.pathmorphing}

\tikzset{->-/.style={decoration={
markings,
mark=at position #1 with {\arrow{>}}},postaction={decorate}}}

\usepackage{fullpage}

\usepackage{graphicx}
\usepackage[T1]{fontenc}
\usepackage{amsmath}
\usepackage{amssymb}
\usepackage{amsfonts}
\usepackage{amsthm}

\usepackage{mathabx}
\usepackage{euscript}

\usepackage[all,2cell]{xy}
\UseAllTwocells
\usepackage{mathrsfs}
\usepackage{enumitem}
\usepackage{subfigure}
\usepackage{float}
\usepackage{caption}

\usepackage[latin1]{inputenc}
\usepackage{amsxtra}
\usepackage{amsmath}
\usepackage{amssymb}
\usepackage{amsfonts}
\usepackage[all,2cell]{xy}
\UseAllTwocells
\usepackage{mathrsfs}
\usepackage{amsthm}
\usepackage{enumitem}
\usepackage{subfigure}
\usepackage{epigraph}

\newtheorem{thm}[equation]{Theorem}

\newtheorem{cor}[equation]{Corollary}
\newtheorem{lem}[equation]{Lemma}
\newtheorem{prop}[equation]{Proposition}
\newtheorem{para}[equation]{}
\newtheorem{conj}[equation]{Conjecture}

\newtheoremstyle{example}{\topsep}{\topsep}%
{}
{}
{\bfseries}
{.}
{2pt}
{\thmname{#1}\thmnumber{ #2}\thmnote{ #3}}

\theoremstyle{example}
\newtheorem{ass}[equation]{Assumption}

\newtheorem{conven}[equation]{Convention}

\newtheorem{Defi}[equation]{Definition}
\newtheorem{defi}[equation]{Definition}

\newtheorem{ex}[equation]{Example}

\newtheorem{exas}[equation]{Examples}

\newtheorem{propdef}[equation]{Proposition-Definition}

\newtheorem{rem}[equation]{Remark}
\newtheorem{rems}[equation]{Remarks}

\numberwithin{equation}{subsection}

\renewcommand{\theparagraph}{ \Alph{paragraph}.} 
\def\Ac{\mathcal{A}}
\def\Bc{\mathcal{B}}
\def\Cc{\mathcal{C}}
\def\Kc{\mathcal{K}}
\def\Dc{\mathcal{D}}
\def\Ec{\mathcal{E}}
\def\Fc{\mathcal{F}}
\def\Gc{\mathcal{G}}

\def\Lc{\mathcal{L}}
\def\Mc{\mathcal{M}}
\def\Nc{\mathcal{N}}
\def\Hc{\mathcal{H}}
\def\Oc{\mathcal{O}}
\def\Pc{\mathcal{P}}
\def\Qc{\mathcal{Q}}

\def\Tc{\mathcal{T}}

\def\Uc{\mathcal{U}}
\def\Vc{\mathcal{V}}
\def\Wc{\mathcal{W}}

\def\gen{\mathfrak{g}}
\def\hen{\mathfrak{h}}

\def\ken {\mathfrak{k}}

\def\wen{\mathfrak{w}}

\def\Cen{\mathfrak{C}}
\def\Fen{\mathfrak{F}}

\def\Hen{\mathfrak{H}}
\def\Len{\mathfrak{L}}
\def\Pen{\mathfrak {P}}
\def\Sen{\mathfrak{S}}

\def\VV{\mathbb{V}}

\def\CC{\mathbb{C}}
\def\DD{\mathbb{D}}
\def\FF{\mathbb{F}}

\def\HH{\mathbb{H}}
\def\II{\mathbb{I}}

\def\LL{\mathbb{L}}

\def\OO{\mathbb{O}}
\def\PP{\mathbb{P}}

\def\RR{\mathbb{R}}

\def\TT{\mathbb{T}}
\def\ZZ{\mathbb{Z}}

\def\E {{\tt{E}}}
\def\M {{\tt{M}}}
\def\R {{\tt{R}}}

\def\Vs{{\EuScript{V}}}

\def\Aff{{\on{Aff}}}
\def\Alg{{\on{Alg}}}
\def\an{{\on{an}}}

\def\ba{{\mathbf{a}}}
\def\Bar{{\on{Bar}}}
\def\BBB{{\on{BB}}}
\def\be{\begin{equation}}
\def\bef{\begin{figure}}
\def\beg{\begin{gathered}}
\def\BF{{\on{BF}}}
\def\bg{\begin{gathered}}
\def\bgen{{\overset\bullet\gen}}
\def\bh{{\mathbf{h}}}
 \def\bl{{\mathbf{l}}}
\def\bPhi{{\mathbf{\Phi}}}
\def\bPsi{{\mathbf{\Psi}}}

\def\Br{{\on{Br}}}
\def\bs{{\bigstar}}
\def\btp {\begin{tikzpicture}}
\def\bv{{\mathbf{v}}}

\def\centerarc[#1](#2)(#3:#4:#5)
{ \draw[#1] ($(#2)+({#5*cos(#3)},{#5*sin(#3)})$) arc (#3:#4:#5); }
\def\Coalg{{\on{Coalg}}}
\def\codim{{\on{codim}}}
\def\Coh{{\on{Coh}}}
\def\Coker{{\on{Coker}}}

\def\Cone{{\on{Cone}}}

\def\cont{{\on{cont}}}
\def\Conv{{\on{Conv}}}
\def\CY{{\on{CY}}}

\def\Def{{\on{Def}}}
\def\del{{\partial}}
\def\dgCat{{\tt{dgCat}}}
\def\dgFun{{\on{dgFun}}}
\def\dgVect{{\on{dgVect}_\k}}
\def\Diffeo{{\on{Diffeo}}}
\def\dirr{{\on{dir}}}

\def\drarr{\draw  [decoration={markings,mark=at position 0.7 with
{\arrow[scale=1.5,>=stealth]{>}}},postaction={decorate},
line width=.2mm]}

\def\ed{{\end{document}}}
\def\ee{\end{equation}}
\def\eg{\end{gathered}}
\def\egen{{\on{Xgen}}}
\def\enf{\end{figure}}
\def\End{{\on{End}}}
\def\eng{\end{gathered}}
\def\eps{{{\varepsilon}}}
\def\etp{\end{tikzpicture}}
\def\exc{{\on{exc}}}

\def\FC{{\on{FC}}}
\def\FM{{\on{FM}}}
\def\Fuk{{\on{Fuk}}}
\def\Fun{{\on{Fun}}}

\def\genoo{{\on{gen},\infty}}
\def\genr{{\on{gen}}}

\def\gr{{\on{gr}}}

\def\hra{{\hookrightarrow}}
\def\ho{{\on{ho}}}

\def\hol{{\on{hol}}}

\def\Hom{{\on{Hom}}}

\def\hr{{{\on h}. \on{r}.}}

\def\Id{{\on{Id}}}
\def\Im{{\on{Im}}}

\def\k{{\mathbf{k}}}
\def\kappa{{\varkappa}}
\def\Ker{{\on{Ker}}}

\def\LFS{{\Phi{\on S}}}
\def\Lie{{L}}
\def\Lin{{\on{Lin}}}
\def\lla{\longleftarrow}

\def\lra{\longrightarrow}
\def\LS{{\on{LS}}}

\def\max{{\on {max}}}
\def\MC{{\on{MC}}}
\def\Mod{{\text{-}\on{Mod}}}
\def\mod{{\text{-}\on{mod}}}
\def\Mor{{\on{Mor}}}

\def\n{n$^\circ$}
\def\naive{{\on{naive}}}
\def\Ndg{{\on{N}_{\on{dg}} }}

\def\Ob{{\on{Ob}}}
\def\ol{\overline}
\def\on{\operatorname}
\def\oo{{\infty}}
\def\op{{\on{op}}}
\def\OR{{\on{or}}}

\def\para{{\parallel}}
\def\PBr{{\on{PBr}}}
\def\Perf{{\on{Perf}}}
\def\perv{{\on{perv}}}
\def\Perv{{\on{Perv}}}
\def\phi{{{\varphi}}}

\def\Pretr{{\on{Pre}\text{-}\on{Tr}}}
\def\pt{{\on{pt}}}

\def\Re{{\on{Re}}}

\def\reg{{\on{reg}}}
\def\res{{\on{res}}}
\def\RHom{{\on{RHom}}}

\def\Schob{{\on{Schob}}}
\def\Sec{{\on{Sec}}}

\def\Sgn{{\on{Sgn}}}
\def\sgn{{\on{sgn}}}
\def\Sol{{\on{Sol}}}

\def\Th{{\on{Th}}}

\def\Tot{{\on{Tot}}}
\def\Tr{{\on{Tr}}}

\def\ul{\underline}

\def\vC{{{\overrightarrow C}}{}}
\def\Vect{{\on{Vect}}}
\def\Vert{{\on{Vert}}}

\def\wc{\widecheck}
\def\wh{\widehat}
\def\wt{\widetilde}

\def\1{{\mathbf{1}}}
\def\<{{\langle }}
\def\>{{\rangle }}
\def\-{{\setminus}}
\def\={{\, \simeq \, }}
\def\!{{ \rotatebox[origin=c]{180}{!} }}

\begin{document}

\title{ Perverse schobers and the Algebra of the Infrared}

\author {Mikhail Kapranov, Yan Soibelman, Lev Soukhanov}

\maketitle

\begin{abstract}
  We relate the Algebra of the Infrared of Gaiotto-Moore-Witten with the theory of perverse schobers which
  are (conjectural, in general) categorical analogs of perverse sheaves. A perverse schober on a complex plane
  $\CC$ can be seen as an algebraic structure that can encode various categories of D-branes of
  a 2-dimensional supersymmetric field theory, as well as the interaction (tunnelling) between such categories.
  We show that many constructions of the Algebra of the Infrared can be developed once we have a schober on
  $\CC$. These constructions can be seen as giving various features of the analog, for
  schobers, of the geometric Fourier transform well known for $\Dc$-modules and perverse sheaves.

 \end{abstract}

\tableofcontents

\addtocounter{section}{-1}

\section{Introduction}
\numberwithin{equation}{section}

\paragraph{The goal of the paper.}

The term ``Algebra of the Infrared'' was coined by D. Gaiotto, G.  Moore and
E. Witten  \cite{GMW1} to describe
a remarkable new chapter of homological algebra dealing with triangulated categories of D-branes
in quantum field theory. In the infrared limit, a physical theory reduces to the data of vacua and of  tunnelling between the
vacua. For a massive $2$-dimensional theory with $(2,2)$ supersymmetry,
they introduced  a remarkable algebraic structure:  a
homotopy Lie ($\Lie_\oo$-) algebra $\gen$   acting by higher derivations on
a homotopy associative ($A_\oo$-)algebra $R$ and a Maurer-Cartan element in $\gen$ giving
a deformation of $R$ describing  the global category of D-branes (Fukaya-Seidel category).
Remarkably, these structures make appeal to rather deep aspects of convex geometry.
In \cite  {KaKoSo} a part of the constructions of \cite{GMW1} was interpreted using the language of secondary
polytopes.

\vskip .2cm

The goal of this paper is to relate the Algebra of the Infrared  with
the  theory of {\em perverse schobers}, which are (conjectural, in general)
categorical analogs of perverse sheaves introduced in  \cite{KS-schobers}   and studied further in
\cite{BKS, donovan-PS-crossing,  {donovan-perv-RS}, {donovan-kuwagaki}}.
More precisely, we show that perverse schobers form the natural framework for this type of algebra:
many of its constructions
can be developed each time we have a schober on $\CC$.

\paragraph{Tunnelling between vacua and perverse sheaves.}
Let us summarize (with some formal generalizations) the approach of \cite{GMW1}. A theory of the above kind
has a discrete set $A$ of vacua; for simplicity we assume $A$ finite. Further,
 $A\= \{w_1,\cdots, w_N\}\subset\CC$  is embedded into the complex plane   parametrizing charges
 of the supersymmetry algebra \cite{witten-olive}. For a Landau-Ginzburg (LG) theory with superpotential
 $W: Y\to\CC$ the vacua correspond to critical points $y_i$ of $W$ (which we assume isolated but
 not necessarily Morse), and $w_i=W(y_i)$ are the corresponding critical values (assumed distinct).

 \vskip .2cm

 In the IR limit the vacua become separated (cease to interact)
and so one can associate to the $i$th vacuum
 its own  ``category of infrared D-branes'' $\Phi_i$. If the vacuum is massive (in the LG example
 this means that $y_i$ is Morse), then $\Phi_i=D^b\Vect_\k$ is ``the triangulated category freely generated
 by one object'' (the geometric vanishing sphere in the LG Morse example). More generally, if $y_i$ is
 an isolated singular point of $W$, then $\Phi_i$ is the
 {\em local Fukaya-Seidel category} of $W$ at $y_i$ (see \S \ref{subsec:lef-schober} below for a discussion). This is a categorical
 analog of the space of vanishing cycles of $W$ at $y_i$, in particular, the rank of its Grothendieck
 group is $\mu_i$, the multiplicity, or Milnor number of $y_i$.

 \vskip .2cm

  Now, switching the interaction, i.e., the  tunnelling between the vacua,
gives rise to {\em transport functors}
\[
M_{ij}: \Phi_i \lra \Phi_j.
\]
 In the LG example, when both critical points $y_i$ and $y_j$ are Morse, $M_{ij}: D^b\Vect_\k \to
  D^b\Vect_\k$ is the functor of tensor multiplication with the Floer complex spanned by $\zeta$-solitons,
  i.e.,  by gradient trajectories of $\Re(\zeta^{-1} W)$ joining $y_i$ and $y_j$, see \cite{GMW1}.
  Here $\zeta$ is the slope of the interval $[w_i, w_j]$.

   \vskip .2cm

   We observe that data $(\Phi_i, m_{ij})$ (formed by vector spaces and linear operators instead of categories
   and functors) appear in the mathematical problem of {\em classification of perverse sheaves} \cite {gelfand-MV}.
   Originally \cite{BBD}, perverse sheaves were introduced as a language for intersection cohomology theory,
   which is a cohomology theory for singular varieties that satisfies Poincar\'e duality
   \cite{goresky-macpherson}. They are ``self-dual'' local objects, situated half-way between sheaves and
   cosheaves, much like half-densities, so important in quantum mechanics
   (they form a natural $L_2$-space), are situated between functions and volume forms. One manifestation
   of this is the relation of intersection  cohomologies and $L_2$-cohomologies, see \cite{goresky} for a review.

    \vskip .2cm

    An additional remarkable property of perverse sheaves is their intimate connection with the concept
    of vanishing cycles: forming vanishing cycles in the most general situation preserves perversity
    \cite{BBD, {Ka-Scha}}. In particular, S. Gelfand, R. MacPherson and K. Vilonen
\cite {gelfand-MV} used vanishing cycles to identify the category of perverse sheaves on $\CC$
or on a disk (modulo constant sheaves, see below),  with the category of data $(\Phi_i, m_{ij})$ subject to
a mild non-degeneracy property.

This parallelism leads us to suggest that perverse schobers should be a natural conceptual tool to
analyze infrared behavior of theories as above.

\paragraph {The concept of a perverse schober.}
Intuitively, perverse schobers should be categorical analogs of perverse sheaves,
where ``vector spaces are replaced by triangulated categories'' \cite{KS-schobers}.  The original
 observation of  \cite{KS-schobers} was that a datum of
 a spherical functor $a$ and its right adjoint $a^*$
\[
\xymatrix{
\Phi \ar@<.4ex>[r]^a&\Psi \ar@<.4ex>[l]^{a^*},
}
\]
see \cite{AL}, can be viewed as a categorical analog of the diagram of \cite{GGM}
  describing a perverse sheaf on $\CC$ with one singularity at $0$.

  \vskip .2cm

  The  term  {\em Schober}  is one of the German analogs of the English word {\em sheaf}:
  it has the same root, also present in the German {\em schieben} (to push), {\em Verschiebung}
  ({\em shift}), {\em Scheibe} (moving disk), as well as
  in the English {\em shove} (as in  ``when push comes to shove''). For economy, we often speak simply
  about ``schobers'', omitting the adjective ``perverse''.

\vskip .2cm

However, a general theory of perverse schobers has not yet been
developed. While perverse sheaves can be defined in terms of the derived category formed by
complexes of sheaves of vector spaces \cite{BBD},  there is no analogous pre-existing concept
of ``complexes of sheaves of categories''.  Nor it is clear what such a concept should look like
(although some important  hints as to categorified homological algebra are provided by
\cite{dyckerhoff:DK}).
Already in the case of Riemann surfaces, a flexible topological treatment of perverse
schobers is quite non-trivial, and is being developed in \cite{DKSS}.  

\vskip .2cm

In the present paper our focus is not on the general theory but on things one can do with perverse
schobers in a physically motivated context.  This provides an additional motivation for future foundational
studies. Accordingly, we adopt the most naive approach to perverse schobers on $\CC$, obtained by categorifying the description of \cite {gelfand-MV} for perverse sheaves.
 That is,  we consider a schober $\Sen$ with singularities at
$A=\{w_1,\cdots, w_N\}$ as represented by a collection
 $(a_i: \Phi_i\to\Psi)_{i=1}^N$ of spherical functors with common target. This description proceeds
 in the presence of a ``spider''  $K$  which is a system of cuts joining some outside point $v$
 with the $w_i$. The consideration of such spiders  goes back to the Arnold school of singularity theory
 in the 1960-70's, see \cite{AGV}.

 \vskip .2cm

 Already with this minimal background, it is possible to capture a lot of beautiful phenomena. Some further
 features are not fully captured but it becomes clear that they would follow from a more systematic theory.
 We indicate them and leave them for future work.

 \paragraph{The role of the Fourier transform.}  To be precise,
  the physical picture corresponds to looking at perverse sheaves and schobers
{\em with the intent of making the Fourier transform}. That is,
there is an important difference between the   above two approaches. In  the physical  tunnelling picture,
it is  the rectilinear  intervals $[w_i, w_j]$ that play a distinguished role, so the transport happens along
such   intervals. But  the approach of
Gelfand-MacPherson-Vilonen involves the transport along curved paths passing
through some  common faraway point (``Vladivostok'', see Remark \ref {rem:vladi}). In fact, the concept of a perverse sheaf being  purely  topological,
the  mathematical reason to  consider rectilinear intervals  may not be clear. The Fourier transform point of view makes  this reason apparent because of the dominance order of
various exponentials $w\to e^{zw}$ in the complex domain.

The relevance of the Fourier transform for the physical picture is clear by looking at
  the Landau-Ginzburg  example:
a typical  oscillatory integral
\be\label{eq:QC-int}
I(\hbar) \,= \int_{\Gamma\subset Y}  e^{{i\over\hbar}W(y)} \Omega(y) dy,
\ee
fundamental  in the theory, can be seen as the result of, first, integrating over the fibers of $W$ and
then taking the single variable Fourier transform.

\vskip .2cm

To formulate our point of view more clearly, {\em a large part of the Algebra of the Infrared can be interpreted
as the study of the Fourier transform for perverse schobers}. This last concept needs,  of course, to
be defined, and we can proceed by analogy with the case of perverse sheaves. Here we have the concept
of the Geometric Fourier Transform
 \cite{br-mal-ver, malgrange, daia} which corresponds to the formal Fourier transform
 ($z\mapsto \del_w, \,\, \del_z\to -w$) applied to a holonomic regular $\Dc$-module $\Mc$ on $\CC$.
 The latter procedure produces a holonomic, typically irregular $\Dc$-module $\wc \Mc$ on the
 dual complex plane with the only possible singularity at $0$. So ignoring the Stokes phenomena, the
 most basic datum provided by $\wc\Mc$, is a local system $\Lc$ on $\CC\- \{0\}$.

 In the categorical case the analog of $\Lc$ is provided by the
Fukaya-Seidel category with its dependence
on the direction at infinity \cite{seidel, GMW1}. We show that this data can be defined for an arbitrary
schober $\Sen$  on $\CC$ and define the further features that one can expect from the Fourier transform $\wc\Sen$,
see \S \ref{subsec:schober-rect}. Many of the remarkable constructions of the Algebra of the Infrared
can be interpreted as tools for finding these features in terms of the data  defining $\Sen$.

\vskip .2cm

 In fact, these constructions  can be found already in their ``baby version''
at the level of perverse sheaves, more specifically in the  (seemingly new) formulas expressing various curvilinear data
in terms of rectilinear ones, using convex geometry.
  For example, we give a formula (Proposition \ref{prop:stokes-infra})  for the Stokes matrix
of the Fourier transform of a perverse sheaf which can be seen as the de-categorifcation
of the  infrared $A_\oo$-algebra $R$, involving alternating sums over convex paths.

\paragraph{Schobers vs. localized schobers.}  It is worth stressing that  a ``tunnelling datum''
$(\Phi_i, m_{ij}: \Phi_i\to\Phi_j)$ of vector spaces and linear operators
represents not exactly a perverse sheaf, but a
{\em localized perverse sheaf}, i.e., an object of the category $\ol\Perv(\CC,A)$ obtained by
quotienting $\Perv(\CC, A)$,  the category of all perverse sheaves with singularities in $A$,  by the subcategory of
constant sheaves \cite{gelfand-MV}. Any $\Fc\in\Perv(\CC,A)$ gives rise to a localized perverse sheaf
$\ol\Fc$, i.e., a tunnelling datum as above. To lift $\ol\Fc$ to an $\Fc$ in terms of tunnelling data
amounts to giving an
extra vector space $\Psi$ (the generic stalk) and maps
$\xymatrix{
\Phi_i \ar@<.4ex>[r]^{a_i}&\Psi \ar@<.4ex>[l]^{b_i}
}$
so that (the Vladivostok version of) $m_{ij}$ factors as $b_ja_i$.

\vskip .2cm

Accordingly, it is ``localized schobers'' that are, strictly speaking, relevant for our problem. In this paper
we do not attempt to define the corresponding ($\oo$-)category, as its meaningful study would require
more foundational work. Instead, we start with a schober represented by a diagram $(\Phi_i\buildrel a_i\over\to\Psi)$
of spherical functors and use some but not all information contained in such a diagram.
The physical meaning of the category $\Psi$ for a schober seems less clear but possibly important as one can
see from the following example.

\vskip .2cm

A LG superpotential $W: Y\to\CC$ gives rise to the {\em Lefschetz perverse sheaf} $\Lc_W$,
see \S \ref {subsec:PL-class} and its categorification, {\em Lefschetz schober} $\Len_W$, see \S \ref{subsec:lef-schober}.
The vector space $\Psi$ for $\Lc_W$ is the middle cohomology of a generic fiber $W^{-1}(b)$;
the category $\Psi$ for $\Len_W$ is the Fukaya category of $W^{-1}(b)$. As such, they are sensitive to
the way $W$ is (or is not) compactified, i.e., made into a proper map.

Note that traditionally, the superpotential mirror dual to a Fano variety $V$ is  written in a bare-bones,
non-proper form such as a Laurent polynomial $(\CC^*)\to\CC$, e.g.,
\[
W(z_1,\cdots, z_m) \,=\, z_1 +\cdots + z_m + {1\over z_1\cdots z_m}
\]
for $V=\PP^m$, so there is in principle some ambiguity as to its compactification. This ambiguity has a counterpart
in the discussion of the  quantum differential equation (QDE) of $V$ which we write in the  1-variable form as
\be\label{eq:QDE}
{d\Psi\over dq} \,=\, {\mu \Psi\over q} + K *_{q^K} \Psi, \quad \Psi: \CC^*\to H^\bullet(X,\CC).
\ee
Here $K\in H^2(X,\ZZ)$ is the canonical class of $V$ and $\mu$ is the grading operator equal to $i-\dim(V)/2$ on $H^i$.
A part of the mirror symmetry package for $V$ can be expressed by saying that the QDE ``is'' the Fourier transform
$\wc\Lc_W$ of $\Lc_W$ (in particular, its solutions are expressed via oscillatory integrals for $W$). But
$\wc\Lc_W$ is a $\Dc$-module on $\CC$, while \eqref{eq:QDE} gives, a priori, a $\Dc$-module only on $\CC^*$.
Extension of it to all of $\CC$  corresponds  precisely to specifying $\Lc_W$ as a perverse sheaf
(and not just a localized one), i.e., to
specifying the above vector space $\Psi$ precisely, i.e., to the question of compactification of $W$.

\vskip .2cm

From this point of view, the Dubrovin Conjecture (\cite{dubrovin} \S 4.2.2),
 involving categorical lifts of the Stokes matrices of
the QDE, should likely admit a formulation in terms of schobers.

\paragraph{Webs vs. secondary polytopes.}\label{par:web-vs-sec}
Our approach to the ``Algebra of the Infrared in the presense of a schober'' is, following  \cite  {KaKoSo},
based on the concept of secondary polytopes \cite{GKZ}. In that, it is dual to the original approach
of Gaiotto-Moore-Witten \cite{GMW1, GMW2} which uses certain planar graphs called webs.
For convenience of the reader, let us
recall (cf. also \cite{KaKoSo} \S 13) this duality, which is made very natural by the Fourier transform
point of view.

\vskip .2cm

We view elements of $A=\{w_1,\cdots, w_N\}$ as points of the complex plane $\CC=\CC_w$ on which a perverse
schober $\Sen$ is defined, and consider the convex polygon $Q=\Conv(A)$. The secondary
polytope $\Sigma(A)$ has faces corresponding to decompositions $\Pc$ of $Q$ into a union $Q=\bigcup Q'_i$
of (suitably marked, see \S \ref{subsec:back-sec}) convex subpolygons with vertices in $A$. These
decompositions must be {\em regular}, i.e., admit convex piecewise-affine-linear functions
$f: Q\to\RR$ that do break along every possible juncture of the $Q'_i$.
Given a schober $\Sen$ with singularities in $A$,
each subpolygon $Q'\subset Q$
gives rise to a natural ``Clebsch-Gordan space'' $I(Q')$, see \S \ref{subsec:L-oo-andR_oo}.
It consists of natural transformations between the compositions of the transport functors $M_{ij}$
along two complementary arcs in the boundary of $Q'$. The infrared $\Lie_\oo$-algebra
$\gen=\gen_\Sen$ is formed out of such $I(Q')$ using the factorization properties of secondary
polytopes.

\vskip .2cm

In comparison, the original approach of \cite{GMW1} proceeds in the dual complex plane $\CC_z$
(on which, in our language, the Fourier transform $\wc\Sen$ is supposed to live). Points of $\CC_z$
represent $\RR$-linear functionals $\CC_w\to\RR$. Given a convex piecewise-affine-linear $f:Q\to\RR$
as above, the linear parts of $f$ on different $Q'_i$ give points $p_i\in\CC_z$. Joining $p_i$ with
$p_j$ by an edge when $Q'_i$ and $Q'_j$ have an edge in common, gives a graph $\wen$ in $\CC_z$
which is precisely a web in the sense of \cite{GMW1, GMW2}.
Topologically, $\wen$ is a graph Poincar\'e dual to the subdivision $\Pc$, see Fig. \ref{fig:web-dec}.
So $\wen$ encodes both the decomposition $\Pc=\{Q'_i\}$ and the function $f$.

\bef[H]
 \btp[scale=0.4]

 \node (i) at (0,6){};
\node (j) at (6,6){};
 \node (k) at (9,3){};
 \node (l) at (6,0){};
 \node (m) at (0,0){};
 \node (n) at (-3,3){};

 \fill (i) circle (0.1);
  \fill (j) circle (0.1);  \fill (k) circle (0.1);  \fill (l) circle (0.1);  \fill (m) circle (0.1);  \fill (n) circle (0.1);

  \draw (0,6) -- (6,6) -- (9,3) -- (6,0) -- (0,0) -- (-3,3) -- (0,6);

  \draw (0,6) -- (0,0) -- (6,6) -- (6,0);

  \draw (-7, -3) -- (11, -3);
\draw (-7, 9) -- (11, 9);
\draw (-7,-3) -- (-7,9);
\draw (11, -3) -- (11,9);

\node at (-1,3){$Q'_1$};
\node at (2,4){$Q'_2$};
\node at (4,2){$Q'_3$};
\node at (7.2,3){$Q'_4$};

\node at (10,8.2){$\CC_w$};
\node at (-4, 0){\large$\Pc$};

 \etp
 \quad\quad\quad\quad
  \btp[scale=0.4]

  \node at (-0.7, 2.2){$p_1$};

  \node (a) at (-1,3){};
    \node at (-0.7, 2.2){$p_1$};

  \node (b) at (2,4){};
  \node at (2,3.2){$p_2$};

  \node (c) at (4,2){};
  \node at (4.2, 2.9){$p_3$};

  \node (d) at (7,3){};
  \node at (8,3){$p_4$};

  \fill (a) circle (0.15);  \fill (b) circle (0.15);
 \fill (c) circle (0.15);
 \fill (d) circle (0.15);

 \node (1) at (-6,6){};
 \node (2) at (2,9){};
 \node (3) at (8,7){};
 \node (4) at (9,-1){};
 \node (5) at (2,-2){};
 \node (6) at (-4,-1){};

 \draw[line width=0.4mm, ->] (-1,3) -- (1);
  \draw[line width=0.4mm, ->] (2,4) -- (2);
   \draw[line width=0.4mm, ->] (7,3) -- (3);
    \draw[line width=0.4mm, ->] (7,3) -- (4);
     \draw[line width=0.4mm, ->] (4,2) -- (5);
      \draw[line width=0.4mm, ->] (-1,3) -- (6);

  \draw[line width=0.4mm] (-1,3) -- (2,4);
   \draw[line width=0.4mm] (2,4) -- (4,2);
    \draw[line width=0.4mm] (4,2) -- (7,3);

\draw (-7, -3) -- (11, -3);
\draw (-7, 9) -- (11, 9);
\draw (-7,-3) -- (-7,9);
\draw (11, -3) -- (11,9);

\node at (10,8.2){$\CC_z$};
\node at (-5,0){\large$\wen$};
 \etp

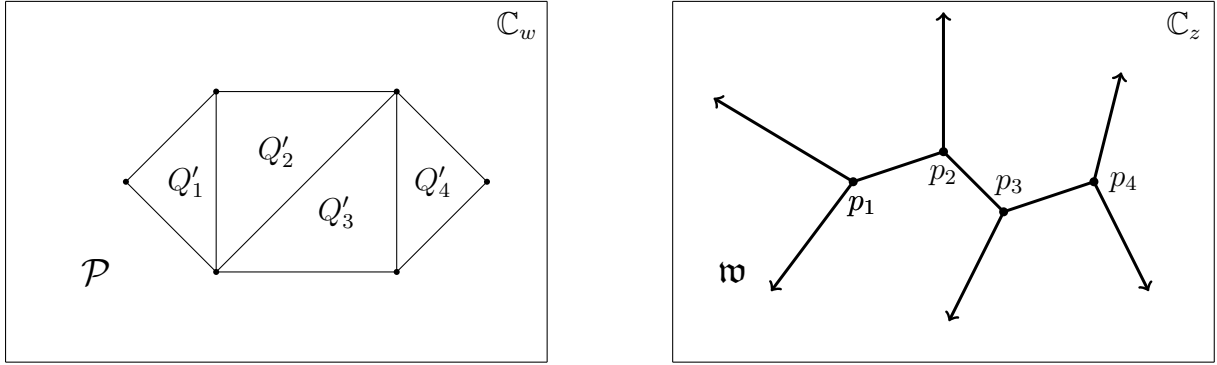
\caption{Webs $\wen$ are dual to subdivisions $\Pc$ together with convex PL-functions.}
\label{fig:web-dec}
\enf

Thus the two approaches are completely equivalent and one can easily translate any argument from one
language to the other.

 \paragraph{The role of higher coherences.}
 An intrinsic difficulty in any upgrade from vector space to triangulated categories is that of
 {\em higher coherences}.

 \vskip .2cm

 First, ``classical'' triangulated categories cannot be used as primary
 objects since their Hom-spaces are not fundamental: they
 appear as the cohomology of certain cochain complexes or
 the sets of connected components of some topological spaces. Keeping track of these
 complexes or spaces leads to {\em enhancements}: pre-triangulated dg-categories
 \cite{BK} (the approach adopted here), or more fundamental stable $\oo$-categories \cite{HA}.

 \vskip .2cm

 Second, once such an approach is adopted, the meaning of any identity among the
 morphisms changes: it should now hold up to some ``homotopy'', and there should be
 a coherence condition for these homotopies, which is itself a next level homotopy and so on
 {\em ad infinitum}. This of course complicates the study since every single step is now, so to say,
 replaced by an infinite routine.

 \vskip .2cm

 Some of these routines are by now standard and do not need special emphasis. Such are, for example, $A_\oo-$ and $\Lie_\oo$-structures
 (higher coherence data for associativity or the Jacobi identity); we use them in several
 places throughout the paper, mostly in the same context as \cite{GMW1, GMW2}.

 Beyond that, in this paper we keep such higher homotopy techniques to a minimum, in line with our
 ``elementary'' approach to schobers. That is, we work at the triangulated category level or,
 implicitly, both at the triangulated and the enhanced levels, when the transition between the two is easy.
 Hopefully, this makes the paper more accessible and helps to keep its size under control.

 \vskip .2cm

 At the very end, we highlight one situation where higher coherences are essential and which will be
 studied in a future paper. This is the Maurer-Cartan element of $\gen$ which gives,
 by deformation, the Fukaya-Seidel category. In  \cite{GMW1}, such an element was defined
 using the count of   solutions of a   $\ol\del$-equation. From our point of view, such an element
 should be present for any schober and appear as a collection of ``higher coherence data for
 Picard-Lefschetz identities'' in the following sense.

 \vskip .2cm

 For a perverse sheaf $\Fc$ on $\CC$ with singularities in $A$, every ``empty triangle''
 (i.e., a rectilinear triangle $\Delta_{ijk}=\Conv \{w_i, w_j, w_k\}$ not containing other $w_l$)
 gives rise to a Picard-Lefschetz (``wall-crossing'')   identity  among three linear transport operators:
 \be\label{eq:PL-intro}
 m_{ik}'\,=\, m_{ik} - m_{jk} m_{ij},
 \ee
see Proposition \ref{prop:PL-form} for details.

 \vskip .2cm

For a perverse schober $\Sen$ on $\CC$ with singularities in $A$, each such identity upgrades  to an exact triangle involving the transport functors
 \[
M'_{ik} \lra M_{ik} \buildrel u_{ ijk}\over \lra M_{jk} \circ M_{ij}
\buildrel v_{ ijk}\over\lra M_{ik} [1],
\]
cf.  Proposition \ref{prop:PL-triang}. It involves new data: the {\em Picard-Lefschetz arrows} (natural transformations
of functors) $u_{ijk}, v_{ijk}$.
Up to easy modifications (dualization etc.) there is just one Picard-Lefschetz arrow for each
geometric empty  triangle $\Delta_{ijk}$. Any subpolygon $Q'\subset Q$ with vertices in $A$,
triangulated into empty triangles, gives then a $2$-dimensional composition (pasting) of the Picard-Lefschetz
arrows corresponding to the triangles. An important fact is that {\em at the level of triangulated categories}
this pasting is independent of the triangulation\footnote{We have here  a clash of two uses of the term  ``triangulated'':
one for the algebraic concept of a triangulated category,  the other for a geometric triangulation of
a plane polygon.}, see Propositions \ref{prop:PL-ass} and \ref{prop:PL-3:1} for the effect of $2\to 2$
and $3\to 1$ moves. But at the dg-level we have a higher coherence problem for such independence.
The condition of the corresponding higher coherences are precisely those of being a Maurer-Cartan
element in $\gen_\Sen$.

\vskip .2cm

We expect  that such an element is canonically, up to gauge equivalence,
determined by the schober $\Sen$. This should follow  from a more ``advanced'' definition of
a schober which would include a sufficient supply of coherence data at the very outset. We leave this
for future work.

\paragraph{Overview of the paper.} The paper consists of five chapters.

 \vskip .2cm

{\bf Chapter \ref{sec:Perv(C,A)}} is devoted to an overview of perverse sheaves on topological surfaces,
especially on the complex plane $\CC$.
This subject is elementary but our presentation emphasizes several subtler points important for the sequel.

 \vskip .2cm

In \S \ref{subsec:abs-PL} we develop, following \cite{KS-schobers},
 an ``abstract'' version of the classical Picard-Lefschetz theory, in the context of a perverse sheaf instead of a
 singular point of a Lefschetz pencil, as is usually done,  see, e.g.,  \cite{AGV} Ch.1.
 The key concept for this is the transport map $m_{ij}(\gamma): \Phi_i\to\Phi_j$
 between spaces of vanishing cycles associated to a path $\gamma$.
Remarkably, the concept of
 a perverse sheaf appears to be intrinsically designed to produce Picard-Lefschetz identities such as
\eqref{eq:PL-intro}, see  Proposition \ref{prop:PL-form}.
 One can further interpret (Proposition \ref{cor:PL-vassiliev})
  these identities in terms of ``Vassiliev derivatives'' of the monodromy  living on various strata
 in the space of paths.
We also recall from  \cite {gelfand-MV} the explicit quiver-like description of perverse sheaves on $\CC$
which we later use as a blueprint for defining schobers.

\vskip .2cm

In \S \ref {subsec:PL-class} we connect the abstract Picard-Lefschetz theory for perverse sheaves with
the classical theory for Lefschetz pencils $W:Y\to X$ where $X$ is a Riemann surface. In fact, we consider
a slightly more general situation when $W$ has isolated but not necessarily Morse critical points.
The connection goes via a particular perverse sheaf $\Lc_W$ on $X$ associated with $W$ which
we call the {\em Lefschetz perverse sheaf}. It seems to be an object of some importance.
In particular, it has a natural categorification, the {\em Lefschetz schober},
which underlies many features of Fukaya-Seidel theory  see \S \ref{subsec:lef-schober}.

\vskip .2cm

In \S \ref{subsec:perv-bil} we extend the classical $(\Phi, \Psi)$-description of perverse sheaves on the disk
\cite{GGM} to a description of perverse sheaves equipped with a  nondegenerate,
not necessarily symmetric, bilinear form. Recall that a triangulated category can be seen as a ``higher'' analog not
just of a vector space but of a vector spaces with a bilinear form, whole role is played by the spaces $\Hom(x,y)$.
This form is, in general, non-symmetric, hence the difference between left and right adjoint functors. Our description is,
in agreement with the original proposal in \cite{KS-schobers}, completely analogous to the data defining a spherical
functor. In particular, we give, in this decategorified context, analogs of the entire spherical functor package
\cite{AL} and of the alternative definition of spherical functors due to Kuznetsov \cite{kuznetsov}.

\vskip .2cm

In \S \ref{subsec:quot-LC} we recall the description, due to \cite{gelfand-MV}, of the category
$\ol\Perv(\CC,A)$ of localized
perverse sheaves  (discussed above) in terms of tunnelling diagrams $(\Phi_i, m_{ij})$.

\vskip .2cm

{\bf Chapter \ref{sec:recti-FT}} is devoted to a study, in the form convenient for us, of the Fourier transform
of a perverse sheaf $\Fc$ on $\CC$ (that is, of the holonomic regular $\Dc$-module $\Mc$
corresponding to $\Fc$). This problem was studied before, first of all by Malgrange
\cite{malgrange} in much greater generality ($\Mc$ is allowed to be irregular). The work of Malgrange
was revisited by D'Agnolo-Kashiwara \cite{dagnolo-microlocal} from the point of view of  their
approach  via enhanced ind-sheaves \cite{dagnolo}.   It was recognized in
\cite {KKP, dagnolo-sabbah}  that  focusing specially on
the case of regular $\Mc$ leads to important insights and further developments.
 We continue on this path, bringing out
the connections with convex geometry and with the Algebra of the Infrared.

\vskip .2cm

In \S \ref{subsec:rec-approach}  we introduce, besides motivation and terminology, the rectilinear
version of  transport maps $m_{ij}$. We formulate Proposition \ref{prop:rect-cat} which says that these maps
give an alternative description of the category  $\ol\Perv(\CC,A)$.

\vskip .2cm

In \S \ref{subsec:iso-recti} we prove Proposition \ref{prop:rect-cat}. Our method is based on
{\em isomonodromic deformations} of perverse sheaves, when we move the allowed
singular points $w_i$ while keeping the local behavior near each $w_i$ unchanged.
This concept is usually considered for linear ODE (the Schlessinger system), not perverse
sheaves. In the perverse sheaf context, isomonodromic deformations have a transparent
``wall-crossing'' effect on the rectilinear transports. The resulting wall-crossing formula
is another manifestation of the Picard-Lefschetz identity. The corresponding
{\em wall of collinearity} (called {\em wall of marginal stability} in \cite{GMW1})
consists of configurations $(w_1, \cdots, w_N)\in\CC^N$ such that three points
$w_i, w_j, w_k$ lie on a real line. We also discuss {\em configuration chambers}
(connected components of the complement to the walls) and their relation with the
combinatorial concept of an oriented matroid
\cite{bjorner, bokowski-sturmfels}. This relation is surprisingly nontrivial.

\vskip .2 cm

In \S \ref{subsec:FTPS-top}, we describe, following \cite{malgrange, dagnolo-sabbah}, the Fourier transform
of (the holonomic regular $\Dc$-module associated to a) perverse sheaf $\Fc$ on $\CC$, at the
topological level, as a perverse sheaf $\wc\Fc$ on the dual plane $\CC_z$. We present a description of the
$(\Phi, \Psi)$-diagram  of $\wc\Fc$ (Proposition \ref{prop:wc-M-top-2}) as well as the representation of the section
of $\wc\Fc$ outside $0$ as actual Fourier integrals (Proposition \ref{prop:FT=integrals}).

\vskip .2cm

In \S \ref {subsec:FT-stokes}, we recall the concept of a Stokes structure (a local system on $S^1$
with a fitration indexed by the sheaf of exponential germs) due to Deligne and Malgrage
\cite{deligne-docs} and focus specifically on the case of {\em exponential Stokes structures}
where the exponential germs are given by a collection of simple exponentials $z\mapsto e^{wz}$,
$w\in A$. We emphasize, in Proposition \ref{prop:exp-stokes-OM} and Remark \ref{rem:manin-schechtman},
 the little noticed  fact that the combinatorial data encoded by an exponential Hodge structure
are related to the concept of oriented matroid.
We  also recall the description of the (exponential) Stokes structure on $\wc\Fc$ due to
\cite{malgrange, dagnolo-sabbah}.

\vskip .2cm

In \S \ref{subsec:baby-IA} we present a ``baby'', or decategorified,
 version of the Algebra of the Infrared, by
which we mean a class of formulas related to perverse sheaves (not schobers)
and involving summation over convex polygonal paths. Recall that the actual Algebra of
the Infrared of  \cite{GMW1, GMW2} involves complexes whose summands are labelled by such paths. We give
two examples of such formulas. One, in Proposition \ref{prop:baby-IA},  expresses the transport along a curvilinear path in terms of iterated transport along convex polygonal paths.
The other, in Proposition \ref{prop:stokes-infra}, gives a similar expression for the
Stokes matrix of  $\wc\Fc$. These results support our main thesis that the right context for the
Algebra of the Infrared is provided by perverse schobers.

\vskip .2cm

{\bf Chapter \ref{sec:schob-surf}} begins an upgrade, to perverse schobers, of the
study of perverse sheaves presented in Chapters  \ref{sec:Perv(C,A)} and
\ref{sec:recti-FT}.

\vskip .2cm

In \S \ref{subsec:schob-surf} we present, following \cite{KS-schobers}, our working
definition of schobers on surfaces in terms of spherical functors at singular points
complemented by a local system on the complement. We are especially interested in the
case when the surface is the disk or the complex plane. In this case a schober
can be defined in terms of a diagram of spherical functors, once we choose a ``spider''.
This is analogous to the Gelfand-MacPherson-Vilonen description \cite{gelfand-MV}
for perverse sheaves.

\vskip .2cm

In \S \ref{subsec:PL-tri} we upgrade the Picard-Lefschetz identities of Proposition
\ref{prop:PL-form} to {\em Picard-Lefschetz triangles} \cite{KS-schobers}.
These  are exact triangles  in the category of (dg-)functors connecting
three  transport functors
which generalize the transport
maps of perverse sheaves. The arrows (natural transformations) in Picard-Lefschetz
triangles, which we call {\em Picard-Lefschetz arrows}, will be fundamental objects
of study in the rest of the paper. We prove, in Propositions \ref{prop:PL-ass} and
\ref{prop:PL-3:1}, geometric identities on the  ``pasting'' of these arrows which have
the meaning of invariance under elementary moves of plane triangulations.
These identities can be seen as strating points of the analysis leading to the
Algebra of the Infrared.

\vskip  .2cm

In \S \ref{subsec:rem-SOD} we review the theory of semi-orthogonal decompositions
(SODs) in triangulated categories and their mutations \cite{BK-serre}. We also
recall the construction of $2$-term SODs via gluing (dg-)functors \cite{kuz-lun}.

\vskip .2cm

In \S \ref{subsec:unitri-mon} we develop a generalization of the technique of gluing functors to
construct $N$-term SODs. Instead of a single functor, we need a uni-triangular {\em monad},
i.e., an associative algebra in the category of functors. We study mutations of uni-triangular monads, the Koszul duality as well
as the straightening procedure for $A_\oo$-monads. In the case of exceptional collections
(which are SOD's into copies of the category $D^b\Vect_\k$),  monads reduce to usual
associative dg-algebras and our analysis reduces to that of  \cite{bondal-polishchuk}.

\vskip .2cm

In \S \ref{subsec:monads-glue} we define the Fukaya-Seidel category associated to
a schober $\Sen$ on $\CC$ in terms of a diagram of spherical functors associated to
$\Sen$ via a spider. Our approach is a modification of the classical Barr-Beck theory
\cite{barr-wells}
that studies a monad associated to a pair of adjoint functors. In  our case we have
 not one but several pairs of adjoint functors with common target and we
 associate to them a  uni-triangular monad which we call the {\em Fukaya-Seidel monad}.
 The Fukaya-Seidel category is defined as the category of algebras over this monad.

 \vskip .2cm

In \S \ref{subsec:lef-schober} we discuss a particular example of a schober on $\CC$,
which we call the {\em Lefschetz schober} $\Len_W$ associated to a generalized Lefschetz pencil. i.e., to a superpotential $W:Y\to\CC$
with isolated critical points. It can be seen as a categorification of the Lefschetz perverse
sheaf $\Lc_W$ of \S \ref{subsec:PL-class}.
As with any Fukaya-style object, the construction of
$\Len_W$ is bound to involve non-trivial technical questions of symplectic geometry. We give references
to the literature where such questions have been addressed at various level of generality
and detail, some of the work still ongoing. In the more ``classical'' case of a Lefschetz pencil
(all critical points non-degenerate) we can rely on the techniques developed by Seidel
\cite{seidel}.

\vskip .2cm

{\bf Chapter \ref{sec:AIR-I}} is the first of two chapters dedicated to the  analysis
of schobers on $\CC$ from the ``infrared''  perspective
In it, we collected the aspects involving  more elementary convex geometry, i.e.,  making use of convex hulls, convex paths etc. but not of the concept of secondary polytopes.

\vskip .2cm

In \S  \ref{subsec:schober-rect}, we explain the rectlinear approach to schobers and spell out the schober
analog of the geometric Fourier transform for perverse sheaves. This analog features the Fukaya-Seidel
category.

\vskip .2cm

In \S \ref{subsec:infra-com-schob}, we develop a categorical analog of the Baby Infrared Algebra
of  \S \ref{subsec:baby-IA}. It appears as a cochain complex in the derived category
whose terms are direct sums of iterated rectilinear transports. We include it into
and a larger diagram  of exact triangles called the
{\em infrared Postnikov system} which categorifies the individual identities used to obtain
the formulas of  \S \ref{subsec:infra-com-schob}.

\vskip .2cm

In \S \ref{subsec:infra-FS-mon}, we  realize the components of the Fukaya-Seidel monad $M$ as
``circumnavigating'' transports of \S \ref{subsec:infra-com-schob}, and define
another monad $R$, which we call the {\em infrared monad}, formed by iterated
rectilinar transports. This monad is the schober analog of the $A_\oo$-algebra of half-plane webs
from \cite{GMW1, GMW2}, see also \cite{KaKoSo}. We formulate a conjecture that
each component $M_{ij}$ of $M$ can be obtained by adding lower order terms to the differential
 in $R_{ij}$. Note that this is somewhat different from the main trust of \cite{GMW1, GMW2}, where
  $A_\oo$-deformations are considered.

 \vskip .2cm

 {\bf Chapter  \ref{sec:AIR-II}} is devoted to the schober analogs of some of the main constructions of
  \cite{GMW1, GMW2}: the $\Lie_\oo$-algebra $\gen$ formed by
  ``closed polygons''  and its action on the monad $R$. Following \cite{KaKoSo}, we use the
  language of secondary polytopes.

  \vskip .2cm

  In \S \ref{subsec:back-sec} we give a brief  overview of  secondary polytopes, the main
  reference being \cite{GKZ}.

  \vskip .2cm

  In \S \ref{subsec:W-II-kind} we discuss the variation of the secondary
  polytope $\Sigma(A)$ under the change of $A$. Here we have a subtle phenomenon:
 the combinatorics of  $\Sigma(A)$ can change without the naive convex geometry of $A$
 (what lies in the convex
  hull of what) undergoing any visible change.  This phenomenon has been noticed
  in \cite{GMW1} at the level of webs and called {\em exceptional wall-crossing}.
  We develop a quantified version of this phenomenon in terms of jumps of
  the cohomology of a certain complex. We prove that the locus of
  ``exceptional'' configurations $A$ (for which we do have a jump) has real codimension
  $\geq 1$ so that we can indeed speak about walls. We also prove that for $A$ outside
  of such walls, $\Sigma(A)$  has the factorization property: each face of it is
  itself a product of other secondary polytopes. In this, we correct the treatment of
  \cite{KaKoSo} where the condition of exceptionally general position was overlooked.

  \vskip .2cm

  In \S \ref{subsec:L-oo-andR_oo} we define the $\Lie_\oo$-algebra $\gen = \gen_A(\Sen)$
  associated to a schober $\Sen\in \Schob(\CC,A)$. It consists of certain Clebsch-Gordan type spaces
  (or spaces of intertwiners) $\II(Q')$ associated to closed convex polygons $Q'$ with vertices in $A$.
  An element of $\II(Q')$  can be seen as a natural transformation between two iterated
  transport functors corresponding to two polygonal paths forming the boundary of $Q'$.
  The $\Lie_\oo$-operations in $\gen$ are induced by pasting of natural transformations.

  \vskip .2cm

  In \S \ref{subsec:def-R} we define a $\Lie_\oo$-morphism from $\gen$ to the deformation
  complex of the infrared operad $R$. In particular, any Maurer-Cartan element $\beta\in \gen^1$
  gives rise to an $A_\oo$-deformation $R(\beta)$ of $R$.

  \vskip .2cm

  In  \S \ref{subsec:MC} we analyze Maurer-Cartan elements $\beta$ in $\gen$ more explicitly.
  We conjecture that
  for any schober $\Sen$ there is a canonical (up to gauge equivalence)  $\beta$
  such that $R(\beta)$ is quasi-isomorphic to the Fukaya-Seidel monad $M$. We indicate
  an approach to proving this conjecture by deforming $A$ from the ``maximally concave''
  position, deforming $\Sen$ accordingly (isomonodromically),
   thereby crossing some walls
   and tracing the effect of the  wall-crossings on $\Sigma(A)$ and $\gen = \gen_A(\Sen)$.

  \vskip .2cm

  Finally,  {\bf Appendix \ref{app:enh}} collects notations and conventions related to
  triangulated categories, their dg-categorical enhancements, as well as adjunctions,
  Serre functors and other related  features.

\vskip 2cm

\paragraph{Acknowledgements.} We  would like to thank M. Abouzaid and K. Hori for useful dicussions
and correspondence. The research of M.K. was supported by World Premier International Research Center Initiative (WPI Initiative),
 MEXT, Japan. The research of Y.S. was partially supported by the Munson-Simu award at Kansas State University and NSF grant. He also thanks to IHES where part of the project was carried out.

\section{Review of $\Perv(\CC, A)$} \label{sec:Perv(C,A)}
\numberwithin{equation}{subsection}

\subsection{Abstract Picard-Lefschetz theory}\label{subsec:abs-PL}

\paragraph {Setup and definitions.}\label{par:setup}

We fix a field $\k$  of characteristic $0$.
Let $X$  be a smooth ($C^\infty$) surface, possibly non-compact and  with boundary $\del X$. Let
$A=\{w_1,\cdots w_N\} \subset X-\del X$ be  a finite set
of interior points. We refer to elements of $A$ as {\em marked points}.
The set $A$ defines a stratification of $X$ with strata being $X\-A$ and $1$-element subsets
$\{w_i\},  i=1,\cdots, N$.
We denote by
$D^b(X,A)$ the derived category formed by bounded complexes of $\k$-vector spaces on $X$  which are constructible
(in particular, have stalks with finite-dimensional cohomology)
with respect to this stratification. 

\vskip .2cm
Let $\Perv(X,A)\subset D^b(X, A)$ be the  abelian category of perverse sheaves with respect to the middle perversity. Explicitly,
$\Fc\in\Perv(X,A)$ iff:

\begin{itemize}
\item[($P^+$)] $\ul H^i(\Fc)$ is supported on a set of real codimension $\geq 2i$. That is,
$\ul H^1$ is supported on $A$ and
$\ul H^{\geq 2}=0$.

\item[($P^-$)] For any stratum $S$ the sheaves $\ul\HH^i_S(\Fc)$ vanish for $i< 2\on{codim}_\RR(S)$.
\end{itemize}
This implies that $\Fc|_{X\-A}$ is a local system in degree $0$ and any local system on $X$ in degree $0$ is perverse.
Note that our normalization differs by shift from that of
\cite{BBD}, in which a local system in degree $(-1)$ is perverse.

\vskip .2cm

The category $\Perv(X, A)$ has a perfect duality
\be
\Fc \,\mapsto \, \Fc^* := \DD(\Fc)[2],
\ee
where $\DD$ is Verdier duality \cite{Ka-Scha}. The shift by $2$ comes from our normalization of perversity so that on local systems in degree $0$
we get the usual (unshifted) duality.

\vskip .2cm

Further,  $\Perv(X,A)$ is the heart of the perverse t-structure on $D^b(X,A)$, in particular, we have the functors of {\em perverse cohomology}
\be
\ul H^i_\perv: D^b(X,A)\lra \Perv(X,A).
\ee

Let $\Diffeo(X,A)$ be the topological group of all self-diffeomorphisms of the pair $(X,A)$
which preserve  $\del X$ pointwise and $A$ as a set. Let also $\Diffeo(X, \Id_A)\subset \Diffeo(X,A)$
be the subgroup
formed by self-diffeomorphisms which preserve $A$ pointwise.
If
$X$ is oriented, we denote by $\Diffeo^+(X,A)\subset \Diffeo(X,A)$, $\Diffeo^+(X, \Id_A\subset \Diffeo(X,\Id_A))$ the subgroups
of orientation preserving diffeomorphisms.    We denote by
\[
\begin{gathered}
\Gamma_{X,A}\,=\, \pi_0\, \,  {\Diffeo}(X,A), \quad  \quad  \Gamma^+_{X,A}\,=\, \pi_0\, \,  {\Diffeo}^+(X,A),
\\
\Gamma_{X,\Id_A}\,=\, \pi_0\, \,  {\Diffeo}(X,\Id_A), \quad  \quad
\Gamma^+_{X, \Id_A}\,=\, \pi_0\, \,  {\Diffeo}^+(X, \Id_A)
\end{gathered}
\]
the corresponding mapping class groups
of $(X,A)$.

\begin{rem}
Alternatively, we can choose one new marked point on each component of $\del X$, denote the set of
these new points by $B$ and consider the group
of diffeomorphisms (resp. orientation preserving diffeomorphisms) of $(X, A)$ which preserve not the full
$\del X$ but just the points from $B$. This will give the same groups of connected components.
\end{rem}

The topological nature of the perversity conditions implies:

\begin {prop}\label{prop:Gamma-acts}
The group $\Gamma_{X,A}$ acts on the category $\Perv(X,A)$. \qed
\end{prop}

\begin{ex}\label{ex:class-group-braid}
If $X= D=\{|z|\leq 1\}$ be the closed unit disk in $\CC$ and $A=\{w_1 < \cdots < w_n\}\subset (-1,1)$
be a set of $N$ points on the real line inside $D$. Then
$\Gamma^+_{D,A}$
is canonically identified with $\Br_N$, and $\Gamma^+(D, \Id_A)$ with $\PBr_N$,  the  braid
group, resp. the pure braid group on $N$ strands,
cf \cite{farb} \S 9.1.3.  Thus $\Br_N$ acts on $\Perv(D,A)$.
\end{ex}

\paragraph{The Classical $(\Phi, \Psi)$ description.}\label{par:phi-psi}
Let $X=D=\{|z| \leq 1\}$ be the unit disk in $\CC$ and $A=\{0\}$. We recall the classical \cite{beil-gluing, GGM} description of
the category $\Perv(D,0)$.

 \begin{prop}\label{prop:ggm}
 $\Perv(D,0)$ is equivalent to the category $\Qc$ of diagrams
\[
\xymatrix{
\Phi \ar@<.4ex>[r]^a&\Psi \ar@<.4ex>[l]^b
}
\]
of finite-dimensional  $\k$-vector spaces
are linear maps such that
$T_\Psi: = \Id_\Psi - ab$   is an isomorphism
(or, what is equivalent, such that  $T_\Phi=\Id_\Phi-ba$ is an isomorphism). \qed
\end{prop}

For future use, we recall from \cite{GGM} the explicit construction of the equivalence
\be\label{eq:eq-phipsi}
\Perv(D, 0) \lra\Qc, \quad \Fc\mapsto \bigl\{
\xymatrix{
\Phi(\Fc) \ar@<.4ex>[r]^{a_\Fc}&\Psi (\Fc) \ar@<.4ex>[l]^{b_\Fc}
}
\bigr\}.
\ee
Consider the radial ``cut'' $K=[0,1)\subset D$ (one can choose any other simple curve starting at $0$ and ending on the boundary of $D$),
see Fig. \ref{fig:disk-1}, where $D$ is depicted as an ellipse.

Consider the sheaves $\ul\HH^i_K(\Fc)$ of hypercohomology with support in $K$. The assumption that $\Fc$ is perverse, implies that
$\ul\HH^i_K(\Fc)=0$ for $i\neq 1$. The sheaf $\ul\HH^1_K(\Fc)$ on $K$ is locally constant (hence constant) on $K -  \{0\}$. Let $\eps\in (0,1)$ be a small nonzero real number.
We define:
\[
\bg
\Phi(\Fc)\, \,:= \, \,\ul \HH^1_K(\Fc)_0  \,\,\simeq\,\,  \HH^1_K(\Delta, \Fc) \quad \text{(the space of {\em vanishing cycles})},\\
\Psi(\Fc)\,\, := \,\,   \ul \HH^1_K(\Fc)_\eps  \quad \text{(the space of {\em nearby cycles})}.
\eg
\]

\bef[H]
\centering
\btp [scale=.8, baseline=(current  bounding  box.center)]
\node (0) at (0,0){};
\fill (0) circle (0.1);
\node (eps) at (3.5,0){};
\fill (eps) circle (0.1);

\draw (0,0) ellipse (5.5cm and 3cm);

\draw (0,0) circle (1.1cm);
\draw (3.5, 0) circle (1.2cm);

\node (b) at (5,0){};
\draw[->, line width = .3mm] (0,0)  --  (5.5,0);
\draw[->, line width = 0.5mm,  dotted] plot [smooth,tension=1.5] coordinates{(3.7, 0.5)     (1 ,0.5)   (-0.8,0)  (1, -0.5) (3.7, -0.5)};

\node at (-0.3, 0) {\small$0$};
\node at (5.7, 0) {\small$1$};
\node at (2.7,1.3) {$V$};
\node at (3.5,0.9) {\small$V_+$};
\node at (3.5,-0.9) {\small$V_-$};
\node at (4,2.5) {$D$};
\node at (5,0.3){$K$};
\node at (-1,1){\small$U$};
\node at (3.6, -0.25){\small$\varepsilon$};

\etp

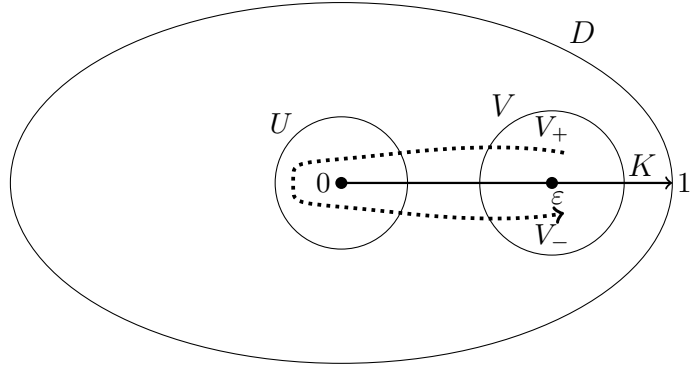
\caption{ The spaces $\Phi(\Fc)$ and $\Psi(\Fc)$.  }
\label{fig:disk-1}
\enf

Since $\Fc|_{D\setminus\{0\}}$ is a local system in degree $0$, after taking a small neighborhood $V$ of $\eps$ and denoting by
$V_\pm$
the two components of $V -  K$, we have
\[
\Psi(\Fc) \,=\, \Fc(V -  K)/\Fc(V)\,=\, {{\Fc(V_+)\oplus\Fc(V_-)}\over {\Fc(V)}} \,=\,{  { \Fc_\eps \oplus \Fc_\eps} \over
\Fc_\eps } \,\simeq \, \Fc_\eps
\]
is just the stalk of $\Fc$ at $\epsilon$, which explains the name ``nearby cycles''.

As for $\Phi(\Fc)$, it includes into the exact sequence (the standard sequence relating sheaves of hypercohomology with and without support)
\be\label{eq:LES-Phi}
0\to \HH^0(U, \Fc) \to \Fc(U-K) \buildrel \beta\over\to \Phi(\Fc) \to\HH^1(U,\Fc) \to 0,
\ee
where $U$ is any small contractible neighborhood of $0$ in $D$ (in fact, $U=D$ is also allowed). We now explain the construction of
the maps $a_\Fc$ and $b_\Fc$.

\vskip .2cm

\noindent \ul{The map $a_\Fc$:} It is just the {\em generalization map} of the constructible sheaf $\ul\HH^1_K(\Fc)$ from the stalk
at $0$ to the stalk at a nearby point $\eps$. More explicitly, restrictions from $D$ to $U$ and $V$ give the diagram
\[
\Phi(\Fc) = \HH^1_{U\cap K}(U,\Fc) \buildrel r_U\over \lla \HH^1_K(D,\Fc) \buildrel r_V\over\lra H^1_{V\cap K}(V,\Fc) = \Psi(\Fc),
\]
with $r_U$ being an isomorphism, and we define $a_\Fc= r_V\circ r_U^{-1}$.

\vskip .2cm

\noindent\ul{The map $b_\Fc$:}  Let $\psi\in \Psi(\Fc)=\Fc_\eps$. Continue $\psi$ to a section of $\Fc$ in $V_+$. After this, extend it
to a section $\wt{\psi}$ of $\Fc$ in $D-K$, so $\wt{\psi}|_{V_+}$ coincides with $\psi$. Then $\wt{\psi}|_{U-K} \in \Fc|_{U-K}$, and we define
\[
b_\Fc(\psi)\,=\, \beta\bigl( \wt{\psi}|_{U-K}\bigr) \,\in \, \Phi(\Fc).
\]
In other words, we continue $\psi$ along the dotted contour in Fig. \ref{fig:disk-1} and look at the ``multivaluedness'' of this continuation
near $0$. The map $b_\Fc$ is known as the {\em variation map}.

\vskip .2cm

The equivalence \eqref{eq:eq-phipsi} depends on the choice of the cut $K$. In particular, $\Phi(\Fc)$ and $\Psi(\Fc)$
are, intrinsically, not vector spaces but local systems on the (infinite-dimensional) space $\Kc$ formed by all possible cuts $K$
(i.e., all simple closed curves joining $0$ with the boundary of $D$). This space is homotopy equivalent to the circle $S^1$ formed
by directions at $0$. So we can think of $\Phi(\Fc)$ and $\Psi(\Fc)$ as being local systems on this circle
and denote these local systems  $\bPhi(\Fc)$ and $\bPsi(\Fc)$ Their monodromies are the
operators $T_\Phi$ and $T_\Psi$ above.

\paragraph{The transport maps $m_{ij}(\gamma)$ and the  Picard-Lefschetz identities.}\label{par:PL-id}
Let $X, A=\{w_1,\cdots, w_N\}$ be  as in \n  \ref{par:setup}
For any $w\in Y$ we have the natural {\em circle of directions}
\[
S^1_w = \bigl(T_w Y -\{0\}\bigr)/\RR_{>0}^*.
\]
Let $\Fc\in\Perv(Y,A)$. For any $i=1,\cdots, N$ we can restrict $\Fc$ to a small disk around $w_i$ and associate to it
local systems $\bPhi_i = \bPhi_{w_i}(\Fc)$, $\bPsi_i = \bPsi_{w_i}(\Fc)$ on $S^1_{w_i}$, as in \n \ref{par:phi-psi}

\vskip .2cm

Let $\alpha$ be a simple, piecewise smooth  arc in $X$ joining  two marked points $w_i$ and $w_j$
and not passing through any other marked points, see  Fig. \ref {fig:transport}.
Denote by $\dirr_i(\alpha)\in S^1_{w_i}$ and $\dirr_j(\alpha)\in S^1_{w_j}$ the tangent directions to $\alpha$ at the
end points
and  equip $\alpha$ with the orientation going from $w_i$ to $w_j$.
We want to study $\Fc$ in a neighborhood of $\alpha$. Considering $\alpha$ as a closed subset in $\Sigma$,
we have the sheaf $\ul\HH^1_\alpha(\Fc)$ on $\alpha$ which is constant on the open arc $\alpha-\{w_i, w_j\}$.

\bef[H]
\centering
\btp [scale=.8, baseline=(current  bounding  box.center)]

\node (i) at (0,0){};
\fill (i) circle (0.1);

\node (j) at (12,0){};
\fill (j) circle (0.1);

\node (k) at (5,-2){};
\fill (k) circle (0.1);

\node (b) at (4,0.3){$$};
\fill (b) circle (0.08);

\draw [dotted, line width =0.3mm]  (0,0) circle (1cm);
\draw [dotted, line width =0.3mm]  (12,0) circle (1cm);

\draw[->, line width = 0.2mm] plot [smooth,tension=1.5] coordinates{
(0,0) (4,0.3)
(8,-0.3)
(12,0)
};
\node at (-0.5,0){$w_i$};
\node at (12.5, 0){$w_j$};
\node at (5.5, -2){$w_k$};
\node at (9,-0.7){\large$\alpha$};
\node at (4,-0.3){$p$};
\node at (4,0.8){$\Psi_\alpha$};
\node at (1.4, 0.7){$\Phi_{i,\alpha}$};
\node at (10.5, 0.3){$\Phi_{j,\alpha}$};

\etp

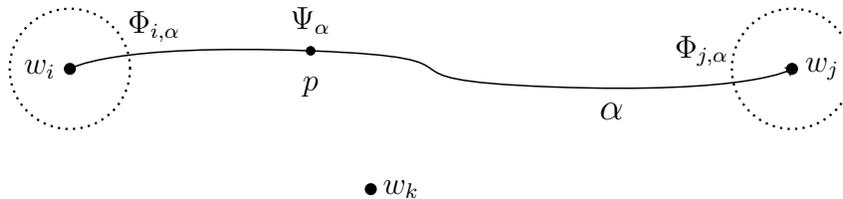
\caption{ The transport map.  }
\label{fig:transport}
\enf
The stalks
\[
\Phi_{i,\alpha}=\ul\HH^1_\alpha(\Fc)_{w_i}, \quad \Phi_{j,\alpha}=\ul\HH^1_\alpha(\Fc)_{w_j}
\]
are just the stalks of the local systems $\bPhi_i, \bPhi_j$ at the tangent directions $\dirr_i(\alpha)$ and $\dirr_j(\alpha)$
respectively.  The stalk
\[
\Psi_\alpha \,=\,\Gamma(\alpha -\{w_i, w_j\}, \ul\HH^1_\alpha(\Fc))
\]
is identified with the stalk $\Fc_p$ at any intermediate point $p\in\alpha$. Note that this identification depends on the
chosen orientations of $X$ and $\alpha$ which give a {\em co-orientation} of $\alpha$, i.e., a particular numeration
of the $2$-element set of ``sides'' of $\Sigma-\alpha$ near $p$ (change of either orientation incurs a minus sign in the identification).  These spaces are connected by the maps
\[
\xymatrix{
\Phi_{i,\alpha}  \ar@<.4ex>[r]^{a_{i,\alpha}}&\Psi_\alpha  \ar@<.4ex>[l]^{b_{i,\alpha}}
\ar@<-.4ex>[r]_{b_{j,\alpha}}
&
\Phi_{j, \alpha}
\ar@<-.4ex>[l]_{a_{j,\alpha}},
}
\]
obtained from the description of $\Fc$ on small disks near $w_i$ and $w_j$ and using $\gamma$ as  the choice for a cut $K$.
We define the {\em transport map}
along $\alpha$ as
\be\label{eq:transport}
m_{ij}(\alpha) = b_{j,\alpha}\circ a_{i,\alpha}: \Phi_{i,\alpha} \lra \Phi_{j,\alpha}.
\ee
By an {\em admissible isotopy} of paths from $w_i$ to $w_j$ we mean a  continuous $1$-parameter family of simple arcs
$(\alpha_t)_{t\in[0,1]}$,
each $\alpha_t$ joining $w_i$ with $w_j$ and not passing through any other $w_k$.
The maps $m_{ij}(\alpha)$ are ``unchanged'' under such an admissible isotopy. More precisely, we have a commutative
diagram
\[
\xymatrix{
\Phi_{i,\alpha_0}
\ar[d]_{t_i}
\ar[r]^{m_{ij}(\alpha_0)} & \Phi_{j, \alpha_0}\ar[d]^{t_j}
\\
\Phi_{i,\alpha_1} \ar[r]^{m_{ij}(\alpha_1)} & \Phi_{j, \alpha_1},
}
\]
where $t_i$ is the monodromy of the local system $\bPhi_i$ on $S^1_{w_i}$ from $\dirr_i(\alpha_0)$ to
$\dirr_i(\alpha_1)$, and similarly for $t_j$.

\vskip .2cm

We now describe what happens when a path crosses a marked point.  That is, we consider a situation as in Fig. \ref{fig:PicLef},
where a path $\gamma'$ from $w_i$ to $w_k$ approaches the composite path formed by $\beta$ from $w_i$ to $w_j$
and $\alpha$ from $w_j$ to $w_k$. After crossing $w_j$, the path $\gamma'$ is changed to $\gamma$.

\bef[H]
\centering
\btp[scale=.8, baseline=(current  bounding  box.center)]

\node (i) at (3.5, -3){};
\fill(i) circle (0.15);

\node (k) at (-3, -3.5){};
\fill (k) circle (0.15);

\node (j) at (.7, -1){};
\fill (j) circle (0.15);

 \centerarc[line width=0.5](0.7, -1)(325:580:0.7)

\draw [ decoration={markings,mark=at position 0.7 with
{\arrow[scale=1.5,>=stealth]{>}}},postaction={decorate},
line width = .3mm]   (-3,-3.5) .. controls (1,2) ..  (3.5, -3) ;

\draw [decoration={markings,mark=at position 0.7 with
{\arrow[scale=1.5,>=stealth]{>}}},postaction={decorate},
line width = .3mm]   (-3,-3.5) .. controls (1,-4) ..   (3.5,-3);

\draw [decoration={markings,mark=at position 0.7 with
{\arrow[scale=1.5,>=stealth]{>}}},postaction={decorate},
line width = .3mm]  (0.7,-1)  .. controls (2,-2) ..  (3.5,-3) ;

\draw [ decoration={markings,mark=at position 0.7 with
{\arrow[scale=1.5,>=stealth]{>}}},postaction={decorate}, line width = .3mm]  (-3,-3.5) .. controls (-2,-3) .. (0.7,-1);

\node at (4.5, -3.5){$w_k$};

\node at ( -4,-3.5){$w_i$};
\node at ( 0.7,-0.6){$w_j$};
\node at (0,-4.2){\large$\gamma'$};
\node at (-.8,-1.7){\large$\beta$};
\node at (2,-1.5){\large$\alpha$};
\node at (1.4, 1.){\large$\gamma$};

\etp

\caption{ The Picard-Lefschetz situation.  }
\label{fig:PicLef}
\enf
In this case we have identifications
\be\label{eq:3-isotop}
\Phi_{i,\gamma} \to \Phi_{i,\beta} \to \Phi_{i, \gamma'}, \quad
\Phi_{k,\gamma'}\to \Phi_{k, \alpha}\to  \Phi_{k,\gamma},
\quad \Phi_{j, \beta} \to \Phi_{j, \alpha},
\ee
given by {\em clockwise} monodromies of the local systems $\bPhi$ around the corresponding arcs in the circles of directions.
So after these identications we can assume that we deal with single vector
spaces denoted by $\Phi_i, \Phi_k$ and $\Phi_j$ respectively.

\begin{conven}\label{conv:clock}
In this paper we will often consider compositions of several transport maps
(and  later,  transport functors) corresponding
to  composable sequences of several oriented paths.
In all such cases our default convention will be to identify
 the intermediate $\Phi$-vector spaces (and later, categories) using the {\em clockwise} monodromies
 from the incoming direction to the outgoing one.
This extends the
   last identification in \eqref{eq:3-isotop},  which is indicated by an arc around $w_j$ in
   Fig. \ref{fig:PicLef}.
\end{conven}

\begin{prop}[Abstract Picard-Lefschetz identity]\label{prop:PL-form}
We have the equality of linear operators $\Phi_i\to\Phi_k$:
\[
m_{ik}(\gamma') = m_{ik}(\gamma) - m_{jk}(\alpha) m_{ij}(\beta).
\]
\end{prop}

This was proved in \cite[Prop. 1.8]{KS-schobers} but with slightly different sign conventions. It can be seen as a disguised version of the identity $T_\Psi=\Id_\Psi-ab$ for the monodromy in the $(\Phi,\Psi)$-description of $\Perv(D,0)$. For convenience of the reader we give a detailed argument.

\vskip .2cm

\noindent{\sl Proof:} Let us deform the situation of Fig. \ref{fig:PicLef} by pinching the
directions of the three paths at $w_i$ into a single direction, so that they all go together from
$w_i$ to a nearby point $\wt w_i$ and then diverge. Let us do the same at $w_k$, see
Fig. \ref{fig:PicLef-pinch}.

\bef[H]
\centering
\btp[scale=.8, baseline=(current  bounding  box.center)]

\node (wtk) at (3.5, -3){$\bullet$};

\node (wti) at (-3, -3.5){$\bullet$};

\node (j) at (.7, -1){$\bullet$};

\draw [ decoration={markings,mark=at position 0.7 with
{\arrow[scale=1.5,>=stealth]{>}}},postaction={decorate},
line width = .3mm] (-3,-3.5)  .. controls (1,2) ..  (3.5, -3) ;

\draw [ decoration={markings,mark=at position 0.7 with
{\arrow[scale=1.5,>=stealth]{>}}},postaction={decorate},
line width = .3mm]  (-3,-3.5) .. controls (1,-4) ..   (3.5,-3) ;

\draw [ decoration={markings,mark=at position 0.7 with
{\arrow[scale=1.5,>=stealth]{>}}},postaction={decorate},
line width = .3mm]   (0.7,-1) .. controls (2,-2) ..   (3.5,-3);

\draw [ decoration={markings,mark=at position 0.7 with
{\arrow[scale=1.5,>=stealth]{>}}},postaction={decorate},
line width = .3mm]  (-3,-3.5) .. controls (-2,-3) .. (0.7,-1);

\node at (4.0, -2.6){$\wt w_k$};
\node at (5.5, -3.5) {$w_k$};
\node at ( -3.5,-3){$\wt w_i$};
\node at (-5, -4){$w_i$};
\node at ( 0.7,-0.6){$w_j$};
\node at (0,-4.2){\large$\wt\gamma'$};
\node at (-.8,-1.7){\large$\wt\beta$};
\node at (2,-1.5){\large$\wt\alpha$};
\node at (1.4, 1.){\large$\wt\gamma$};

\node (i) at (-4.5,-4){$\bullet$};
\node (k) at (5, -3.5){$\bullet$};

\draw  [ decoration={markings,mark=at position 0.7 with
{\arrow[scale=1.5,>=stealth]{>}}},postaction={decorate},
line width = .3mm]  (i.center) -- (wti.center) ;

\draw  [  decoration={markings,mark=at position 0.7 with
{\arrow[scale=1.5,>=stealth]{>}}},postaction={decorate},
line width = .3mm] (wtk.center) -- (k.center);

\etp

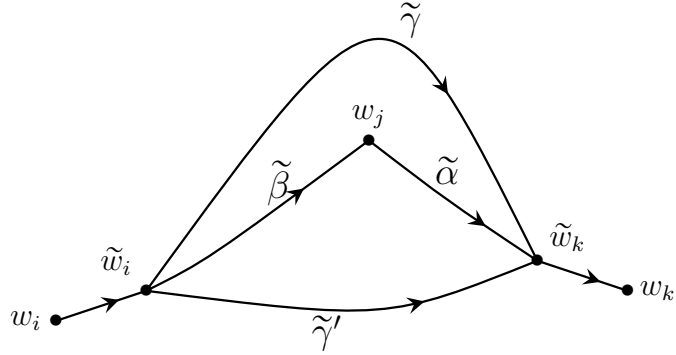
\caption{ The Picard-Lefschetz situation, pinched.  }
\label{fig:PicLef-pinch}
\enf

It is enough to prove the statement for this deformed situation which we now assume.
Thus, $\gamma$ now consists of a short segment $[w_i, \wt w_i]$
(not necessarily straight) followed by a path $\wt\gamma$ followed by
a similar segment $[\wt w_k, w_k]$. Similarly,
\[
\gamma' \,=\, [w_i, \wt w_i]\cup \wt\gamma'\cup [\wt w_k, w_k],
\quad \alpha \,=\,\wt\alpha  \cup [\wt w_k, w_k], \quad \beta = [w_i, \wt w_i]\cup\wt\beta.
\]
We view $\wt w_i$ and $\wt w_k$ as the midpoints to define the transport maps.
Let $\Psi_i, \Psi_k$ be the stalks of $\Fc$ at $\wt w_i, \wt w_k$ respectively.
Then, we have maps
\[
a_i: \Phi_i\lra \Psi_j, \quad b_k: \Psi_k\lra\Phi_k
\]
corresponding to the intervals $[w_i, \wt w_i]$ and $[w_k, \wt w_k]$, the monodromy maps
\[
T_{\wt\gamma}, T_{\wt\gamma'}: \Psi_i\lra \Psi_k
\]
corresponding to $\wt\gamma,\wt\gamma'$  and the maps
\[
a_{j,\wt\beta}: \Phi_{j,\beta}\lra \Psi_i, \quad b_{j,\wt\beta}: \Psi_i \lra \Phi_{j,\beta}
\]
\[
a_{j,\wt\alpha}: \Phi_{j,\alpha}\lra \Psi_k, \quad b_{j,\wt\alpha}: \Psi_k \lra \Phi_{j,\alpha}
\]
corresponding to the segments near $w_j$. So by definition,
\be\label{eq:m_{ij}-deformed-noncat}
\begin{gathered}
m_{ik}(\gamma) \,=\, b_k \circ T_{\wt\gamma}\circ a_i, \quad
m_{ik}(\gamma') \,=\, b_k\circ T_{\wt\gamma'}\circ a_i,
\\
m_{ij}(\beta)\,=\, b_{j,\wt\beta}\circ a_i, \quad m_{jk}(\alpha) \,=\, b_k\circ a_{j,\wt\alpha}.
\end{gathered}
\ee

Note that the transformation $T = T_{\wt\gamma}^{-1}\circ T_{\wt\gamma'}: \Psi_i\to\Psi_i$ is a full counterclockwise monodromy of $\Psi_i$ around $w_j$. Hence, \[T = 1-a_{j, \wt\beta}b_{j, \wt\beta}\]

Now, we also have \[T_{\wt\gamma} a_{j, \wt\beta} = a_{j, \wt\alpha}\]
due to the clockwise identification $\Phi_{j, \beta} \simeq \Phi_{j, \alpha}$. Combining previous equalities, we obtain
\[T_{\wt\gamma'} = T_{\wt\gamma}T = T_{\wt\gamma}(1-a_{j, \wt\beta}b_{j, \wt\beta}) = T_{\wt\gamma} - a_{j, \wt\alpha}b_{j, \wt\beta}\]

Multiplying this equality by $b_k$ on the left and $a_i$ on the right we obtain the desired statement. \qed

\begin{figure}[H]
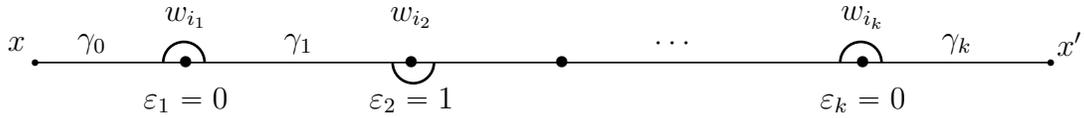

\centering
\btp[scale=0.5]

\node at (0,0){$\bullet$};
\node at (8,0){$\bullet$};
\node at (-4,0){$\bullet$};
\node at (-10,0){$\bullet$};
\draw[line width = 0.6] (-14,0) -- (13,0);

\node at (-14.5,0.5){$x$};
\node at (13.5,0.5){$x'$};

\node at (-10, -1){$\eps_1=0$};
\node at (-4, -1){$\eps_2=1$};
\node at (8, -1){$\eps_k=0$};

\node at (-10,1.2){$w_{i_1}$};
\node at (-4,1.2){$w_{i_2}$};
\node at (8,1.2){$w_{i_k}$};

 \draw[line width=1] ($(-10,0)+({0.5*cos(0)},{0.5*sin(0)})$) arc (0:180:0.55);

  \draw[line width=1] ($(-4,0)+({0.5*cos(180)},{0.5*sin(1800)})$) arc (180:360:0.55);

   \draw[line width=1] ($(8,0)+({0.5*cos(0)},{0.5*sin(0)})$) arc (0:180:0.55);

   \node at (-12.5, 0.5){$\gamma_0$};
   \node at (-7, 0.5){$\gamma_1$};
   \node at (10.5, 0.5) {$\gamma_k$};

   \node at (3, 0.5){$\cdots$};

\node at (-14,0){\tiny$\bullet$};
\node at (13,0){\tiny$\bullet$};

\etp
\caption{Iterated transport situation.} \label{fig:iter-trans}
\end{figure}

\paragraph{The iterated transport map as the ``Vassiliev derivative'' of
the monodromy.}\label{par:vassiliev}
Consider the following {\em interated transport situation} depicted in Fig. \ref{fig:iter-trans}.
That is, let $x,x'\in X\- A$ and $\gamma$ be a (piecewise smooth) simple path joining $x$ with $x'$ and
passing, possibly, through some $w_{i_1}, \cdots, w_{i_k}\in A$, in this order.

Let $\gamma_0$ be the part of $\gamma$ between $x$ and $w_{i_1}$, let $\gamma_1$ be the part
of $\gamma$ between $w_{i_1}$ and $w_{i_2}$, and so on, until $\gamma_k$ which is defined to
be the part of $\gamma$ between $w_{i_k}$ and $x'$.
Let us define the {\em iterated transport map}
\be\label{eq:m(gamma)}
m(\gamma) \,=\, a_{i_k, \gamma_k}\circ m_{i_{k-1}, i_k}(\gamma_{k-1})\circ\cdots
\circ m_{i_1, i_2}(\gamma_1)\circ b_{i_1,\gamma_0}: \, \Fc_x \lra\Fc_{x'}.
\ee
Here the identification of the $\Phi$-categories at the intermediate points is made using the upward
(clockwise)
half-monodromies understood with respect to the direction of $\gamma$ from $x$ to $x'$ and to the
chosen orientation of $X$, see Convention \ref{conv:clock}.

If $\gamma$ does not pass through any $w_i$, then $k=0$ and $m(\gamma)$ is defined
to be the monodromy of $\Fc$ along $\gamma$.

\vskip .2cm

Let $\eps_1,\cdots, \eps_k\in \{0,1\}$. Denote by $\gamma_{\eps_1,\cdots, \eps_k}$ the deformation
of $\gamma$ obtained by bending it around $w_{i_\nu}$ upwards (in the above sense),
if $\eps_\nu=0$ and downwards, if $\eps_\nu=1$, see Fig. \ref{fig:iter-trans}. Thus
$\gamma_{\eps_1,\cdots,\eps_k}$ does not pass through any $w_i$.

The following is obtained by iterating the argument in the proof of Proposition
\ref{prop:PL-form}.

\begin{cor}\label{cor:PL-vassiliev}
We have an equality of maps $\Fc_x\to\Fc_{x'}$:
\[
m(\gamma) \,=\,\sum_{\eps_1,\cdots,\eps_k\in\{0,1\}} (-1)^{\eps_1+\cdots + \eps_k}
 m(\gamma_{\eps_1,\cdots, \eps_k}). \qed
\]
\end{cor}
This has an attractive interpretation in terms of Vassiliev's approach to invariants of knots, cf.
\cite{vassiliev} Ch. V. Let $\Pi$ be the (infinite-dimensional) space of all simple paths joining $x$
with $x'$. For each $i=1,\cdots, N$ we have a hypersurface $H_i\subset \Pi$ consisting of paths
passing through $w_i$. The hypersurface $\Hc=\bigcup_{i=1}^N H_i$ (``discriminant'')
cuts $\Pi$ into various {\em chambers} (defined as connected components of $\Pi\- \Hc$).
For each chamber $C$ we have a well defined monodromy
\[
m(C) \,=\, m(\gamma): \Fc_x\lra\Fc_{x'}, \quad \forall\,  \gamma\in C.
\]
The correspondence $C\mapsto m(C)$ can be seen as an ``invariant'' in the spirit of
\cite{vassiliev}. Further, the discriminant itself is further stratified into ``faces'' by multiplicity of intersections
of different components $H_i$.

 Corollary \ref{cor:PL-vassiliev} means that $m(\gamma)$ for $\gamma$ lying on such a face
 $F\subset H_{i_1}\cap \cdots H_{i_k}$ is obtained as the alternating sum of the $m(C)$,
 where $C$ runs over the $2^k$ chambers surrounding $F$. In other words, $m(\gamma)$
 is realized as the $k$th Vassiliev derivative of the invariant $C\mapsto m(C)$.

 We will take up this point of view in a slightly different form in \S \ref{subsec:infra-com-schob}.

\paragraph{The spider descriptions of $\Perv(D,A)$ and  $\Perv(\CC, A)$.}\label{par:spider}
Let  $D$ be the unit disk and $A=\{w_1,\cdots, w_N\}\subset D-\del D$. A description of the
category $\Perv(D,A)$ was given by Gelfand-MacPherson-Vilonen \cite{gelfand-MV}.
To recall it, it is convenient to
introduce the following terminology. Fix a point $\bv\in \del D$. By a {\em $\bv$-spider} for $(D,A)$ we will mean a system
$K = \{\gamma_1,\cdots, \gamma_N\}$ of simple closed  piecewise smoorth arcs in $D$ so that $\gamma$ joins $v$ with $w_i$ and
different $\gamma_i$ do not meet except at $v$. See the left part of Fig. \ref{fig:disk-2}.

Note that a $\bv$-spider for $(D,A)$ defines a total order on $A$, by {\em clockwise}
ordering of the slopes of the $\gamma_i$ at $\bv$ (assumed distinct). Let us assume that the numbering
$A=\{w_1, \cdots, w_N\}$ is chosen in this order, see again the left part of Fig. \ref{fig:disk-2}.

\bef[H]
\centering
\btp[scale=.4, baseline=(current  bounding  box.center)]

\node (1)  at (-1, -2.5){};
\fill (1) circle (0.15);

\node (2)  at (-2.5, -1){};
\fill (2) circle (0.15);

\node (n)  at (-2,3){};
\fill (n) circle (0.15);

\draw (0,0) circle (7cm);

\node (b) at (7,0){};
\fill (b) circle (0.15);

\draw[line width = .3mm] (7,0) .. controls (2, -1) .. (-1, -2.5);

\draw[line width = .3mm] (7,0) .. controls (0,1) .. (-2.5,-1);

\draw[line width = .3mm] (7,0) .. controls (2,-6) and (-9,-6)   .. (-2,3);

\node at (7.7, -.5) {$\bv$};

\node at (6,-6) {$D$};

\node at (-1.6, -2.8) {$w_2$};

\node at (-2,-2){$\cdots$};

\node at (-3.3, -1.5){$w_N$};

\node at (-1, +3){$w_1$};

\node at (-4,+2){$\gamma_1$};

\node at (2,-2){$\gamma_2$};

\node at (0.5,0){$\gamma_N$};

\etp
\quad\quad\quad
\btp[scale=.4, baseline=(current  bounding  box.center)]

\node (1)  at (4, -3){};
\fill (1) circle (0.15);
\node at (5.5, -3) {$w_1$};
\draw[line width = .3mm] (0,0) .. controls (3, -1) .. (4, -3);
\node at (2, -1.5){$\gamma_1$};

\node (2)  at (-2.5,-4){};
\fill (2) circle (0.15);
\node at (-4,-4){$w_2$};
\draw[line width = .3mm] (0,0) .. controls (-1,-1) ..  (-2.4,-4);
\node at (-2.5,-2){$\gamma_2$};

\node (3) at (-4,+4){};
\fill (3) circle (0.15);
\node at (-4.7, +4.7){$w_3$};
\draw[line width = .3mm] (0,0) .. controls (-2,+3) ..  (-4,+4);
\node at (-2,+1.5){$\gamma_3$};

\node (n)  at (3,+2){};
\fill (n) circle (0.15);
\node at (3,+3){$w_N$};
\draw[line width = .3mm] (0,0) .. controls (1,+1) .. (3,+2);
\node at (1,+1.7){$\gamma_N$};

\node (v) at (0,0){};
\fill (v) circle (0.15);
\node at (-1,0){$\bv$};

\node at (0,+3.5){$\cdots$};

\node at (6,6) {$\CC$};

\draw (-7,-7) -- (-7,7) -- (7,7) -- (7,-7) -- (-7,-7);

\etp

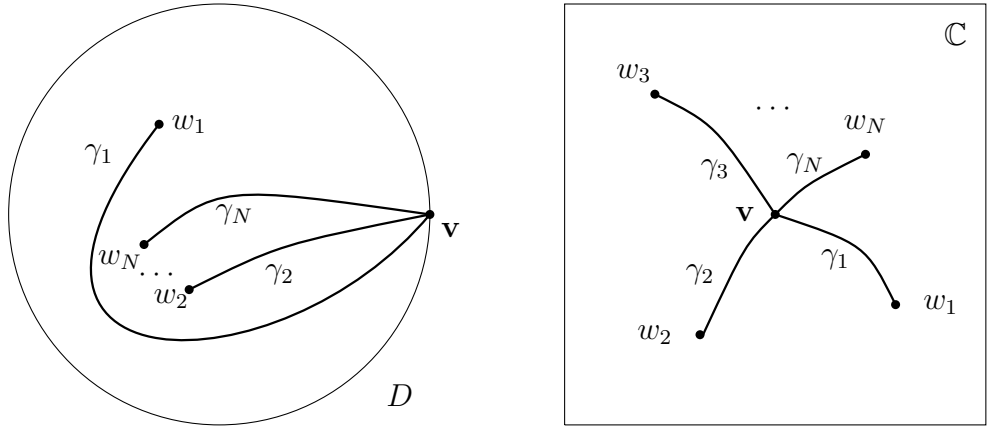
\caption{ Equivalence depending on a spider.  }
\label{fig:disk-2}
\enf

\begin{Defi} Let  $\Qc_N$ be the category formed by diagrams (quivers) of finite-dimensional $\k$-vector spaces
\[
\xymatrix{\Phi_N
\ar@<.4ex>[rd]^{a_N}
&
\\
\vdots & \Psi \ar@<.4ex>[ul]^{b_N}
\ar@<.4ex>[dl]^{b_1}
\\
\Phi_1 \ar@<.4ex>
[ur]^{a_1}
}
\]
such that   $T_{\Phi_i}:= \Id_{\Phi_i}-b_ia_i$ is invertible for each $i$. This implies that
each $T_{i, \Psi} = \Id_\Psi-a_ib_i$ is  invertible as well.
\end{Defi}

Let now $K$ be a $\bv$-spider for $(D,A)$ and $\Fc\in\Perv(D,A)$.
For any $i=1,\cdots, N$ we  denote $\Phi_{i, K}(\Fc)= \Phi_{i,\gamma_i}(\Fc)$
the stalk of the local system $\bPhi_i(\Fc)$ at the point in $S^1_{w_i}$ represented by the direction of
$\gamma_i$.
Let us also identify the stalk $\Psi_{\gamma_i}(\Fc)$ with the stalk $\Fc_v$,
as in \S \ref{subsec:abs-PL}  \ref{par:PL-id}
This gives a diagram
\[
\Theta_K(\Fc) \,=\, ( \Psi(\Fc),  \Phi_{i, K}(\Fc),  a_{i, \gamma_i} , b_{i, \gamma_i}  ) \,\in \, \Qc_N.
\]
Here $a_{i, \gamma_i}$ and $b_{i, \gamma_i}$ are the canonical maps along $\gamma_i$ as in
\eqref{eq:transport}.

\begin{prop} [\cite {gelfand-MV}] \label{prop:GMV1}
Let $K$ be a $\bv$-spider for $(D,A)$.
The functor $\Theta_K: \Perv(D,A)\to\Qc_N$ is an equivalence.
\qed
\end{prop}

\begin{figure}[H]
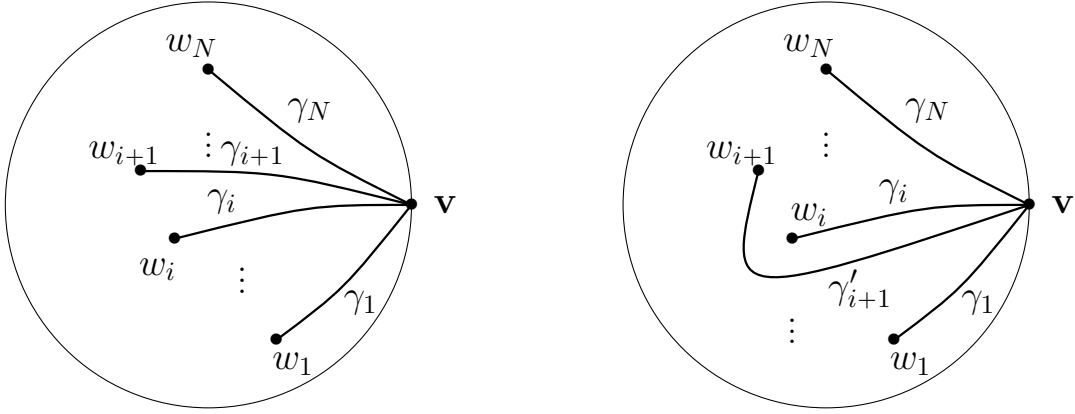

\centering
\btp[scale=.45, baseline=(current  bounding  box.center)]
\draw (0,0) circle (6cm);

\node at (6,0){$\bullet$};
\node at (2,-4){$\bullet$};
\node at (-1,-1){$\bullet$};
\node at (-2,+1){$\bullet$};
\node at (0, +4){$\bullet$};
\node at (7,0){\large$\bv$};

\draw[line width = .3mm] (6,0) .. controls (4, -2.5) ..  (2,-4);
  \draw[line width = .3mm] (6,0) .. controls (3,0) ..  (-1,-1);
   \draw[line width = .3mm] (6,0) .. controls (2,+1) ..  (-2,+1);
 \draw[line width = .3mm] (6,0) .. controls (3,+1.5) ..  (0,+4);

\node at (2.5, -4.7){\large $w_1$};
\node at (4.5, -2.9){\large$\gamma_1$};
\node at (1, -2){\large$\vdots$};
\node at (-1.5, -1.9){\large$w_i$};
\node at (0.4, +0.1){\large$\gamma_i$};
\node at (-2.5, +1.5){\large$w_{i+1}$};
\node at (0, +2){\large$\vdots$};
\node at (1.3, +1.5){\large$\gamma_{i+1}$};
\node at (-0.5, +4.7){\large$w_N$};
\node at (3, +2.8){\large$\gamma_N$};

\etp
\quad\quad\quad\quad\quad
\btp[scale=.45, baseline=(current  bounding  box.center)]
 \draw (0,0) circle (6cm);
 \node at (6,0){$\bullet$};
\node at (2,-4){$\bullet$};
\node at (-1,-1){$\bullet$};
\node at (-2,+1){$\bullet$};
\node at (0, +4){$\bullet$};
\node at (7,0){\large$\bv$};

\draw[line width = .3mm] (6,0) .. controls (4, -2.5) ..  (2, -4);
\draw[line width = .3mm] (6,0) .. controls (3,0) ..  (-1, -1);
\draw[line width = .3mm] (6,0) .. controls (-3, -3) ..  (-2,+1);
\draw[line width = .3mm] (6,0) .. controls (3,+1.5) ..  (0,+4);

\node at (2.5, -4.7){\large $w_1$};
\node at (4.5, -2.9){\large$\gamma_1$};
\node at (-0.5, -0.3){\large$w_i$};
\node at (2, +0.4){\large$\gamma_i$};
\node at (1, -2.5){\large$\gamma'_{i+1}$};
\node at (0, +2){\large$\vdots$};
\node at (-2.5, +1.5){\large$w_{i+1}$};
\node at (-0.5, +4.7){\large$w_N$};
\node at (3, +2.8){\large$\gamma_N$};
\node at (-1, -3.5){\large$\vdots$};

 \etp
\caption{Action of $\tau_i\in\Br_N$ on spiders.}\label{fig:br-spid}
\end{figure}

\begin{rems}\label{rem:vladi}
(a) We can think of the point $\bv\in\del D$ as being  being far away at the infinity (``Vladivostok'').
For this reason we will sometimes refer to the description of $\Perv(D,A)$ given by Proposition
\ref{prop:GMV1} as the {\em Vladivostok description}.

\vskip .2cm

(b) The significance of the choice of  position of $\bv$ becomes clearer when,
instead of the disk $D$, we can consider the complex plane $\CC$ so that the boundary of $\CC$
is not given. Then we can take for $\bv$ any point in $\CC\- A$ and consider
spiders spreading out from that point, see the right part of Fig. \ref{fig:disk-2}. This will
also lead to an identification $\Perv(\CC,A)\to\Qc_N$.

\end{rems}
We further recall the dependence of the equivalence $\Theta_K$ on the choice of a spider $K$.
Note that isotopic spiders give equivalences which are {\em canonically isomorphic} as functors
because the space of spiders in a given isotopy class (i.e., each connected component of the space of
spiders) is contractible.
Let $\Kc(D,\bv,A)$ be the set of isotopy classes of $\bv$-spiders for $(D,A)$. It is classical, see, e.g.,
\cite{AGV} \S 2.6, that  $\Kc(D,\bv,A)$ forms a torsor over the braid group $\Br_N$.
One way of seeing this is to use the realization
$\Br_N=\Gamma^+(D,A)$ as the mapping class group, see Example \ref{ex:class-group-braid}. That is,
every diffeomorphism $f: D\to D$ preserving $A$ as a set and $\del D$ pointwise, takes
a spider $K$  to a new spider $f(K)$. More explicitly, recall that $\Br_N$ is generated by
elementary twists $\tau_1, \cdots, \tau_{N-1}$ subject to the relations
\[
\tau_i\tau_j=\tau_j \tau_i, \,\,\, |i-j|\geq 2, \quad \tau_i\tau_{i+1}\tau_i = \tau_{i+1}\tau_i\tau_{i+1}.
\]
In this presentation, the action of $\tau_i$ on a spider $K = \{\gamma_1,\cdots, \gamma_N\}$
leaves $\gamma_j$, $j\neq i+1$, unchanged and replaces $\gamma_{j+1}$ by
the path $\gamma'_{j+1}$ obtained by moving $\gamma_{j+1}$ counterclockwise around $\gamma_j$,
see Fig. \ref{fig:br-spid}.    Thus the homomorphism $\Br_N\to S_N$ comes from the change of the order of the slopes of
the $\gamma_j$ at $\bv$.

\begin{prop}[\cite{gelfand-MV}, Prop. 1.3] \label{prop:ch-spid-1}
(a) Replacing $K$ by $\tau_i(K)$ changes a diagram
$(\Psi, \Phi_j, a_j, b_j)$ to  $(\Psi, \tau_i(\Phi_j), \tau_i(a_j), \tau_i(b_j))$, where:
 \[
\begin{gathered}
\tau_i(\Phi_1,\cdots, \Phi_N)\,=\, (\Phi_1, \cdots, \Phi_{i-1}, \Phi_{i+1}, \Phi_i, \Phi_{i+2}, \cdots, \Phi_N),
\\
\tau_i(a_1,\cdots, a_N) \,=\,
(a_1, \cdots, a_{i-1}, \, T_{i,\Psi} \, a_{i+1}, \, a_i, a_{i+2}, \cdots, a_N),
\\
\tau_i(b_1,\cdots, b_N) \,=\,
(b_1,\cdots, b_{i-1}, \, b_{i+1} \, T_{i, \Psi}^{-1}, \, b_i, b_{i+2}, \cdots, b_N).
\end{gathered}
\]

(b) Put differently, the formulas in (a) define an action of $\Br_N$ on the category $\Qc_N$.
For any spider $K\in\Kc(D,\bv,A)$ (and the corresponding numbering of $A$ by slopes of the $\gamma_i$ at $\bv$),
the equivalence $\Theta_K$
takes the geometric action of $\Br_N=\Gamma^+_{D,A}$ on $\Perv(D,A)$ (Proposition \ref {prop:Gamma-acts}) to the above
action   of $\Br_N$ on $\Qc_N$.  \qed
\end{prop}


\subsection {Relation to the classical Picard-Lefschetz theory.}\label{subsec:PL-class}

\paragraph{The Lefschetz perverse sheaf.}\label{par:PL-perv}

When speaking of complex manifolds, we allow them to have a boundary.
By a {\em Riemann surface} we mean a $1$-dimensional complex manifold (possibly with boundary).

\begin{Defi}\label{def:GLP}
Let  $X$ be a Riemann surface.
An $X$-valued {\em  generalized Lefschetz pencil} of relative dimension $n$ is a datum $(Y,W)$, where:
\begin{itemize}
\item[(1)]    $Y$ is an
$(n+1)$-dimensional complex manifold.

\item[(2)] $W: Y \to X$ is a holomorphic,  proper (fibers are compact without boundary)
 map  with only finitely many isolated
singular points, lying in the interior of $Y$.
 We denote  the critical points $y_1,\cdots, y_N$, and put $w_i=W(y_i)$ to be the corresponding critical values.
We assume that all the $w_i$ are distinct and put $A=\{w_1,\cdots, w_N\}\subset X$.
\end{itemize}

A {\em Lefschetz pencil} is a generalized Lefschetz pencil with only
quadratic (Morse) singular points.

\end{Defi}

Given a generalized Lefschetz pencil  $(Y,W)$, we can form  a constructible complex
$RW_* \,  \ul\k_Y \in D^b(X ,A)$. It gives rise to:

\begin{itemize}
\item Constructible sheaves $R^q W_*\, \ul\k_Y = \ul H^q (RW_* \,\ul\k_Y)$.

\item Perverse sheaves
$R^q_\perv W_* \ul\k_Y = \ul H^q_\perv (RW_*\,\ul\k_Y) \in \Perv(\Sigma,A)$.  In particular, we have the {\em Lefschetz perverse sheaf}
\[
\Lc = \Lc_W = R^{n}_\perv W_*\, \ul\k_Y,
\]
corresponding to the middle dimension.
\end{itemize}

\begin{prop} [Picard-Lefschetz theory I]
Let $(Y,W)$ be a generalized Lefschetz pencil. Then:

\begin{itemize}
\item[(a)] For $q\neq n$ the sheaf $R^qW_*\, \ul\k_Y$ is locally constant
(i.e., it does not have a singularity at any point of $A$) and is identifed with
$R^q_\perv W_* \,\ul \k_Y$.

\item[((b)] For $q=n$, the perverse sheaf $\Lc$ is identified, outside $A$, with $R^nW_* \, \ul\k_Y$, the local system of
middle cohomology of the smooth fibers $W^{-1}(b)$.

\item[(c)] Suppose that $(Y,W)$ is a Lefschetz pencil. Then, for any $w_i\in A$,
the local system of vanishing cycles
$\bPhi_{w_i}(\Lc)$ on $S^1_{w_i}$ is $1$-dimensional, with monodromy $(-1)^{n+1}$.

\end{itemize}
\end{prop}

\bef[H]

\centering
\btp[scale=0.5]

\node (0) at (0,0){};
\fill (0) circle (0.2);

\draw[ line width=.4mm] (0,3) ellipse (2cm and 1cm);
\draw[ line width=.4mm] (0,3) ellipse (2.7cm and 1.7cm);

 \draw[ line width=.4mm] (0,-3) ellipse (2cm and 1cm);
  \draw[ line width=0.4mm] (0,-3) ellipse (2.7cm and 1.7cm);
\draw [ line width=.4mm] ( -2.035, -2.8) -- ( 2.035, 2.8);
\draw [ line width=.4mm] (2.035, -2.8) -- ( -2.035, 2.8);

\pgfmathsetmacro{\e}{1. 77}   
\pgfmathsetmacro{\a}{2}
\pgfmathsetmacro{\b}{(\a*sqrt((\e)^2-1)}
\draw[line width =0.4mm] plot[domain=-0.8:0.8] ({\a*cosh(\x)},{\b*sinh(\x)});
\draw [line width =0.4mm]  plot[domain=-0.8:0.8] ({-\a*cosh(\x)},{\b*sinh(\x)});

\draw[ color=red, line width=.7mm] (0,0) ellipse (2cm and 1cm);

\node at (4.6,3){\large$W^{-1}(\eps)$};
\node (N0) at (9,-5){};
\fill  (N0) circle (0.1);
\draw[line width=0.3mm] (9,-5) circle (1.5);
\node at (8.5, -4.5){$w_i$};
\draw [line width=0.3mm]   (9,-5)  .. controls (9.5, -4.8) .. (10.5, -5);
\node at (10.9, -4.5){$\eps$};
\node (eps) at (10.5, -5){};
\fill (eps) circle (0.1);
\node at (-3.5,0) {\large$\Delta_i$};
\draw [->, dotted, line width = 0.6mm] (3,-1) -- (6.7, -3.7);
\node at (6,-2){\large$W$};

\etp
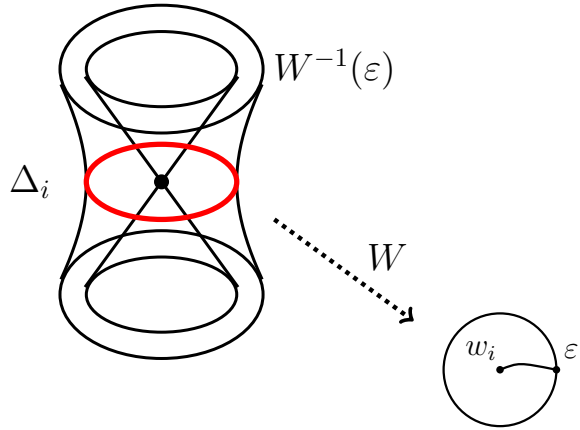
\caption{The vanishing cycle.}\label{fig:vanish}
\enf

Further, in the classical theory, we can think of a radial direction $\eps\in S^1_{w_i}$   as a nearby point of $Y$
connected to $w_i$ by a path, see Fig. \ref{fig:vanish}.  Such a datum
gives rise to
the {\em geometric  (Lefschetz) vanishing cycle} $\Delta_i = \Delta_i(\eps)\subset W^{-1}(\epsilon)$
which is a submanifold homeomorphic to the $n$-sphere $S^n$ and defined canonically up to isotopy.
See Fig. \ref{fig:vanish} and  \cite{AGV}, Ch. 1.
The stalk of $\bPhi_{w_i}(\Lc)$ at $\eps$ can be thought of (after identifying cohomology with
homology by Poincar\'e duality) as spanned by $[\Delta_i]$,  the homology class of $\Delta_i$. Note that there are two choices
of $[\Delta_i]$ corresponding to two orientations of $\Delta_i \simeq S^n$, of which none is preferred.
The monodromy around $w_i$ can be seen geometrically (see \n \ref{par:symcon}) as induced
by the antipodal map of $S^n$ which preserves orientation for
$n+1$ even and changes it for $n+1$ odd, which explains (b).

\begin{prop} [Picard-Lefschetz theory II]\label{prop:PLII}
Let $(Y,W)$ be a Lefschetz pencil.
The maps
\[
\xymatrix{
\bPhi_{w_i}(\Lc)_\eps \ar@<.4ex>[rr]^{\hskip -1.8cm a_i}&&\bPsi_{w_i}(\Lc)_\eps = H^n(W^{-1}(\eps), \k) \ar@<.4ex>[ll]^{\hskip- 1.8cm b_i}
}
\]
describing $\Lc$ on a small disk near $w_i$, have the following form:
\begin{itemize}
\item [(a)] The map $a_i$ considers (any choice of)  $[\Delta_i]$ as an element of the  middle (co)homology of $W^{-1}(\eps)$.

\item[(b)] The map  $b_i$ takes any $\xi\in H^n(W^{-1}(\eps),\k)$ into $([\Delta_i],\xi )\cdot [\Delta_i]$. Here
$[\Delta_i],\xi)\in\k$
is the intersection index (or, if we consider $[\Delta_i]$ as a homology class, just the pairing of homology and cohomology).
\end{itemize}
\end{prop}

Combining Proposition \ref {prop:PLII} with the  identity $T_\Psi=\Id-ab$ gives the classical {\em Picard-Lefschetz theorem}: the monodromy of $H^n(W^{-1}(\eps))$
around $w_i$ is the transvection associated to the vanishing cycle.

\vskip .2cm

Let now $\alpha$ be a path joining $w_i$ and $w_j$ as in \n \ref{subsec:abs-PL} \ref{par:PL-id}
After choosing orientations of the two vanishing cycles $\Delta_i$ and $\Delta_j$, we identify
the vector spaces $\Phi_i$ and $\Phi_j$ with $\k$, so the transport map $m_{ij}(\alpha): \Phi_i\to\Phi_j$
becomes an element of $\k$, in fact, an integer, the intersection index in $W^{-1}(b)$:
\be\label{eq:mij-inter-I}
m_{ij} \,=\, \bigl( [\Delta_i]_b, [\Delta_j]_b\bigr)_{W^{-1}(b)} \,\in \, \ZZ \,\subset \, \k
\ee
Here $b\in\alpha$ is an intermediate point, as in Fig. \ref {fig:PicLef}, while
$[\Delta_i]_b, [\Delta_j]_b\in H^n(W^{-1}(b),\k)$
are the parallel transports of $[\Delta_i]$ and $[\Delta_j]$ along $\alpha$ (with respect to the Gauss-Manin connection)
from the vicinity of $w_i$ or $w_j$ to $b$.

\vskip .2cm

More generally, let $(Y,W)$ be a generalized Lefschetz pencil. Then the space
$\Phi_{w_i}(\Lc_W)_\eps$, with its action of the monodromy operator,
is the classical {\em space of vanishing cyclies}
of the isolated singularity $y_i$ of $W$, see \cite{AGV} \S 2.2.

\paragraph{ The symplectic connection.} \label{par:symcon}
In the situation of \n  \ref{par:PL-perv},
assume that $Y$ is equipped with a K\"ahler metric $h=g+\sqrt{-1}\cdot \omega$ so that $g$ is a Riemannian metric and
$\omega$ is a symplectic form on $Y$ (considered as a $C^\infty$-manifold of dimension $2n+2$).

\vskip .2cm

Let $y\in Y-\{y_1,\cdots, y_N\}$ be a non-critical point for $W$. Then the fiber $W^{-1}(W(y))$ is, near $y$, a smooth
complex manifold of dimension $n$. We define
\be\label{eq:symp-conn}
\Xi_y \,=\, \bigl(T_y (W^{-1}(W(y)))\bigr)^\perp_h\,\subset \, T_y Y
\ee
to be the orthogonal complement (with respect to $h$)   to its tangent space. Thus $\Xi_y$
is a $1$-dimensional complex submanifold in $T_y Y$ with the property that the differential
$d_y W: \Xi_y\to T_{W(y)} X$
is an isomorphism. We get a   $C^\infty$ complex rank $1$ subbundle $\Xi$ of
$T_{Y-\{y_1,\cdots, y_N\}}$.

\vskip .2cm

Let $X_0 = X\-A$ be the set of non-critical values. The map
\[
Y_0 \, := \, W^{-1}(X_0) \,=\, Y\,- \bigcup_{i=1}^N W^{-1}(w_i) \buildrel W\over\lra X_0
\]
is a $C^\infty$-locally trivial fibration, and $\Xi$ can be considered as a nonlinear connection in this fibration,
known as the {\em symplectic connection} \cite{macduff} \S 6.3. Denote this connection
$\nabla_W$.
Since we assumed $W$ to be proper,  any (piecewise) smooth path   $\gamma\subset X_0$,
joining two points $x_1$ and $x_2$, gives a well defined diffeomorphism
\[
T_\gamma: W^{-1}(x) \to W^{-1}(x').
\]
It is further known (see \cite{macduff} Lemma 6.3.5) that $T_\gamma$ is a symplectomorphism with respect to the symplectic
forms on $W^{-1}(x)$ and $W^{-1}(x')$ obtained as the restrictions of $\omega$.

The connection $\nabla_W$ is not flat but it can be seen as
a geometrization of the Gauss-Manin connection
on the cohomology of the $W^{-1}(x)$, which it induces. In particular, the ``monodromy'' of the
geometric vanishing cycle $\Delta_i\simeq S^n$ with respect to an appropriate closed contour $\gamma$
around $w_i$  (and appropriate K\"ahler metric)
can be seen as a self-diffeomorphism of $\Delta_i$ isotopic to the antipodal map (\cite{macduff} Ex. 6.3.7).

\paragraph{ $m_{ij}(\alpha)$ as the number of curved gradient trajectories.}
We keep the notations of \n \ref{par:symcon}

\begin{Defi}
Let $\gamma$ be a piecewise smooth path in $X_\genr$. By  \em{ $\gamma$-gradient trajectories}
we will mean  integral trajectories of $\Xi$ lifting $\gamma$.
\end{Defi}

Thus, if $\gamma$ begins at some $b\in X_\genr$, then there is a unique $\gamma$-gradient
trajectory passing through any $y\in W^{-1}(b)$.  The name ``$\gamma$-gradient trajectories is
explained as follows.

\begin{prop}\label{prop:gamma=line}
Let $X=\CC$, so $W$ is a ``holomorphic Morse function''.  Let $p\in\CC$ be a non-critical value.

\vskip .2cm

(a)  Let
$\gamma$ be the horizontal half-line $(p-\infty, p]$ going from $p$ to the left. We assume that this half-line
consists of non-critical vales. Then $\gamma$-gradient trajectories coincide (as unparametrized curves)  with the gradient trajectories
of the real  function $\Re(W): Y\to\RR$ with respect to the Riemannian metric $g$, beginning in $W^{-1}(p)$.

\vskip .2cm

(b) More generally, if $|\zeta|=1$ and $\gamma = (p-\zeta\cdot \infty, p]$ is the half-line of slope $\zeta$ emanating
from $p$, then $\gamma$-gradient trajectories coincide with   gradient trajectories of the real  function
$\Re(\zeta^{-1} W)$.
\end{prop}

\noindent {\sl Proof:} Part (b) follows from (a). Part (a) follows from the next well known fact.

\begin{prop}
The gradient flow of $\Re(W)$ with respect to $g$ coincides with the Hamiltonian flow of $\Im(W)$
with respect to $\omega$. In particular, gradient trajectories of $\Re(W)$ are mapped, by $W$, to
horizontal lines  $\{ \Im(W)=\on{const}\}$  in $\CC$. \qed
\end{prop}

Recall that $\Re(W)$ is a Morse function with the  critical points $y_i$ (same as for $W$),
 each of the same index
(number of negative eigenvalues) equal to
$n+1$. The same holds for $\Re(\zeta^{-1} W)$, $|\zeta| = 1$. Thus any critical point $y_i$ gives
rise to the {\em unstable manifold}, or {\em Lefschetz thimble} $\Th_\zeta(y_i)\simeq \RR^{n+1}$
which is the union of all descending gradient trajectories of  $\Re(\zeta^{-1} W)$ originating from $y_i$.
By Proposition \ref {prop:gamma=line}, the function $W:Y\to\CC$ projects $\Th_\zeta(y_i)$ to the
half-line $(w_i-\zeta\infty)$ with slope $\zeta$ originating at $w_i$.

\vskip .2cm

Let now $\alpha$ be a piecewise smooth path in $X$ joining $w_i$ and $w_j$ and not passing through
any other $w_k$, as in \S \ref{subsec:abs-PL} \ref{par:PL-id}
Let $\zeta$ be the slope of $\alpha$ at $w_i$, so that the half-line $(w_i-\zeta\infty)$ is
tangent to $\alpha$ at $w_i$. We can  compare $\alpha$-gradient trajectories with the $(w_i-\zeta\infty)$-gradient
trajectories, i.e., with gradient trajectories of $\Re(\zeta^{-1}W)$. This shows that $\gamma$-gradient trajectories
which approach $y_i$ over $w_i$, also form an $(n+1)$-dimensional manifold which we call the
{\em $\alpha$-Lefschetz thimble} originating at $y_i$ and denote $\Th_\alpha(y_i)$, see Fig. \ref{fig:cur-thim}.
By construction, $W$ projects $\Th_\alpha(y_i)$ to $\alpha$, so that the preimage of any intermediate
point $p\in\alpha$ is an $n$-sphere
\[
\Delta_{i,\alpha}(p) \,=\, W^{-1}(p) \cap \Th_\alpha(y_i) \,\simeq \, S^n.
\]
This is nothing but the translation to $p$ (along $\alpha$, with respect to the symplectic connection)
of the geometric vanishing cycle $\Delta_i$.

Starting from the other end, we have a similar Lefschetz thimble $\Th_\alpha(y_j)$ originating
at $y_j$ and the sphere $\Delta_{j,\alpha}(p)\subset
W^{-1}(p)$.

\vskip .2cm

We now make the following

\begin{ass}
The submanifolds $\Delta_{i,\alpha}(p), \Delta_{j,\alpha}(p)$ in $W^{-1}(p)$ intersect transversely.
In particular, being of middle dimension, they intersect in finitely many points.
\end{ass}

If we choose orientations of the $(n+1)$-dimensional manifolds $\Th_\alpha(y_i)$,
$\Th_\alpha(y_j)$,
then we get orientations of the $n$-spheres $\Delta_{i\alpha}(p)$ and $\Delta_{j,\alpha}(p)$.
In particular, such a choice identifies the $1$-dimensional vector spces $\Phi_{i,\alpha}$ and $\Phi_{j,\alpha}$
with $\k$, and \eqref{eq:mij-inter-I} implies:

\begin{prop}
$m_{ij}(\alpha)$ equals the signed number of intersection points  of  $\Delta_{i\alpha}(p)$ and $\Delta_{j,\alpha}(p)$
(with the signs given by the orientations).
\end{prop}

Since points of, say  $\Delta_{i, \alpha}(p)$ correspond to $\alpha$-gradient trajectories beginning at $y_i$ and
ending in $W^{-1}(p)$, we can reformulate the above as follows.

\begin{cor}\label{cor:mij=grad-traj}
$m_{ij}(\alpha)$ equals the signed number of $\alpha$-gradient trajectories joining $y_i$ and $y_j$.
\end{cor}

\bef[H]
\centering
\btp[scale=0.5]
\node (i) at (-9,0){};
\node (j) at (11,0){};
\node (p) at (0,0){};
\fill (i) circle (0.2);
\fill (j) circle (0.2);
\fill (p) circle (0.2);
\node (xi) at (-9,9){};
\node(xj) at (11,7){};
\fill (xi) circle (0.2);
\fill (xj) circle (0.2);

\draw[ line width=.4mm] (0,7) ellipse (1cm and 2cm);
\draw[ line width=.4mm] (0,10) ellipse (1cm and 2cm);

\draw[line width=0.3mm] plot[smooth, tension=1]  coordinates{ (-9,0) (-4.5,-0.5)  (0,0)  (5.5,0.5) (11,0)
};

\draw (-3,2) -- (-3,12) -- (3,16) -- (3,6) -- (-3,2);

\draw plot [smooth, tension=1]  coordinates {(0,12) (-4,10.5) (-9,9) (-4.5, 8) (0,8)};

\draw plot [smooth, tension=1]  coordinates {(0,9) (6,8.5) (11,7) (6,6.5)  (0,5)};

\node at (-11,0) {\huge$w_i$};
\node at (13,0) {\huge$w_j$};
\node at (0,-2){\huge$p$};
\node at (-11,9) {\huge$y_i$};
\node at (13,7) {\huge$y_j$};

\node at (2,17.5){\large$W^{-1}(p)$};

\node at (-7,12){\large$\Th_\gamma(y_i)$};
\node at (9,9.5){\large$\Th_\gamma(y_j)$};
\node at (1,13){\large$\Delta_{i,\alpha}(b)$};
\node at (-2.3,4.7){\large$\Delta_{j,\alpha}(b)$};
\node at (8,-0.5){\huge$\alpha$};

\etp
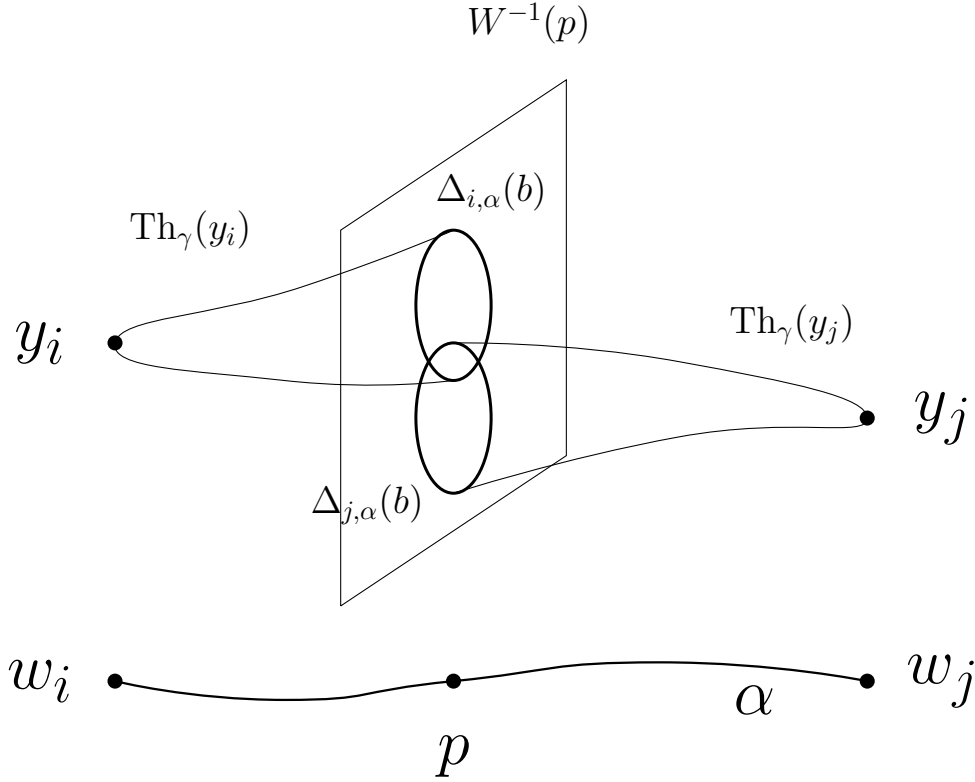
\caption{ Curved Lefschetz thimbles over a path $\alpha$. }
\label{fig:cur-thim}
\enf

\begin{ex}\label{ex:straight-lines}
Let $X=\CC$ and $\alpha=[w_i,w_j]$ be the straight line interval joining $w_i$ and $w_j$. We assume, as before, that
there are no other critical values on this interval. Let
\[
\zeta_{ij} = \frac{w_i-w_j}{|w_i-w_j|}
\]
be the slope of $\alpha=[w_i,w_j]$. Corollary \ref{cor:mij=grad-traj} implies that $m_{ij}(\alpha)$ is the number of
gradient trajectories of $\Re(\zeta_{ij}^{-1}W)$ joining $y_i$ and $y_j$.
We will return to this point of view in \S \ref{subsec:rec-approach}.

\vskip .2cm

Note that $y_i$ and $y_j$ are critical points of $\Re(\zeta_{ij}^{-1}W)$ of the same index. As well known,
for a {\em generic} Morse function $f$ there are no gradient trajectories joining two critical points of the same
index. In fact, functions that do have such trajectories, form a subvariety of real codimension $1$
in the space of all Morse functions, see \cite{Kap-Saito} and references therein.
So a {\em generic $1$-parameter family} of Morse functions will
contain finitely many such ``exceptional'' functions. This is what happens in our situation: the $1$-parameter
family is $\{\Re(\zeta^{-1}W)\}_{|\zeta|=1}$, and exceptional functions correspond to $\zeta=\zeta_{ij}$
for all pairs $(i,j)$.
\end{ex}


\subsection{Perverse sheaves with bilinear forms}\label{subsec:perv-bil}

\paragraph{Seifert duality for local systems on $S^1$.}
Let $\LS(S^1)$ be the category of local systems of finite-dimensional
$\k$-vector spaces on $S^1$. We introduce the   {\em Seifert duality} functor $\bs$
on this category by defining, for a local system $\Lc$,
\[
(\Lc^\bigstar)_\zeta \,=\, (\Lc_{-\zeta})^*.
\]
That is, the stalk of $\Lc^\bs$ at the point $\zeta\in S^1$ is  set to be the dual to the stalk of $\Lc$ at the opposite point $(-\zeta)$.

It is clear that $\bs$ is a perfect duality on $\LS(S^1)$.
For a local system $\Lc\in\LS(S^1)$, we will call a homomorphism
$q: \Lc \lra \Lc^\bs$ a {\em Seifert bilinear form} on $\Lc$, called {\em non-degenerate}, if it is an isomorphism.
A local system  equipped with   a nondegenerate symmetric (resp. antisymmetric) Seifert form, will be called
a {\em Seifert symmetric} resp. {\em Seifert antisymmetric} local system.

\begin{ex}
Let $\Sgn$ be the $1$-dimensional local system on $S^1$ with monodromy $(-1)$. Then $\Sgn$ carries a natural {\em antisymmetric}
scalar product $q:\Sgn\to\Sgn^\bs$.  To describe it,
let us  view $S^1$ as the unit circle $|z|=1$
in $\CC$, and
let $p: S^1\to S^1$ be the $2$-fold covering (with Galois group $\ZZ/2$)  taking $\zeta$ to $\zeta^2$. Then $\Sgn\subset p_*\ul\k_{S^1}$ consists of
$\ZZ/2$-anti-invariant sections.  Given $\zeta\in S^1$, the preimage $p^{-1}(\zeta)$ consists of two opposite points which generate
a $\ZZ$-sublattice $\ZZ\cdot p^{-1}(\zeta)\subset \CC$ of rank $1$. Anti-invariant functions $p^{-1}(\zeta)\to\k$ are the same as
group homomorphisms $\ZZ\cdot p^{-1}(\zeta)\to\k$, that is,
\[
\Sgn_\zeta \,=\,\Hom_\ZZ (\ZZ\cdot p^{-1}(\zeta), \k).
\]
Now consider the standard symplectic form $\omega = dx\wedge dy$ on $\CC$. Given two opposite points $\zeta, -\zeta\in S^1$,
any $\xi\in p^{-1}(\zeta)$ and $\eta\in p^{-1}(-\zeta)$ are orthogonal unit vectors and so $\omega(\xi,\eta)=\pm 1$. Therefore $\omega$
defines a nondegenerate pairing
\[
(\ZZ\cdot p^{-1}(\zeta)) \otimes_\ZZ  (\ZZ\cdot p^{-1}(-\zeta)) \lra \ZZ
\]
which by the above, induces a pairing $q_\zeta: \Sgn_\zeta\otimes\Sgn_{-\zeta}\to \k$. Since $\omega$ is anti-symmetric, $q_{-\zeta} = -q_\zeta^t$.
\end{ex}

Let now $(\Lc,q)$ be a local system on $S^1$  with a Seifert bilinear form
(not necessarily symmetric or antisymmetric).
Choose a point $\zeta\in S^1$ and let $V=\Lc_\zeta$ be the stalk at $\zeta$.
By the monodromy $M$ along the positive (counterclockwise) half-circle we can identify $V=\Lc_\zeta$ with $\Lc_{-\zeta}$. Then $q$ gives
a bilinear form
\be
B=B_q: V=\Lc_\zeta \buildrel q_\zeta \over\lra  \Lc_{-\zeta}^* \buildrel M^t\over\lra  \Lc_\zeta^*=V^*
\ee
which we call the {\em Seifert form} on $V$ associated to $q$. Compare with the definition of the Seifert form
for an isolated singularity in \cite{AGV} \S 2.3.

Note that if $q$ is symmetric (as a Seifert form on $\Lc$), then $B_q$ does not have to be
symmetric as a bilinear form on a vector space. In fact, we have the following fact that is
proved by direct comparison of definitions.

\begin{prop}

(a) The category of Seifert symmetric local systems on $S^1$ is equivalent to the category of pairs $(V, B: V\to V^*)$
where $V$ is a finite-dimensional $\k$-vector space (stalk of $\Lc$ at some point)
and $B$ is a (not necessarily symmetric) nondegenerate bilinear form on $V$. Under this equivalence,  the monodromy of $\Lc$ corresponds to the
{\em Serre operator} $S= B^{-1}\circ B^t$.

\vskip .2cm

(b) Tensoring with $\Sgn$ identifies the categories of Seifert symmetric and Seifert antisymmetric local systems
on $S^1$.
In particular, the category of Seifert antisymmetric local systems is again equivalent to the category of
pairs $(V, B; V\to V^*)$ as in (a). Under this equivalence, the monodromy of $\Lc$ corresponds to
the operator $-S=-B^{-1}\circ B^t$.
\qed
\end{prop}

\paragraph{Duality for vanishing cycles.}
Recall the functor of vanishing cycles $\bPhi: \Perv(D,0) \to \LS(S^1)$.

\begin{prop}\label{prop:phi-dual-seifert}
The functor $\Phi$ takes  (shifted) Verdier duality on perverse sheaves to  Seifert duality of local systems,
i.e., we have
natural isomorphisms
\[
\bPhi(\Fc^*) \,\simeq \, (\bPhi(\Fc))^\bs.
\]
\end{prop}

\noindent {\sl Proof:}  For simplicity we assume that $\Fc$ is given on the entire complex plane $\CC$, i.e., $\Fc\in \Perv(\CC,0)$.
Let $\zeta\in S^1$. We have the closed half-line $\RR_+\cdot \zeta\subset \CC$ and the embeddings
\[
\{0\} \buildrel k_\zeta\over\lra \RR_+\cdot\zeta \buildrel i_\zeta\over\lra \CC.
\]
By definition,
\[
\bPhi(\Fc)_\zeta \,=\, k_\zeta^* \,i_\zeta^!\,  (\Fc)[1].
\]
Since the  Verdier duality functor $\DD$ interchanges $i^!$ with $i^*$, we have
\[
\bPhi(\Fc^*)_\zeta \,\simeq \,( k_\zeta^! \,i_\zeta^* \, (\Fc) [1])^*,
\]
where the last $*$ in the RHS is the usual duality of vector spaces.
We now consider the commutative diagram
\[
\xymatrix{
\{0\}
\ar[d]_{k_{-\zeta}}
\ar[r]^{k_\zeta} & \RR_+\cdot \zeta
\ar[d]^{i_\zeta}
\\
\RR_+\cdot(-\zeta) \ar[r]_{\hskip 0.5cm i_{-\zeta}}&\CC.
}
\]
It is a Cartesian square of closed embeddings, and so gives a natural  morphism
\be\label{eq:iden-zeta}
k_{-\zeta}^! \, i_{-\zeta}^* \,\Fc \lra k_\zeta^* \, i_\zeta^! \Fc
\ee
and it is immediate that in our case ($\Fc\in \Perv(D,0)$) it is a quasi-isomorphism.  This gives an identification
$\bPhi(\Fc^*) \simeq (\bPhi(\Fc))^\bs$. \qed

\begin{rem}
The identification \eqref {eq:iden-zeta} would of course hold if we take, instead of $(-\zeta)$, any
$\zeta'\neq\zeta$ and so we have an identification
$\bPhi(\Fc^*)_{\zeta} \simeq \bPhi(\Fc)_{\zeta'}^*$. The role of the
exact opposite point $(-\zeta)$ is that it provides a consistent way to choose such a $\zeta'$ for
any $\zeta$. Any other consistent way would do just as well, and  would be homotopy equivalent to
$\zeta\mapsto (-\zeta)$.
\end{rem}

The identification induced by  \eqref {eq:iden-zeta} can be rewritten as a non-degenerate pairing (note that $(\zeta\cdot\RR_+)\cap (-\zeta\cdot\RR_+)=\{0\}$)
\be\label{eq:pair-beta2}
\bPhi(\Fc)_\zeta \otimes \bPhi(\Fc^*_{-\zeta} ) \,=\,\HH^1_{\zeta\cdot \RR_+}(\CC,\Fc) \otimes \HH^1_{-\zeta\cdot \RR_+}(\CC, \Fc^*)
\buildrel{\beta_\zeta} \over \lra
H^2_{\{0\}} (\CC, \ul \k) = \k.
\ee
We note that this pairing, having the form $H^1\otimes H^1\to H^2$, is anti-symmetric, i.e., $b_{-\zeta}$ becomes, after interchanging the arguments,
equal to the negative of $b_\zeta$. This implies the following fact.

\begin{prop}\label{prop:e-eps}
Let $\Fc\in\Perv(D,0)$ and $e: \Fc\to \Fc^{**}$ and $\eps: \bPhi(\Fc) \lra \bPhi(\Fc)^{\bs\bs}$  be the canonical isomorphisms. Then the identification
of Proposition \ref{prop:phi-dual-seifert} takes $e$ to $(-\eps)$. \qed
\end{prop}

\paragraph {Duality in the $(\Phi, \Psi)$-description.}

\begin{prop}\label{prop:phi-psi-dual}
Let $\Fc\in\Perv(D,0)$ be represented, via Proposition \ref{prop:ggm},
 by a  diagram
$
\xymatrix{
\Phi \ar@<.4ex>[r]^a&\Psi \ar@<.4ex>[l]^b
}
$.
Then $\Fc^*$ is represented by the diagram
$
\xymatrix{
\Phi^* \ar@<.4ex>[r]^{a'}&\Psi^* \ar@<.4ex>[l]^{b'}
},
$
where
\[
\begin{cases}
a'=b^t (1-a^tb^t)^{-1}
\\
b' = - a^t.
\end{cases}
\]
\end{prop}

\begin{rems}
(a)  Sometimes a simpler formula is asserted,  namely $a'=b^t, b'=a^t$.
However, this cannot be true for the following reason.
The monodromy of $\Fc$ on $D\- \{0\}$ (where it is a local system) is $T=1-ab$. Now, for any local system $\Lc$
on any space $X$ and
any path $\gamma$ joining points $x$ and $y$, the monodromies of $\Lc$ and $\Lc^*$ along $\gamma$ are related by
\[
M(\Lc^*)_x^y = \bigl( (M(\Lc)_x^y)^t\bigr)^{-1}: \Lc^*_x\lra \Lc^*_y, \quad M(\Lc)_x^y: \Lc_z\lra\Lc_y
\]
(this is the only possible expression that makes sense for any $x$ and $y$). So the monodromy of $\Fc^*$ should be $(T^t)^{-1} = (1-b^ta^t)^{-1}$,
while the ``simpler formula'' above would give $1-b^ta^t=T^t$. The formula in Proposition \ref{prop:phi-psi-dual}
does give the correct monodromy due to the ``Jacobson identity''
\be\label{eq:jacobson}
(1-uv)^{-1} \,=\, 1+ u (1-vu)^{-1} v.
\ee
See the discussion around Eq. (1.1.6) in  \cite{KS-shuffle} for the significance of this identity for  the $(\Phi, \Psi)$-description.

\vskip .2cm

 (b)
 Applying  Proposition \ref{prop:phi-psi-dual}  twice, we get   the $(\Phi, \Psi)$ description of $\Fc^{**}$
as the diagram  $
\xymatrix{
\Phi \ar@<.4ex>[r]^{a''}&\Psi \ar@<.4ex>[l]^{b''}
},
$
with
\[
\begin{cases}
b'' = (a')^t = -(1-ba)^{-1}b
\\
a'' = -(b')^t (1-(a')^t(b')^t)^{-1} = -a(1 + (1-ba)^{-1} ba)^{-1} = -a(1-ba)^{-1}\end{cases}
\]
This diagram is isomorphic to  $
\xymatrix{
\Phi \ar@<.4ex>[r]^{a}&\Psi \ar@<.4ex>[l]^{b}
},
$ via the isomorphisms
\[
e_\Psi = \Id: \Psi\lra\Psi, \quad e_\Phi =  -(1-ba)^{-1}: \Phi\lra\Phi.
\]
The equality $e_\Psi=\Id$ reflects the fact that dualization commutes with taking general fiber.
The appearance of the minus sign in the formula for $e_\Phi$ is a manifestation of Proposition \ref{prop:e-eps}.

\vskip .2cm

(c) The  formula of Proposition \ref{prop:phi-psi-dual} can be compared with the formula given
in \cite{bezr-kapr} (Prop. 4.5)
for the effect of  the  Fourier-Sato transform on the
$(\Phi, \Psi)$-description of
perverse sheaves on $(\CC,0)$. Both formulas are remindful of cluster transformations.
\end{rems}

\paragraph{Duality for transport maps.}
Given $\Fc\in\Perv (X,A)$ and a path  $\alpha$ joining $w_i$ with $ w_j$. Then the transport map
$m_{ij}^\Fc(\alpha): \bPhi_i(\Fc)_\alpha \to \bPhi_j(\Fc)_\alpha$
acts between the stalks of $\bPhi_i(\Fc)$ and $ \bPhi_j(\Fc)$ in the directions given by the
path $\alpha$, see
\eqref{eq:transport}. By Proposition \ref{prop:phi-psi-dual}, the dual map can be seen as acting
can be seen as a map
$m_{ij}^\Fc(\alpha)^t:  \bPhi_j(\Fc^*)_{-\alpha} \lra  \bPhi_i(\Fc^*)_{-\alpha}$
acting between the stalks in the directions, opposite to the directions given by $\alpha$.
On the other hand, we have the dual perverse sheaf $\Fc^*$ and its own transport map
$m_{ji}^{\Fc^*} (\alpha^{-1})$ associated to the  opposite path $\alpha^{-1}$, i.e., to
$\alpha$ run in the opposite direction.
Proposition  \ref{prop:phi-psi-dual} implies at once:

\begin{cor}
We have an equality
 \[
 m_{ji}^{\Fc^*}(\alpha^{-1}) \,=\,
-T_{+,i} \,  m_{ij}^\Fc(\alpha)^t
\, T_{+,j}^{-1},
\]
where $T_{i,+}$ is the counterclockwise  half-monodromy from the opposite to direction of $\alpha$ to
the direction of $\alpha$, and similarly for $T_{+,j}$. \qed
\end{cor}

\paragraph{Proof of Proposition  \ref{prop:phi-psi-dual}.}
For simplicity we assume that $\Fc$ is given on the entire complex plane $\CC$, i.e., $\Fc\in \Perv(\CC,0)$.

\vskip .2cm

The $(\Phi, \Psi)$-diagram describing $\Fc$ is the stalk at $\zeta=1$ of a local system of similar diagrams on the unit circle $S^1$:
\be\label{eq:LS-phi-psi}
\HH^1_{\zeta\cdot\RR_+}(\CC,\Fc) =    \xymatrix{
\Phi_\zeta(\Fc)  \ar@<.4ex>[r]^{a_{\zeta,\Fc}} &\Psi_\zeta(\Fc)  \ar@<.4ex>[l]^{b_{\zeta, \Fc}}
}\buildrel (I)\over\simeq \Fc_\zeta\buildrel (II)\over\simeq \HH^0(\CC\-(\zeta\cdot\RR_+), \Fc).
\ee
Here $\Fc_\zeta$ is the stalk of $\Fc$ at $\zeta\in S^1$, and the identifications (I) and (II) go as follows:
\begin{itemize}
\item[(I)]  (Local Poincar\'e duality)  As in Fig. \ref{fig:disk-1}, choose a small disk $V$ around $\zeta$ and write $V\-\zeta\cdot \RR_+ = V_+\sqcup V_-$
with $V_+$ following $\zeta\cdot \RR_+$ in the counterclockwise direction and $V_-$ preceding $\zeta\cdot\RR_+$. Then by definition
\[
\Psi_\zeta(\Fc) = \HH^1_{V\cap(\zeta\cdot\RR_+)}(V,\Fc) =   { \Fc(V_+)\oplus \Fc(V_-)\over \Fc(V)},
\]
and we send the class  $[s_+, s_-]$ of a pair of sections $s_\pm\in \Fc(V_\pm)$, to the image of  $s_+-s_-$ in $\Fc_\zeta$.

\item[(II)] Given an element $s$ of $\Fc_\zeta$, we extend it to a section in $\CC\-(\zeta\cdot\RR_+)$ by continuing counterclockwise from $\zeta$.
\end{itemize}
With these identification, the map $b_{\zeta,\Fc}$ is identified with the coboundary map
\[
\delta = \delta_\zeta: \HH^0(\CC\-(\zeta\cdot\RR_+), \Fc)\lra \HH^1_{\zeta.\cdot\RR_+}(\CC,\Fc).
\]
We focus in particular on $\zeta=\pm 1$ and write $\Phi_\pm(\Fc)$ for $\Phi_{\pm 1}(\Fc)$ etc. We use similar notations for $\Fc^*$.

By Proposition
\ref{prop:phi-dual-seifert},  $\Phi_-(\Fc^*)$  is identified with $\Phi_+(\Fc)^*$.   The non-degenerate pairing
that  Proposition
\ref{prop:phi-dual-seifert} gives, can be written as  the cup-product pairing (note that $\{0\}=\RR_+\cap \RR_-$)
\be\label{eq:pair-beta}
\Phi_+(\Fc)\otimes \Phi_-(\Fc^*) \,=\,\HH^1_{\RR_+}(\CC,\Fc) \otimes \HH^1_{\RR_-}(\CC, \Fc^*)
\buildrel\beta\over \lra
H^2_{\{0\}} (\CC, \ul \k) = \k.
\ee
At the same time, $\Psi_+(\Fc^*)$ is identified with $\Psi_+(\Fc)^*$.

\begin{lem}\label{lem:beta}
Let $\sigma\in \HH^0(\CC\-\RR_-, \Fc^*)$ and $x\in \HH^1_{\RR_+}(\CC, \Fc)$.
Then the pairing $\beta(x, \delta_-(\sigma))$ can be found as follows:
\begin{itemize}

\item Choose  a small disk   $V$  around $1$, and let  $\gamma(x)\in\HH^1_{V\cap\RR_+}(V, \Fc)$ be the restriction of $x$
(i.e.,  the result of applying the generalization map
of the sheaf $\ul\HH^1_{\RR_+}(\Fc)$ to the image of $x$ in $\ul\HH^1_{\RR_+}(\Fc)_0$).

\item Represent $\gamma(x)$ as the class $[s_+, s_-]$  of a pair  $s_\pm\in\HH^0(V_\pm,\Fc)$ (as in the identification (I) above).

\item Then $\beta(x, \delta_-(\sigma)) = (s_+, \sigma)-(s_-,\sigma)$, where $(s_\pm, \sigma)$ is the canonical pairing of sections of the local systems
$\Fc$ and $\Fc^*$ on $V_\pm$ .
\end{itemize}
\end{lem}

\noindent{\sl Proof:} Consider the diagram
\[
\xymatrix{
\HH^1_{\RR_+}(\CC, \Fc) \otimes \HH^0(\CC\-\RR_-, \Fc^*)
\ar[dd]_{  \Id\otimes\delta}
\ar[rr]^{\alpha = \cup} && \HH^1_{(\CC\-\RR_-)\cap \RR_+} (\CC\-\RR_-, \ul\k)
\ar@{-->}[llddd]
\ar[d]^{=}
\\
&& \HH^0(\CC\-\RR_-, \ul\HH^1_{\RR_+}(\ul\k))
\ar[d]^{\delta}_{\simeq}
\\
\HH^1_{\RR_+}(\CC, \Fc)\otimes\HH^1_{\RR_-}(\CC, \Fc^*)
\ar[d]_{\beta = \cup}
&& \HH^1_{\RR_-} (\ul \HH^1_{\RR_+}(\ul \k))
\ar[dll]_{=}
\\
\HH^2_{\RR_+\cap \RR_-}(\CC, \ul\k) = \k &&
}
\]
whose outer rim is commutative by naturality of the coboundary maps $\delta$ (Note that both $\alpha$ and $\beta$ are given by
the cup-product.)  Now, after the identification of the target of
$\alpha$ with $\k$ given by the  maps on the right side (the dotted arrow), the map $\alpha$ is immediately found to have the form
claimed in the lemma.  \qed

\vskip .2cm

We now describe the standard $(\Phi, \Psi)$-diagram (i.e., the stalk at $\zeta=1$ of the local system \eqref{eq:LS-phi-psi} of diagrams)
but for $\Fc^*$ instead of $\Fc$. We use the  isomorphisms
\[
\begin{gathered}
\Phi_+(\Fc^*) \buildrel\beta\over   \lra \Phi_-(\Fc)^* \buildrel (T_+^t)^{-1}\over\lra \Phi_+(\Fc)^*,
\\
\Psi_+(\Fc^*) \buildrel (I_{\Fc^*}) \over\lra \Fc_1^* \buildrel (I_\Fc)^t\over\lra \Psi_+(\Fc)^*,
\end{gathered}
\]
where:
\begin{itemize}
\item $T_+: \Phi_+(\Fc)\to\Phi_-(\Fc)$ is the anti-clockwise half-monodromy (along the upper semi-circle);

\item $\beta=\beta_{\Fc^*}$ is the pairing \eqref{eq:pair-beta} for $\Fc^*$;

\item $(I_\Fc)$ and $(I_{\Fc^*})$ are the identifications in \eqref{eq:LS-phi-psi} for $\Fc$ and $\Fc^*$.
\end{itemize}
After this, we define $a'$ and $b'$  by commutativity of the diagram
\[
\xymatrix{
\Phi_+(\Fc^*)
\ar[d]_{(T_+^t)^{-1}\circ\beta}
\ar@<.4ex>[r]^{a_{+,\Fc^*}} &\Psi_+(\Fc^*)  \ar@<.4ex>[l]^{b_{+, \Fc^*}}
\ar[d]^{ (I_\Fc)^t\circ (I_{\Fc^*})}
\\
\Phi _+(\Fc)^*\ar@<.4ex>[r]^{a'}&\Psi_+(\Fc)^*, \ar@<.4ex>[l]^{b'}
}
\]
that is, by
\[
\begin{gathered}
a'\,=\, (I_\Fc)^* \circ (I_{\Fc^*}) \circ a_{+, \Fc^*} \circ\beta^{-1} \circ T_+^t,
\\
b'\,=\, (T_+^t)^{-1} \circ\beta\circ b_{+,\Fc^*}\circ (I_\Fc^*)^{-1} \circ((I_\Fc)^t)^{-1}.
\end{gathered}
\]
To find $a'$, let $\phi\in \Phi_+(\Fc)^*$ and $\psi\in \Psi_+(\Fc)=\Fc_1$ and find $(a'\phi, \psi)$.

\vskip .2cm

Let us denote  the counterclockwise  half-monodromies from
$+1$ to $(-1)$ in all the local systems involved,   by  the same symbol $T_+$, the counterclockwise half-monodromies from $(-1)$ to $1$
by $T_-$, and the full monodromies from $1$ to $1$ by $T= T_-\circ T_+$.

\vskip .2cm

In this notation set
\[
\phi'=T_+(\phi)\in\Phi_-(\Fc)^*, \quad f=\beta^{-1}(\phi')\in\Phi_+(\Fc^*).
\]
By Lemma \ref{lem:beta} applied to $\Fc^*$ instead of $\Fc$,
\[
(a_{+,\Fc^*}(f),\psi) \,=\, \beta(f, \delta_- s_-(\psi)),
\]
where $s_-(\psi)$ is the section of $\Fc$ on $\CC\- \RR_-$ equal to $\psi$ at $1$ and $\beta$ is viewed as a pairing.
Viewing $\beta$ as an isomorphism $\Phi_+(\Fc^*)\to\Phi_-(\Fc)^*$, we can rewrite the RHS of the previous equality as
\[
(\beta(f), \delta_- s_-(\psi)) \,=\, (T_+(\phi), \delta_- s_-(\psi)).
\]
Now, recall that $b_\Fc(\psi)\in\Phi_+(\Fc)$ is obtained as $\delta_+ s'(\psi)$, where $s'(\psi)$ is the section of $\Fc$
on $\CC\-\RR_+$ obtained by extending $\psi$ counterclockwose from $\RR_+$.

Note that $s_-(\psi) = T_-^{-1}(s'(\psi))$.  So we write the RHS of the previous displayed equality as
\[
\begin{gathered}
(T_+(\phi), \delta_- T_-^{-1} s'(\psi)) \,=\, (T_+(\phi), T_-^{-1}\delta_+s'(\psi))\,=
\\
=\/ (T_-T_+(\phi), \delta_+s'(\psi)) \,=\, (T(\phi),  b_\Fc(\psi)).
\end{gathered}
\]
This means that $a'=b^t\circ T$, where $T$ is the monodromy of $\bPhi(\Fc)^*$, i.e.,
\[
T= ((1-ba)^t)^{-1} =   (1-a^tb^t)^{-1}
\]
as claimed.

\vskip .2cm

Let us now find $b'$. For this, fix $x\in \Psi_+(\Fc)^*=\Fc_1^*$ and $y\in\Phi_+(\Fc)$ and find
\[
(b'(x),y) \,=\,  \bigl( (T_+^t)^{-1}(\beta(b_{+,\Fc^*}(x))), y\bigr)\,=\, \bigl(\beta(b_{+,\Fc^*}(x)), T_+(y)\bigr).
\]
Now, $b_{+,\Fc^*}(x)=\delta_- s'(x)$, where $s'(x)$ is the section of $\Fc^*$ in $\CC\-\RR_+$ obtained by continuing $x$ counterclockwise from
$\RR_+$.  So by Lemma \ref{lem:beta} in which the roles of $\RR_+$ and $\RR_-$ are interchanged,  the RHS of the previous displayed equality
is equal to  {\em minus} the pairing of the values of $s'(x)$ and $a_{-, \Fc}(T_+(y))$ at the point $(-1)$.
The minus sign comes from the fact that interchanging  roles of $\RR_+$ and $\RR_-$ amounts to change of
orientation or, put differently, it comes from the anticommutativity of the product $H^1\otimes H^1\to H^2$.
But the value of $s'(x)$ at $(-1)$ is equal to $T_+(x)$,
we we get
\[
(T_+(x), a_{-,\Fc}(T_+(y))) \,=\, (T_+(x), T_+(a_{+,\Fc}(y)) \,=\, (x, a_{+,\Fc}(y)).
\]
This means that $b'=-a^t$, as claimed.

\paragraph{Perverse sheaves with bilinear forms. The baby spherical functor package. }
Denote by $\Perv(D,0)^\BF$ the category formed by pairs
$(\Fc, q)$, where $\Fc\in\Perv(D,0)$ is a perverse sheaf and  $q$ is a nondegenerate bilinear form on $\Fc$, i.e., an
isomorphism $q: \Fc\to \Fc^*$.
No  symmetry conditions on $q$ are assumed.
Morphisms are morphisms of perverse sheaves compatible with bilinear forms. We will describe $\Perv(D,0)^\BF$
in terms of
$(\Phi, \Psi)$-diagrams with additional structure.
The answer is completely analogous to  (is   a ``decategorification of'') the concept of a spherical functor
\cite{AL}.
We start with some generalities.

\vskip .2cm

Let $V$ be a finite-dimensional $\k$-vector space equipped with a non-degenerate bilinear form $B$. We view $B$ alternately as
an isomorphism $B: V\to V^*$ and as a pairing
\[
(v, v') \,\mapsto \, B(v, v') \, :=\,\<v, B(v')\>,
\]
where on the right hand side we have the canonical pairing of $V^*$ and $V$.

\vskip .2cm

The change of arguments in the pairing $B(v, v')$ corresponds to  passing to the dual map
$B^t: V=(V^*)^*\to V^*$.
The form $B$ gives rise to the {\em Serre operator} $S_B= B^{-1}\circ B^t: V\to V$,
so that we have
\[
B(v_1, v_2) = B(v_2, S_B(v_1)).
\]

Further,  let $V, W$ be two finite-dimensional $\k$-vector spaces equipped with non-degenerate bilinear forms $B_V, B_W$ and $f: V\to W$ be a linear map.
Then the {\em right} and {\em left adjoints} to $f$ are the maps $  f^*,\, ^*f: W\to V$ defined by
\[
B_W(f(v), w) = B_V(v, f^*(w)), \quad B_V(^*f(w), v) = B_W(v, f(w)).
\]
Viewing $B_V, B_W$     as isomorphisms $B_V: V\to V^*$ and $B_W: W\to W^*$, we have
\be\label{eq:adj-invar}
f^* = B_{V}^{-1} f^t B_{W},
\quad
^*f = (B_{V}^t)^{-1} f^t B_{W}^t
\ee
where $f^t: W^*\to V^*$ is the map dual to $f$.
One can further define iterated   right and left  adjoints such as $f^{**} = (f^*)^*,\,\, ^{**}f$ etc.

Specializing to $(W,B_W)=(V,B_V)$, let us note the following.

\begin{lem}\label{lem:isometry}
Let $T: V\to V$ be an endomorphism such that $T^*\circ T=1$. Then $^*T = T^* = T^{-1}$ and $T$ is an isometry, i.e., $B_V(v, v')= B_V(Tv, Tv')$. \qed
\end{lem}

\begin{Defi}\label{def:sph-map-baby}
Let $(\Phi, B_\Phi)$ and $(\Psi, B_\Psi)$ be finite-dimensional  $\k$-vector spaces equipped with non-degenerate bilinear forms. A linear map
$a: \Phi\to\Psi$ is called {\em spherical}, if the following conditions are satisfied:
\begin{itemize}
\item[(S1)]  The map $T_\Phi = 1-(a^*)a: \Phi\to\Phi$ is invertible.

\item[(S2)] We have $a^* + (^*a) - (^*a)aa^*=0$.
\end{itemize}
\end{Defi}

Compare with \cite{AL}. The condition (S2) is analogous to Kuznetsov's definition of a spherical functor, see
\cite{kuznetsov}, Def. 2.12 . Futher, the next proposition is a de-categorified version of the full spherical functor
package, see  \cite{AL}, Th. 1.1.

\begin{prop}\label{prop:baby-sph-pack}
If $a$ is spherical, then:

\begin{itemize}
\item[(a)] The map  $T_\Psi = 1-aa^*: \Psi\to\Psi$ is invertible.

\item[(b)] Both $T_\Phi$ and $T_\Psi$ are isometries (for $B_\Phi$ and $B_\Psi$ respectively).

\item[(c)] The left multiplication by $T_\Phi^{-1}$ takes the right adjoint to minus the left adjoint:
\[
(^*a) = -(1-(a^*)a)^{-1}(a^*).
\]

\item[(d)] We have $(^*a)aa^* = (a^*) a (^*a)$.
\end{itemize}

\end{prop}

Note that the requirement that $T_{\Phi}$ and $T_{\Psi}$  are isometries is strictly weaker than the requirement that $a$ is a spherical map,
while in categorical version it is enough to require that $T_{\Phi}$ and $T_{\Psi}$ are autoequivalences.

\vskip .2cm

\noindent{\sl Proof:}   Multiplying (S2) by $a$ on the left and adding $1$ we get
\[
(^* T_\Psi) \circ T_\Psi = (1-a(^*a))(1-aa^*)  = 1-a(a^*+(^*a) - (^*a)aa^*) = 1,
\]
so $T_\Psi$ is invertible, giving (a) and is an isometry by Lemma \ref{lem:isometry}.
Similarly, multiplying (S2) by $a$ on the right and adding $1$, we get
\[
(^* T_\Phi) \circ T_\Phi =   (1-(^*a)a)(1-a^*a) = 1 - (a^* + (^*a) + (^*a)aa^*)a = 1
\]
so $T_\Phi$ is an isometry, giving (b). Next,  (c) is reduced to (S2) by substituting $T_{\Phi}^{-1} = 1-(^*a)a$.
Finally,  let us rewrite the equality claimed in (d) by using (c). It will read:
\[
(^*a) a (1-(a^*) a) (^*a)   \, \buildrel ?\over = \,  (1-(a^*)a)(^*a) a(^*a).
\]
But this is obviously true since $(^*a) a $ commutes with $(1-(^*a) a)^{-1} = (1-a^* a)$. \qed

\begin{prop}
The category $\Perv(D,0)^\BF$ is equivalent to the category formed by data consisting of:
\begin{itemize}
\item Finite-dimensional $\k$-vector spaces $\Phi,\Psi$ equipped with non-degenerate bilinear forms $B_\Phi, B_\Psi$;
\item A spherical map $a: \Phi\to\Psi$.
\end{itemize}
\end{prop}

\noindent{\sl Proof:}
Suppose  $\Fc \in \Perv$ is represented by a diagram
$
\xymatrix{
\Phi \ar@<.4ex>[r]^a&\Psi \ar@<.4ex>[l]^b
}
$.
Note that the diagram representing, by Proposition \ref{prop:phi-psi-dual}, the dual  $\Fc^*$ is naturally isomorphic to
    $\xymatrix{
\Phi^* \ar@<.4ex>[r]^{a^{\dagger}}&\Psi^* \ar@<.4ex>[l]^{b^{\dagger}} }$
with
\[
\begin{cases}
a^{\dagger} = -b^t (1-a^tb^t)^{-1}
\\
b^\dagger =  a^t.
\end{cases}
\]
The isomorphism is given by $+1$ on $\Phi^*$ and $(-1)$ on $\Psi^*$.  Let us use this diagram for $\Fc^*$.
 The bilinear form is then  a map $q: \Fc \rightarrow \Fc^*$, which we represent as  the following commutative diagram:
\[
\xymatrix{
\Phi
\ar[d]_{B_{\Phi}}
\ar@<.4ex>[r]^{a} &\Psi \ar@<.4ex>[l]^{b}
\ar[d]^{B_{\Psi}}
\\
\Phi^*\ar@<.4ex>[r]^{a^\dagger}&\Psi^*, \ar@<.4ex>[l]^{b\dagger}
}
\]
Using \eqref{eq:adj-invar}, we write
\[
\begin{cases}
a^* = B_{\Phi}^{-1} a^t B_{\Psi} = B_{\Phi}^{-1} b' B_{\Psi} = b
\\
(^*a) = (B_{\Psi} a B_{\Phi}^{-1})^t = (a')^t = -(1-ba)^{-1} b =  -(1-(a^*)a)^{-1}a^* = -T_{\Phi}^{-1} (a^*),
\end{cases}
\]
whence the statement. \qed


\paragraph{Calabi-Yau perverse sheaves.}

We now focus on symmetry conditions for $q:\Fc\to\Fc^*$.

\begin{Defi} We call an {\em even} (resp. {\em odd}) {\em Calabi-Yau perverse sheaf} a perverse sheaf
$\Fc$ with a symmetric (resp. antisymmetric) nondegenerate bilinear form $q: \Fc\to\Fc^*$.
\end{Defi}

Let $\Perv(D,0)^{\CY+}$ and $\Perv(D,0)^{\CY-}$ be the categories formed by even and odd Calabi-Yau perverse sheaves on $(D,0)$.
They are full subcategories in $\Perv(D,0)^\BF$.

\begin{prop}\label{prop:CY-decat}
\begin{itemize}
\item[(a)] A spherical map $a: \Phi\to\Psi$ corresponds to an even Calabi-Yau perverse sheaf, if and only if:
\begin{itemize}
\item[(a1)] $B_\Psi$ is symmetric.

\item[(a2)] $T_\Phi= 1-(a^*)a$ is equal to $(-S_{B_\Phi})$, the negative of the Serre operator for $B_\Phi$.
\end{itemize}
\item[(b)]
A spherical map $a: \Phi\to\Psi$ corresponds to an odd Calabi-Yau perverse sheaf, if and only if:
\begin{itemize}
\item[(b1)] $B_\Psi$ is anti-symmetric.

\item[b2)]  $T_\Phi= 1-(a^*)a$ is equal to $S_{B_\Phi}$,   the Serre operator for $B_\Phi$.
\end{itemize}
\end{itemize}
\end{prop}

The proof is an immediate consequence of the next corollary, which, in turn, follows at once from Proposition \ref{prop:e-eps}.

\begin{cor}
If $\Fc$ is  even (resp. odd) Calabi-Yau perverse sheaf  on $(D,0))$, then $\bPhi(\Fc)$ is
a Seifert anti-symmetric (resp.  symmetric)  local system. \qed

\end{cor}

\begin{ex}
Let $W: X\to \Sigma$ be a generalized Lefschetz pencil (Definition \ref{def:GLP})
of relative dimension $n$ and $\Lc_W\in\Perv(\Sigma, A)$ be the corresponding
Lefschetz perverse sheaf. Poincar\'e duality, formulated in the language of perverse
sheaves, implies that $\Lc_W$ is  Calabi-Yau.

More precisely, since the constant sheaf $\ul\k_X$ is Verdier self-dual
(up to a shift by $n+1$), the complex $RW_*(\ul\k_X)$ is Verdier self-dual up to a shift
by $n$ and so  $\Lc_W = \ul H^n_\perv (RW_*(\ul\k_X))$
is self-dual,
i.e., equipped with a non-degenerate bilinear form $q: \Lc_W\to\Lc_W^*$ which is
symmetric for even $n$ and anti-symmetric for odd $n$.

In particular, let us restrict $\Lc_W$ to a small disk around a singular value $w_i=W(x_i)$.
Corollary \ref{prop:CY-decat} gives then the following classical result: {\em the monodromy of
an isolated singularity transposes (up to a sign) the Seifert form}. See \cite{AGV} Th. 2.5,
where the operator $\on{Var}$ differs from the Seifert form by a sign.

\end{ex}


\subsection{Quotient by (locally) constant sheaves}\label{subsec:quot-LC}

\paragraph{The category $\ol\Perv(X, A)$.}
Let $\Ac$ be an abelian category. We recall  that a full subcategory $\Bc\subset\Ac$ is called a {\em Serre subcategory},
if  it is  closed under subobjects, quotient objects and extensions. In particular,  $\Bc$ is abelian itself
and a strictly full subcategory in $\Ac$ (i.e., any object of $\Ac$ isomorphic to an object in $\Bc$,
is in $\Bc$).
In this case we have  a new abelian category $\Ac/\Bc$, called the {\em quotient category} of $\Ac$ by $\Bc$.
By definition, $\Ob(\Ac/\Bc)=\Ob(\Ac)$, and
\[
\Hom_{\Ac/\Bc}(A,A') = \Hom_\Ac(A,A')/I_{A,A'},
\]
where $I_{A,A'}$ consists of all morphisms $A\to A'$ which factor through an object of $B$. See
\cite  {gelfand-MV} and references therein for more details. An alternative definition of $\Ac/\Bc$
is as follows.

\begin{prop}
Let $Q_\Bc\subset\Mor(\Ac)$ consist of morphisms $f: A\to A'$ such that both $\Ker(f)$ and $\Coker(f)$ lie in $\Bc$.
Then $\Ac/\Bc$ is identified with $\Ac[Q_\Bc^{-1}]$, the localization of $\Ac$ with respect to $Q_\Bc$. \qed
\end{prop}

We now specialize to $\Ac=\Perv(X, A)$.
Let $\LS(X) = \Perv(X, \emptyset) \subset \Perv(X, A)$ be the full subcategory formed by local systems on
$X$.
Following \cite  {gelfand-MV}, we note that $\Bc=\LS(X)$ is a Serre subcategory in the abelian category $\Perv(X, A)$,
The corresponding quotient category will be denoted
\be
\ol\Perv(X, A) \,=\, \Perv(X, A) \bigl/ \LS(X).
\ee
Objects of $\ol\Perv(X,A)$ will be called {\em localized perverse sheaves}.

\vskip .2cm

Denoting, as before, $A=\{w_1,\cdots, w_N\}$, we note that the vanishing cycle functors $\bPhi_i$  descend to the quotient category,
giving the functors that we still denote
\[
\bPhi_i: \ol\Perv(X, A) \lra \LS(S^1_{w_i}).
\]
Further, for any path $\alpha$ joining $w_i$ and $w_j$ and not passing through any other $w_k$, we have the
natural transformation of these new functors
\[
m_{ij}(\alpha): \bPhi_{i,\alpha} \lra \bPhi_{j,\alpha},
\]
where $\bPhi_{i,\alpha}$ takes, as in \eqref{eq:transport}, a perverse sheaf $\Fc$ to the stalk of $\bPhi_i(\Fc)$ at the direction of
$\alpha$ at $w_i$. These transformations satisfy the formal Picard-Lefschetz identities in Proposition \ref{prop:PL-form}.

\paragraph {The spider description of $\ol\Perv(D, A)$.} We now specialize to $X=D$ being a disk.
In this case $\LS(D)$
consists of constant sheaves. The spider description of $\Perv(D,A)$
(Proposition \ref{prop:GMV1})  implies, as pointed out in \cite{gelfand-MV},
a description of $\ol\Perv(D,A)$ as well. To formulate it, let $K = \{\gamma_i,\cdots, \gamma_N\}$ be a $v$-spider for $(D,A)$,
see again Fig. \ref{fig:disk-2}.  For any $i,j=1,\cdots, N$
let $\alpha_{ij}$ be the path obtained by going first from $w_i$ to $v$ via $\gamma_i$ and then back to $w_j$ via
the reverse of $\gamma_j$. Let also $T_i(\Fc): \Phi_{i,\gamma_i}(\Fc) \to  \Phi_{i,\gamma_i}(\Fc)$ be the monodromy
of the local system $\bPhi_i(\Fc)$.

\begin{Defi}\label{def:cat-M_A}

Let $\Mc_N$ be the category whose objects are   diagrams consisting of:

\begin{itemize}
\item [(0)] Vector spaces $\Phi_i$, $i=1,\cdots, N$.

\item[(1)] Linear operators $m_{ij}: \Phi_i\to \Phi_j$ given for all $i,j$ and such that $ \Id_{\Phi_i}-m_{ii}$ is invertible.
\end{itemize}
\end{Defi}


\begin{prop}[\cite{gelfand-MV}]\label{prop:GMV2}
Let $K$ be a $v$-spider for $(D,A)$. The functor $\Xi_K:  \ol\Perv(D,A)\to \Mc_N$ taking
$\Fc$ to $\Phi_i = \Phi_{i, \gamma_i}(\Fc)$ and
\[
m_{ij} \,=\begin{cases}
m_{ij}(\alpha_{ij}), &\text { if } i\neq j;
\\
\Id-T_i(\Fc), & \text { if } i=j,
\end{cases}
\]
is an equivalence of categories. In other words, we have the following commutative diagram of functors
\[
\xymatrix{
\Perv(D,A)\ar[d]_{\lambda}
\ar[r]^{\hskip .6cm \Theta_K}_{\hskip .6cm \sim} & \Qc_N\ar[d]^{\mu}
\\
\ol\Perv(D,A) \ar[r]_{\hskip .6cm \Xi_K}^{\hskip .6cm \sim} &\Mc_N
}
\]
with horizontal arrows being equivalences, $\lambda$ being the canonical localization functor and $\mu$
given
explicitly by
\[
(\Psi, \Phi_i, a_i, b_i) \,\mapsto ( \Phi_i, m_{ij} = b_j a_i).
\]
\qed
\end{prop}

\begin{rem}
In fact, $\ol\Perv(D,A)$ can be canonically embedded back into $\Perv(D,A)$ as the full subcategory
of perverse sheaves $\Fc$ s.t. $\HH^\bullet(D,\Fc)=0$, see \cite{gelfand-MV},  Prop. 3.1
and \cite{KKP}, \S 2.3.2. At the level of quivers, this corresponds to the full embedding
$t: \Mc_N \hra \Qc_N$ given by
\[
t(\Phi_i, m_{ij}) \,=\, \biggl( \Psi = \bigoplus_{j=1}^N \Phi_j, \,\, \Phi_i = \Phi_i, \,\, a_i: \Phi_i \buildrel \bigoplus_j m_{ij}\over\lra \bigoplus _{j=1}^N \Phi_j, \,\,\, b_i: \bigoplus_{j=1}^N \Phi_j \buildrel \on{proj.}_i\over\lra \Phi_i
\biggr),
\]
and which is right adjoint to $\mu$, see \cite{gelfand-MV}, proof of Prop. 2.3.
Given a diagram $(\Phi_i, m_{ij})\in\Mc_N$, we will refer to the perverse sheaf $\Fc\in \Perv(D,A)$
corresponding to the quiver $t(\Phi, m_{ij})\in\Qc_N$
as the {\em Gelfand-MacPherson-Vilonen representative}
of $(\Phi_i, m_{ij})$.

\vskip .2cm

 Let us also note the relation of the GMV description  with the Katzarkov-Kontsevich-Pantev description
of nc Betti structures, see Sect. 2.3.2 of  \cite{KKP}. 
\end{rem}

\paragraph{  Dependence  on the spider.}
For later use let us recall, following \cite{gelfand-MV}, the effect of the change of spider $K$ on
the equivalence $\Xi_K: \ol\Perv(D,A)\to\Mc_N$. We use the action of the braid group $\Br_N$
on the set $\Kc(D,v,A)$ of spiders, see Fig. \ref{fig:br-spid}.  Proposition \ref{prop:ch-spid-1}
is complemented by the following fact which is its direct consequence
(it is also a consequence of the Picard-Lefschetz identities).

\begin{prop}[\cite{gelfand-MV}, Prop. 2.4]\label{prop:GMV-mut}
(a) Replacing $K$ by $\tau_i(K)$ changes a diagram $(\Phi_j, m_{\nu,j})$ to the  diagram
$(\Phi'_j, m'_{\nu, j})$, where;
\[
 (\Phi'_1,\cdots, \Phi'_N) \,=\, (\Phi_1, \cdots, \Phi_{i-1}, \Phi_{i+1}, \Phi_i, \Phi_{i+2},
\cdots, \Phi_N)\}
\]
\[
 m'_{\nu, j} \,=\, \begin{cases}
 m_{\nu, j}, &\text{ if }\,  i \notin \{j, j+1, \nu, \nu+1\};
 \\
 m_{\nu,i},&\text{ if }\,  j=i+1,\, \nu \notin \{i, i+1\};
 \\
 m_{i,j}, &\text{ if }\,  \nu=i+1, \, j\notin \{i, i+1\};
 \\
 m_{\nu, i+1} + m_{i, i+1} (1-m_{ii})^{-1} m_{\nu,i}, & \text{ if }\,  j=i, \, \nu\notin\{i, i+1\};
 \\
 m_{i+1,j} - m_{ij}m_{i+1,i},& \text{ if }\,  \nu=i, \, j\notin \{i, i+1\};
 \\
 (1-m_{i,i})m_{i+1,i}, & \text{ if } \, \nu-i,\, j=i+1;
 \\
 m_{i, i+1}(1-m_{i,i}), & \text{ if } \, \nu=i+1,\, j=i;
 \\
 m_{i+1, i+1}, & \text{ if } \, \nu=j=i;
 \\
 m_{i,i}, & \text{ if } \, \nu=j=i+1.
 \end{cases}
 \]

(b) Put differently, the formulas in (a) define an action of $\Br_N$ on the category $\Mc_N$.
For any spider $K\in\Kc(D,v,A)$ (and the corresponding numbering of $A$ by slopes of the $\gamma_i$ at $v$),
the equivalence $\Xi_K$
takes the geometric action of $\Br_N=\Gamma^+_{D,A}$ on $\ol\Perv(D,A)$ (coming from
Proposition \ref {prop:Gamma-acts}) to the above
action   of $\Br_N$ on $\Mc_N$.  \qed
\end{prop}


\section {Rectilinear point of view and the Fourier transform}\label{sec:recti-FT}

\subsection {The rectilinear approach}\label{subsec:rec-approach}

\paragraph{Notation and terminology.}
Let $A=\{w_1, \cdots, w_N\}\subset\CC$. The following notation and terminology,
pertaining to convex geometry of $\CC=\RR^2$,  will be
used throughout  the rest of the paper.

\vskip .2cm

For any points $x,y\in\CC$ we denote by $[x,y]$ the straight line interval between $x$ and $y$.
We denote by $Q=\Conv(A)$ the convex hull of $A$, a convex polygon in $\CC$.

\begin{Defi}\label{def:conv-pos}
Let $\zeta\in \CC$, $|\zeta|=1$.
We say that $A$ is:
\begin{itemize}

\item[(a)]   {\em in convex position}, if each element of $A$ is a vertex of $Q$.

\item[(b)]
{\em in linearly general position},   if no three elements of $A$ lie on a real line in $\CC$.

\item[(c)] in {\em strong linearly general position}, if no two intervals $[w_i, w_j]$ and $[w_k, w_l]$
are parallel.

\item[(d)]  {\em in linearly general position including $\zeta$-infinity},  if $A$ is in linearly general position and all the $\Im(\zeta^{-1}w_i)$
are distinct. When $\zeta=1$, we simply speak about {\em linear position including infinity}.

\end{itemize}

\end{Defi}

Condition (d) can be understood as $A\cup \{R\zeta\}$ being in linearly general position,
where $R \gg 0$
is a very large positive number.

\vskip .2cm

The {\em rectilinear approach} to the study of perverse sheaves $\Fc\in\Perv(\CC,A)$, which we develop in
this section,  consists in using the
tools from convex geometry, in particular, the above concepts, to analyze $\Fc$.

\paragraph{ Motivation: Fourier transform brings convex geometry.}
The rectilinear approach may seem unnatural, since
the concept of a perverse sheaf is purely topological.  However, convex geometry appears naturally,
{\em  if we want to do the Fourier transform}. In other words, assume that our
perverse sheaves are formed by some analytic  functions (solutions of linear differential equations)
and we are interested in Fourier transforms of these functions.

\vskip .2cm

More precisely, let $w$ be the coordinate in $\CC$, and  write $\CC=\CC_w$ to signify
that we deal with  the plane with this coordinate.
Let $z$ be the ``dual variable'' to $w$, so that the complex plane $\CC_z$ with coordinate $z$ is the dual
$\CC$- vector space to $\CC_w$ via the pairing $(z,w)=z\cdot w$.
The naive formula for the Fourier (Laplace) transform
\be\label{eq:FT}
\wc f(z) \,=\,\int_\gamma f(w) e^{-zw} dw
\ee
needs, in the context of complex analytic  functions,  several well known
qualifications.

\vskip .2cm

First,  $f$ being possibly multi-valued,
the countour $\gamma$ must be a $1$-cycle with coefficients in the local
system of determinations of $f$.

\vskip .2cm

Further, $\gamma$ must be an unbounded (Borel-Moore) cycle
such that the integrand decays appropriately at its infinity.
The case relevant to our situation is when (each branch of) $f$ has at most polynomial growth: this means that the
differential equation has regular singularities.
Then, for $\wc f(z_0)$ to be defined, the function $w\mapsto e^{-z_0w}$ must decay at the infinity of $\gamma$,
i.e., $\gamma$ can go to infinity only in the half-plane $\Re(z_0w)>0$.

\vskip .2cm

In this way,  studying  the relative growth of  exponentials $e^{-z_0w}$  for different $z_0$,  leads naturally to  questions of
convexity, to  emphasis  on straight lines etc. The following elementary example
(interchanging the roles of $z$ and $w$)
may be useful.

\begin{ex}\label{ex:exp-sums}
Let $A=\{w_1,\cdots, w_N\}\subset \CC_w$ be as before  and consider an exponential sum
\[
g(z) =\sum_{i=1}^N c_i e^{w_i z},
\]
with all $c_i\neq 0$. The  asymptotic geometry of the  distribution of zeroes of $g(z)$ is governed by the convex polygon
$Q=\Conv(A)$. More precisely,  the zeroes concentrate,  asymptotically, on the rays which are normal
to the edges of $Q$.  Indeed,  if we approach the infinity in any other direction, there will be a single summand
in $g(z)$  that will dominate all the other summands, so $g(z)\neq 0$. This observation, due to  H. Weyl
\cite{weyl}, is relevant for understanding
of integrals like \eqref{eq:FT}  and \eqref{eq:QC-int}  by the stationary phase approximation.
\end{ex}

\paragraph{Rectilinear transport maps.} \label{par:rect-trans}
Let us assume that $A$ is in linearly general position including infinity, i.e., $\zeta=1$
in Definition \ref{def:conv-pos}(d).
For any $w_i\in A$, we
identify $S^1_{w_i}$, the circle of directions at  $w_i$,  with the unit circle
$S^1=\{|\zeta |=1\}\subset\CC$.  As in Example \ref{ex:straight-lines}, denote by
\be\label{eq:zeta-ij}
\zeta_{ij} \,=\, {w_i -w_j \over |w_i-w_j|} \,\in \, S^1
\ee
the slope of the interval $[w_i, w_j]$. Thus $\zeta_{ij}\notin \RR$ by our assumption.

\vskip .2cm

\begin{figure}[H]
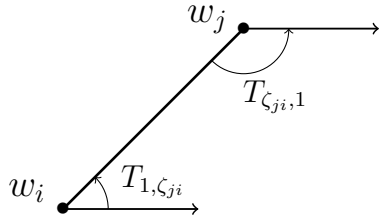

\centering

\btp[scale=0.6]

\node at (2,2){$\bullet$};
\draw[->, line width=0.7] (2,2) -- (5,2);
\node at (-2,-2){$\bullet$};
\draw[->, line width = 0.7] (-2,-2) -- (1, -2);
\draw[line width=1] (2,2) -- (-2,-2);

\node at (1.2,2.2){\large$w_j$};
\node at (-2.8 ,-1.6){\large$w_i$};

\draw [<-] (3,2) arc (0: -135: 1) ;

\draw [->]  (-1,-2) arc (0:45:1);

\node at (0,-1.4){$T_{1, \zeta_{ji}}$};

\node at (2.7, 0.5) {$T_{\zeta_{ji}, 1}$};

\etp
\caption{Rectilinear transport.} \label{fig:rect-trans}

\end{figure}

Let $\Fc\in\Perv(\CC, A)$. We denote by $\Phi_i(\Fc)$ the stalk of the local system $\bPhi_i(\Fc)$
at  the horizontal direction
$1\in S^1_{w_i}\simeq S^1$. For any $i\neq j$ we define the {\em rectilinear transport map}
$m_{ij}=m_{ij}(\Fc): \Phi_i(\Fc)\to\Phi_j(\Fc)$ as the composition
\be\label{eq:rect-trans-def}
\Phi_i(\Fc) = \bPhi_i(\Fc)_1 \buildrel T_{1,\zeta_{ji}} \over \lra \bPhi_i(\Fc)_{\zeta_{ji}}
\buildrel m_{ij}([i,j])\over\lra \bPhi_j(\Fc)_{\zeta_{ij}} \buildrel T_{\zeta_{ji}, 1}\over\lra
\bPhi_j(\Fc)_1 = \Phi_j(\Fc),
\ee
where:
\begin{itemize}
\item $m_{ij}([w_i,w_j])$ is the transport map along the rectilinear interval $[w_i, w_j]$.

\item $T_{1,\zeta_{ji}}$, resp. $T_{\zeta_{ij}, 1}$ is the monodromy from $1$ to $\zeta_{ji}$,
resp. from $\zeta_{ij}$ to $1$, taken in the counterclockwise direction, if
$\Im(w_i)<\Im(w_j)$,  and in the clockwise direction, if $\Im(w_i)>\Im(w_j)$, see
Fig. \ref{fig:rect-trans}.
\end{itemize}

\vskip .2cm

We also define
\[
m_{ii} = b_i  a_i: \Phi_i(\Fc)\lra\Phi_i(\Fc),
\]
where $a_i$ and $b_i$ are the map in the standard $(\Phi,\Psi)$-diagram
$
\xymatrix{
\Phi_i(\Fc) \ar@<.4ex>[r]^{a_i}&\Psi_i(\Fc) \ar@<.4ex>[l] ^{b_i}
},
$
representing $\Fc$ near $w_i$ in the horizontal direction. Thus
$1-m_{ii}= T_i(\Fc)$ is the counterclockwise monodromy of $\bPhi_i(\Fc)$,
in particular, it is invertible.

\vskip .2cm

Thus the data $(\Phi_i(\Fc), m_{ij})$ form an object of the category $\Mc_N$,
see Definition \ref{def:cat-M_A} and we have the {\em  functor of rectilinear transport data}
\be\label{eq:phi-para}
\Phi^\para = \Phi^\para_A: \ol\Perv(\CC,A) \lra \Mc_N, \quad \Fc\mapsto (\Phi_i(\Fc), m_{ij}).
\ee

\begin{prop}\label{prop:rect-cat}
The functor $\Phi^\para$ is an equivalence of categories.
\end{prop}

\begin{figure}[H]
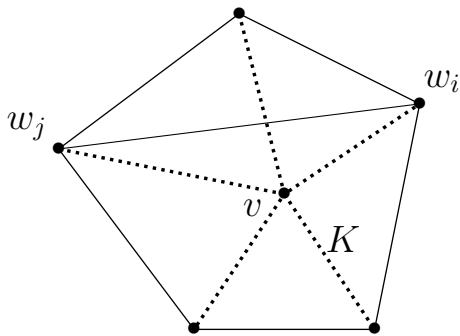

\centering

\btp[scale=0.6]

\node at (0,0){$\bullet$};
\node at (3,2){$\bullet$};
\node at (-1,4){$\bullet$};
\node at (-5,1){$\bullet$};
\node at (-2,-3){$\bullet$};
\node at (2,-3){$\bullet$};

\draw[line width = 0.5] (3,2) -- (-1,4) -- (-5,1) -- (-2, -3) -- (2, -3) -- (3,2);

\draw (-5,1) -- (3,2);

\draw[dotted, line width=1.3]  (-5,1) -- (0,0) -- (3,2);
\draw[dotted, line width=1.3] (-1,4) -- (0,0) -- (-2, -3);
\draw[dotted, line width=1.3] (0,0) -- (2, -3);

\node at (3.5, 2.5) {\large$w_i$};
\node at (-5.7, 1.5) {\large$w_j$};
\node at (-0.7, -0.3){\large$v$};

\node at (1.3,-1){\large$K$};

\etp
\caption{ For a convex configuration, rectilinear = GMV.} \label{fig:rect-GMV}

\end{figure}

\begin{rem}\label{rem:Phi-convexpos}
This statement is similar but not identical to Proposition \ref{prop:GMV2},
because the $[w_i, w_j]$ do not all pass through a single point $v$
(``Vladivostok''). Moreover, they cannot, in general, be deformed (without crossing other
$w_k$) to paths which have this property. However, if $A$ is   in   convex position,
then
we can take $v$ inside $Q$ and form a spider $K$ consisting of the intervals $[v, w_i]$
(dotted lines in Fig. \ref {fig:rect-GMV}). Then each interval $[w_i, w_j]$   will be isotopic to the path
$[w_i,v]\cup [v,w_j]$ used to define the Gelfand-MacPherson-Vilonen (GMV) transport map,
so in this case the statement does reduce to  Proposition \ref{prop:GMV2}.
\end{rem}

We postpone the proof  of Proposition
\ref {prop:rect-cat} until the discussion of isomonodromic deformations in the rectilinear approach.


\subsection {Isomonodromic deformations in the rectilinear approach}\label{subsec:iso-recti}

\paragraph{ Isomonodromic deformations of perverse sheaves. }

Denote by  $\CC^N_\neq\subset \CC^N$   the space of all numbered tuples
(``configurations'')  of $N$ distinct points $A= (w_1,\cdots, w_N)$.
We refer to $\CC^n_\neq$ as the {\em configuration space}. It is a classifying space of the pure braid group $\PBr_N$.

\vskip .2cm

Let $A(t) = (w_1(t),\cdots, w_N(t))_{0\leq t\leq 1}$ be a continuous path in $\CC^N_\neq$ joining two
configurations  $A(0)$ and $A(1)$. Since the concept of a perverse sheaf is purely topological, we have a canonical
(defined uniquely up to a unique isomorphism) equivalence of categories
\[
I_{A(t))} : \Perv(\CC, A(0))\lra\Perv(\CC, A(1)),
\]
which we call the {\em isomonodromic deformation along the path} $A(t)$.
This elementary but important concept  can be understood in either of three folllowing ways.

\vskip .2cm

(1)
Intuitively, we ``move the singularities'' $w_i(t)$ of a perverse sheaf $\Fc$ while
keeping the data
near the singularities, i.e., the local systems of diagrams
$
\xymatrix{
\bPhi_i \ar@<.4ex>[r]^{a_i}&\bPsi_i \ar@<.4ex>[l]^{b_i}
}
$
on $S^1_{w_i}=S^1$, the same.

\vskip .2cm

(2)  Note that
$A(t)$ gives a unique isotopy class of diffeomorphisms
$f: (\CC, A(0))\to (\CC, A(1))$ sending $w_i(0)$ to $w_i(1)$ and identical  near infinity.
The equivalence $I_{A(t)}$ is just the direct image $f_*$.
It is clear that  $I_{A(t)}$ depends
(up to a unique isomorphism) only on the homotopy class of the path $A(t)$ in $\CC^N_\neq$.
Therefore the $I_{A(t)}$ for all possible $A(t)$ (with arbitrary beginning and end configurations),
define a local system of abelian categories on $\wt\CC^N_\neq$.
The corresponding action of $\pi_1(\CC^N_\neq, A(0)) = \PBr_N$ on $\Perv(\CC,A(0))$ is the one
described in  Example \ref{ex:class-group-braid}.

\vskip .2cm

(3) Let  $p: \wt\CC^N_\neq\to\CC^N_\neq $  be  the  universal covering. Explicitly,
we take $\wt\CC^N_\neq$ to consist of homotopy classes of paths originating
at some given point $A(0)$, so (the constant path at) $A(0)$ is a distinguished point of
$\wt\CC^N_\neq$. Consider the trivial bundle
\[
\pi:  \wt\CC^N_\neq\times\CC\lra\wt\CC^N_\neq
\]
with fiber $\CC$ and let $\wt w_i, i=1, \cdots, N$ be its tautological section whose value
at $B$ with $p(B)=(w_1,\cdots, w_N)$, is $w_i$. Denote by $\Gamma_i \subset
\wt\CC^N_\neq\times\CC$ be the graph of $\wt w_i$.

\vskip .2cm

Let $\wt A\subset  \wt\CC^N_\neq\times\CC$
be the union of the $\Gamma_i$  and $\Perv( \wt\CC^N_\neq\times\CC, \wt A)$
be the category of perverse sheaves on $ \wt\CC^N_\neq\times\CC$
whose microsupport \cite{Ka-Scha} is contained in the union of the conormal bundles to
the $\Gamma_i$.
 Since such sheaves are smooth in directions transversal
to the fibers of $\pi$ and since $\wt\CC^N_\neq$ is homeomorphic to a ball, for any
$B\in\wt\CC^N_\neq$ with $p(B)=A$, the restriction on $\pi^{-1}(B)$ gives an
equivalence
\[
r_B: \Perv( \wt\CC^N_\neq\times\CC, \wt A) \lra \Perv(\CC,A).
\]
In particular, if a path $A(t)$ joins $A(0)$ with some $A(1)$, it lifts to a path $\wt A(t)$ in
$\wt\CC^N_\neq$
joining $A(0)$ with some point $B$ over $A(1)$, and we define $I_{A(t)} = r_B\circ r_{A(0)}^{-1}$.

\vskip .2cm

Our goal is to analyze the effect of isomonodromic deformations on rectilinear transport maps.

\paragraph{Walls of collinearity/marginal stability and configuration chambers.}
Let  $\CC^N_\genr\subset \CC^n_\neq$ consist   of $A$ in linearly general position.
Explicitly, for $A=(w_1,\cdots, w_n)$ consider the  real matrix
\[
P(A)\,= \,
\left[
\begin{matrix}
1&1& \cdots & 1
\\
\Re(w_1)&\Re(w_2)&\cdots & \Re(w_N)
\\
\Im(w_1)&\Im(w_2)&\cdots & \Im(w_N)
\end{matrix}
\right].
\]
 For $1\leq i<j<k\leq N$ let $p_{ijk}(A)$
be the $3\times 3$-minor of $P(A)$ on the columns numbered $i,j,k$.
 Then
$\CC^N_\genr = \CC^N_\neq \- \bigcup_{i,j,k} D_{ijk}$,
where $D_{ijk}$ is the  real  (of real codimension $1$) quadratic hypersurface
\[
D_{ijk} \,=\, \bigl\{p_{ijk}(A)=0\bigr\}\quad (w_i,w_j, w_k \text{ lie on a straight line})
\]
which we call the
{\em wall of collinearity}.  (In   \cite{GMW1} it is called a wall of marginal stability.)

\vskip .2cm

The connected components of $\CC^N_\genr$ (into which the walls $D_{ijk}$ subdivide $\CC^N_\neq$) are
open sets which we call
{\em configuration chambers}.
Describing these chambers in a more combinatorial fashion is a classical but  very difficult
problem of geometry. An obvious discrete invariant is as follows.

\begin{Defi}\label{def:or-datum}
(a)  The {\em  orientation datum} of $A=(w_1,\cdots, w_N)\in\CC^N_\genr$ is
the
collection $\OO(A) = (e_{ijk}(A))$ of signs
\[
e_{ijk}(A) \,=\, \sgn\bigl(p_{ijk}(A)\bigr)\,\in \, \{\pm 1\}, \quad 1\leq i<j<k\leq N.
\]

(b) The orientation datum $\OO(C)$ of a configuration chamber $C$ is defined as
$\OO(A)$ for any $A\in C$.

\end{Defi}

\paragraph{Oriented matroids.}
Note that the signs $e_{ijk}(A)$  are subject to restrictions coming from the Pl\"ucker relations among the minors $p_{ijk}$.
To write them down, let us extend the $p_{ijk}(A)$ (and $e_{ijk})(A)$, in a unique way, to alternating functions
defined for all triples  of  $i,j,k\in \{1,\cdots, N\}$ (in particular, make them equal to $0$ when two of the indexes coincide).
Then the relations are as follows.

\begin{prop}
For any  $i_1,\cdots, i_4$ and $j_1, j_2$ between $1$ and $N$, we have
\[
\sum_{\nu=1}^4 (-1)^\nu p_{i_1, \cdots, \wh i_\nu, \cdots, i_4} \cdot p_{j_1, j_2, i_\nu} \,=\,0.  \qed
\]
\end{prop}

This leads to the following definition, see  \cite{bjorner, bokowski-sturmfels}
for more information.

\begin{Defi} \label{def:or-matr}
(a) A structure of an {\em oriented matroid}
on $\{1,\cdots, N\}$  is a collection $\OO$ of
signs $\eps_{ijk} = \eps_{ijk}^\OO \in\{\pm 1\}$ given for any triple $i,j,k\in \{1,\cdots, N\}$ of distinct integers, which satisfy the axioms:
\begin{itemize}
\item[(1)] $\eps_{ijk}$ depends on $i,j,k$ in an alternating way: $\eps_{ijk}=-\eps_{jik}$ etc.

\item[(2)] Let us extend $\eps_{ijk}$ to all triples by putting $\eps_{ijk}=0$ if $i,j,k$ are not all distinct.
Then, for any  $i_1,\cdots, i_4$ and $j_1, j_2$ between $1$ and $N$, the set
\[
\bigl\{ (-1)^\nu \eps_{i_1, \cdots, \wh i_\nu, \cdots, i_4} \cdot \eps_{j_1, j_2, i_\nu}, \bigl|\,  \nu = 1\cdots, 4\bigr\}\,\subset \, \{0, \pm 1\}
\]
\end{itemize}
either contains $\{\pm 1\}$ or is equal to $\{0\}$.

\vskip .2cm

(b) For an  oriented matroid $\OO = (\eps_{ijk})$  let $\CC^N_\OO\subset\CC^N_\genr$ be the
set of $A$ such that $e_{ijk}(A)=\eps_{ijk}$ for all $i,j,k$. We say that $\OO$ is {\em realizable},
if $\CC^N_\OO\neq\emptyset$.
\end{Defi}

Thus,  the orientation datum $\OO(A)$ for any $A\in\CC^N_\genr$ is  an  oriented matroid.

\begin{rems}\label{rem:or-matr-real}
(a) More precisely, the concept we defined is known as a {\em uniform oriented matroid}
(or {\em simplicial chirotope}) {\em of rank $3$}, see  \cite{bjorner, bokowski-sturmfels}.
Since in this paper we consider only this specific type of oriented matroids,
we use the simplified terminology.

\vskip .2cm

(b) By construction, each $\CC^N_\OO$ is an open and closed subset, i.e., a union of
possibly several connected components (configuration chambers) in $\CC^N_\genr$. In fact,
$\CC^N_\OO$ can be disconnected.
In other words, it
is possible that more than one configuration chamber  has the same orientation datum.
For an  example (with $N=14$)  and estimates  for large $N$ see, respectively, \cite{suvorov} and
\cite{bokowski-sturmfels}, Cor. 6.3.

\vskip .2cm

(c) For large $N$, most   oriented matroids are not realizable, and it is a difficult problem to decide
whether a given oriented matroid is realizable or not.
See \cite{bokowski-sturmfels}, Fig. 1-1 for an example of non-realizable  oriented matroid with $N=9$.

\end{rems}

\paragraph{ Walls of horizontality and refined configuration chambers.}
Let $\CC^N_\genoo\subset\CC^N_\genr$ be the open part formed by confugurations in linearly general position including infinity.
It is obtained by removing the additional walls
\[
D_{ij} = D_{ij\oo} \,=\, \bigl\{ (w_1,\cdots, w_N)\, \bigl|  \, \Im(w_i)=\Im(w_j)\bigr\}, \quad 1\leq i<j\leq N,
\]
which we call the {\em walls of horizontality}.
The connected components of $\CC^N_\genoo$ will be called {\em refined  configuration chambers}. They
give rise to uniform oriented matroids on $\{1,\cdots, N, \oo\}$.

\paragraph{Isomonodromic deformations and rectilinear transport.}

Let $A(t)$ be a path in $\CC^N_\neq$ whose end points $A(0)$ and $A(1)$ are in linearly
general position including infinity, so the functors $\Phi^\para_{A(0)}$ and
$\Phi^\para_{A(1)}$ are defined, see \eqref{eq:phi-para}. We would like to describe
the effect of the isomonodromic deformation $I_{A(t)}$ on these functors by
constructing an explicit equivalence $J_{A(t)}$ fitting into the commutative diagram
\be\label{eq:I-and-J}
\xymatrix{
\Perv(\CC, A(0))
\ar[d]_{I_{A(t)}}^\sim
\ar[r]^{\hskip 1cm\Phi^\para_{A(0)}} & \Mc_N
\ar[d]^{J_{A(t)}}_\sim
\\
\Perv(\CC, A(1)) \ar[r]^{\hskip 1cm \Phi^\para_{A(1)}} & \Mc_N
}
\ee
and compatible with composition of paths. For this, it suffices to consider the
elementary situation of ``wall-crossing''. That is, we assume that:
\begin{itemize}
\item $A(0)$ and $A(1)$ lie in neighboring refined configuration chambers $C_0$ and $C_1$,
separated by a wall of collinearity $D_{ijk}$ or of horizontality $D_{ij}$.

\item $A(t)$ crosses the wall once, transversally at a smooth point not lying on other
walls.
\end{itemize}
The corresponding {\em wall-crossing formulas} for the change of $\Phi^\para$
can be easily extracted from the general Picard-Lefschetz formula of
Proposition \ref{prop:PL-form}.

\vskip .2cm

First of all, crossing any wall does not change the vector spaces $\Phi_i$ themselves,
as $\Phi_i(\Fc)$ is determined by the behavious of $\Fc$ on a small neighborhood of $w_i$.
So the effect of $J_{A(t)}$ will be the change of the maps $m_{ij}$. We will
denote by $m_{ij}$ the ``old'' map, corresponding to $A(0)$ and by $m'_{ij}$ the
``new'' map corresponding to $A(1)$. We write $T_i - 1-m_{ii}$, the invertible
transformation corresponding, under $\Phi^\para$, to the counterclockwise monodromy
of $\Fc$ around $w_i$.

\vskip .2cm

Consider first  the crossing of a wall of collinearity $D_{ij}$. That is, at some $t_0\in (0,1)$,
we have $\Im(w_i(t_0))=\Im(w_j(t_0))$. The shape of $J_{A(t)}$ will depend on whether
$w_j$ moves above $w_i$ or vice versa and also on the relative position of
$\Re(w_i)$ and $\Re(w_j)$ (at $t=t_0$). Following through various monodromies, we get:

\begin{prop}[Wall crossing: horizontality]\label{prop:wc-horiz}
Crossing $D_{ij}$  leaves unchanged  all $m_{pq}$
other than $m_{ij}$
and $m_{ji}$. Further, we have
\[
m'_{ij} = \begin{cases}
m_{ij}T_i, & \text{ if  $w_j$ moves above $w_i$ and } \Re(w_j)<\Re(w_i);
\\
T_j m_{ij}, & \text{ if  $w_j$ moves above $w_i$ and } \Re(w_j)>\Re(w_i);
\\
m_{ij}T_i^{-1}, & \text{ if  $w_j$ moves below $w_i$ and } \Re(w_j)<\Re(w_i);
\\
T^{-1}_j m_{ij}, & \text{ if  $w_j$ moves below $w_i$ and } \Re(w_j)>\Re(w_i).
\end{cases}
\]
and similarly for $m'_{ij}$ (interchange the roles of $i$ and $j$ above). \qed

\end{prop}

\begin{figure}[H]
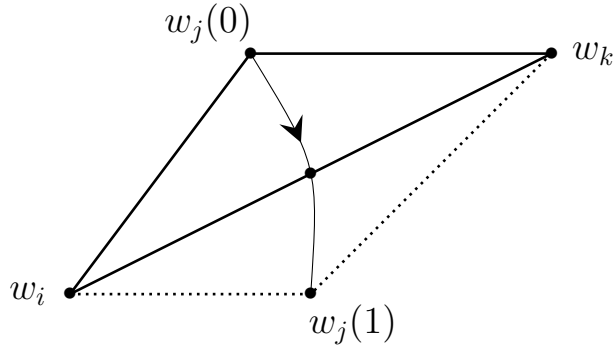

\centering
\btp[scale=0.8]
\node at (0,0){$\bullet$};
\node at (4,2){$\bullet$};
\node at (-1,2){$\bullet$};
\node at (-4,-2){$\bullet$};
\node at (0,-2){$\bullet$};
\draw[line width = 01]  (4,2) -- (-1,2) -- (-4, -2);
\draw[dotted, line width = 1]  (-4, -2) -- (0, -2) -- (4,2);
\draw[line width=1]  (-4, -2) -- (4,2);
\draw [decoration={markings,mark=at position 0.4 with
{\arrow[scale=3,>=stealth]{>}}},postaction={decorate}]
(-1,2) .. controls (0.14,0.14) .. (0, -2);

\node at (-4.7, -2){\large$w_i$};
  \node at (4.7, 2){\large$w_k$};
\node at (-1.7, 2.5){\large$w_j(0)$};
  \node at (.7, -2.5){\large$w_j(1)$};

\etp
\caption{Crossing a wall of collinearity. }\label{fig:wc-coll}
\end{figure}

Consider now the crossing of a wall of collinearity $D_{ijk}$. We choose the order
of indices $i,j,k$ so that for some $t=t_0\in(0,1)$ the point $w_j(t)$ crosses
the interval $[w_i(t), w_k(t)]$ in its interior point (so not necessarily $i<j<k$). See
Fig. \ref{fig:wc-coll}.
We denote by $\eps_{ijk}\in\{\pm 1\}$
the sign of Definition \ref {def:or-datum}
for $A(0)$, i.e., the
orientation of the triangle
$\bigtriangleup [w_i(0), w_j(0), w_k(0)]$ and by $\eps'_{ijk}$ the similar sign for
$A(1)$. Clearly, $\eps'_{ijk}=-\eps_{ijk}$: the orientation changes as we cross the wall.

\begin{prop}[Wall-crossing:collinearity]\label{prop:wc-coll}
Crossing $D_{ijk}$ (so that $w_j$ crosses $[w_i, w_k]$)
does not change any $m_{pq}$ other than $m_{ik}$ and $m_{ki}$. Further,
\[
m'_{ik} \, = \,
m_{ik} + \eps_{ijk} m_{jk}\circ m_{ij},
\]
and similarly for $m_{ki}$ (interchange the roles of $i$ and $k$ above). \qed
\end{prop}

\begin{rem}
Propositions \ref{prop:wc-horiz} and \ref{prop:wc-coll} can be seen as giving an action,
on the category $\Mc_N$, of the {\em fundamental groupoid} of $\CC^N_{\neq,\infty}$,
in which we take one base point in each refined configuration chamber.
It is a convex geometry counterpart of the (purely topological)
action of the {\em fundamental group}
$\pi_1(\CC^N_\neq)=\PBr_N$ on $\Mc_N$ described in \cite{gelfand-MV}, Prop. 2.4.

\end{rem}

\paragraph{Proof of Proposition \ref{prop:rect-cat}.}

We now go back to the proof of the Proposition \ref{prop:rect-cat}.
We have already seen in Remark \ref{rem:Phi-convexpos} that for $A$ in
convex position, $\Phi^\para_A$ is an equivalence. Now,
any configuration $A\in \CC^N_\genoo$ can be connected to a configuration
in convex position by a chain of wall-crossings. For each wall-crossing
we have a diagram
\eqref{eq:I-and-J} with $I_{A(t)}$ and $J_{A(t)}$ being equivalences.
Therefore,  if $\Phi^\para_{A(0)}$ is an equivalence, then
$\Phi^\para_{A(1)}$ is an equivalence as well. In this way we inductively deduce our statement
for arbitrary $A\in \CC^N_\genoo$ from the case of convex position.


\subsection {Fourier transform of perverse sheaves: topology}\label{subsec:FTPS-top}

\paragraph{ Reminder on $\Dc$-modules and the Riemann-Hilbert correspondence.}
We recall some basic material, see \cite{hotta, malgrange} for more detail.

Let $X$ be a smooth algebraic variety over $\CC$ We denote by $\Dc_X$ the sheaf of algebraic differential operators
on $X$ in Zariski topology. Let $\Dc_X\Mod$ be the category of sheaves of left $\Dc_X$-modules $\Mc$ on $X$ which
are coherent (=locally finitely presented) over $\Dc_X$. Such an $\Mc$ can be regarded, in a standard way (see {\em loc. cit.})
as a system of algebraic linear partial differential equations on $X$.

Let $X^\an$ the complex analytic manifold corresponding to $X$,
with its sheaf $\Oc_X^\an$ of holomorphic functions. Clearly, $\Oc_{X^\an}$ is a left $\Dc_X$-module (not coherent).
The {\em solution complex} of $M\in\Dc_X\Mod$ is
the complex of sheaves
\be
\Sol(M) \,=\, \ul\RHom_{\Dc_X}(\Mc, \Oc_X^\an)
\ee
on $X^\an$. Its $0$th cohomology sheaf, i.e., $\ul\Hom_{\Dc_X}(\Mc, \Oc_X^\an)$,  is the sheaf of holomorphic solutions, in the  usual  sense,  of the system of PDE corresponding to $\Mc$. Note that $\Sol$ is a contravariant functor.

Let
\[
\Dc_X\Mod^\hol \,\, \supset \, \,\Dc_X\Mod^\hr
\]
be the full subcategories in $\Dc_X\Mod$ formed by modules which are {\em holonomic} (``maximally overdetermined'')
and {\em holonomic regular} (``maximally overdetermined with regular singularities''). We refer to
\cite{hotta, malgrange} for precise definitions. Let us just mention that the condition of being regular involves,  for
a non-compact variety $X$, also the behavior at infinity, i.e., on a compactification of $X$.

Let $\Perv(X)$ be the category of perverse sheaves on $X^\an$ (constructible with respect to some algebraic
stratification). The base field $\k$ is assumed to be $\CC$.

The following is
a summary of the basic results of the theory;  part (b) is known as the {\em Riemann-Hilbert correspondence}
(for $\Dc$-modules).

\begin{thm}\label{thm:RH}
(a) If $\Mc$ is any holonomic $\Dc_X$-module (regular or not), then $\Sol(M)$ is a perverse sheaf.

\vskip .2cm

(b)  The restriction of the functor $\Sol$ to $\Dc_X\Mod^\hr$ is an equivalence of categories
\[
\Dc_X\Mod^\hr \buildrel\sim\over \lra \Perv(X)^\op. \qed
\]
\end{thm}

\begin{rem}
The theorem implies that for any holonomic (possibly irregular) $\Mc$,  there is a canonical holonomic regular
module $\Mc^\reg$ (the ``regularization'' of $\Mc$) such that $\Sol(\Mc^\reg)\simeq \Sol(\Mc)$. See \cite{KK}
for the general discussion and \cite{malgrange} (I.4.6) for an example which shows that $\Mc^\reg$
cannot, in general, be defined in a purely algebraic way.
\end{rem}

\paragraph {Fourier transform for $\Dc$-modules and perverse sheaves.}
Let $V$ be a finite-dimensional $\CC$-vector space and $V^*$ is the dual space.
Considering $V$ as an affine algebraic variety, we see that $\Dc_V\Mod$ is equivalent to the category
of finitely generated modules $M$ over the Weyl algebra $\Dc(V)$ of global polynomial differential operators
on $E$, via the {\em localizattion}
\[
M\,\mapsto \Mc \,=\, M\otimes_{\CC[V]} \Oc_{V}, \quad M=H^0(V, \Mc).
\]
Here $\CC[V]$ is the algebra of polynomial functions on $V$.
Explicitly, $\Dc(V)$ is generated by linear functions $l\in V^*$ and constant vector fields $\del_v, v\in V$,
subject to $[\del_v, l]=l(v)$. This leads to the identification
\[
f: \Dc(V^*) \lra\Dc(V), \quad v\mapsto -\del_v, \quad \del_f \mapsto f, \quad v\in V=(V^*)^*, f\in V^*.
\]
Given a $\Dc(V)$-module $M$, its {\em Fourier transform} is the $\Dc(V^*)$-module $\wc M$ which,
as a vector space, is the same as $M$ but with action of $\Dc(V^*)$ given via $f$. Note that the
{\em sheaves} $\Mc$ and $\wc \Mc$ of $\Dc_V$ and $\Dc_{V^*}$-modules associated to $M$ and $\wc M$,
can look quite different, since it is the  action of the $\del_v$  on $M$ that is  localized in forming $\wc\Mc$.

\vskip .2cm

Let now $\Fc\in\Perv(V)$ be a perverse sheaf.
By  Theorem \ref{thm:RH}(b), $\Fc=\Sol(\Mc)$ for a unique holonomic regular  sheaf $\Mc$
of $\Dc_V$-modules, which is encoded by the $\Dc(V)$-module $M=H^0(V,\Mc)$.
Let us  the Fourier transform $\wc M$ and the corresponding sheaf
$\wc\Mc$ of $\Dc_{V^*}$-modules. Note that $\wc\Mc$ may be irregular, but by Theorem \ref{thm:RH}(a),
its solution complex is perverse, so
we have a new perverse sheaf, now on $V^*$,
\be\label{eq:FT-of-F}
\wc\Fc \,=\,\Sol(\wc\Mc) \,=\,
\ul{\RHom}_{\Dc(V^*)}(\wc M, \Oc^\an_{V^*}),
\ee
which we call the {\em Fourier transform} of $\Fc$.

\paragraph {Description of $\wc\Fc$.}

In this paper we will be interested in the $1$-dimensional case only. That is, we take $V=\CC_w$, $V^*=\CC_z$
to be the complex lines with ``dual'' coordinates $w$ and $z$. As before, we fix $A=\{w_1,\cdots, w_N\}$
and take $\Fc\in\Perv(\CC_w, A)$.

\vskip .2cm

A description of $\wc\Fc$ in this case
be deduced from  results of Malgrange \cite{malgrange}.
We present this description in an explicit form which highlights the analogy with the
Fukaya-Seidel
category. First, we have:

\begin{prop}\label{prop:wc-Phi}\label{prop:wc-M-top-1}
The only possible singularity of $\wc M$ in the finite part of  $\CC_z$ is at $z=0$, and this singularity is regular.
In particular, $\wc \Fc \in \Perv(\CC_z, 0)$.
\end{prop}

\noindent{\sl Proof:} This is a particular case of more general statements of \cite{malgrange}
involving not necessarily regular holonomic $\Dc_\CC$-modules. More precisely, the determination of
singularities follows from  \cite{malgrange}  (XII.2.2). In fact, both this determination and
the claim that the singularity at $0$ is regular follow from the ``stationary phase theorem'' of
\cite{malgrange} (see \cite{dagnolo-microlocal} Th. 1.3 for a concise formulation) which describes
the exponential Puiseux factors of $\wc M$ as the Legendre transforms of those of $M$. \qed

\vskip .2cm

Next, as $\wc\Fc$ lies in $\Perv(\CC_z,0)$, it is completely determined by its $(\Phi, \Psi)$-diagram
(Proposition \ref{prop:ggm})
\be\label{eq:Phi-Psi-FT}
\xymatrix{
\wc\Phi := \bPhi(\wc\Fc)_\zeta \ar@<.4ex>[rr]^{\wc a=\wc a_\zeta }&&\  \bPsi(\wc\Fc)_\zeta =: \wc\Psi \ar@<.4ex>[ll]^{\wc b=\wc b_\zeta}
}
 \ee
for some direction $\zeta$.  To describe it, we introduce some notation.

\vskip .2cm

\begin{figure}[H]
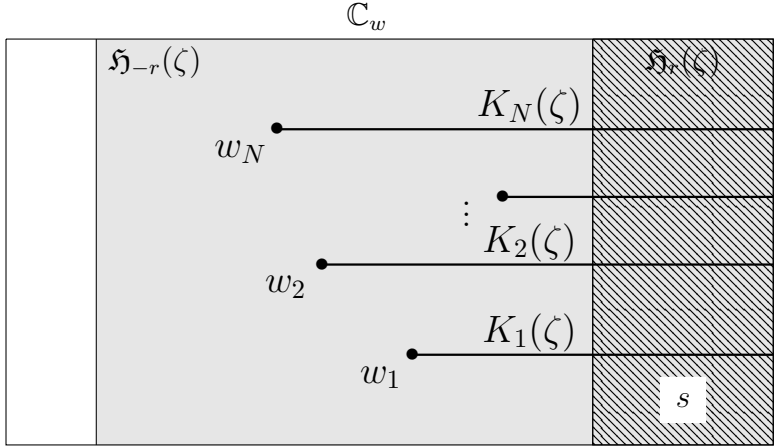

\centering
\btp[scale=.6]

\draw (-8, 4) -- (-8,-5) -- (9,-5) -- (9, 4) -- (-8,4);
\node at (1,-3){$\bullet$};
\node at (-1,-1) {$\bullet$};
\node at (3,0.5) {$\bullet$};
\node at (-2,2) {$\bullet$};

\node at (0,4.5){$\CC_w$};

\draw [line width=.3mm] (1,-3) -- (9,-3);
\draw [line width=.3mm] (-1,-1) -- (9,-1);
\draw [line width=.3mm] (3,0.5) -- (9,0.5);
\draw [line width=.3mm] (-2,2) -- (9,2);

\node at (0.3, -3.5){\large$w_1$};  \node at (3.6,-2.5){\large$K_1(\zeta)$};
\node at (-1.8,  -1.5){\large$w_2$};  \node at (3.6, -0.5){\large$K_2(\zeta)$};
\node at (2.2,  0.3){\large$\vdots$};
\node at (-2.8,  1.5){\large$w_N$};    \node at (3.6, 2.5){\large$K_N(\zeta)$};

\filldraw[opacity=0.1] (-6,4) -- (-6,-5) -- (9,-5) -- (9, 4) -- (-6, 4);
\draw[line width= 0.1mm] (-6, 4) -- (-6,-5);

\filldraw[ pattern=north west lines] (5,4) -- (5,-5) -- (9,-5) -- (9, 4) -- (5, 4);
\draw[line width= 0.1mm] (5, 4) -- (5,-5);

\node at (-4.7, 3.5){$\Hen_{-r}(\zeta)$};
\node at (7, 3.5){$\Hen_r(\zeta)$};

\filldraw[color=white] (7.5, -4.5) -- (6.5, -4.5) -- (6.5, -3.5) -- (7.5, -3.5) -- (7.5, -4.5);
\node at (7,-4) {$s$};
\etp

\caption{The plane $\CC_w$ with the  cuts $K_i(\zeta)$  at  $w_i$.  Here $\zeta=-1$. } \label{fig:horiz-cuts}
\end{figure}

The circles of directions around $w_i$ in $\CC_w$ and
around $0$ or  $\oo$ in $\CC_z$ are both parametrized by the standard unit circle $|\zeta|=1$
but we will use the following  {\em orientation reversing} identification:
\be\label{eq:sigma-i}
\sigma_i:  S^1_0(\CC_z) = S^1_\oo(\CC_z) \lra S^1_{w_i}(\CC_z), \quad \zeta \,\mapsto \, -\ol\zeta.
\ee
For any $\zeta\in S^1$ and any $r\in\RR$ we consider the  following half-plane,
half-rays of slope $-\ol \zeta$ originating
from the $w_i$ and their union:
\be\label{eq:half-plane(zeta)}
\begin{gathered}
\Hen_r (\zeta) \,=\,\bigl\{ w\in\CC
\bigl|\,
\,  \Re(- \zeta\cdot  w) \geq r\bigr\},
\\
K_i(\zeta) \,=\, w_i - \ol \zeta\cdot\RR_+, \quad K(\zeta) \,=\,  \bigcup K_i(\zeta).
\end{gathered}
\ee
see Fig.  \ref{fig:horiz-cuts}  which corresponds to $\zeta=-1$.
The choice of $(-\zeta)$ vs.  $(-\ol\zeta)$ in these definitions is explained by the following
obvious fact.

\begin{prop}\label{prop:role-zeta}
(a) The function $w\mapsto e^{\zeta\cdot w}$ decays in $\Hen_r(\zeta)$.

(b) If $r\geq |w_i|$ for all $i$, then $K(\zeta)\subset \Hen_{-r}(\zeta)$. In particular,
$w\mapsto e^{\zeta\cdot w}$ decays on each $K_i(\zeta)$.
  \qed
\end{prop}

We assume that elements of $A$ are in linearly general position including $(-\ol\zeta)$-infinity and numbered so that
$\Im(-\zeta w_1) < \cdots < \Im(-\zeta w_N)$.
Then the  rays $K_i(\zeta)$  do not intersect.
The union $K(\zeta)$ can be viewed as a  ``spider'' in which the ``Vladivostok'' point $\bv$ has been moved to
$(-\ol\zeta\oo)$, see \S \ref{subsec:abs-PL}\ref{par:spider}

\noindent
Let  $ \bPhi_i(\Fc)\in\LS(S^1_{w_i})$ be the local system of vanishing cycles of $\Fc$ at $w_i$
and $\bPsi(\Fc)\in \LS(S^1_\oo)$ be the local system of nearby cycles of $\Fc$ at $\oo$.
Let us denote for short
\[
\Phi_ i\,  = \, \bPhi_i(\wc\Fc)_{ -\ol\zeta}  \,=\,\ul\HH^1_{K_i(\zeta) }(\Fc)_{w_i}, \quad \Psi =  \bPsi(\wc\Fc)i_{-\ol\zeta}  \,=\, \HH^0(\Hen_r (\zeta), \Fc),
\,\, r\gg 0.
\]
 The restriction of sections from $\Hen_r(\zeta) $ to the sufficiently far part of
$K_i(\zeta)$
identifies $\Psi$ with the space $\Psi_i(\Fc)_{-\ol\zeta} $ of the standard $(\Phi, \Psi)$-description of $\Fc$ near $w_i$ given
by $K_i(\zeta)$.
 Thus we have, for each $i$,  the diagram
\[
\xymatrix{
\Phi_{i}   \ar@<.4ex>[r]^{a_{i}}& \Psi\ar@<.4ex>[l]^{ b_{i}}
},\quad a_i = a_{i, -\ol\zeta},\,\, b_i =  b_{i,-\ol\zeta},
\]
describing $\Fc$ near $w_i$ in the direction $(-\ol\zeta)$.
 In particular, $T_{i, \Phi} =1-b_{i} a_{i}$
  is the full counterclockwise monodromy on $\Phi_{i}$.
We also denote $T_{i,\Psi} = 1-a_{i}  b_{i}$ the monodromy on $\Psi$ obtained by going around
$w_i$ along $K_i(\zeta)$ (i.e., from $\oo$ to near $w_i$ along $K_i(\zeta)$, then around $w_i$, then back to
$\oo$ along $K_i(\zeta)$).

\begin{prop}\label{prop:wc-Phi}\label{prop:wc-M-top-2}
 We have identifications of the  spaces in the diagram \eqref{eq:Phi-Psi-FT} as
 \[
 \wc \Psi \,\,\simeq \,\,   \, \bigoplus_{i=1}^N \Phi_{i} = \bigoplus_{i=1}^N \bPhi_i(\Fc)_{-\ol\zeta},\quad \quad
\wc \Phi \, \, \simeq \, \, \Psi =  \bPsi(\Fc)_{-\ol\zeta},
\]
so that the maps are identified as
\[
\wc{b} = -\sum_{i=1}^N a_i T_{i,\Phi}^{-1}, \quad
\wc{a} =  (\wc a_1, \cdots, \wc a_N), \quad  \wc a_i \,=\, b_i \, T_{i-1, \Psi}^{-1} ... T_{1, \Psi}^{-1}.
\]

\end{prop}

\noindent{\sl Proof:} This is a reformulation of  \cite{dagnolo-sabbah}, Prop. 6.1.4.
Alternatively, the statements
  are particular cases of more general
results (involving not necessarily regular $\Dc(\CC_w)$-modules) given in \cite{malgrange} \S XIII.2,
see especially   formulas (XIII.2.4)  and (XIII.2.5).
 In comparing with  \cite{dagnolo-sabbah, malgrange}, one need to be mindful of different conventions.
 For instance, in \cite{malgrange} p. 199,  the plane $\CC_z$ is equipped with the  orientation opposite
 to the one given by the complex structure, which replaces various monodromies by their inverses.
 This corresponds to our identifications \eqref{eq:sigma-i}.
 \qed

\begin{rems}\label{rems:one-shot}
(a) A one-shot description of the map $\wc a$ is as follows. Take a section $s\in \HH^0(\Hen_r, \Fc)$ and
look at its values in the area below  $K=\bigcup K_i$ (i.e., in the lower right corner of Fig. \ref{fig:horiz-cuts}).
Extend  $s$ from that area  to a section $ s'$ in $\CC\- K$. Then $\wc a(s)$ is obtained by   taking
the image of $s'$
in each $\ul\HH^1_{K_i}(\Fc)_{w_i}$ under the coboundary map.

\vskip .2cm

(b)
Proposition \ref {prop:wc-M-top-2}  implies that the monodromy operator on $\wc \Phi$
has the form
\[
T_{\wc \Phi}= T_{N,\Psi}^{-1} T_{N-1,\Psi}^{-1} \cdots T_{1,\Psi}^{-1}
\]
the full {\em clockwise} monodromy on $\Psi(\Fc)$, cf. \cite{malgrange} (XIII.2.6).
For example, for $N=2$ we have, applying the Jacobson identity \eqref{eq:jacobson},
\[
\begin{gathered}
T_{\wc\Psi}= 1-\wc b\wc a = 1 + T_{1,\Phi}^{-1} b_1 + a_2 T_{2,\Phi}^{-1} b_2 T_{1,\Phi}^{-1}=
\\
= 1 + a_2(1-b_1 a_1)^{-1} b_1 + a_2 (1-b_2 a_2)^{-1} b_2 (1-b_1a_1)^{-1} =
\\
= (1-a_1b_1)^{-1}  -  (1-a_1b_1)^{-1} + (1+a_2 (1-b_2 a_2)^{-1}b_2) (1-a_1b_1)^{-1} =
\\
= (1-a_2b_2)^{-1} (1-a_1b_1)^{-1},
\end{gathered}
\]
and similarly for any $N$. This relation is clear from the first principles since moving $\zeta$ counterclockwise
amounts to rotating the half-planes $\Hen_r(\zeta)$ and the cuts $K_i(\zeta)$ clockwise.

\vskip .2cm

(c) Let $A=\{0\}$. In this case  the $\Dc_w$-module $M$ corresponding to $\Fc$ is
 monodromic, and so $\wc M$ is also regular,
see \cite{brylinski-ast, br-mal-ver}.  Proposition \ref{prop:wc-M-top-2}   in this case reduces
 to Proposition 4.5 of \cite{bezr-kapr}.

\end{rems}

 \paragraph{$\wc\Fc$ and Fourier integrals.} The first identification of  Proposition \ref{prop:wc-M-top-2}
 which we write as
 \be\label{eq:g-zeta-def}
 g_{\zeta} \,=\,\sum_{i=1}^N g_{i,\zeta}:\,\,  \bigoplus_{i=1}^N  \bPhi_i(\Fc)_{-\ol\zeta}
 \buildrel \sim\over  \lra \bPsi (\wc\Fc)_{\zeta},
 \quad g_{i,\zeta}: \bPhi_i(\Fc)_{-\ol\zeta} \to \bPsi(\wc\Fc)_\zeta,
 \ee
is in fact given in terms of actual Fourier integrals of (distribution) solutions of $\Mc$, the regular
holonomic $\Dc$-module corresponding to $\Fc$. Here we present this interpretation, following
\cite{malgrange} Ch. XII.

\vskip .2cm

For this, we   take
$\Mc = \Dc_{\CC_w}/\Dc_{\CC_w}\cdot p$  to correspond to a single polynomial differential operator
$p\in \Dc(\CC_w)$. This is known to be possible for any holonomic $\Dc$-module in one variable,
see \cite{malgrange} \S V.1.  (It is probably possible to adjust the argument to avoid this assumption, by working
with matrix differential operators.)
Thus $\Fc = \Sol(\Mc)$ consists, outside of $A$, of holomorphic solutions of $p(f)=0$.

\vskip .2cm

Let $\Gamma\subset \CC$ be an arc, i.e., an $1$-dimensional real analytic submanifold, possibly  with boundary.
Recall \cite {harvey, morimoto} that the sheaf of {\em hyperfunctions} on $\Gamma$ is defined as
\[
\Bc_\Gamma \,=\, \ul H^1_\Gamma(\Oc_\CC) \quad \text{(and we have $H^{\neq 1}_\Gamma(\Oc_\CC)=0$).}
\]
 Sometimes sections of $\Bc_\Gamma$ near
the endpoints of $\Gamma$ are referred to as {\em microfunctions}. The sheaf $\Bc_\Gamma$ is acted
upon by $\Dc_\CC$ via its action on $\Oc_\CC$. According to Sato, see {\em loc. cit.}, sections of
$\Bc_\Gamma$ are interpreted as a version of distributions (generalized functions)
supported on $\Gamma$. For $\Fc=\Sol(\Mc)$ with $\Mc$ holonomic regular as above, we have identifications
of sheaves
\be\label{eq:Phi=hyperf}
\ul\HH^1_\Gamma \Sol(\Mc) \, =\, \ul\HH^1_\Gamma\, \ul\RHom_\Dc(\Mc, \Oc)\,\=\,
\ul\Hom_\Dc(\Mc, \ul H^1_\Gamma(\Oc)) \,=\, \ul\Hom_\Dc(\Mc,\Bc_\Gamma)
\ee
with the sheaf of solutions of $p(f)=0$ in $\Bc_\Gamma$. Such solutions can be viewed as distributions
in any other classical sense.

\vskip .2cm

We take $\Gamma=K_i(\zeta)$, the half-line from \eqref{eq:half-plane(zeta)}, and let
\[
\phi \,\in\, \bPhi_i(\Fc)_{-\ol\zeta} \,=\, \ul\HH^1_{K_i(\zeta)}(\Fc)_{w_i} \,=\, \HH^1_{K_i(\zeta)}(\CC, \Fc),
\]
 the second identification coming from the fact that $\ul\HH^1_{K_i(\zeta)}(\Fc)$ is (locally) constant
 on the open ray $K_i(\zeta)\-\{w_i\}$. Let
 \be\label{eq:alpha-i-phi}
 \alpha_i(\phi) \,\in\, \Hom_\Dc(\Mc, \Bc_{K_i(\zeta)}) \,=\,\bigl\{ f\in\Bc_{K_i(\zeta)}(K_i(\zeta)) \,\bigl| \,
 p(f)=0\bigr\}
 \ee
be the global hyperfunction solution of $\Mc$ on $K_i(\zeta)$ corresponding to $\phi$ by
\eqref{eq:Phi=hyperf}. As. distribution, $\alpha_i(\phi)$ has at most polynomial growth at infinity
(since $\Mc$ is regular). Therefore the Fourier integral
\be\label{eq:fourier-integral}
G_{i,\zeta} (\phi) (z) \,=\,\int_{w\in K_i(\zeta)} \alpha_i(\phi)(w)  e^{zw} dw
\ee
converges for $z\in\zeta\cdot \RR_{>0}$ and extends to a holomorphic solution of $\wc\Mc$ in an open sector
around $\zeta \cdot \RR_{>0}$. Thus we have a map
$G_{i,\zeta}: \bPhi_i(\Fc)_{-\ol\zeta}\to\bPsi(\wc\Fc)_\zeta$.

\begin{prop}\label{prop:FT=integrals}
The map $G_{i,\zeta}$ coincides with $g_{i,\zeta}$, i.e., the latter is given by the integral
\eqref{eq:fourier-integral}.
\end{prop}

\noindent{\sl Proof:} This is \cite{malgrange} Thm. XII.1.2. \qed


\subsection{  Fourier transform of perverse sheaves: Stokes structure}\label{subsec:FT-stokes}

\paragraph{ Setup.}
As before, let $A \subset\CC_w$ and
 $\Fc\in\Perv(\CC_w, A)$.
Let $\Mc$ be the holonomic regular $\Dc_{\CC_w}$-module with $\Sol(\Mc)=\Fc$, let $\wc\Mc$ be its Fourier transform
and $\Fc=\Sol(\Mc)$. By Proposition \ref {prop:wc-M-top-1}, $\wc \Mc$ has only one possible singularity in the finite part
of $\CC_z$, namely $0$, which is regular.
But the singularity at $\oo$ may be irregular, so $\wc\Mc$ contains more information
than $\wc\Fc$, encoded in its   Stokes structure  near $\oo$.
This structure has been studied,   rather explicitly,  in    \cite {malgrange, dagnolo-sabbah}.
The goal of this section and the next one is to express it in a particularly explicit and striking for rently, not noticed before)
which can be seen as a de-categorification
of the $A_\oo$-Algebra of the Infrared.

\paragraph{Stokes structures: general case.}

We start with recalling the necessary concepts.
The basic references are \cite{babbitt, deligne-docs, KKP, sabbah}.

\vskip .2cm

We consider the case of only one irregular singular point, at $\oo$. Let
\[
\CC((z^{-1}))   \, \supset \, \CC\{\{z^{-1}\}\}
\]
be the fields of formal Laurent series vs. of germs of meromorphic functions near $\oo$.
Both fields are equipped
with the action of the derivative $\del_z$.

A {\em formal} (resp. {\em convergent}) {\em meromorphic connection}  near $\oo$
is a  datum $(E, \nabla)$, where $E$ is a
finite-dimensional vector space over   $\CC((z^{-1}))$ (resp. $\CC\{\{z^{-1}\}\}$)
and $\nabla:  E\to E$ is a $\CC$-linear map satisfying the Leibniz rule
\[
\wh\nabla (f\cdot e) \,=\, (\del_zf)\cdot e\, + \, f\cdot \wh\nabla (e), \quad f\in \CC((z^{-1}))
\text{ or } \CC\{\{z^{-1}\}\}, \,\, e\in E.
\]
A holonomic  $\Dc_{\CC_z}$-module $\Nc$ gives a
convergent meromorphic connection $(E,\nabla)$ near $\oo$,
where $E=\Nc\otimes_\Oc \CC\{\{z^{-1}\}\}$ and $\nabla$ is induced by the action of
$\del_z$ on $\Nc$.

\vskip .2cm

Let $S^1=S^1_\oo$ be the circle of directions at $\oo$, identified with the unit circle
$\{|\zeta|=1\}$ in $\CC$. For an open  $U\subset S^1$ let $C_U=\{z| \arg(z)\in U\}$
be the corresponding domain in $\CC$ (so it is an open  sector if $U$ is connected).
\vskip .2cm

The circle $S^1$ carries a locally constant sheaf $\Ac$ called {\em the sheaf of
exponential germs}.  By definition, for an open arc $U$ as above, $\Ac(U)$ consists
of germs of $1$-forms $\alpha$ defined near the infinity of $C_U$ and having
a Puiseux expansion of the form (for some  $p\in\ZZ$ depending on $\alpha$)
\[
\alpha \,=\, \sum_{k\in {1\over p} \ZZ,  \,\, k >-1} a_k z^k dz.
\]
Further, $\Ac$ is made into a
{\em sheaf of partially ordered sets}. That is, for any open $U\subset S^1$
we have a partial order $\leq_U$ on $\Ac(U)$. It is defined by
\[
\alpha \leq_U \beta \quad \Leftrightarrow \quad \exp \int (\alpha-\beta)
\quad \text {is of at most polynomial growth a the infinity of $C_U$}.
\]
This partial order  is compatible with
restrictions and gluing of sections. In particular, it gives an order $\leq_\zeta$
on the stalk $\Ac_\zeta$ at any $\zeta\in S^1$.

\vskip .2cm

Note that for two given sections $\alpha\neq \beta$ on some $U$, there are finitely
many points  $\zeta_1,\cdots, \zeta_n \in U$ at which the order $\leq_\zeta$ switches,
so  $\alpha$ and $\beta$ become incomparable at $\zeta=\zeta_j$,
while for all other $\zeta\neq\zeta_j$ we have either $\alpha <_\zeta\beta$
or $\beta <_\zeta\alpha$.

\vskip .2cm

As $(\Ac, \leq)$ is a sheaf of posets, one can speak about $\Ac$-{\em filtered} and
$\Ac$-{\em graded}  sheaves   on $S^1$
in a straightforward way, see \cite{sabbah}  for precise definitions.

\vskip .2cm

Thus, an $\Ac$-filtration $F$ on a sheaf $\Lc$ gives,  for each $\zeta$,
a filtration $F_\alpha \Lc_\zeta$ in the stalk $\Lc_\zeta$ labelled by  elements $\alpha$
of the poset $\Ac_\zeta$.
An $\Ac$-grading on $\Lc$ gives, for each $\zeta$, a grading, i.e., a direct
sum decomposition $\Lc_\zeta =\bigoplus_{\alpha\in \Ac_\zeta} \Lc_{\zeta, \alpha}$
on the stalk. These local filtrations and gradings should vary with $\zeta$
in an appropriately continuous way, see {\em loc. cit.}

\vskip .2cm

Any grading gives a tautological filtration called {\em split}:
\[
F_\alpha \Lc_\zeta \,=\,\bigoplus_{\beta\leq\alpha} \Lc_{\zeta,\alpha}.
\]
Any filtration $F$ in $\Lc$ gives  the {\em associated graded sheaf}
$\gr^F\Lc$  with grading summands being the quotients of the filtration:
\be\label{eq:gr-alpha-F}
\gr_\alpha^F\Lc \,=\, F_\alpha \Lc
\biggl/\sum_{\beta<\alpha} F_\beta\Lc.
\ee
If $F$ is the  split filtration corresponding to a grading, then $\gr^F\Lc$ recovers
the grading.

\begin{Defi}
Let $\Lc$ be a  finite rank local system of $\CC$-vector spaces on $S^1$.

\begin{itemize}

\item[(a)] A {\em Stokes structure} on $\Lc$ is
an  $\Ac$-filtration $F$ which is {\em locally split}, i.e., locally isomorphic to the  filtration
coming from a grading. A local system with a Stokes structure is called
a {\em Stokes local system}.

\item[(b)] The {\em Stokes sheaf} $\Uc(\Lc)$  of a Stokes local system $(\Lc,F)$
is the sheaf of automorphisms of $\Lc$ preserving the filtration and inducing
the identity on $\gr^F \Lc$.

\end{itemize}
\end{Defi}

Thus $\Uc(\Lc)$ is a constructible sheaf of (non-abelian, in general) unipotent groups on $S^1$.
In particular, any $\Ac$-graded local system $L$ is a Stokes local system
via the split filtration, and we can speak about the sheaf $\Uc(L)$.
The following is then straightforward.

\begin{prop}
Let $L$ be an $\Ac$-graded local system on $S^1$.
The following   are in natural bijection:
\begin{itemize}
\item[(i)] The set of isomorphism classes of  Stokes local systems $(\Lc,F)$ together with
an isomorphism $\gr^F\Lc\to L $ of $\Ac$-graded local systems.

\item[(ii)] The first non-abelian cohomology set $H^1(S^1, \Uc(L))$.
\qed
\end{itemize}

\end{prop}

Given a convergent meromorphic connection $(E,\nabla)$ near $\oo$,
we have the local system $\Lc=E^\nabla$ of flat sections on $S^1$. That is,
$E^\nabla(U)$,  for an open $U\subset S^1$, consists of flat  holomorphic sections
of $E$ near the infinity of $C_U$. The sheaf $E^\nabla$ comes with a Stokes structure
given by the filtrations
\[
F_\alpha E^\nabla (U) \,=\,\biggl\{ s\in E^\nabla (U)\biggl|\, \, \exp\bigl(-\int \alpha)\cdot s \text
{ has at most polynomial growth at the infiity of } C_U\biggr\},
\]
where $\alpha\in \Ac(U)$. The  ``irregular Riemann-Hilbert  correspondence''
due to Deligne and Malgrange, is as follows.

\begin{thm}\label{thm:irr-rh1}
The above construction gives an equivalence of categories $\on{RH}$ fitting into
a commutative diagram whose lower row is also an equivalence:
\[
\xymatrix{
\text{Convergent meromorphic connections near } \oo
\ar[d]_{\text{Formalization}}
\ar@{<->}[r]^{\hskip 1.4cm \on{RH}}
&  \text{Stokes local systems }
\text{ on } S^1
\ar[d]^{\on{gr}}
\\
\text{Formal meromorphic connections near } \oo
\ar@{<->} [r]^{\hskip 1cm  \widehat{ \on{RH}}} & \Ac\text{-graded local systems }
\text{ on } S^1
}
\]
If $\Lc$ is a filtered local system corresponding, under $\on{RH}$,
to a connection $(E,\nabla)$, then $\Uc(\Lc)$ is identified
with the sheaf of flat sections of $\ul\End(E)$ whose asymptotic expansion at $\oo$ is identically $1$. \qed
\end{thm}

Let $\Gc\in\Perv(\CC)$ be a perverse sheaf of $\CC$-vector spaces (with some
finite set of singularities). The behavior of  $\Gc$ near $\oo$ gives a local system
on $S^1_\oo=S^1$ which we denote $\Psi_\oo(\Gc)$.

Theorem \ref{thm:irr-rh1} is complemented, for $\Dc$-modules, as follows,
see \cite{sabbah}

\begin{thm}\label{thm:Irr-RH-Perv}
The following categories are equivalent:
\begin{itemize}
\item[(i)] Holonomic $\Dc_\CC$-modules $\Nc$ which are regular everywhere except,
possibly, $\oo$.

\item[(ii)] Data $(\Gc, F)$, where $\Gc\in\Perv(\CC)$ is a perverse sheaf and $F$ is
a Stokes structure on the local system $\bPsi_\oo(\Gc)$. \qed
\end{itemize}

\end{thm}

\paragraph{Exponential Stokes structures. }
Let $A=\{w_1,\cdots, w_N\}\subset \CC_w$.  For simplicity we assume that $A$ is strong linearly
general position, see Definition \ref{def:conv-pos}(c). One can weaken this assumption by more careful reasoning.

\vskip .2cm

The constant sheaf $\ul A$ on $S^1$ is embedded into $\Ac$
as  consisting of constant forms $\alpha_i = w_i dz$. Thus $\exp\int \alpha_i = e^{w_i z}$
is a simple exponential.

\begin{figure}[H]
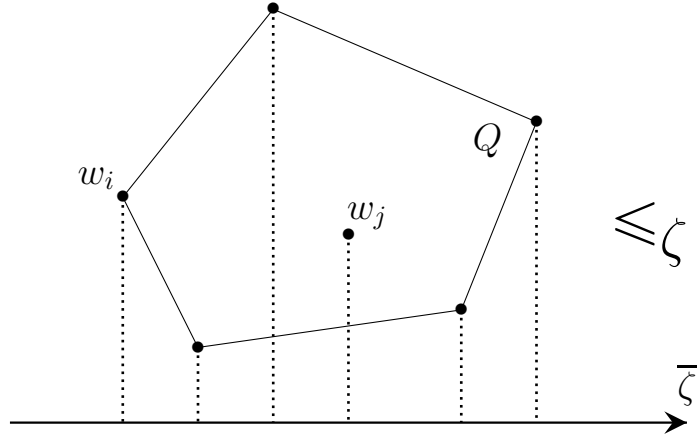

\centering
\btp[scale=0.5]

\node at (0,-1){$\bullet$};
\node at (5,2){$\bullet$};
\node at (-2,5){$\bullet$};
\node at (-6,0){$\bullet$};
\node at (-4, -4){$\bullet$};
\node at (3,-3){$\bullet$};

\draw (5,2) -- (-2,5) -- (-6,0) -- (-4,-4) -- (3, -3) -- (5,2);

\draw   [decoration={markings,mark=at position 1 with
{\arrow[scale=2,>=stealth]{>}}},postaction={decorate},
line width=.3mm] (-9, -6) -- (9, -6);

\node at (9, -5){\large$\ol\zeta$};

\draw [dotted, line width=1] (5,2) -- (5, -6);
\draw [dotted, line width=1] (-2,5) -- (-2, -6);
\draw [dotted, line width=1] (-6,0) -- (-6, -6);
\draw [dotted, line width=1] (-4,-4) -- (-4, -6);
\draw [dotted, line width=1] (3, -3) -- (3, -6);
\draw [dotted, line width=1] (0,-1) -- (0, -6);

\node at (8, -1){\huge$\leq_\zeta$};

\node at (3.7,1.5){\large$Q$};

\node at (-6.7, 0.5){\large $w_i$};

\node at (0.5, -0.5){\large $w_j$};
\etp

\caption{Convex geometry of exponential Stokes structures.}\label{fig:exp-geo}
\end{figure}

While the sheaf $\ul A$ is constant, the partial order  $\leq_\zeta$ on it induced from $\Ac$,
changes with the direction of $\zeta\in S^1$. Explcitly,
\be\label{eq:leq-zeta}
w_i\leq_\zeta w_j \quad \Leftrightarrow \quad \Re(\zeta w_i) \leq \Re(\zeta w_j)
\ee
(this means that $e^{(w_i-w_j)z} $ does not grow exponentially on $\zeta\cdot\RR_+$).
In other words,  after the orthogonal projection to the line $\ol\zeta\cdot \RR$
(note the complex conjugation!)
the image of $w_i$ precedes the image of $w_j$ in the orientation given by
$\zeta\cdot\RR_+$, see Fig. \ref{fig:exp-geo}.

\begin{figure}[H]
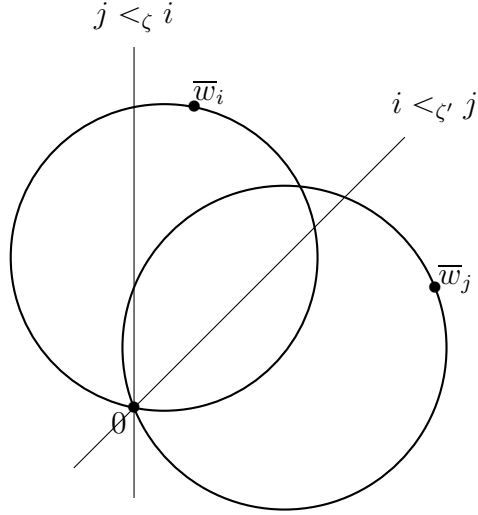

\centering
\btp[scale=0.4]

\node at (0,0){$\bullet$};
\node at (2,10){$\bullet$};
\node at (10,4){$\bullet$};

\draw [line width = 0.8] (1,5) circle (5.099);
\draw  [line width = 0.8] (5,2) circle (5.38516);


\node at (-0.5,-0.5){$0$};

\node at (2.5, 10.5){$\ol w_i$};
\node at (10.7, 4.3){$\ol w_j$};

\draw (0,-3) -- (0,0) -- (0,12);
\draw (-2,-2) -- (0,0) -- (9,9);

\node at (0,13){$j <_\zeta i $};
\node at (10,10){$i <_{\zeta'} j $};
\etp
\caption{The Lissajoux figure for an exponential Stokes structure is a union of circles. }
\label{fig:lissajoux}
\end{figure}

 \begin{rem}
 It is convenient, following Kontsevich,  to depict a local system of finite posets on $S^1$ by a ``Lissajoux figure'',
 i.e., a curve $L$ encircling $0$ so that cutting $L$ with an oriented line $\RR\zeta$ gives the order
 $\leq_\zeta$ on the finite set $L\cap \RR\zeta$ (the stalk of the local system at $\zeta$).
 The above implies the following  version of this description in our case. Note that the locus of orthogonal projections
 of $w_i$ to all the lines through $0$ is $C_i$, the circle with diameter $[0,w_i]$.
 Thus $i<_\zeta j$, if $\RR\ol\zeta$ meets $C_i$ before $C_j$ or, equivalently,  if $\RR\zeta$ meets
 $\ol C_i$ before $\ol C_j$.
 This means that $L$ in our case can be taken as the union of the circles $\ol C_i$ with diameters
  $[0,\ol w_i]$,
 see Fig. \ref{fig:lissajoux}.
 \end{rem}

Clearly, the data of all the $\leq_\zeta$ on $A$   contains the   information about
the convex geometry
of   $A$, in particular, about the polygon $Q=\Conv(A)$.  To illustrate this, let us
assume for simplicity
hat $A$ is in linearly general position and recover, from the data of the $\leq_\zeta$,
the uniform oriented matroid structure on $A$, i.e., the orientation datum
$(\eps_{ijk}) $, see
Definition \ref{def:or-datum}.

\begin{prop}\label{prop:exp-stokes-OM}
Let $i,j,k\in \{1,\cdots, N\}$  and let us choose $\zeta_0$ such that $w_i <_{\zeta_0}  w_j <
_{\zeta_0}  w_k$.
After we rotate $\zeta$  clockwise from $\zeta_0$ to $(-\zeta_0)$, the order will change
to the opposite, so we have three changes of order.
\begin {itemize}
\item[(a)] If the sequence of orders is
$(ijk) \to (ikj) \to (kij)\to (kji)$, then $\eps_{ijk}=+1$, i.e., the triangle $[w_i, w_j, w_k]$ is
positively oriented.

\item[(b)] If the sequence is $(ijk)\to (jik)\to (jki)\to (kji)$, then $\eps_{ijk}=-1$, i.e., the triangle
is negatively oriented. \qed
\end{itemize}
\end{prop}

More generally,
the order $\leq_\zeta$
switches  at $N(N-1)$ values  $\zeta=\sqrt{-1}\cdot\ol  \zeta_{ij}$
where $\zeta_{ij}$ are the slopes  defined  in \eqref{eq:zeta-ij}. These slopes are distinct by our assumption on $A$.
We will call the $\sqrt{-1}\cdot\ol \zeta_{ij}$ the {\em anti-Stokes directions}.

Passing through each $\sqrt{-1}\cdot \ol \zeta_{ij}$ switches the order of just two elements,
$w_i$ and $w_j$, leaving all the other order relations intact.
If  $\zeta$ is   not anti-Stokes, then
$\leq_\zeta$ is a total order on $A$.

Starting from a non-anti-Stokes $\zeta_0$ and going counterclockwise to $(-\zeta_0)$, we reach the opposite
order by
a sequence of $N(N-1)/2$ elementary switches, i.e., we get  a reduced decomposition of
the maximal permutation $s_\max\in S_N$ into elementary tranpositions:
\[
s_\max \,=\, s_{i_1} \cdots s_{i_l}, \quad l=N(N-1)/2, \,\, s_i = (i, i+1), \, \, i=1,\cdots, N-1.
\]

\begin{rem}\label{rem:manin-schechtman}
Not every reduced decomposition $R$ of $s_\max$ comes from some $A$.  In fact,    $R$
 gives rise to a uniform oriented matroid structure
(Definition \ref{def:or-matr}) on $\{1,\cdots, N\}$. The correspondence goes through
 $B(N,2)$, the second higher Bruhat order of Manin-Schechtman \cite{manin-schechtman}.
 More precisely, a datum of $R$, in the form of  a sequence of  switches as above,
  is  the same an admissible, in the sense of  \cite[Def. 2.2]{manin-schechtman}
 order on the set of  $N\choose 2$ pairs $1\leq i < j\leq N$. The set $B(N,2)$ is the quotient of the set of
 such orders by a natural equivalence relation (moves between reduced decompositions by
 interchanging neighboring pairs of  elementary transpositions that commute).   So $R$
 gives an element of $B(N,2)$ and all elements of $B(N,2)$ appear from some $R$.
  Further, any element of $B(N,2)$ gives a uniform oriented matroid  of rank $3$
 and each such matroid, after renumbering $\{1,\cdots, N\}$,
 appears from $B(N,2)$.  This is shown in \cite{KV1, KV2, ziegler}.

 Now, decompositions that  come  from some $A$,
give,  in this process,  realizable oriented matroids while most oriented matroids are not realizable,
see Remark \ref{rem:or-matr-real}(b).
\end{rem}

\begin{Defi}
 An {\em $A$-Stokes structure}
on a local system $\Lc$ on $S^1$ is a Stokes structure $F$
such that $\gr^F_\alpha\Lc=0$ unless $\alpha = wdz$ is a constant differential form and $w\in A$.
An {\em exponential Stokes structure} is an $A$-Stokes structure for some $A\subset\CC$.
\end{Defi}

Theorem \ref{thm:Irr-RH-Perv} together with Proposition \ref {prop:wc-M-top-1},
specialize to the following.

\begin{cor}\label{cor:stokes-ft-classify}
 Passing from a perverse sheaf to a holonomic regular $\Dc_{\CC_w}$-module and then to the Fourier
transform, establishes an equivalence between the following categories:
\begin{itemize}
\item[(i)] $\Perv(\CC_w, A)$.

\item[(ii)] The category of pairs $(\Gc,F)$, where $\Gc\in\Perv(\CC_z, 0)$ and $F$
is an $A$-Stokes structure on $\bPsi_\oo(\Gc)$. \qed
\end{itemize}
\end{cor}

\paragraph{The Stokes structure on the Fourier transform $\wc\Fc$.} \label{par:stokes-str-fourier}
We now specialize to  our original situation, that is:

\vskip .2cm

$A=\{w_1, \cdots, w_N\}\subset\CC_w$: a finite set, assumed  in strong linearly general position;

$\Fc\in\Perv(\CC,A)$: a perverse sheaf;

$\bPhi_i(\Fc)$: the local system on $S^1_{w_i}(\CC_w) \simeq S^1$  formed by
vanishing cycles of $\Fc$ near $w_i$;

$\Mc$: the holonomic regular $\Dc_{\CC_w}$-module with $\Sol(\Mc)=\Fc$;

$\wc\Mc$:  the Fourier transform of $\Mc$;

$\wc\Fc = \Sol(\wc\Mc)\in\Perv(\CC_z, 0)$ (Proposition \ref {prop:wc-M-top-1},);

$F$: the $A$-Stokes structure on the local system $\bPsi_\oo(\wc\Fc)$
classifying   $\wc\Mc$
(Corollary \ref{cor:stokes-ft-classify}).

\vskip .2cm

\noindent Our goal is to describe $F$ explicitly.
We  will call {\em Stokes directions} the conjugates $\ol \zeta_{ij}$, and
 {\em Stokes intervals} the $N(N-1)$ open arcs into which the Stokes directions
cut $S^1$. Thus the anti-Stokes directions   $\sqrt{-1}\cdot \ol\zeta_{ij}$ are the $90^\circ$ rotations of the Stokes
directions.
 We further assume for simplicity that no Stokes direction is anti-Stokes and vice versa
(this assumption  can  also be eliminated by more careful reasoning). Recall the union $K(\zeta)$
of the half-lines $K_i(\zeta)$
from \eqref{eq:half-plane(zeta)}. In terms of it,
  the role of Stokes vs.  anti-Stokes directions is as follows:

\begin{itemize}
\item At Stokes directions, the topology of $K(\zeta)$ (with respect to $A$) changes.

\item An anti-Stokes directions, the order $\leq_\zeta$ of dominance of exponentials (but not the topology)
changes.
\end{itemize}

\noindent We further recall the maps $g_\zeta$ from \eqref{eq:g-zeta-def}, which,
by Proposition \ref{prop:FT=integrals}, are given by
  Fourier integrals \eqref{eq:fourier-integral}. Recall also the identifications $\sigma_i$ from \eqref {eq:sigma-i}.

\begin{prop}\label{prop:g-zeta-varies}
(a) When $\zeta$ varies within a Stokes interval,  $g_\zeta$ is covariantly constant with respect
to the local system structures on the source and target, and
\[
F_i \bPsi_\oo(\wc\Fc)_\zeta \,=\, g_\zeta \biggl(\bigoplus_{w_j\leq_\zeta w_i} \bPhi_j(\Fc)_{-\ol\zeta}\biggr)
\]
becomes the split filtration associated with the grading in the RHS.

\vskip .2cm

(b) The $g_\zeta$ for non-Stokes $\zeta$ glue together into an identification of  $A$-graded
local systems
\[
\ol g: \, \bigoplus_{i=1}^N  \, \sigma_i^* \bPhi_i(\Fc)  \buildrel\sim\over\lra \gr^F \bPsi_\oo(\wc\Fc).
\]

\end{prop}

\noindent {\sl Proof:} The ``covariantly constant'' part in (a) is clear, as the topology of
the picture with  $K(\zeta)$ as well as $A$, does not change
within a Stokes interval.

\vskip .2cm

The identification of $F_i$ in (a) follows from    Proposition \ref{prop:FT=integrals}
and  the fact that  the leading term in the growth of $g_{i,\zeta} (\phi)(z)$ along $\zeta\cdot\RR_+$
is $\on{const}\cdot e^{w_iz}$, the contribution  to the integral from near the origin of the ray $K_i(\zeta)$.

\vskip .2cm

\begin{figure}[H]
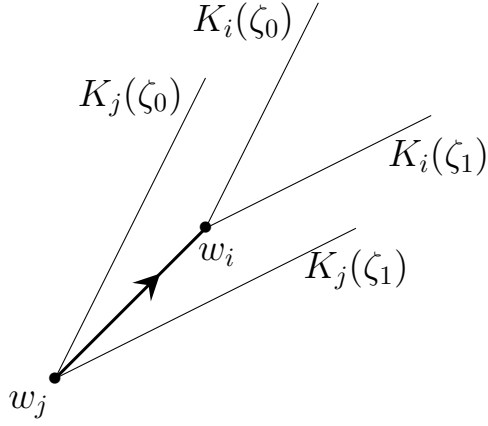

\centering
\btp[scale=0.5]
\node at (0,0){$\bullet$};
\draw (0,0) -- (8,4);
\draw (0,0) -- (4,8);
\node at (4,4){$\bullet$};
\draw (4,4)-- (10,7);
\draw (4,4)-- (7, 10);
\draw  [decoration={markings,mark=at position 0.7 with
{\arrow[scale=2,>=stealth]{>}}},postaction={decorate},
line width=.4mm] (0,0) -- (4,4);

\node at (-0.7,-0.7) {\large$w_j$};
\node at (4.3, 3.3){\large$w_i$};

\node at (8,2.9){\large$K_j(\zeta_1)$};
\node at (10.2, 5.9){\large$K_i(\zeta_1)$};
\node at (2,7.5){\large$K_j(\zeta_0)$};
\node at (5, 9.5){\large$K_i(\zeta_0)$};
\etp
\caption {Crossing a Stokes direction.}\label{fig:cross-stokes}
\end{figure}

To show (b), consider $\zeta$ crossing a single Stokes direction $\ol\zeta_{ij}$ counterclockwise,
starting from $\zeta=\zeta_0$ just before the crossing  and arriving at $\zeta=\zeta_1$ just after.
Then the $K_\nu(\zeta)$, $\nu=1,\cdots, N$, rotate clockwise, see Fig. \ref{fig:cross-stokes}.

For each $\nu$ let us identify the $\bPhi_\nu(\Fc)_{-\ol\zeta}$, $\zeta\in [\zeta_0, \zeta_1]$ with
each other by monodromy and denote the resulting single space by $\Phi_i$.
Now, for $\nu \neq j$, the integration contour $K_\nu(\zeta_2)$ is isotopic (relatively to $A$
and to the allowable ways to approach $\oo$) to
$K_\nu(\zeta_0)$. Therefore for $\phi\in\Phi_\nu$ the Fourier integral $g_{\nu,\zeta_1}(\phi)(z)$
is the analytic continuation of $g_{\nu, \zeta_0}(\phi)(z)$. We simply write
\[
g_{\nu,\zeta_1}(\phi)(z)\,=\, g_{\nu,\zeta_0}(\phi)(z), \quad \nu\neq j
\]
(as analytic functions in a sector near $\zeta\cdot\RR_+$).

But $K_j(\zeta_1)$ is not isotopic (rel. $A$ and the allowable infinity) to $K_j(\zeta_0)$: it differs by crossing $w_i$.
Therefore, instead of the equality, we have the Picard-Lefschetz relation
\be\label{eq:stokes-PL}
g_{j, \zeta_1}(\phi)(z) \,=\, g_{j, \zeta_0}(\phi)(z) -  g_{i, \zeta_1}\bigl(m_{ji}(\phi)\bigr)(z), \quad \phi\in\Phi_j,
\ee
where $m_{ji}: \Phi_j\to\Phi_i$ is the rectilinear transport.

It remains to notice that for $\zeta\in[\zeta_0, \zeta_1]$ we have $w_i \leq_\zeta w_j$, so the
second summand in the RHS of the above equality lies in the lower level of  filtration and
so the filtration remains unchanged. \qed


\subsection{The baby  Algebra of the Infrared}\label{subsec:baby-IA}

\paragraph{Meaning of the term.} By   {\em baby Algebra of the  Infrared } we mean a class of   identities
for  perverse sheaves on $\CC$ that express the transport or variation map along a curved path as a
(possibly alternating)
alternating sum of compositions of rectilinear transport maps corresponding to certain convex polygonal paths.
These identities are elementary and are obtained by iterated application of the Picard-Lefschetz identities
(Proposition \ref{prop:PL-form}). But the sums or alternating sums that appear, can be seen as de-categorifications
(baby versions) of filtrations or complexes\footnote
{There is no essential difference between the two concepts.  A filtration can  be
seen as a complex in derived category  whose terms are  shifted quotients of the filtration.
}
  whose summands are labelled by such convex polygonal paths. Such complexes constitute an
essential part of the Algebra of the Infrared of \cite{GMW1}.

In this section we collect several examples starting from the simplest one which establishes the pattern.

\begin{figure}[H]
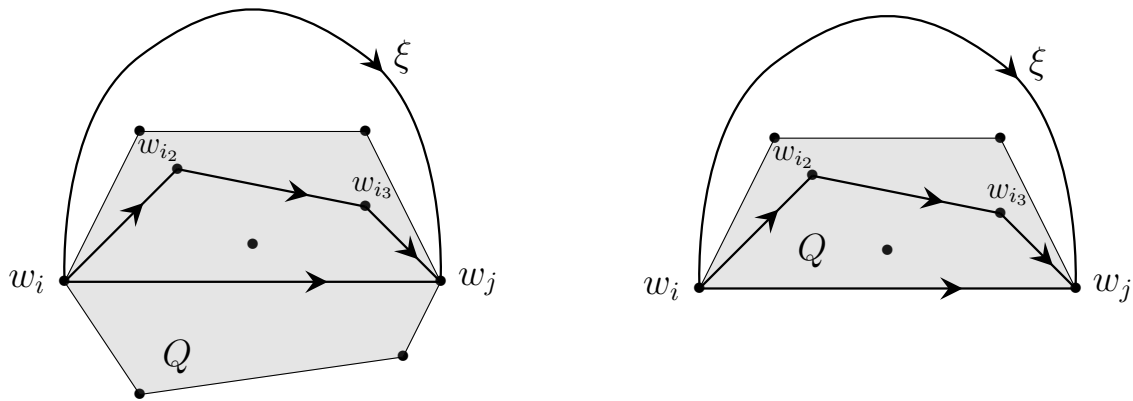

\centering
\btp[scale=0.5]

\node at (0,0){\small$\bullet$};
\node (wi) at (-5,-1){\small$\bullet$};
\node (wi2) at (-2, 2){\small$\bullet$};
\node (wi3) at (3,1){\small$\bullet$};
\node (wj) at (5,-1){\small$\bullet$};

\node at (-3,3){\small$\bullet$};
\node at (3,3){\small$\bullet$};
\node at (4, -3){\small$\bullet$};
\node at (-3,-4){\small$\bullet$};

\filldraw [color=gray, opacity=0.2] (-5, -1) -- (-3,3) --(3,3) -- (5, -1) -- (4, -3) -- (-3, -4) -- (-5, -1);

\draw  (-5, -1) -- (-3,3) --(3,3) -- (5, -1) -- (4, -3) -- (-3, -4) -- (-5, -1);

\draw  [decoration={markings,mark=at position 0.7 with
{\arrow[scale=2,>=stealth]{>}}},postaction={decorate},
line width=.3mm] (-5, -1) -- (5, -1);

\draw  [decoration={markings,mark=at position 0.7 with
{\arrow[scale=2,>=stealth]{>}}},postaction={decorate},
line width=.3mm] (-5, -1) -- (-2,2);

\draw  [decoration={markings,mark=at position 0.7 with
{\arrow[scale=2,>=stealth]{>}}},postaction={decorate},
line width=.3mm] (-2,2) -- (3,1);

\draw  [decoration={markings,mark=at position 0.7 with
{\arrow[scale=2,>=stealth]{>}}},postaction={decorate},
line width=.3mm] (3,1) -- (5, -1);

\draw  [decoration={markings,mark=at position 0.7 with
{\arrow[scale=2,>=stealth]{>}}},postaction={decorate},
line width=.3mm]   plot  [ smooth, tension=1] coordinates{   (-5, -1)    (-3, 5)  (3, 5)    (5, -1)   } ;

\node at (-2, -3){\large$Q$};

\node at (-6, -1){\large$w_i$};
\node at (6, -1){\large$w_j$};
\node at (4, 5) {\large$\xi$};

\node at (-2.5, 2.5){$w_{i_2}$};
\node at (3.2, 1.5){$w_{i_3}$};

\etp
\quad\quad\quad\quad
\btp[scale=0.5]

\node at (0,0){\small$\bullet$};
\node (wi) at (-5,-1){\small$\bullet$};
\node (wi2) at (-2, 2){\small$\bullet$};
\node (wi3) at (3,1){\small$\bullet$};
\node (wj) at (5,-1){\small$\bullet$};

\node at (-3,3){\small$\bullet$};
\node at (3,3){\small$\bullet$};
\node at (4, -3){};
\node at (-3,-4){};

\filldraw [color=gray, opacity=0.2] (-5, -1) -- (-3,3) --(3,3) -- (5, -1);

\draw  (-5, -1) -- (-3,3) --(3,3) -- (5, -1);


\draw  [decoration={markings,mark=at position 0.7 with
{\arrow[scale=2,>=stealth]{>}}},postaction={decorate},
line width=.3mm] (-5, -1) -- (5, -1);

\draw  [decoration={markings,mark=at position 0.7 with
{\arrow[scale=2,>=stealth]{>}}},postaction={decorate},
line width=.3mm] (-5, -1) -- (-2,2);

\draw  [decoration={markings,mark=at position 0.7 with
{\arrow[scale=2,>=stealth]{>}}},postaction={decorate},
line width=.3mm] (-2,2) -- (3,1);

\draw  [decoration={markings,mark=at position 0.7 with
{\arrow[scale=2,>=stealth]{>}}},postaction={decorate},
line width=.3mm] (3,1) -- (5, -1);

\draw  [decoration={markings,mark=at position 0.7 with
{\arrow[scale=2,>=stealth]{>}}},postaction={decorate},
line width=.3mm]   plot  [ smooth, tension=1] coordinates{   (-5, -1)    (-3, 5)  (3, 5)    (5, -1)   } ;

\node at (-2, 0){\large$Q$};
\node at (-2, -3){};

\node at (-6, -1){\large$w_i$};
\node at (6, -1){\large$w_j$};
\node at (4, 5) {\large$\xi$};

\node at (-2.5, 2.5){$w_{i_2}$};
\node at (3.2, 1.5){$w_{i_3}$};

\etp
\caption{Circumnavigation vs.  inland  convex paths.} \label{fig:circum}
\end{figure}

\paragraph{ Circumnavigation around a polygon.  } Let $A=\{w_1,\cdots, w_N\}\subset\CC$ be in linearly general position,
and $Q=\Conv(A)$. Note that not every $w_k$ is a vertex of $Q$, some of the $w_k$ may lie inside.
Let $w_i, w_j$ be two vertices of $Q$ and $\xi$ be a path going around $Q$ (say, clockwise)
and joining  $w_i$ and $w_j$, see the left part of Fig. \ref{fig:circum}.

Without loss of generality, we can assume that $[w_i, w_j]$ is actually an edge of $Q$;
if not, we can replace $Q$ with one
of the two parts on which $[w_i, w_j]$ dissects it and replace $A$ by its intersection
with that part, see the right part of  Fig. \ref{fig:circum}.

Suppose $\Fc\in\Perv(\CC,A)$ is given, so we have the local system of
vanishing cycles $\bPhi_k(\Fc)$ on $S^1_{w_k}$ for each $k$.
Let us identify, using the monodromy,  the stalks $\bPhi_i(\Fc)_\zeta$ for all directions $\zeta\in S^1_{w_i}$
between  $\zeta_{ji}$ (the direction towards $w_j$) and the tangent direction to $\gamma$ and denote the resulting
single space $\Phi_i$. Let us define $\Phi_j$ similarly on the other end. Then we have
the transport map $m_{ij}(\xi): \Phi_i\to\Phi_j$.

\vskip .2cm

In this section we will consider polygonal paths
\be\label{eq:polyg-path}
\gamma \,=\, [w_{1_1}, \cdots, w_{i_n}] \, := \,
 [w_{i_1}, w_{i_2}] \cup [w_{i_2}, w_{i_3}]\cup\cdots\cup [w_{i_{p-1}}, w_{i_p}].
\ee
with vertices in $A$.
Let $\Lambda(i,j)$ be the set of such  paths
joining $w_i$ and $w_j$  (i.e., with $i_1=i, i_p=j$) and satisfying the following {\em convexity property}:

\begin{itemize}
\item   $\gamma\cup [w_i, w_j]$ is the boundary of a convex polygon.
\end{itemize}
See  Fig. \ref{fig:circum}. For each $\gamma\in \Lambda(i,j)$ we define the {\em iterated transport map}
$m(\gamma): \Phi_i\to\Phi_j$ as the composition
\be\label{eq:iter-trans-2}
m(\gamma) \,=\, m_{i_{p-1}, j}\circ \cdots\circ m_{i_2, i_3} \circ m_{i,i_2}
\ee
where $m_{kl}:\Phi_k\to\Phi_l$ is the rectilinear transport, and the identification of the stalks of the intermediate local
systems  $\bPhi$ is done  using clockwise monodromies,  per Convention \ref{conv:clock}.

\begin{prop}\label{prop:baby-IA}
We have the identity
\[
m_{ij}(\xi) \,=\,\ \sum_{\gamma\in \Lambda (i,j)} m(\gamma).
\]
\end{prop}

\begin{figure}[H]
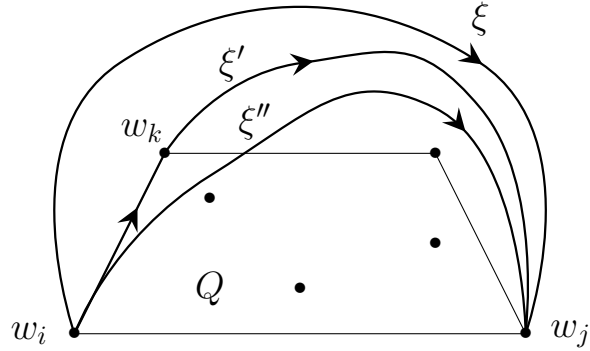

\centering
\btp[scale=0.6]

\node at (0,0){\small$\bullet$};
\node (wi) at (-5,-1){\small$\bullet$};
\node (wi2) at (-2, 2){\small$\bullet$};
\node (wi3) at (3,1){\small$\bullet$};
\node (wj) at (5,-1){\small$\bullet$};

\node at (-3,3){\small$\bullet$};
\node at (3,3){\small$\bullet$};

\draw  (-5, -1) -- (-3,3) --(3,3) -- (5, -1) -- (-5, -1);

\draw  [decoration={markings,mark=at position 0.7 with
{\arrow[scale=2,>=stealth]{>}}},postaction={decorate},
line width=.3mm]   plot  [ smooth, tension=1] coordinates{   (-5, -1)    (-4, 5)  (4,  5)    (5, -1)   } ;

  \draw  [decoration={markings,mark=at position 0.3 with
{\arrow[scale=2,>=stealth]{>}}},postaction={decorate},
line width=.3mm]   plot  [ smooth, tension=1] coordinates{   (-3,3)    (0, 5)  (4,4)    (5, -1)   } ;

  \draw  [decoration={markings,mark=at position 0.7 with
{\arrow[scale=2,>=stealth]{>}}},postaction={decorate},
line width=.3mm]  (-5, -1) -- (-3,3);

  \draw  [decoration={markings,mark=at position 0.7 with
{\arrow[scale=2,>=stealth]{>}}},postaction={decorate},
line width=.3mm]   plot  [ smooth, tension=1] coordinates{   (-5, -1)    (-2,2.5)  (3,  4)    (5, -1)   } ;

\node at (-2, 0){\large$Q$};

\node at (-6, -1){\large$w_i$};
\node at (6, -1){\large$w_j$};
\node at (4, 6) {\large$\xi$};

\node at (-3.5, 3.5){\large$w_k$};

\node at (-1.5,5){\large$\xi'$};

\node at (-1, 3.7){\large$\xi''$};

\etp
\caption{ Tightening the cord.} \label{fig:cord}
\end{figure}

\noindent{\sl Proof:} We view the path $\xi$ as a cord and  proceed by ``tightening the cord'',
straightening it first so that it contains the segment $[w_i, w_k]$, where $w_k$ is the next after $w_i$ vertex of $Q$
on the way to $w_j$, see
Fig. \ref{fig:cord}. More precisely, we get a composite path consisting of $[w_i, w_k]$ and a curved path
$\xi'$ going around $Q$ and connecting $w_k$ with $w_j$.
Moving the cord further to cross the vertex $w_k$ and no other vertex, we get a path $\xi''$.

\noindent The Picard-Lefschetz identity (Proposition  \ref{prop:PL-form}) gives:
\[
m_{ij}(\xi) \,=\, m_{ij}(\xi'') + m_{kj}(\xi') m_{ik}.
\]
Then we proceed to tighten similarly the paths $\xi'$ and $\xi''$,
keeping the straight part  $[w_i, w_k]$  intact (i.e., viewing $w_k$
as a fixed ``eye''  through which we pass the rope). Continuing like this,
we represent  $m_{ij}(\xi)$  as the sum of
the effects of  polygonal paths consisting entirely of straight line segments.
These polygonal paths will have precisely the convex nature claimed. \qed

\paragraph{Counterclockwise circumnavigation: alternating sums.}
 For future reference, we note
a version of Proposition \ref{prop:baby-IA} which concerns the path $\xi$ run
  counterclockwise, i.e., in our  general notation,   the path $\xi^{-1}$
running  from $w_j$ to $w_i$.

Let $\Lambda(j,i)$ be the set of the same polygonal convex paths as in $\Lambda(i,j)$ but
oriented from $w_j$ to $w_i$. As before, each $\gamma\in \Lambda(j,i)$ gives
the iterated transport map $m(\gamma):\Phi_j\to\Phi_i$. Here, too,  we identify the
stalks of the  intermediate local systems $\bPhi$
  using clockwise monodromies,  per Convention \ref{conv:clock}.

 \begin{prop}\label{prop:baby-IA-counter}
 We have an equality in $\Hom(\Phi_j, \Phi_i)$
 \[
 m_{ji}(\xi^{-1}) \,=\,\sum_{\gamma\in \Lambda(j,i)} (-1)^{l(\gamma)} m(\gamma).
 \]
Here $l(\gamma)= |\gamma\cup A|-2$ is the number of points from $A$ other than $w_i$ and $w_j$ lying on $\gamma$.

 \end{prop}

\noindent{\sl Proof:} The argument is similar to that of Proposition  \ref{prop:baby-IA}
using Picard-Lefschetz identities which in this case will lead to signs in the RHS due to difference
in orientation. For
example, the very first identity above will
  take the form
 \[
 m_{ji}(\xi^{-1}) \,=\, m_{ji}((\xi'')^{-1}) - m_{ki} m_{jk}((\xi')^{-1}).
 \]
Therefore, in the process of tightening the cord, the sign of a new summand will change each
 time we acquire a new element from $A$ on the path.
 Collecting the signs, we get the statement. \qed

\paragraph{The Stokes matrices of an exponential Stokes structure.}
Next, we want to give a "baby infrared'' formula for the Stokes matrices of the Fourier transform of a perverse sheaf.
As a first step, we
explain, following \cite{deligne-var},  the correspondence between ``abstract''  exponential Stokes {\em structures}
and  ``concrete''  Stokes {\em matrices}.

\vskip .2cm
Let $A=\{w_1,\cdots, w_N\}$. Assume that  $A$  is in strong
linearly general position. Note that an $A$-Stokes structure
on a local system $\Lc$ can be equivalently seen as a filtration indexed just by the sheaf of posets $(\ul A, \leq)$.
That is, we have globally defined subsheaves $F_i\Lc\subset \Lc$  (corresponding to $w_i\in A$)
such that, for a given
direction $\zeta$, the inequality
$w_i \leq_\zeta w_j$ implies that $F_i\Lc\subset F_j\Lc$ at tle level of stalks at $\zeta$.
Note that $F_i\Lc$ is not locally constant, in general.

\vskip .2cm

Further, an $\ul A$-grading on a local system $L$ is the same as an $A$-grading, i.e., a direct
sum decomposition $L=\bigoplus_{i=1}^N L_i$.

\begin{prop}\label{prop:H-of stokes-sheaf}
Let $L$ be an $A$-graded local system and $\Uc(L)$  be the corresponding Stokes sheaf.
Let $I\subset S^1$ is an open (connected)  arc.
\begin{itemize}
\item[(a)] If $I$   contains at least  $N(N-1)/2$ consecutive
anti-Stokes directions $\sqrt{-1}\zeta_{ij}$, then  we have $H^0(I, \Uc(L))=\{1\}$.

\item[(b)]  If $I$   contains at most $N(N-1)/2$ consecutive.
anti-Stokes directions, then  we have
$H^1(I, \Uc(L))=\{1\}$.

\end{itemize}
\end{prop}

\begin{figure}[H]
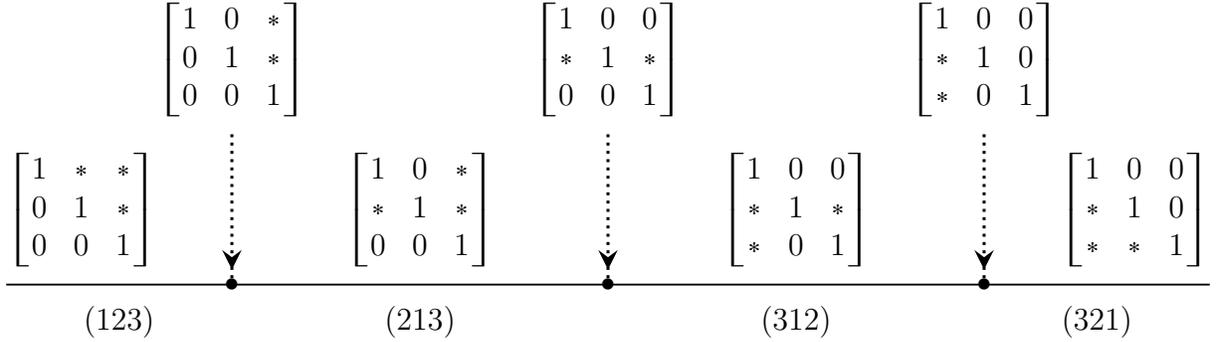

\centering
\btp

\draw[line width = 0.7]  (-8,0) -- (8,0);
\node at (-5,0){$\bullet$};
\node at (5,0){$\bullet$};
\node at (0,0){$\bullet$};

\node at (-6.5, -0.5){$(123)$};
  \node at (-2.5, -0.5){$(213)$};
   \node at (2.5, -0.5){$(312)$};
    \node at (6.5, -0.5){$(321)$};

\node at (-7, 1){
$
\begin{bmatrix}
1&*&*
\\0&1&*
\\0&0&1
\end{bmatrix}
$
};

\node at (-2.5, 1){
$
\begin{bmatrix}
1&0&*
\\
*&1&*
\\0&0&1
\end{bmatrix}
$
};

  \node at (2.5, 1){
$
\begin{bmatrix}
1&0&0
\\
*&1&*
\\ *&0&1
\end{bmatrix}
$
};

    \node at (7, 1){
$
\begin{bmatrix}
1&0&0
\\
*&1&0
\\ *&*&1
\end{bmatrix}
$
};

    \node at (-5,3){
$
\begin{bmatrix}
1&0&*
\\
0&1&*
\\ 0&0&1
\end{bmatrix}
$
};

      \node at (-0,3){
$
\begin{bmatrix}
1&0& 0
\\
*&1&*
\\ 0&0&1
\end{bmatrix}
$
};

        \node at (5,3){
$
\begin{bmatrix}
1&0& 0
\\
*&1&0
\\ *&0&1
\end{bmatrix}
$
};

\draw   [dotted, decoration={markings,mark=at position 0.9 with
{\arrow[scale=1.5,>=stealth]{>}}},postaction={decorate},
line width=.4mm] (5,2) -- (5,0);

\draw   [dotted, decoration={markings,mark=at position 0.9 with
{\arrow[scale=1.5,>=stealth]{>}}},postaction={decorate},
line width=.4mm] (0,2) -- (0,0);

\draw   [dotted, decoration={markings,mark=at position 0.9 with
{\arrow[scale=1.5,>=stealth]{>}}},postaction={decorate},
line width=.4mm] (-5,2) -- (-5,0);

\etp
\caption{A part of the Stokes sheaf in the matrix form for $N=3$.}\label{fig:stokes-sheaf}
\end{figure}

\noindent{\sl Proof:} We can trivialize each $L_i$ on $I$, that is, assume $L_i=\ul V_i$
where $V_i$ is some finite-dimensional vector space. Then $\Uc(L)$ is identified with
the subsheaf in the constant sheaf with stalk $\prod_{i,j} \Hom(V_i, V_j)$, consisting
of block-matrices with $1$'s on diagonals which are upper-triangular with respect to
switching orders, see Fig. \ref{fig:stokes-sheaf}. The orders switch at the anti-Stokes directions.

In the situation (a), since the orders corresponding to $\zeta$ and $(-\zeta)$ are opposite,
the only possible matrix which will be upper-triangular with respect to all the orders $\leq_\zeta$, $\zeta\in I$
and  with $1$'s on diagonal,
will be the identity, so $H^0=\{1\}$.

\vskip .2cm

To prove (b), we proceed by induction. Let us write $I$ as an open interval $(\theta_1, \theta_2)$
in the counterclockwise direction. Let $\phi\in I$ be the last anti-Stokes direction on the way from
$\theta_1$ to $\theta_2$.

Let $I' = (\theta_1, \phi)$  and $I_2=(\psi, \theta_2)$, where $\psi$ is the direction just before $\phi$.
Then $I'\cup I'' =I$ form a covering of $I$ and by the  inductive assumption $\Uc(L)$ has
no higher cohomology (i.e., no $H^1$) on $I', I''$ and $I'\cap I''$. Therefore
the first cohomology in question is identified with the set of double cosets
\be\label{eqref:double-coset}
H^1(I, \Uc(L)) \,=\, \Uc(L)(I')\biggl \backslash \Uc (L)(I'\cap I'') \biggr/ \Uc(L)(I'').
\ee
Since $I'\cap I''$ is a single arc which
does not contain anti-Stokes directions, $ \Uc (L)(I'\cap I'')$
consists of block-unitriangular matrices with respect to some total order $\leq_\zeta$.
The two subgroups on the left and right in \eqref {eqref:double-coset} consist of
matrices unitriangular  with respect to two partial orders $\leq', \leq ''$ of which
$\leq_\zeta$ is a minimal common refinement. So it is immediate to see that the
any matrix $g$ uni-triangular with respect to $\leq_\zeta$ is the product $g' g''$ of matrices
unitriangular with respect to $\leq'$ and $\leq''$ respectively. \qed

\vskip .2cm

Note that the cases (a) and (b) of Proposition \ref{prop:H-of stokes-sheaf} have nontrivial intersection:
when $I$ contains {\em exactly }  $N(N-1)/2$  anti-Stokes directions (so the width of $I$
is slightly greater than $180^\circ$). Such arcs will be called {\em balanced}. Thus, on
a balanced arc,  $\Uc(L)$ has neither $H^0$ nor $H^1$.

\vskip .2cm

\begin{figure}[H]
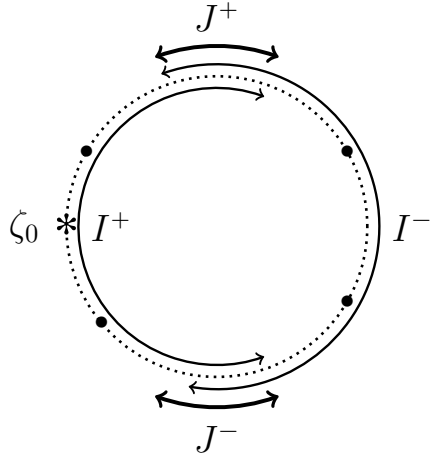

\centering
\btp[scale=0.4]

\draw [dotted, line width=1]  (0,0) circle (5cm);

\node at (30:5){$\bullet$};
\node at (150:5){$\bullet$};
\node at (220:5){$\bullet$};
\node at (330:5){$\bullet$};


\node (A) at (70:4.6){};
\draw [<->, line width=0.9]
(A.center) arc (70: 290: 4.6);

\node (B) at (260:5.4){};
\draw [<->, line width=0.9]
(B.center) arc (260: 470: 5.4) ;

\node (C) at (70:6){};
\draw [<->, line width=1.5]
(C.center) arc (70: 110: 6) ;

\node (D) at (250:6){};
\draw [<->, line width=1.5]
(D.center) arc (250: 290: 6) ;

\node at (-3.5,0){\large$I^+$};
\node at (6.5,0){\large$I^-$};

\node at (0,7){\large$J^+$};
\node at (0,-7){\large$J^-$};

\node at (-5,0){\huge$\ast$};
\node at (-6.4,0){\large$\zeta_0$};


\etp
\caption{A balanced covering. }\label{fig:balanced}
\end{figure}

Note further that we can form a covering of $S^1$ by two balanced arcs
(a {\em balanced covering}, for short) $S^1=I^+\cup I^-$ whose  intersection must consist of two
intervals: $I^+\cap I^- = J^+\sqcup J^-$, containing no anti-Stokes directions and giving
rise to opposite orders $\leq_\zeta$, see Fig. \ref{fig:balanced}.

\begin{cor}\label{cor:bal-cover}
For a balanced covering $S^1=I_+\cup I_-$ we have a canonical identification
\[
H^1(S^1, \Uc(L)) \,\simeq \, \Uc(L)(J^+) \times \Uc(L)(J^-).
\]
\end{cor}

\noindent{\sl Proof:}
Indeed, because on the $I^\pm$ the sheaf $\Uc(L)$ does not have $H^1$, we can write
$H^1(S^1, \Uc(L))$ as the double coset  space as in \eqref{eqref:double-coset}. We have
\[
\Uc(L)(I^+\cap I^-) = \Uc(L)(J^+) \times \Uc(L)(J^-),
\]
while the two subgroups
by which it is factorized, reduce to $\{1\}$. \qed

\begin{Defi}
Let $(\Lc,F)$ be an $A$-Stokes structure and $L=\gr^F(\Lc)$. Let $S^1=I^+\cup I^-$ be
a balanced covering.
The {\em Stokes matrices}
associated to $F$ (and the  covering)
are the sections $C_\pm \in \Uc(L)(J_{\pm})$ so that
$(C_+, C_-)$ corresponds to $F$ via  Corollary \ref{eqref:double-coset}.
\end{Defi}

More explicitly, Proposition \ref{prop:H-of stokes-sheaf} means that there exist unique splittings,
i.e., isomorphisms
\[
g_{\pm}: L|_{I^\pm} \lra \Lc|_{I^\pm}
\]
of $\ul A$-filtered sheaves. The Stokes matrices are defined as
\be
\begin{gathered}
C^+ \,=\, \bigl(g_+^{-1}|_{J^+}\bigr)  \circ \bigl(g_-|_{J^+}\bigr): L(J^+)\lra L(J^+),
\\
C^- \,=\, \bigl(g_+^{-1}|_{J^-}\bigr)  \circ \bigl(g_-|_{J^-}\bigr): L(J^-)\lra L(J^-).
\end{gathered}
\ee
Note that
each $\zeta\in J^\pm$ gives the same order $\leq_\zeta = \leq_\pm$ on $A$.
Thes $\leq_\pm$ are total orders, opposite to each
other. By construction, $C^\pm$ being sections of $\Uc(L)$, are unitriangular with respect to these orders,
i.e., if we write $C^\pm = \|C^\pm_{ij}\|$ with $C^\pm_{ij}: L_i \to L_j$ , then $C^\pm_{ij}\neq 0$ only for
$w_i \leq_\pm w_j$ and $C^\pm_{ii}=1$.

\vskip .2cm

Note further that $C^+$ and $C^-$ act in different spaces. To bring them together, let
\[
\begin{gathered}
T_{+/-}^L: L(J^+) \lra L(J^-), \quad T_{-/+}^L: L(J_-) \lra L(J^+), \\
T^L \,= \, T_{-/+} ^L\circ T_{+/-}^L:
L(J^+) \lra L(J^+)
 \end{gathered}
\]
be the counterclockwise monodromies for $L$. As $L$ is graded, they are represented by diagonal matrices.
Put
\[
\wt C^- \,=\, (T^L_{+/-})^{-1}\circ  C^- \circ T^L_{+/-} : L(J^+) \to L(J^+).
\]
Let also $T^\Lc: \Lc(J^+)\to\Lc(J^+)$ be the full counterclockwise monodromy of $\Lc$.
Using the  identification $g_+: L(J^+) \to \Lc(J^+)$, we can view $T^\Lc$ as an automorphism
(normalized monodromy)
\[
T^\Lc_{\on{norm}} = g_+^{-1} T^\Lc g_+: L(J^+)\lra L(J^+).
\]

\begin{prop}\label{prop:mono=2Stokes}
 We have an equality
 $
 T^\Lc_{\on{norm}} \, = \, C^+ T^L (\wt C^-)^{-1}
 $ which realizes $T^\Lc_{\on{norm}} $ as the product of an upper triangular, diagonal and lower triangular
 matrices.
\end{prop}

\noindent{\sl Proof:}
This is straightforward  by analyzing the commutative diagram
\[
\xymatrix{
\Lc(J^+) \ar[r]^{\res^{-1}}&\Lc(I^+)
\ar[rr]^\res && \Lc(J^-)
 \ar[r]^{\res^{-1}}& \Lc(I^-) \ar[r]^\res& \Lc(J^+) &
\\
L(J^+) \ar[u]^{g_+}
\ar[r]_{\res^{-1}}& L(I^+) \ar[r]_\res\ar[u]^{g_+}
& L(J^-) \ar[r]_{(C^-)^{-1}}\ar[ul]_{g_+}
& L(J^-) \ar[r]_{\res^{-1}}\ar[u]^{g_-}
& L(I^-) \ar[r]_\res\ar[u]^{g_-}
& L(J^+) \ar[r]_{C^+}\ar[u]^{g_-}
& L(J^+), \ar[ul]_{g_+}
}
\]
where $\res$ means the restriction isomorphisms for the local systems $\Lc$ and $L$.
The composition of the top row is $T^\Lc$.  The composition of the bottom row is
$C^+ T^L_{-/+} (C^-)^{-1} T^L_{+/-}$, and our statement follows from this. \qed

\paragraph{  $\zeta$-convex   paths.}
 Next, we introduce some terminology.
Let $\zeta\in\CC$ be  such that $|\zeta|=1$. For a subset $B\subset A $ we call the {\em $\zeta$-convex hull}
and denote $\Conv^\zeta(B)$ the set
\be\label{eq:zeta-conv}
\Conv^\zeta(B) \,=\, \Conv\biggl( \, \bigcup_{w\in B} (w+\zeta\cdot \RR_{\geq 0}) \biggr).
\ee
 Thus $\Conv^\zeta(B)$ is an unbounded polygon in $\CC$ with  (one or) two infinite edges in the
$\zeta$-direction.

\begin{figure}[H]
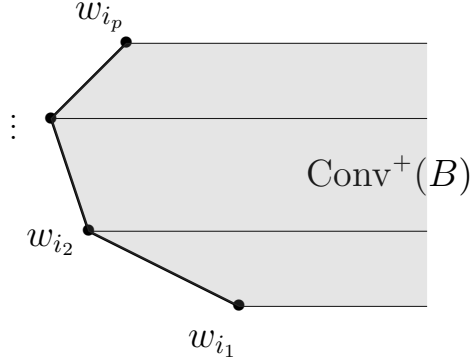

\centering
\btp[scale=0.5]
\node at (-4,3){$\bullet$};
\node at (-6,1){$\bullet$};
\node at (-5, -2){$\bullet$};
\node at (-1, -4){$\bullet$};
\draw (-4,3) -- (4,3);
\draw (-6,1) -- (4,1);
\draw (-5, -2) -- (4, -2);
\draw (-1,-4) -- (4, -4);

\draw [line width= 1] (-4,3) -- (-6,1) -- (-5,-2) -- (-1, -4);

\node at((-4.7, 3.7){\large$w_{i_p}$};
\node at((-7, 1){$\vdots$};
\node at((-6, -2.3){\large $w_{i_2}$};
\node at((-1.7, -5){\large$w_{i_1}$};

\node at (3, -0.5){\large$\Conv^+(B)$};
\filldraw [color=gray, opacity=0.2] (4,3) -- (-4,3) -- (-6,1) -- (-5, -2) -- (-1, -4)-- (4, -4);
\etp

\caption{A left convex path.}\label{fig:right-convex}
\end{figure}

A polygonal path $\gamma$ as in \eqref{eq:polyg-path}  will be called
  {\em $\zeta$-convex}, if
  $\gamma$ is the finite  part of the boundary of the polygon
  $\Conv^\zeta (\gamma)$.
 For $\zeta=1$ we will use the term {\em left convex} instead of ``$1$-convex''
and  for  $\zeta= -1$ we will use the term {\em right convex} instead of ``$(-1)$-convex''.
We also
write $\Conv^\pm$ for $\Conv^{\pm 1}$.
See Fig. \ref{fig:right-convex}.

\paragraph{The Stokes matrix of $\widecheck{\Fc}$.}
We now specialize  to the situation and notations of \S \ref{subsec:FT-stokes} \ref{par:stokes-str-fourier}
Let
\[
\Lc = \bPsi_\oo(\wc\Fc),\quad L=\gr^F(\Lc)=\bigoplus_{i=1}^N \sigma_i^* \bPhi_i(\Fc).
\]
Given a balanced covering $S^1=I^+\cup I^-$ as above, take a generic
(non-Stokes, nor anti-Stokes) point $\zeta_0$ in the middlle of $I^+$,
so that $\zeta_+=-\sqrt{-1}\cdot\zeta_0\in J_+$, see Fig. \ref{fig:balanced},
where $\zeta_0=-1$. Let $\leq_+ =\leq_{\zeta_+}$ and let us number $A=\{w_1,\cdots, w_N\}$
so that $w_1 >_+ \cdots >_+ w_N$. For the choice $\zeta_0=-1$ this means ordering
$w_i$ so that $\Im(w_1) < \cdots <  \Im(w_N)$, as in Fig. \ref{fig:horiz-cuts}.
The Stokes matrix $C^+$ describing the Stokes structure $F$ on $\Lc$, has
components
\[
C^+_{ij}: L_{i, \zeta_+} = \bPhi_i(\Fc)_{-\ol{\zeta}_+}\lra  \bPhi_j(\Fc)_{-\ol{\zeta}_+} =
L_{j, \zeta_+}, \quad i\geq j, \,\, C^+_{ii}=1.
\]
For each $i$ let us identify, by the quarter-circle monodromy,
\[
\Phi_i(\Fc)_{-\ol{\zeta}_+} \buildrel\sim\over\lra  \Phi_i := \Phi_i(\Fc)_{-\ol{\zeta}_0}.
\]
For the choice $\zeta_0=-1$ this corresponds to considering the stalks of the $\bPhi_i$
in the real positive direction, as in \S \ref{subsec:rec-approach} \ref{par:rect-trans}
So we have the rectilinear transport maps $m_{ij}:\Phi_i\to\Phi_j$, see
\eqref{eq:rect-trans-def}.

\vskip .2cm

For $i<j$ let $\Lambda(i,j)$ be the set of $(-\ol\zeta_0)$-convex  (left convex, for $\zeta_0=-1$)
paths starting from $w_i$ and ending
at $w_j$. For any $\gamma\in \Lambda (i,j)$
we have the iterated transport map $m(\gamma):\Phi_i\to\Phi_j$,
as in \eqref{eq:iter-trans-2}.

\begin{prop}\label{prop:stokes-infra}
For $i < j$ we have
\[
C^+_{ij} \,=\sum_{\gamma\in \Lambda(i,j)} m(\gamma).
\]

\end{prop}

In a similar way, $C^-_{ij}$ is interpreted as a sum over $(+\ol\zeta_0)$-convex
(right convex, for $\zeta_0=-1$) paths.

\begin{rem}
Another formula for the Stokes matrices of $\wc\Fc$ was given
in \cite {dagnolo-sabbah} (Th. 5.2.2 and Rem. 5.3.3) using a different approach
(enhanced perverse sheaves). This formula is essentially equivalent to
Proposition
\ref {prop:stokes-infra}.  The ``infrared'' shape of our  formula
is due to the rectilinear approach adopted.
\end{rem}

\noindent {\sl Proof of Proposition \ref{prop:stokes-infra}:}
Recall from \eqref {eq:g-zeta-def} and Proposition \ref  {prop:FT=integrals} the isomorphism
$g_\zeta: L_\zeta\to\Lc_\zeta$ betweent the stalks of $L$ and $\Lc$ over $\zeta$.
Here $\zeta$ is any non-Stokes direction.
Let $g_\pm: L\to\Lc$ be the isomorphism of local systems over $I^\pm$ which coincides
with $g_{\pm\zeta_0}$ at the point $\pm\zeta_0$.

\vskip .2cm

To analyze $g_\pm$, we extend the construction of solutions $G_{i,\zeta}(\phi)(z)$,
$\phi\in\bPhi_i(\Fc)_{-\ol\zeta}$, see \eqref{eq:fourier-integral}, to a more flexible setup
which allows us to deform the path of integration.
Let $z_0\in\CC_z$ be given and
$\gamma$ be any simple path, starting at $w_i$ in direction $-\ol\zeta$, avoiding other $w_j$ and eventually going to
infinity along a straight line on which the function $w\mapsto e^{z_0w}$ decays. Then, as in
\eqref {eq:alpha-i-phi}, we construct, out of $\alpha_{i,\zeta}(\phi)\in \bPsi_i(\Fc)_{-\ol\zeta}$,
a distribution solution $\alpha_{i,\gamma}(\phi)(w)$ of $p(f)=0$ defined on the path $\gamma$.
  We then put
\[
G_{i,\gamma}(\phi)(z) \,=\, \int_\gamma \alpha_{i,\gamma}(\phi)(w) e^{zw} dw
\]
It is a solution of $\wc\Mc$, i.e., a section of $\wc\Fc$,  defined for $z$ with $\arg(z)$ close to $\arg(z_0)$.
Thus the
$G_{i,\zeta}(\phi)(z)$  of  \eqref{eq:alpha-i-phi} corresponds to the case  when
$\gamma=K_i(\zeta)$ is  the straight path in the direction
$-\ol\zeta$ and $z_0=\zeta$.

\vskip .2cm

Now, to analytically continue $G_{i,\gamma}(\phi)(z)$ to $z$ with other values of
$\arg(z)$, all we need to do is to simply bend the eventual direction of the path $\gamma$
so that $w\mapsto e^{zw}$ continues to decay on it. With these notation, we prove:

\begin{lem}\label{lem:g-compat}
$g_{\pm}$ is compatible with the filtration everywhere on $I^\pm$
(split filtration in the source, Stokes filtration in the target).
\end{lem}

\noindent{\sl Proof of the Lemma:}  This an extension of the argument in the proof
of Proposition \ref{prop:g-zeta-varies}(b).  More precisely, look at the specialization of
$g_+$ at some $\zeta\in I^+$. By the above, it is given by the collection of the
$G_{i, \gamma_i}$, where $\gamma_i$ is obtained by bending the eventual
direction of the ray $K_i(\zeta_0)$ so that it approaches the infinity in the direction
$(-\ol\zeta)$, not $(-\ol\zeta_0)$.

To see that $g_+$ is compatible, over $\zeta$, with the filtrations, we compare it with
$G_\zeta$ which is obtained by integration over straight rays $K_i(\zeta)$
which go in the direction $(-\ol\zeta)$ right away, without bend.
As we move the direction from $\zeta_0$ to $\zeta$,  we may cross some Stokes
directions and so the  $G_{i,\gamma_i}$ will be expressed through the $G_{i,\zeta}$
via a series of Picard-Lefschetz formulas, like \eqref{eq:stokes-PL}.
Now,  by keeping track of the dominance of exponentials, we see the following:
\begin{itemize}
\item[(!)] {\em As long as $\zeta$ does not deviate from $\zeta_0$ by more than
$90^\circ$ in either direction},   the correction term in each such Picard-Lefschetz formula
will belong to the lower layer of the filtration!
\end{itemize}
This means that  compatibilty with
filtrations  will  hold all the way
{\em until the first Stokes directions differing from $\zeta_0$ by more than
$\pm 90^\circ$}, i.e., everywhere on $I^+$. Similarly for $I^-$. \qed

\begin{figure}[H]
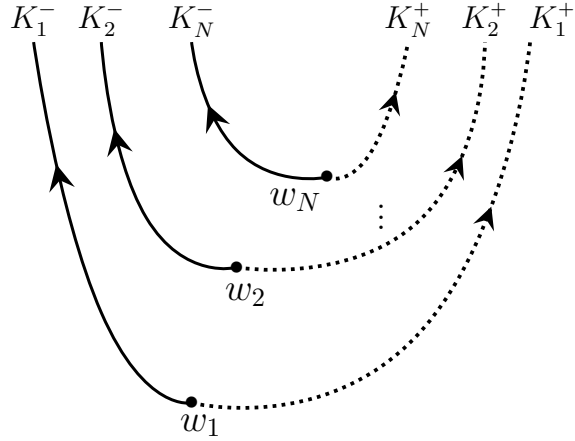

\centering

\btp[scale=.6]

\node at (1,3){$\bullet$};
\node at (-1,1) {$\bullet$};

\node at (-2,-2) {$\bullet$};

\draw
[  dotted,  decoration={markings,mark=at position 0.7 with
{\arrow[scale=1.5,>=stealth]{>}}},postaction={decorate},
line width=.5mm]   plot  [ smooth, tension=1] coordinates
{  (1,3) (2,3.6) (2.8,6)     };

    \node at (2.8, 6.5){$ K^+_N$};

\draw
[  dotted,  decoration={markings,mark=at position 0.7 with
{\arrow[scale=1.5,>=stealth]{>}}},postaction={decorate},
line width=.5mm]
plot  [ smooth, tension=1] coordinates {  (-1,1) (3,2) (4.5,6)    };

      \node at (4.5, 6.5){$ K^+_2$};

\draw
[  dotted,  decoration={markings,mark=at position 0.7 with
{\arrow[scale=1.5,>=stealth]{>}}},postaction={decorate},
line width=.5mm]
plot  [ smooth, tension=1] coordinates {   (-2,-2) (3,-0.5) (5.5,6)     };

        \node at (6, 6.5){$ K^+_1$};

\node at (0.3, 2.5){\large$w_N$};
\node at (-0.8,  0.4){\large$w_2$};
\node at (2.2,  2.3){\large$\vdots$};
\node at (-1.8,  -2.5){\large$w_1$};

\draw
[   decoration={markings,mark=at position 0.7 with
{\arrow[scale=2,>=stealth]{>}}},postaction={decorate},
line width=0.4mm]
plot  [ smooth, tension=1] coordinates{    (-2, -2)   (-4,-0) (-5.5,6)     } ;

  \node at (-5.5, 6.5){$ K^-_1$};

\draw
[decoration={markings,mark=at position 0.7 with
{\arrow[scale=2,>=stealth]{>}}},postaction={decorate},
line width=.4mm]
plot  [ smooth, tension=1] coordinates{   (-1,1)    (-3, 2.1) (-4,6)      } ;

   \node at (-4, 6.5){$ K^-_2$};

\draw
[decoration={markings,mark=at position 0.7 with
{\arrow[scale=2,>=stealth]{>}}},postaction={decorate},
line width=.4mm]
 plot  [ smooth, tension=1] coordinates{    (1,3)   (-1, 3.6)
(-2,6) };

 \node at (-2, 6.5){$ K^-_N$};

\etp
\caption{The Stokes matrix: combing at $\pm 90^\circ$.  } \label{fig:two-sys}

\end{figure}

Now, by definition, $C^+ = g_+^{-1} g_-: L\to L$, computed anywhere on $J^+$, for
example at $\zeta_+ = -\sqrt{-1}\cdot \zeta_0$. That is
\[
(g_-)_i \,=\sum_{j\geq i} (g_+)_{j} \, C^+_{ij}.
\]
In other words, $C^+$ expresses the solutions $G_{i, K_i^-}(\phi)(z)$,
$\phi\in\Phi_i$, through the solutions $G_{i, K_i^+}(\phi)(z)$.
Here $K_i^+$ is the path going from $w_i$ in the direction $(-\ol\zeta_0)$
``all the way out of $Q=\Conv(A)$'' and then turning $90^\circ$ left.
Similarly, $K_i^-$ goes from $w_i$ in the opposite direction $\zeta_0$ all the way out and
then turns $90^\circ$ right, see Fig.  \ref{fig:two-sys}.

\begin{figure}[H]
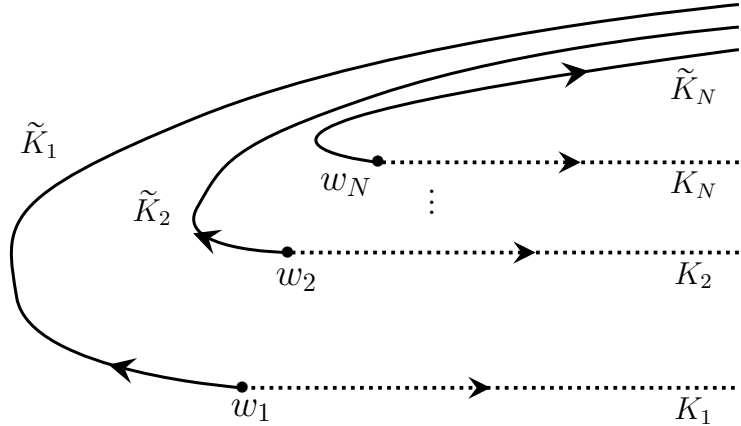

\centering

\btp[scale=.6]

\node at (1,3){$\bullet$};
\node at (-1,1) {$\bullet$};

\node at (-2,-2) {$\bullet$};

\draw
[  dotted,  decoration={markings,mark=at position 0.5 with
{\arrow[scale=1.5,>=stealth]{<}}},postaction={decorate},
line width=.5mm] (9,3) -- (1,3);

\draw
[  dotted,  decoration={markings,mark=at position 0.5 with
{\arrow[scale=1.5,>=stealth]{<}}},postaction={decorate},
line width=.5mm]
(9,1) -- (-1,1);

\draw
[  dotted,  decoration={markings,mark=at position 0.5 with
{\arrow[scale=1.5,>=stealth]{>}}},postaction={decorate},
line width=.5mm]
(-2,-2) -- (9,-2);

\node at (0.3, 2.5){\large$w_N$};
\node at (-0.8,  0.4){\large$w_2$};
\node at (2.2,  2.3){\large$\vdots$};
\node at (-1.8,  -2.5){\large$w_1$};

\draw
[   decoration={markings,mark=at position 0.9 with
{\arrow[scale=2,>=stealth]{<}}},postaction={decorate},
line width=0.4mm]
plot  [ smooth, tension=1] coordinates{   (9,6.5)    (-3,4) ( -7, 0) (-2, -2)} ;

\draw
[decoration={markings,mark=at position 0.9 with
{\arrow[scale=2,>=stealth]{<}}},postaction={decorate},
line width=.4mm]
plot  [ smooth, tension=1] coordinates{   (9, 6)  (1,4.5) ( -3, 2) (-1,1)} ;

\draw
[decoration={markings,mark=at position 0.7 with
{\arrow[scale=2,>=stealth]{>}}},postaction={decorate},
line width=.4mm]
 plot  [ smooth, tension=1] coordinates{    (1,3)   (0.5, 4)
(9,5.5) };

 \node at (8, 4.7){$\wt K_N$};

\node at (-6.5,3.5){$\wt K_1$};
\node at (-4,2){$\wt K_2$};

\node at (8,-2.5){$K_1$};
\node at (8,2.5){$K_N$};
\node at (8, 0.5){$K_2$};
\etp

\caption{ The Stokes matrix:  straight paths vs. combing at $180^\circ$. } \label{fig:St-matrix}
\end{figure}

Now, to find $C^+_{ij}$, it is enough to rotate the eventual directions of the $K_i^\pm$
by $90^\circ$ clockwise so that $K_i^+$ becomes unwound, i.e., becomes the straight
path $K_i=K_i(\zeta_0)$. The paths $K_i^-$ will then become twisted by $180^\circ$
and we denote the resulting paths $\wt K_i$, see Fig.  \ref{fig:St-matrix}
which corresponds to the case $\zeta_0=-1$.
So we need to express the solutions $g_{\wt K_i}$ in terms of the $g_{K_i}$.
Iterated application of the Picard-Lefschetz formulas gives, in the way
completely similar to the proof of Proposition \ref{prop:baby-IA} (``tightening the %
 a sum over left convex paths as claimed. Proposition \ref{prop:stokes-infra} is proved. \qed


\section{Review of schobers on surfaces with emphasis on $\CC$}\label{sec:schob-surf}

We now turn to the categorical analogs of the above constructions.
We use the conventions of Appendix \ref{app:enh}.

\subsection {Schobers on surfaces}\label{subsec:schob-surf}

\paragraph{Schobers on $(D,0)$: spherical functors. }\label{par:schober-D}
Let $\Phi, \Psi$ be  (enhanced) triangulated categories and $a:\Phi\to\Psi$ be an
exact (dg-) functor.
Assume that $a$ has left and right adjoints $^*a, a^*: \Psi\to\Phi$ with the corresponding
unit and counit maps
\[
\begin{gathered}
e: \Id_\Phi\lra a^*\circ a, \quad \eta: a\circ a^* \lra\Id_\Psi,
\\
e': \Id_\Psi \lra a\circ (^*a), \quad \eta': (^*a)\circ a \lra\Id_\Phi
\end{gathered}
\]
Using the dg-ehnancement, we define the functors
\be\label{eq:def-T-Phi}
T_\Phi = T_\Phi^{(a)} = \Cone(e)[-1], \quad T_\Psi = T_\Psi^{(a)} = \Cone (\eta).
\ee
They are called the {\em spherical cotwist} and {\em spherical twist} associated to $a$. The notation
$T_\Phi^{(a)}, T_\Psi^{(a)}$ will be helpful in situations when there is more than one spherical
functor connecting $\Phi$ and $\Psi$.
 By construction,  we have exact triangles of functors
\[
T_\Phi\buildrel \lambda \over \to \Id_\Phi \buildrel e\over\to a^*\circ a \buildrel f \over\to T_\Phi [1], \quad
T_\Psi[-1] \buildrel g \over\to a \circ a^* \buildrel \eta\over\to \Id_\Psi \buildrel \mu \over\to T_\Psi.
\]
and dual triangles
\[
^*T_\Phi[-1] \buildrel ^*f \over\to {}^*a \circ a \buildrel \eta'\over\to\Id_\Phi\buildrel {}^*\lambda \over \to  {}^*T_\Phi, \quad
{}^*T_\Psi \buildrel {}^*\mu \over\to \Id_\Psi \buildrel e'\over\to  a \circ {}^*a  \buildrel {}^*g \over\to {}^*T_\Psi[1].
\]
\begin{prop}[Spherical functor package \cite{AL} Th. 5.1] \label{prop:sph-pack}
In the following list, any two properties imply the two others:
\begin{itemize}
\item[(i)] $T_\Psi$ is an equivalence  (quasi-equivalence of dg-categories).

\item[(ii)] $T_\Phi$ is an equivalence.

\item[(iii)] The composite map
\[
{}^*a \circ T_{\Psi} [-1] \buildrel ({}^*a) \circ g \over\lra {^*a}\circ a\circ a^* \buildrel
\eta' \circ (a^*) \over\lra{ a^*}
\]
is an isomorphism (quasi-isomorphism of dg-functors).

\item[(iv)] The composite map
\[
a^* \buildrel (a^*) \circ e'\over\lra a^*\circ a \circ {}^*a
\buildrel f \circ (^*a) \over\lra T_\Phi \circ {}^*a [1]
\]
is an isomorphism. \qed
\end{itemize}
\end{prop}

The functor $a$ is called {\em spherical}, if the conditions of Proposition
\ref{prop:sph-pack} are satisfied.

\begin{prop}[\cite{kuznetsov} Prop. 2.13]\label{prop:kuznetsov}
A functor $a$ is spherical if and only if the following
two conditions are satisfied:
\begin{itemize}
\item[(1)] The map
\[
a^*\oplus (^*a) \buildrel a^*\circ e' + e\circ (^*a) \over\lra (a^*)\circ a \circ(^*a)
\]
is an isomoprhism.

\item[(2)] The map
\[
^*a\circ a\circ a^* \buildrel ((^*a)\circ\eta, \eta'\circ (^*a))\over\lra {^*a}\oplus a^*
\]
is an isomorphism  .\qed
\end{itemize}
\end{prop}

These conditions can be seen as categorical analogs of the conditions of
Definition \ref{def:sph-map-baby} and Proposition \ref{prop:baby-sph-pack}
which describe perverse sheaves $\Fc\in\Perv(D,0)$ equipped with
a non-degenerate bilinear form.   This motivates the  following \cite{KS-schobers}.

\begin{Defi}
A {\em perverse schober} $\Sen$  on $(D,0)$ is a datum of a spherical functor $a: \Phi\to\Psi$
between enhanced triangulated categories.

\end{Defi}

With a view towards  future more intrinsic theory, we
will say that a perverse schober $\Sen$ {\em is represented} by a spherical functor  $a: \Phi\to\Psi$ and write
$\Phi=\Phi(\Sen)$, $\Psi=\Psi(\Sen)$. Further, $T_\Phi$ and $T_\Psi$ being equivalences,
we can use them to construct
local systems of (emhanced) triangulated categories on $S^1$.
More precisely, we define $\bPhi(\Sen)$  to be the local system with stalk $\Phi$ and
monodromy $T_\Phi[2]$ and $\bPsi(\Sen)$ to be the local system with stalk $\Psi$
and monodromy $T_\Psi$. The introduction of the shift by $2$ is explained as follows.

\begin{prop}\label{prop:sph-loc-sys}
(a) We have canonical identifications of functors
\[
\begin{gathered}
a\circ T^{-1}_\Phi[-1] \,\simeq \, a^{**} \,\simeq \, T^{-1}_\Psi\circ a[1],
\\
a\circ T_{\Phi}[1]\, \simeq \,  {^{**}}a \,\simeq\, T_\Psi\circ a[-1],
\end{gathered}
\]
in particular, $a^{**}$ and $^{**}a$ exist.

\vskip .2cm

(b) The functor $a$ gives a morphism of local systems
\[
\ba = \ba(\Sen): \bPhi(\Sen)\lra\bPsi(\Sen).
\]
\end{prop}

\noindent {\sl Proof:} (a) Let us write the identifications of Proposition \ref{prop:sph-pack} (iii)-(iv)
\[
{}^*a \circ T_\Psi[-1] \,\simeq \, a^* \,\simeq \, T_\Phi \circ {}^*a [1].
\]
Note that for an equivalence of categories, the inverse equivalence is both left and right adjoint.
Therefore, taking the right adjoints of the above identifications, we get
the first line of identifications in (a). The second line is proved similarly.
\vskip .2cm

(b)  The second line of identifications in (a) implies
$a\circ T_\Phi[2]\simeq T_\Psi\circ a$ which means that $a$ is a morphism
of local systems. \qed

\vskip .2cm

We will refer to a morphism of local systems $\ba: \bPhi\to\bPsi$ (or, put differently,
a ``local system of spherical functors'') obtained from a spherical functor $a: \Phi\to\Psi$,
a {\em spherical local system}. This gives an equivalent, more geometric point of view
on a spherical functor.

\vskip .2cm

\begin{cor}\label{cor:spher-adj}
For   a spherical functor $a: \Phi\to\Psi$, all iterated adjoints
\[
\cdots, ^{**}a, {^*a}, a, a^*, a^{**}, \cdots
\]
exist and are spherical functors. \qed
\end{cor}

\begin{rem}\label{rem:spher-adj}
The twist and cotwist associated by \eqref{eq:def-T-Phi} to an iterated adjoint of $a$ will differ,
in general, from the one for $a$.
For example, for the spherical functor ${^*}a: \Psi\to\Phi$ we have
\[
\begin{gathered}
T_\Psi^{({^*}a)}\,=\,\Cone \{\Id_\Psi\buildrel e'\over \lra a\circ {^*a}\}[-1]\,=\, (T_\Psi^{(a)})^{-1},
\\
T_\Phi^{({^*}a)} \,=\, \Cone \{ {^*}a\circ a \buildrel\eta'\over\lra \Id_\Phi\} \,=\, (T_\Phi^{(a)})^{-1},
\end{gathered}
\]
see \cite{AL} \S 1.
\end{rem}

\vskip .2cm

The following is a categorical analog
of the concept of a perverse sheaf with a symmetric or antisymmetric bilinear form.

\begin{Defi}
Let $n\in\ZZ$.
A perverse schober $\Sen$ on $(D,0)$ represented by a spherical functor $a: \Phi\to\Psi$ is called
a {\em Calabi-Yau schober of dimension $n$} (or, shortly, an {\em $n$-CY-schober}), if the following conditions hold:
\item[(1)] $\Psi$ is a Calabi-Yau category of dimension $n$.

\item[(2)] The functor $T_\Phi[n+1]$ is the  Serre functor for $\Phi$.

\item[(3)] The data in (1) and (2) are compatible in the following sense. For any objects
$\phi\in\Phi, \psi\in\Psi$, the diagram of identifications
\[
\xymatrix{
\Hom_\Phi({}^*a(\psi), \phi)^* \ar@{<->}[rr]^{\on{adj}(^*a, a)^t}
\ar@{<->}[d]_{\text{Serre data}}&& \Hom_\Psi(\psi, a(\phi))^*
\ar@{<->}[d]^{\text{CY data}}
\\
\Hom_\Phi\bigl(\phi, \, T_\Phi({}^*a(\psi))[n+1]\bigr)
\ar@{<->}[dr]_{a^*\simeq T_\Phi\circ {}^*a[1]\hskip 1cm}
&& \Hom_\Psi(a(\phi), \psi[n])
\ar@{<->}[dl]^{\on{adj}(a, a^*)}
\\
& \Hom_\Phi(\phi, a^*(\psi)[n])
}
\]
is commutative.
\end{Defi}

\begin{ex}[(Spherical objects)]\label{ex:spher-ob}
The simplest class of examples of spherical functors has
$\Phi=D^b\Vect_\k$. A (dg-)-functor $a: D^b\Vect_\k\to\Psi$ is uniquely determined
by the datum of a single object $E= a(\k)\in\Psi$, by $a(V) = V\otimes_\k E$. The adjoint $a^*: \Psi\to D^b\Vect_\k$
is given by $F\mapsto \Hom^\bullet_\Psi(E,F)$.

Suppose $\Psi$ has a Serre functor $S$.
Recall, cf. \cite[Ch.8] {huybrechts} \cite{seidel-thomas}, that $E$ is called a {\em spherical object
of dimension $n$}, if:
\begin{itemize}
\item[(1)] $S(E)\simeq E[n]$.
\item[(2)] $H^j\Hom^\bullet (E,E) \,= \,
\begin{cases}
\k,& \text{ if } j=0,n;
\\
0,& \text{ otherwise. }
\end{cases}$
\end{itemize}
In this case $T_\Phi: D^b\Vect_\k\to D^b\Vect_\k$ is easily identified, using (2),  with the shift by $-(n+1)$, so it is an
equivalence.  It was shown in  \cite{seidel-thomas}, see also  \cite[Prop. 8.6] {huybrechts} that condition (1)
implies that $T_\Psi$ is also an equivalence, i.e., $a$ is a spherical functor.
Note that if  $\Psi$ is a Calabi-Yau
category  of dimension
$n$, then (1) is automatic, so if (2) holds, then  $a$ gives a Calabi-Yau schober.
\end{ex}

\paragraph{Schobers on surfaces. Transport functors. }
Let $(X,A=\{w_1,\cdots, w_N\})$ be a stratified surface and $S^1_{w_i}$ be the circle
of directions at $w_i$. Let $D_i$ be a small disk around $w_i$ and $D_i^\circ = D_i\- \{w_i\}$
be the punctured disk. Thus we have a homotopy equvalence $D_i^\circ\to S^1_{w_i}$
which is ``homotopy canonical'', i.e., defined uniquely up to a contractible space of
choices. In this paper we adopt the following naive definition
\cite{KS-schobers, {donovan-perv-RS}}.
\begin{Defi}\label{def:schober-surf}
A perverse schober on $(X,A)$ si a datum consisting of:
\begin{itemize}
\item A local system $\Sen_0=\Sen|_{X\- A}$ of pre-triangulated categories on $X\- A$.

\item For each $i$, a spherical local system $\ba_i: \bPhi_i (\Sen)\to\bPsi_i(\Sen)$ on $S^1_{w_i}$
(or, equivalently, on $D_i^\circ$).

\item An identification $\Sen_0|_{D_i^\circ} \,\simeq \, \bPsi_i(\Sen)$ for each $i$.
\end{itemize}
\end{Defi}

In fact, we will really need only the case when $X=D$ (a disk) or $\CC$ (complex plane)
in which an even easier definition is possible (see below).

A more systematic definition, using the language of $\oo$-categories, will be given in
\cite{DKSS}. Such a definition leads naturally to a construction of an $\oo$-category
$\Schob(X,A)$ of perverse schobers for $(X,A)$. In this paper we will not use this
$\oo$-category but will write $\Sen\in\Schob(X,A)$ just to signify that $\Sen$ is
a schober for $(X,A)$.

\vskip .2cm

As in \S \ref{subsec:abs-PL}\ref{par:PL-id}, let $\alpha$ be a piecewise smooth
oriented path in $X$,
joining $w_i$ with $w_j$ and avoiding other $w_k$. For a schober $\Sen\in\Schob(X,A)$
we have the category $\Phi_{i,\alpha}=\bPhi_i(\Sen)_{\dirr_i(\alpha)}$ (the stalk of
the local system $\bPhi_i(\Sen)$ in the direction of $\alpha$) and the similar category
$\Phi_{j,\alpha}$. Denoting $\Psi_\alpha$ the stalk of the local system $\Sen_0$
at any intermediate point of $\alpha$, we have two spherical functors
\[
\Phi_{i,\alpha} \buildrel a_{i,\alpha}\over\lra \Psi_\alpha \buildrel a_{j,\alpha}\over\lla \Phi_{j,\alpha}
\]
and we define the {\em transport functor}
\be\label{eq:transport-Mij}
M_{ij}(\alpha) \,=\, a_{j,\alpha}^* \circ a_{i,\alpha}: \Phi_{i,\alpha} \lra \Phi_{j,\alpha},
\ee
analogous to the transport map \eqref{eq:transport} for perverse sheaves.

\vskip .2cm

Let $\alpha^{-1}$ be the path obtained by reversing the direction of $\alpha$, so $\alpha^{-1}$
goes from $w_j$ to $w_i$.
We have the transport functor
\be\label{eq:transport-Mji}
M_{ji}(\alpha^{-1} ) \,=\, a^*_{i,\alpha}\circ a_{j,\alpha}: \Phi_{j,\alpha} \lra\Phi_{i,\alpha}.
\ee

\begin{prop}\label{prop:adj-trans}
We have identifications of functors
\[
M_{ij}(\alpha) \,\simeq \, T_{\Phi_j}\circ {}^*M_{ji}(\alpha^{-1})\, [1] \,\simeq \,
M_{ji}(\alpha^{-1})^* \circ T_{\Phi_i}[1].
\]
\end{prop}

\noindent {\sl Proof:} Taking the left adjoint of \eqref{eq:transport-Mji}, we get
\[
{}^*M_{ji}(\alpha^{-1}) \,\simeq \, {}^*a_{j,\alpha}\circ a_{i,\alpha} \,\simeq \, T_{\Phi_j}^{-1} \circ a^*_{j,\alpha}\circ
a_{i,\alpha}[-1] \, = \, T_{\Phi_j}^{-1} M_{ij}(\alpha)[-1],
\]
where we used the relation between $a^*_{j,\alpha}$ and $^*a_{j,\alpha}$ from Proposition
\ref{prop:sph-pack}(iv). This gives the first claim of the proposition.
The second claim is proved similarly by considering the right adjoint, using the identification $a^{**}_{i,\alpha}\simeq a_i\circ T_{\Phi_{i,\alpha}}^{-1}[1]$,
see Proposition \ref{prop:sph-loc-sys}. \qed

\begin{Defi}
A schober $\Sen\in\Schob(X,A)$ is called an $n$-Calabi-Yau schober, if :
\begin{itemize}
\item[(1)]  $\Sen_0$ consists of
$n$-Calabi-Yau categories and the monodromies preserve the Calabi-Yau structure.

\item [(2)]  For each $i$, the spherical local system $\ba_i: \bPhi_i(\Sen)\to\bPsi_i(\Sen)$ represents
an $n$-Calabi-Yau schober on the disk $D_i$ around $w_i$.
\end{itemize}
\end{Defi}

\paragraph{Schobers on a disk via spiders.}\label{par:schob-disk-spid}

We now specialize to the case $X=D$ being a closed disk. As in \S \ref{subsec:abs-PL}\ref{par:spider},
we choose $v\in\del D$ and a $v$-spider $K=(\gamma_1,\cdots, \gamma_N)$ for $(,A)$.
Choose the numbering $A=\{w_1,\cdots, w_N\}$ according to the slopes of the $\gamma_i$ at $v$, see
Fig. \ref{fig:br-spid}.
Given a schober $\Sen\in\Schob(D,A)$, we associate to it a diagram
\be\label{sch-disk-diag}
\xymatrix{
\Phi_1
\ar[dr]_{a_1}& \cdots & \Phi_N\ar[dl]^{a_N}
\\
& \Psi &
}
\ee
of $N$ spherical functors with common target as follows:

\vskip .2cm

$\Phi_i=\bPhi_i(\Sen)_{\dirr_{w_i}(\gamma_i)}$ is the stalk of $\bPhi_i(\Sen)$ in the direction of $\gamma_i$.

$\Psi = (\Sen_0)_v$ is the stalk of $\Sen_0$ at $v$.

$a_i$ is the stalk at $\dirr_{w_i}(\gamma_i)$ of the spherical local system
$\ba_i: \bPhi_i(\Sen)\to\bPsi_i(\Sen)$ followed by the monodromy of $\Sen_0$ from $w_i$ to $v$ along
$\gamma_i$.

\vskip .2cm

It is cleat that the data \eqref{sch-disk-diag} determine $\Sen$ completely, in analogy with
Proposition \ref {prop:GMV1}. Indeed, a choice of a spider $K$ identifies $\pi_1(D\- A, v)$
with the free group is generators $\wh\gamma_i$, going from $v$ to $w_i$ along $\gamma_i$,
then making a counterclockwise loop around $\gamma_i$ and then returning back along $\gamma_i$.
Given a diagram  \eqref{sch-disk-diag}, we define a local system $\Sen_0$ on $D\- A$
to correspond to the action of $\pi_1(D\- A, v)$ on $\Psi$ such that $\wh\gamma_i$ acts by
\[
T_{i,\Psi}=\Cone\{a_i \circ a^*_i\buildrel \eta\over\lra \Id_\Psi\}.
\]
The spherical local system $\ba_i: \bPhi_i(\Sen)\to\bPsi_i(\Sen)$ comes from the spherical functor
$a_i$.

\begin{Defi}
We say that $\Sen\in\Schob(D,A)$ {\em  is represented by the diagram  \eqref{sch-disk-diag}
in the presense of the spider} $K$, if it corresponds to the diagram by the procedure just described.
\end{Defi}

Given $\Sen$, the diagram   \eqref{sch-disk-diag}  representing it, depends on the choice of a spider $K$.
The change of the diagram corresponding to the change of spider under an action of the braid group
$\Br_N$, is described by the formulas identical to those of Proposition \ref{prop:ch-spid-1}. More precisely:

\begin {prop}\label{prop:schober-K-change}
Changing $K$ to $\tau_ i (K)$ changes the diagram   \eqref{sch-disk-diag} to
the diagram $(\ol\Phi_\nu\buildrel \ol a_\nu\over\to \Psi)_{\nu=1,\cdots, N}$, where
\[
\begin{gathered}
(\ol \Phi_1,\cdots, \ol\Phi_N) \,=\, (\Phi_1,\cdots, \Phi_{i-1}, \Phi_{i+1},\Phi_i, \Phi_{i+2},\cdots, \Phi_N),
\\
(\ol a_1, \cdots, \ol a_N) \,=\, (a_1,\cdots, a_{i-1}, \,T_{i,\Psi} \, a_{i+1},\,  a_i, a_{i+2}, \cdots, a_N).
\end{gathered}
\]
(Note that the composition of a spherical functor and an equivalence is again spherical).
 \qed
\end{prop}


\subsection{Picard-Lefschetz triangles and Picard-Lefschetz arrows}\label{subsec:PL-tri}

\paragraph{Picard-Lefschetz triangles.}
Let $(X,A)$ be a stratified surface, and $\Sen\in\Schob(X,A)$ be a schober with singularities in $A$.
The Picard-Lefschetz
identity of Proposition \ref{prop:PL-form} is now categorified into a {\em Picard-Lefschetz triangle}
(\cite{KS-schobers}, Prop. 2.13).
More precisely, consider the situation of Fig. \ref{fig:PicLef}, where a path $\gamma$, joining
$w_i$ and $w_k$, is moved across $w_j$. Let us define the categories $\Phi_i, \Phi_j$ and $\Phi_k$
using the
use the same identification of the appropriate  stalks of
$\bPhi_i(\Sen), \bPhi_j(\Sen)$ and $\bPhi_k(\Sen)$ as described in \eqref{eq:3-isotop}

\begin{prop}\label{prop:PL-triang}
In the situation described we have an exact triangle of functors $\Phi_i\to\Phi_k$
\[
M_{ik}(\gamma')[-1] \buildrel u_{\alpha,\beta,\gamma'}\over \lra M_{jk}(\alpha) \circ M_{ij}(\beta)
\buildrel v_{\alpha,\beta, \gamma}\over\lra M_{ik}(\gamma)\lra M_{ik}(\gamma').
\]
\end{prop}

\noindent{\sl Proof:}  we give a detailed argument by revisiting the proof of Proposition \ref{prop:PL-form}, now in the categorical context.  As in that proof, we deform the situation of Fig.  \ref{fig:PicLef} to that of Fig. \ref{fig:PicLef-pinch}, and prove the statement for that deformed situation. We use the same notations for the intermediate points $\wt w_i, \wt w_j, \wt w_k$ and
paths $\wt\alpha, \wt\beta, \wt\gamma, \wt\gamma'$.

\vskip .2cm

We view $\wt w_i$ and $\wt w_k$ as the midpoints to define the transport functors.
Let $\Psi_i, \Psi_k$ be the stalks of $\Sen_0$ at $\wt w_i, \wt w_k$ respectively.
We then have the spherical functors
\[
a_i: \Phi_i\lra \Psi_j, \quad a_k: \Phi_k\lra\Psi_k
\]
corresponding to the intervals $[w_i, \wt w_i]$ and $[w_k, \wt w_k]$, the monodromy functors
\[
T_{\wt\gamma}, T_{\wt\gamma'}: \Psi_i\lra \Psi_k
\]
for $\Sen_0$ corresponding to $\wt\gamma,\wt\gamma'$  and the spherical functors
\[
a_{j,\wt\beta}: \Phi_{j,\beta}\lra \Psi_i, \quad a_{j,\wt\alpha}: \Phi_{j,\alpha} \lra \Psi_k
\]
corresponding to the segments near $w_j$. So by definition we have (similar to \eqref{eq:m_{ij}-deformed-noncat})
\be\label{eq:M_{ij}-deformed}
\begin{gathered}
M_{ik}(\gamma) \,=\, a^*_k\circ T_{\wt\gamma}\circ a_i, \quad
M_{ik}(\gamma') \,=\, a^*_k\circ T_{\wt\gamma'}\circ a_i,
\\
M_{ij}(\beta)\,=\, a^*_{j,\wt\beta}\circ a_i, \quad M_{jk}(\alpha) \,=\, a^*_k\circ a_{j,\wt\alpha}.
\end{gathered}
\ee
We proceed analogously. Consider the area between $\wt\gamma$ and $\wt\gamma'$ as a disk with one marked
point $w_j$ and the radial cut along $\wt\beta$. The transformation
$T= T_{\wt\gamma}^{-1}\circ T_{\wt\gamma'}: \Psi_i\to\Psi_i$ is the full counterclockwise
monodromy of $\Sen_0$, i.e., of $\bPsi_i(\Sen)$, around $w_j$. Thus $T=T_{\Psi_j}$
is the spherical twist for the spherical functor $a_{j,\wt\beta}: \Phi_{j, \beta}\to \Psi_i$ and
we have the defining exact triangle
\[
T[-1]\lra  a_{j,\wt\beta}\circ a^*_{j, \wt\beta} \buildrel \eta\over\lra \Id_{\Psi_i} \lra T
\]
of functors $\Psi_i\to\Psi_i$. Now,  taking the composition $T_{\wt\gamma}\circ (-)$
is an exact functor $\Fun(\Psi_i, \Psi_i)\to \Fun(\Psi_i, \Psi_k)$ so it transforms the
above triangle to
\[
T_{\wt\gamma'}[-1]\lra a_{j, \wt\alpha}\circ {a^*}_{j,\wt\beta} \lra T_{\wt\gamma}
\lra T_{\wt\gamma'}.
\]
Note that here the equality $a_{j, \wt\alpha} \simeq T_{\wt\gamma}\circ a_{j,\wt\beta}$ is used, which is due to the clockwise identification of $\Phi_{j, \beta}, \Phi_{j, \alpha}$.

Our desired triangle is obtained from this by the exact functor
\[
a^*_k\circ (-) \circ a_i: \Fun(\Psi_i, \Psi_k) \lra \Fun(\Phi_i, \Phi_k),
\]
in virtue of \eqref{eq:M_{ij}-deformed}. \qed



\paragraph{The Picard-Lefschetz arrows: duality.}  We now focus on the arrows (natural transformations)
\be\label{eq:PL-arrows}
u_{\alpha\beta\gamma'}: M_{ik}(\gamma')[-1] \lra M_{jk}(\alpha)\circ M_{ij}(\beta), \quad
v_{\alpha \beta \gamma}:  M_{jk}(\alpha)\circ M_{ij}(\beta)\lra M_{ik}(\gamma)
\ee
in the triangle of Proposition \ref{prop:PL-triang}, which we call the {\em Picard-Lefschetz arrows}.
Note that the data of $\alpha,\beta, \gamma$ define $\gamma'$ and $\alpha,\beta, \gamma'$
define $\gamma$ (uniquely up to isotopy), so the notations $u_{\alpha\beta\gamma'}$
and $v_{\alpha\beta\gamma}$ are unambiguous.

\vskip .2cm

We want to express the $v$-arrow in terms of the $u$-arrows  corresponding to the
paths running in the opposite direction.
Recall \eqref{eq:transposed-functors} that a natural transformation $k: F\to G$ between two
functors gives rise to the left and right transposes $^t k: {^*}G \to {^*}F$ and
$k^t: G^*\to F^*$ between the corresponding adjoints (when the adjoints exist).

\begin{prop}\label{prop:v=u^t}
With respect to the identifications of the adjoints of the transport functors given by
Proposition \ref{prop:adj-trans}, we have natural identifications
\[
v_{\alpha\beta \gamma} \,\simeq \, u^t_{\beta, \alpha, \gamma} \circ T_{\Phi_i}
\,\simeq \, T_{\Phi_i} \circ {^t}u_{\beta, \alpha, \gamma}.
\]
\end{prop}

Let us explain the meaning of these identifications more precisely. We have
\[
u_{\beta \alpha\gamma}: M_{ki}(\gamma)[-1]
\lra M_{ji}(\beta) \circ T_{\Phi_j}^{-1}[-2] \circ M_{kj}(\alpha).
\]
Here, the appearance of $T_{\Phi_j}^{-1}[-2]$ in the middle comes from the fact that
we identify the stalks of $\bPhi_j(\Sen)$ at the directions of $\beta$ and $\alpha$ using the
monodromy
{\em along the arc opposite to that used in defining} $u_{\alpha\beta\gamma}$, and the counterclockwise monodromy of $\Phi(\Sen)$ is $T_{\Phi}[2]$. This identification
is {\em different from} that of
Convention \ref{conv:clock}.

\vskip .2cm

So the right transpose acts as
\[
u^t_{\beta, \alpha, \gamma}:
M_{kj}(\alpha)^*  \circ T_{\Phi_j}[2] \circ M_{ji}(\beta)^* \lra M_{ki}(\gamma)^*[1].
\]
Using Proposition \ref{prop:adj-trans} to transform the source and target of this, we write
\[
u^t_{\beta, \alpha, \gamma}: M_{jk}(\alpha) \circ T_{\Phi_j}^{-1}[-1] \circ T_{\Phi_j}[2] \circ M_{ij}(\beta) \circ T_{\Phi_i}^{-1}[-1] \lra M_{ik}(\gamma) \circ T_{\Phi_i}^{-1}
\]
so
\[
u^t_{\beta, \alpha, \gamma} \circ T_{\Phi_i}^{-1}: M_{jk}(\alpha) \circ M_{ij}(\beta)
\lra M_{ik}(\gamma),
\]
just as $v_{\alpha\beta\gamma}$. Similarly for the second identification involving $^t u$.

\vskip .2cm

The proof of Proposition \ref{prop:v=u^t} reduces to a general statement involving a
single spherical functor $a: \Phi\to\Psi$ and the corresponding exact triangle of functors
\be\label{eq:sph-triangle-self}
\xymatrix{
\Id_\Psi \ar[rr]^{\mu} && T_\Psi
\ar[dl]^{g}_{+1}
\\
& a \circ a^* \ar[ul]^\eta
&
}
\ee

\begin{prop}\label{prop:sph-tri-selfd}
The triangle \eqref{eq:sph-triangle-self} is self-dual, i.e., we have natural identifications
\[
\begin{gathered}
\eta^t \,\simeq \, T_{\Psi}^{-1}\circ\eta, \quad {^t}\eta \,\simeq \, g \circ T_{\Psi}^{-1}[1],
\\
\mu^t \,\simeq \, \mu\circ T_\Psi^{-1}, \quad {^t}\mu \,\simeq \, T_{\Psi}^{-1}\circ \mu.
\end{gathered}
\]
\end{prop}

\vskip .2cm

The meaning of these identifications is as follows. By definition, ${^t}\eta: \Id_\Psi \to a\circ {^*}a$.
By Proposition \ref {prop:sph-pack}(iii), we have ${^*}a \simeq a^* \circ T_{\Psi}^{-1}[1]$,
so ${^t}\eta$ acts as ${^t}\eta: \Id_\Psi \to a \circ a^* \circ T_{\Psi}^{-1}[1]$, and so does
$g \circ T_{\Psi}^{-1}[1]$. The following technical proposition will be used in the proof:

\vskip .2cm

\begin{prop}\label{prop:FG-dual}
Suppose we are given a pair of functors $F, G: C \lra D$, admitting adjoints, and the natural transform $\alpha: F \lra G$. Denote by $\eta_F, \eta_G$ the canonical maps from ${^*}F\circ F$ and ${^*}G \circ G$ to $\Id_C$. Then, the following diagram is commutative:
\[
\xymatrix{
{^*}G\circ F \ar[r]^{{^*}F\circ\alpha}\ar[d]^{{^t}\alpha\circ F} & {^*}G\circ G\ar[d]^{\eta_G}\\
{^*}F\circ F \ar[r]^{\eta_F} & \Id_C}
\]
\end{prop}
\noindent{\sl Proof:} Denote by $e_F$ and $e_G$ the canonical maps from $\Id_D$ to $F\circ{^*}F$ and $G\circ{^*}G$ respectively. Then, the map ${^t}\alpha: {^*}G \lra {^*}F$ can be described as the following composition, cf. \eqref{eq:transposed-functors}:
 \[
 {^*}G \buildrel {^*}G\circ e_F \over\lra {^*G}\circ F\circ {^*}F \buildrel{{^*}G \circ \alpha \circ {^*}F}\over\lra {^*G}\circ G \circ {^*}F \buildrel{\eta_G \circ {^*}F}\over\lra {^*}F
 \]
 The statement then follows reduces to the adjointness axiom which states that the following composition is identity:
 \[
 F \buildrel{F \circ e_F} \over\lra F \circ {^*}F \circ F \buildrel \eta_F \circ F \over\lra F
 \]
and the obvious commutativity of the following square:
\[
\xymatrix{
{^*}G\circ G \circ {*}F \circ F \ar[rr]^{\hskip 0.4cm{^*}G\circ G \circ \eta_F}\ar[d]^{\eta_G \circ {^*F}\circ F} && {^*}G\circ G\ar[d]^{\eta_G}\\
{^*}F\circ F \ar[rr]^{\eta_F} && \Id_C}
\]
\qed

\vskip .2cm

\noindent{\sl Proof of Proposition \ref{prop:sph-tri-selfd}:}  First of all, note that for any functor
$a$ with left and right adjoints $^*a, a^*$ and units and counits $e,\eta$ (for $(a, a^*)$) and
$e',\eta'$ (for $({^*}a, a)$), we have identifications
\[
{^t}\eta \simeq e', \,\, {^t}e \,\simeq\, \eta'.
\]
This is because, on one hand, ${^t}\eta$ and ${^t}e$ {\em can serve as} a unit and counit
{\em exhibiting}  $({^*}a,a)$ as an adjoint pair and, on the other hand, the adjoint functor is
defined uniquely up to a unique isomorphism (up to a contractible space of choices in
the dg-setting), see Appendix \ref{app:enh}\ref{par:app-dgadj}
Therefore (Proposition \ref{prop:cone-kt}), the triangle for the cone of $e'$, which
we denote

\[
\Delta_{e'} \,=\,\bigl\{ ... \buildrel {}^t\mu \over\to \Id_\Psi \buildrel e'\over\to  a \circ {}^*a  \buildrel {}^tg \over\to T^{-1}_\Psi[1] \buildrel {}^t\mu[1] \over\to [1] \}
\]

is identified with the left adjoint/transpose of the triangle for the cone of $\eta$, which we denote
\[
\Delta_{\eta} \,=\,\bigl\{T_\Psi[-1] \buildrel g \over\to a \circ a^* \buildrel \eta\over\to \Id_\Psi \buildrel \mu \over\to T_\Psi\}.
\]
So it is enough to construct an isomorphism of triangles $\Delta_{e'}\circ T_\Psi[-1]\to \Delta_{\eta}$.
We consider the vertical arrow given as the composition in Proposition \ref{prop:sph-pack}(iii),
composed with $T_\Psi[1]$ and the other two arrows being $\Id$ with respect to our identifications:

\[
\xymatrix{
a\circ {^*}a\circ T_\Psi[-1]
\ar[rr]_(.6){\hskip 0.5cm {^t}g \circ T_\Psi[-1]}
\ar[d]_{a\circ{^*}a\circ g}
&&
\Id_\Psi
\ar[r]_{{^t}\mu\circ T_\Psi}
\ar[dd]^{\Id}
&
T_\Psi
\ar[r]
\ar[dd]^{\Id}
&
a\circ {^*}a\circ T_\Psi \ar[d]_{a\circ {^*}a \circ g[1]}
\\
a\circ {^*}a \circ a\circ a^*
\ar[d]_{a\circ \eta' \circ a^*}
&&&& a\circ {^*}a \circ a \circ a^*[1] \ar[d]_{a\circ\eta'\circ a^*[1]}
\\
a\circ a^*
\ar[rr]^{\hskip 0.5cm\eta}
&&
\Id_\Psi
\ar[r]^\mu
&
T_\Psi
\ar[r]
&
a\circ a^* [1]
}
\]

First, the middle square is clearly commutative from the Proposition \ref{prop:FG-dual}, applied to $G = T_\Psi$, $F = \Id$. The leftmost square is commutative from the application of the same proposition to $F = T_{\Psi}[-1]$, $G = a \circ a^*$. The rightmost square is the application of the dual statement and also can be reduced to the case of the leftmost square by dualizing on the right and unraveling the $a^{**} \simeq T_{\Psi}^{-1}\circ a[1]$. \qed

\paragraph{Picard-Lefschetz arrows: rotation invariance.}
Consider a curvilinear geometric triangle $\Delta\subset \RR^2$ consisting
of three distinguished points $w_i, w_j, w_k$ and unoriented simple paths
$\alpha$ (joining $w_j, w_k$), $\beta$ (joining $w_i, w_j$) and $\gamma$
(joining $w_i, w_k$) and containing no other distinguished points inside or on the boundary.
We assume that the triple $(w_i, w_j, w_k)$ has the counterclockwise
orientation.

\vskip .2cm

The triangle  $\Delta$ gives rise to several Picard-Lefschetz situations and therefore
to several Picard-Lefschetz $u$- and $v$-arrows. For example (our earlier choice),
moving $\gamma$ past $w_j$ to a path $\gamma'$, we get the arrows
\eqref{eq:PL-arrows}. But we can also move any side of $\Delta$ past the opposite
vertex. In this way, for any cyclic rotation (i.e., even permutation) $(a,b,c)$ of
$(i,j,k)$ we get arrows
\[
u_{c,b,a}: M_{ca}[-1]\lra M_{cb}T_{\Phi_b}^{-1}[-2] M_{ba}, \quad v_{a,b,c}: M_{bc} M_{ab} \lra M_{ac},
\]
where we omitted the notation for the path joining the points, this path being
always one of the sides of $\Delta$.

\begin{figure}[H]
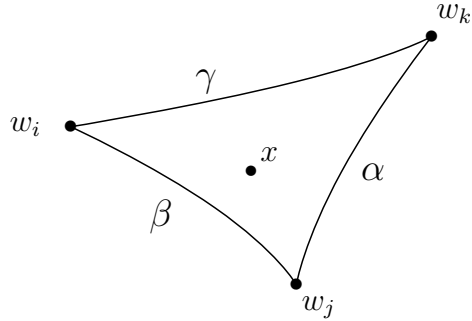

\centering
\btp[scale=.6]
\node at (1, -2.5){$\bullet$};
\node at (4,3){$\bullet$};
\node at (-4,1){$\bullet$};

\node at (0,0){\small$\bullet$};
\node at (0.4, 0.4){$x$};
\draw
[ line width=.2mm]
plot  [ smooth, tension=1] coordinates{   (1, -2.5) (2,0) (4,3)} ;

\draw
[line width=.2mm]
plot  [ smooth, tension=1] coordinates{    (-4,1)  (-0.8,-0.8) (1, -2.5)} ;

\draw
[line width=.2mm]
plot  [ smooth, tension=1] coordinates{    (-4,1)  (1,2) (4,3)} ;

\node at (-5,1){$w_i$};
\node at (4.5, 3.5){$w_k$};
\node at (1.5, -3){$w_j$};

\node at (2.7,0){\large$\alpha$};
\node at (-1,2){\large$\gamma$};
\node at (-2, -1){\large$\beta$};

\etp
\caption{A geometric triangle $\Delta$. }
\end{figure}

The appearance of $T_{\Phi_b}^{-1}[-2] $ in the first but not the second arrow above means the following. In forming $v_{abc}$ we identify the fibers of the local system $\bPhi_b$
using clockwise monodromy (going inside the triangle), in accordance with
Convention \ref{conv:clock}.  In forming $u_{c,b,a}$
we  use the identification using the monodromy along the arc which still goes inside
the triangle but which is now anticlockwise, thus making the choice
opposite to Convention \ref{conv:clock}. This is done
so that the target  of $u_{b,a,c}$ can be more easily related to the
source of $v_{a,b,c}$ by adjunction. After simplifying the shifts, we can write
the $u$-arrow as
\[
u_{c,b,a}: M_{ca}\lra M_{cb}T_{\Phi_b}^{-1}[-1] M_{ba}.
\]
In this way we get three $u$-arrows and three $v$-arrows.

\begin{prop}\label{prop:cyc-rot}
(a) The tensor rotation isomorphisms  (Corollary \ref{cor:tens-rot})
give canonical  identifications
\[
\begin{gathered}
\Hom(M_{jk} M_{ij}, M_{ik}) \, \= \,  \Hom( M_{ki} M_{jk}, M_{ji})\,  \=\,  \Hom( M_{ij}M_{ki},
M_{kj})
\\
\Hom(M_{ki}, M_{kj}T_{\Phi_j}^{-1}M_{ji}[-1]) \, \= \,  \Hom(M_{ij}, M_{ik}T_{\Phi_j}^{-1}M_{kj}[-1])\,  \=\,  \Hom(M_{jk}, M_{ji}T_{\Phi_j}^{-1}M_{ik}[-1]).
\end{gathered}
\]
(b) The isomorphisms in (a) take the $v_{abc}$, resp. $u_{abc}$,  into each other.

\end{prop}

\begin{rem}
Note that the $v$-maps are dual to the $u$-maps by Proposition \ref{prop:v=u^t}.
Proposition \ref{prop:cyc-rot} therefore means that $\Delta$ gives rise to
what is essentially a single
natural transformation (``Picard-Lefschetz plaquette'')
whose behavior is similar to that of  Clebsch-Gordan spaces in representation theory.
\end{rem}

\begin{rem}
The second line of identification is dual to the first up to the application of $T_\Phi$, using the formula $M^*_{ab} = M_{ba}T_{\Phi_a}^{-1}[-1]$ from \ref{prop:v=u^t}:

\[Hom(M_{ik}^*T_{\Phi_i}, M_{jk}^*M_{ij}^*T_{\Phi_i}) \, \= \,  \Hom(M_{ji}^*T_{\Phi_j}, M_{ki}^*M_{jk}^*T_{\Phi_j})\,  \=\,  \Hom(M_{kj}^*T_{\Phi_k}, M_{ij}^*M_{ki}T_{\Phi_k})\]
\end{rem}

\noindent{\sl Proof of  Proposition \ref{prop:cyc-rot}:}
(a) We prove the very first identification, the others are similar.
Let  $A,B,C$  be  three functors such that $\Hom(AB, C)$ makes sense.
By composing two consecutive tensor rotation isomorphisms of Corollary \ref{cor:tens-rot},
we get an identification
\be
\rho_{ABC}: \Hom(AB, C)\lra  \Hom(BC^*, A^*).
\ee
Let now
$
C=M_{ik},  A=M_{jk}, B=M_{ij}.
$
By Proposition \ref{prop:adj-trans},
$A^*= M_{kj}T_{\Phi_i}^{-1}[-1]$ and $C^*= M_{ki}T_{\Phi_i}^{-1}[-1]$, so
\[
R_{ijk} := \rho_{M_{jk}, M_{ij}, M_{ik}} \circ T_{\Phi_i}[1]: \Hom(M_{jk}M_{ij}, M_{ik}) \lra \Hom(M_{ij}M_{ki}, M_{kj})
\]
is the desired identification.

\vskip .2cm

(b) Let us prove that $R_{ijk}(v_{ijk})=v_{jki}$, the other such statements being similar.
Recalling the construction of the tensor rotation isomorphisms,
cf. the proof of Proposition
\ref{prop:FGH-iso}, we see that  for general $A,B,C$ as above, $\rho_{ABC}$ is the composition
\be\label{eq:tens-rot-twice}
\Hom(AB, C) \lra \Hom(A^*ABC^*, \, A^*CC^*) \lra \Hom(BC^*, A^*)
\ee
with the first arrow given by composition with $A^*$ on the left and $C^*$ on the right,
and  the second arrow given by  insertion of the unit $e_A: \Id\to A^*A$ and
the counit $\eta_C: CC^*\to\Id$. This means that for  any $v: AB\to C$,
the arrow $\rho_{ABC}(v)$
is the composition
\be\label{eq:2-ten-rot-expl}
\xymatrix{
BC^* \ar[r] ^{\hskip -0.3cm e_A\circ BC^*} & A^*ABC^*
\ar[rr]^{A^*\circ u\circ C^*} && A^*CC^* \ar[rr]^{\hskip 0.3cm A^*\circ\eta_C}&&A^*.
}
\ee
Now,  specialize this  to
$C=M_{ik}$ etc. as in (a).
Since we are dealing with a single triangle $\Delta$, we can represent our schober $\Sen$
in terms of three spherical functors
$a_i: \Phi_i\to\Psi$, $a_j: \Phi_j\to\Psi$, $a_k: \Phi_k\to\Psi$, cf. \eqref{sch-disk-diag}
where $\Psi$ corresponds to a point $x$ positioned inside $\Delta$.
The transport functors, which represent passing through the midpoints of the edges of $\Delta$,
can be just as well taken to represent passing through $x$ and back, that is
\[
M_{ik}= a^*_k a_i, \,\, M_{jk} = a^*_k a_j,\,\, M_{ij} = a^*_j a_i, \quad\text{etc.}
\]

Then
\[
v = v_{ijk}: a^*_k a_j a^*_j a_i \lra a^*_k a_i, \quad v = {a^*_k}\circ \eta_{j} \circ a_i
\]
is given simply by insertion of the counit $\eta_{j}: a_j a^*_j \to \Id_{\Psi}$.
Now, applying \eqref{eq:tens-rot-twice} to $v$,
i.e., implementing \eqref{eq:2-ten-rot-expl},
we get the  transformation
which is the composition $c$ (the dotted arrow below):
\[
\xymatrix{
(a_j^*)a_i (a_i^*) (a_k^{**}) \ar[r]^{ \hskip -1.3cm (1)}
\ar@{.>}[rrd]_{c}
& (a_j^*) (a_k^{**})(a_k^*) a_j (a_j^*)a_i (a_i^*) (a_k^{**})
\ar[r]^{ \hskip 0.3cm (2)} &
(a_j^*) (a_k^{**})(a_k^*) a_i (a_i^*) (a_k^{**})
\ar[d]^{(3)}
\\
&& (a_j^*) (a_k^{**})
}
\]
Here  (1) is the insertion of two units on the left and
(2) is the contraction via the counit $a_j(a^*_j) \to \Id$ in the middle of
the $8$-term composition, while (3) is the contraction via, first,
the counit $a_i (a_i^*)\to\Id$ and then via the counit $(a^*_k) (a^{**}_k)\to\Id$.

Now,  the axioms of units and counits of adjunctions
imply that $c$ is the equal to just the result of contraction via the counit $a_i(a^*_i) \to \Id$ in the $4$-term composition $(a_j^*)a_i (a_i^*) (a_k^{**})$.
At the same time,  $a_k^{**}=({a_k})T_{\Phi_k}^{-1}[-1]$, so
after composing $c$ on the right with $T_{\Phi_k}[1]$, we get  the
Picard-Lefschetz arrow
$v_{jki}: M_{ij} M_{ki} \to M_{kj}$. \qed

\begin{figure}[H]
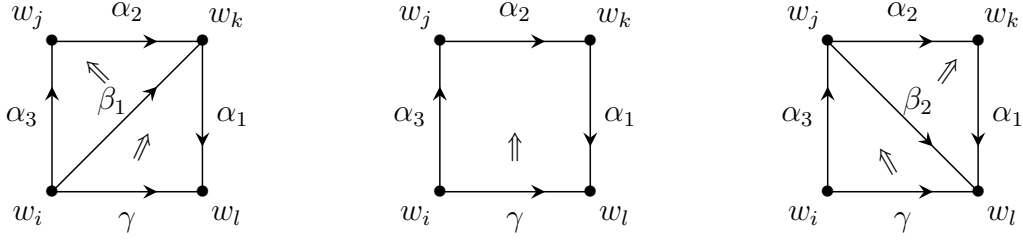

\centering

\btp[scale=0.4]
\node at (0,0){$\bullet$};  \node at (0,5){$\bullet$}; \node at (5,0){$\bullet$}; \node at (5,5){$\bullet$};

\draw  [decoration={markings,mark=at position 0.7 with
{\arrow[scale=1.5,>=stealth]{>}}},postaction={decorate},
line width=.2mm] (0,0) -- (5,0);

\draw  [decoration={markings,mark=at position 0.7 with
{\arrow[scale=1.5,>=stealth]{>}}},postaction={decorate},
line width=.2mm] (0,0) -- (0,5);

\draw  [decoration={markings,mark=at position 0.7 with
{\arrow[scale=1.5,>=stealth]{>}}},postaction={decorate},
line width=.2mm] (0,5) -- (5,5);

\draw  [decoration={markings,mark=at position 0.7 with
{\arrow[scale=1.5,>=stealth]{>}}},postaction={decorate},
line width=.2mm] (5,5) -- (5,0);

\draw  [decoration={markings,mark=at position 0.7 with
{\arrow[scale=1.5,>=stealth]{>}}},postaction={decorate},
line width=.2mm] (0,0) -- (5,5);

\node at (2.5, -1){$\gamma$};
\node at (6, 2.5){$\alpha_1$};
\node at (-1, 2.5){$\alpha_3$};
\node at (2.5, 6){$\alpha_2$};
\node at (2,3){$\beta_1$};

\node at (-0.8, -0.8){$w_i$};
\node at (5.8, 5.8){$w_k$};
\node at (5.8, -0.8){$w_l$};
  \node at (-0.8, 5.8){$w_j$};

\node[rotate=335] at (3,1.5){$\Uparrow$};
\node[rotate=45] at (1.5,4){$\Uparrow$};

\etp
\quad\quad\quad\quad
\btp[scale=0.4]
\node at (0,0){$\bullet$};  \node at (0,5){$\bullet$}; \node at (5,0){$\bullet$}; \node at (5,5){$\bullet$};

\draw  [decoration={markings,mark=at position 0.7 with
{\arrow[scale=1.5,>=stealth]{>}}},postaction={decorate},
line width=.2mm] (0,0) -- (5,0);

\draw  [decoration={markings,mark=at position 0.7 with
{\arrow[scale=1.5,>=stealth]{>}}},postaction={decorate},
line width=.2mm] (0,0) -- (0,5);

\draw  [decoration={markings,mark=at position 0.7 with
{\arrow[scale=1.5,>=stealth]{>}}},postaction={decorate},
line width=.2mm] (0,5) -- (5,5);

\draw  [decoration={markings,mark=at position 0.7 with
{\arrow[scale=1.5,>=stealth]{>}}},postaction={decorate},
line width=.2mm] (5,5) -- (5,0);

\node at (2.5, -1){$\gamma$};
\node at (6, 2.5){$\alpha_1$};
\node at (-1, 2.5){$\alpha_3$};
\node at (2.5, 6){$\alpha_2$};

\node at (-0.8, -0.8){$w_i$};
\node at (5.8, 5.8){$w_k$};
\node at (5.8, -0.8){$w_l$};
  \node at (-0.8, 5.8){$w_j$};

  \node at (2.5, 1.5){$\Uparrow$};

\etp
\quad\quad\quad\quad
\btp[scale=0.4]
\node at (0,0){$\bullet$};  \node at (0,5){$\bullet$}; \node at (5,0){$\bullet$}; \node at (5,5){$\bullet$};

\draw  [decoration={markings,mark=at position 0.7 with
{\arrow[scale=1.5,>=stealth]{>}}},postaction={decorate},
line width=.2mm] (0,0) -- (5,0);

\draw  [decoration={markings,mark=at position 0.7 with
{\arrow[scale=1.5,>=stealth]{>}}},postaction={decorate},
line width=.2mm] (0,0) -- (0,5);

\draw  [decoration={markings,mark=at position 0.7 with
{\arrow[scale=1.5,>=stealth]{>}}},postaction={decorate},
line width=.2mm] (0,5) -- (5,5);

\draw  [decoration={markings,mark=at position 0.7 with
{\arrow[scale=1.5,>=stealth]{>}}},postaction={decorate},
line width=.2mm] (5,5) -- (5,0);

\draw  [decoration={markings,mark=at position 0.7 with
{\arrow[scale=1.5,>=stealth]{>}}},postaction={decorate},
line width=.2mm] (0,5) -- (5,0);

\node at (2.5, -1){$\gamma$};
\node at (6, 2.5){$\alpha_1$};
\node at (-1, 2.5){$\alpha_3$};
\node at (2.5, 6){$\alpha_2$};
\node at (3,3) {$\beta_2$};

\node at (-0.8, -0.8){$w_i$};
\node at (5.8, 5.8){$w_k$};
\node at (5.8, -0.8){$w_l$};
  \node at (-0.8, 5.8){$w_j$};

\node[rotate=30] at (2,1){$\Uparrow$};
\node[rotate=325] at (4,4){$\Uparrow$};

\etp

\caption{Coassociativity of Picard-Lefschetz arrows.} \label{fig:part-comonad}
\end{figure}

\paragraph{Picard-Lefschetz arrows: (co)associativity.}  Among the Picard-Lefschetz arrows
\eqref{eq:PL-arrows}, the arrow $v_{\alpha,\beta, \gamma}$ looks vaguely like a multiplication,
while $u_{\alpha,\beta,\gamma'}$ looks like a comultiplication. To pursue this analogy, consider
a curved $4$-gon formed by $4$ points $w_i,w_j,w_k,w_k\in A$ and paths
$\alpha_1,\alpha_2,\alpha_3,\gamma$,
as in the middle for Fig. \ref{fig:part-comonad}. Let $\beta_1$ and $\beta_2$ be the two diagonals, as in the left and right path of
Fig. \ref{fig:part-comonad}.

\begin{prop}\label{prop:PL-ass}
(a) The $u$-arrows are coassociative, i.e., the diagram
\[
\xymatrix{
M_{il}(\gamma) \ar[rrr]^{u_{\alpha_1, \beta_1, \gamma}[1]}
\ar[d]_{u_{\beta_2,\alpha_3,\gamma}[1]}
&& &M_{kl}(\alpha_1) \circ M_{ik}(\beta_1)[1]
\ar[d]^{M_{kl}(\alpha_1)\circ u_{\alpha_2, \alpha_3,\beta_1}[2]}
\\
M_{jl}(\beta_2) \circ M_{ij}(\alpha_3)[1] \ar[rrr]_{\hskip -1cm u_{\alpha_1, \alpha_2,\beta_2}\circ M_{ij}(\alpha_3)[2]}
&&& M_{kl}(\alpha_1) \circ M_{jk}(\alpha_2) \circ M_{ij}(\alpha_3)[2]
}
\]
is commutative.

\vskip .2cm

(b) The $v$-arrows are associative, i.e., the diagram
\[
\xymatrix{
M_{kl}(\alpha_1) \circ M_{jk}(\alpha_2) \circ M_{ij}(\alpha_3) \ar[rrr]^{\hskip 1cm v_{\alpha_1,\alpha_2, \beta_2}\circ M_{ij}(\alpha_3)}
\ar[d]_{M_{kl}(\alpha_1) \circ v_{\alpha_2,\alpha_3,\beta_1}}
&&& M_{kl}(\alpha_1) \circ M_{ik}(\beta_1)
\ar[d]^{v_{\beta_2,\alpha_3,\gamma}}
\\
M_{jl}(\beta_2) \circ M_{ij}(\alpha_3) \ar[rrr]_{v_{\alpha_1,\beta_1,\gamma}}
&&& M_{il}(\gamma)
}
\]
is commutative.
\end{prop}

Here (as elsewhere) the composition of rectilinear transports is defined following Convention \ref{conv:clock}.
That is, the stalk of the local systems $\bPhi$ at the intermediate points are identified via the
clockwise monodromy from the incoming direction to the outgoing one.

\vskip .2cm

\begin{figure}[H]
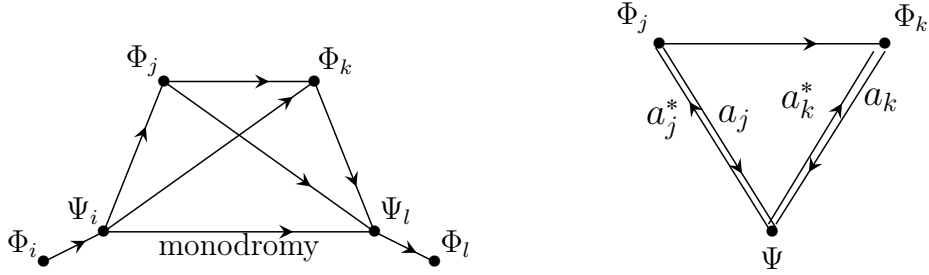

\centering
\btp[scale=0.4]

\node at (0,0){$\bullet$};
\node at (9,0){$\bullet$};
\node at (2,5){$\bullet$};
\node at (7,5){$\bullet$};

\draw  [decoration={markings,mark=at position 0.7 with
{\arrow[scale=1.5,>=stealth]{>}}},postaction={decorate},
line width=.2mm] (0,0) -- (9,0);

\draw  [decoration={markings,mark=at position 0.7 with
{\arrow[scale=1.5,>=stealth]{>}}},postaction={decorate},
line width=.2mm] (0,0) -- (2,5);

\draw  [decoration={markings,mark=at position 0.9 with
{\arrow[scale=1.5,>=stealth]{>}}},postaction={decorate},
line width=.2mm] (0,0) -- (7,5);

\draw  [decoration={markings,mark=at position 0.7 with
{\arrow[scale=1.5,>=stealth]{>}}},postaction={decorate},
line width=.2mm] (2,5) -- (9,0);

\draw  [decoration={markings,mark=at position 0.7 with
{\arrow[scale=1.5,>=stealth]{>}}},postaction={decorate},
line width=.2mm] (2,5) -- (7,5);

\draw  [decoration={markings,mark=at position 0.7 with
{\arrow[scale=1.5,>=stealth]{>}}},postaction={decorate},
line width=.2mm] (7,5) -- (9,0);

\draw  [decoration={markings,mark=at position 0.7 with
{\arrow[scale=1.5,>=stealth]{>}}},postaction={decorate},
line width=.2mm]  (-2,-1) -- (0,0);

\node at (-2,-1){$\bullet$};
\node at (-2.7, -0.4){$\Phi_i$};

\draw   [decoration={markings,mark=at position 0.7 with
{\arrow[scale=1.5,>=stealth]{>}}},postaction={decorate},
line width=.2mm] (9,0) -- (11,-1);

\node at (11, -1){$\bullet$};
\node at  (11.7, -0.4){$\Phi_l$};

\node at (4.5, -0.5){monodromy};

\node at (-0.7,0.7) {$\Psi_i$};
\node at (9.7, 0.7) {$\Psi_l$};
\node at (1.3, 5.7){$\Phi_j$};
\node at (7.7, 5.7){$\Phi_k$};

\etp
\quad\quad\quad\quad
\btp[scale=0.5]
\node (j)  at (0,0){$\bullet$};
\node (k)  at (6,0){$\bullet$};
\node (psi) at (3, -5){$\bullet$};

\draw   [decoration={markings,mark=at position 0.7 with
{\arrow[scale=1.5,>=stealth]{>}}},postaction={decorate},
line width=.2mm] (0,0) -- (6,0);

\draw   [decoration={markings,mark=at position 0.7 with
{\arrow[scale=1.5,>=stealth]{>}}},postaction={decorate},
line width=.2mm] (2.9, -4.8) -- (5.7, -0.2);

  \draw   [decoration={markings,mark=at position 0.7 with
{\arrow[scale=1.5,>=stealth]{>}}},postaction={decorate},
line width=.2mm] (6.0, -0.2) -- (3.1, -4.8);

  \draw   [decoration={markings,mark=at position 0.7 with
{\arrow[scale=1.5,>=stealth]{>}}},postaction={decorate},
line width=.2mm] (0.1, -0.1) -- (3.1, -4.9);

 \draw   [decoration={markings,mark=at position 0.7 with
{\arrow[scale=1.5,>=stealth]{>}}},postaction={decorate},
line width=.2mm]  ( 2.9, -4.9) -- (-0.1, -0.1);

\node at (-0.7, 0.7){$\Phi_j$};
\node at (6.7, 0.7){$\Phi_k$};
\node at (3, -5.7){$\Psi$};

\node at (0.1, -2){\large$a^*_j$};
\node at (2, -2){\large$a_j$};

\node at (3.7,-1.5){\large$a^*_k$};
\node at (5.9, -1.5){\large$a_k$};
\etp
\caption{Curved quadrilateral, pinched... and then contracted.} \label{fig:ijkl-pinch}

\end{figure}

\noindent{\sl Proof:} We prove (b), since (a) then follows by duality (Proposition   \ref{prop:v=u^t}).
As in the construction of the Picard-Lefschetz triangle (Proposition \ref{prop:PL-triang}),
we ``pinch'' the diagram at $w_i$ and $w_l$, after which all the functors in the diagram have a common initial part $\Phi_i\buildrel a_i\over\to \Psi_i$
and a common final part $\Psi_l \buildrel {a^*}_l \over\to \Psi_l$, see the left part of
Fig. \ref{fig:ijkl-pinch}.

\vskip .2cm

So it suffices to prove the analogous statement for the ``stripped''
functors acting from $\Psi_i$ to $\Psi_l$.
Now, after such  stripping off the common beginning and end, all that remains from $M_{il}(\gamma)$
is the monodromy from $\Psi_i$ to $\Psi_l$. So we can eliminate this monodromy also,
contracting the part of the path $\gamma$ joining  the points carrying $\Psi_i$ and $\Psi_j$, to a point
and replacing the spaces $\Psi_i$ and $\Psi_j$ by a common space $\Psi$ (identified with them both
via the  above monodromy).
This gives a diagram in the
right part of Fig. \ref{fig:ijkl-pinch}.

\noindent
The composition $M_{jl}(\beta_2)  \circ M_{ij}(\alpha_3)$, stripped of its initial and final parts, now becomes
$a_j\circ a^*_j: \Psi\to\Psi$, while the composition $M_{kl}(\alpha_1) \circ M_{ik}(\beta_1)$, similarly stripped, becomes
$a_k\circ a^*_k: \Psi\to\Psi$. So we are reduced to the commutativity of the diagram
\[
\xymatrix{
 a_k\circ {a^*_k}\circ a_j\circ a^*_j
\ar[d] \ar[r]& a_k\circ {a^*_k}\ar[d]
\\
a_j \circ {a^*_j} \ar[r] & \Id_\Psi
}
\]
obtained by two separate contractions of the counits for the adjunctions. This diagram clearly commutes.
\qed

\paragraph{Pasting of natural transformations.}\label{par:pasting-NT}
The compositions of the paths in the squares of Proposition \ref{prop:PL-ass} can be seen as
instances of $2$-dimensional composition, or {\em pasting} of natural transformations
(or, more generally, of $2$-morphisms
in a $2$-category),  see \cite{maclane, power-pasting}.
To perform pasting, one needs a {\em pasting diagram} which can be seen as a plane
region decomposed into (curvilinear) polygons so that:
\begin{itemize}

\item[(0)] To each vertex of each polygon there is associated a category.

\item[(1)] Each edge is oriented and carries a functor between the categories
associated to its vertices.

\item[(2)] Each polygon is oriented, which (together with orientation of the edges)
subdivides its edges into ``source''  edges and ``target''  edges. These two groups of
edges should form two composable chains with common beginning and end vertices.
Further, the polygon carries a natural transformation  (indicated by a double arrow)
between the functors
given by these composable chains.  It is indicated by a double arrow inside the polygon.

\item[(3)] Each internal edge appears  once in a source chain of a $2$-cell and
once in the target chain of another $2$-cell.
\end{itemize}

Under these conditions, one can compose the natural transformations associated
to all the $2$-cells to a single natural transformation (``pasting theorem'').
See  \cite{power-pasting} for more details.

\begin{figure}[H]
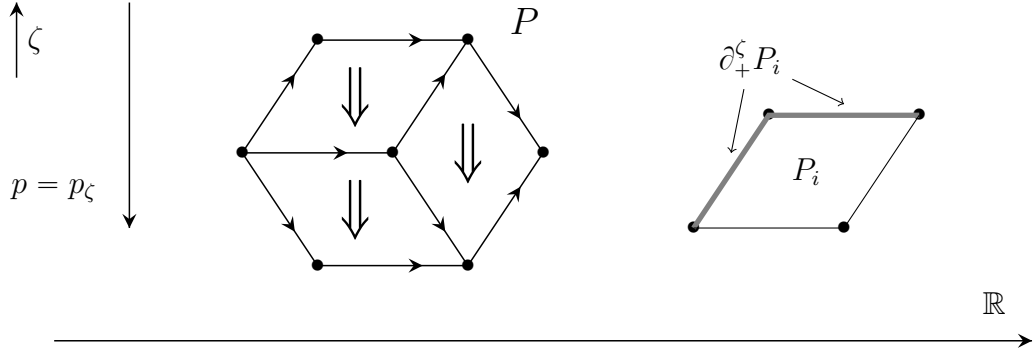

\centering
\btp[scale=0.5]

\node at (0,0){$\bullet$};
\node at (-4,0){$\bullet$};
\node at (-8,0){$\bullet$};
\node at (-6,3){$\bullet$};
\node at (-2,3){$\bullet$};
\node at (-2,-3){$\bullet$};
\node at (-6,-3){$\bullet$};

 \draw  [decoration={markings,mark=at position 0.7 with
{\arrow[scale=1.5,>=stealth]{>}}},postaction={decorate},
line width=.2mm] (-8,0) -- (-6,3);

 \draw  [decoration={markings,mark=at position 0.7 with
{\arrow[scale=1.5,>=stealth]{>}}},postaction={decorate},
line width=.2mm] (-6,3) -- (-2,3);

 \draw  [decoration={markings,mark=at position 0.7 with
{\arrow[scale=1.5,>=stealth]{>}}},postaction={decorate},
line width=.2mm] (-2,3) -- (0,0);

 \draw  [decoration={markings,mark=at position 0.7 with
{\arrow[scale=1.5,>=stealth]{>}}},postaction={decorate},
line width=.2mm] (-8,0) -- (-4,0);

 \draw  [decoration={markings,mark=at position 0.7 with
{\arrow[scale=1.5,>=stealth]{>}}},postaction={decorate},
line width=.2mm] (-8,0) -- (-6,-3);

  \draw  [decoration={markings,mark=at position 0.7 with
{\arrow[scale=1.5,>=stealth]{>}}},postaction={decorate},
line width=.2mm] (-6,-3) -- (-2,-3);

 \draw  [decoration={markings,mark=at position 0.7 with
{\arrow[scale=1.5,>=stealth]{>}}},postaction={decorate},
line width=.2mm] (-2,-3) -- (0,0);

  \draw  [decoration={markings,mark=at position 0.7 with
{\arrow[scale=1.5,>=stealth]{>}}},postaction={decorate},
line width=.2mm] (-4,0) -- (-2,-3);

  \draw  [decoration={markings,mark=at position 0.7 with
{\arrow[scale=1.5,>=stealth]{>}}},postaction={decorate},
line width=.2mm] (-4,0) -- (-2,3);

\node at (-5, 1.5){\huge$\Downarrow$};
\node at (-2,0){\huge$\Downarrow$};
\node at (-5, -1.5){\huge$\Downarrow$};

\node at (-0.5,3.5){\large$P$};

\node at (4,-2){$\bullet$};
\node at (8,-2){$\bullet$};
\node at (6,1){$\bullet$};
\node at (10,1){$\bullet$};

\draw[line width=2, color=gray] (4,-2) -- (6,1) -- (10,1);
\draw (4,-2) -- (8,-2) -- (10,1);

\node at  (7,-0.5){$P_i$};

\node (D+) at (5.5,2.5){$\del^\zeta_+P_i$};

\draw [->] (D+) -- (8,1.3);
\draw [->] (D+) -- (5,0);

\draw  [decoration={markings,mark=at position 1.0 with
{\arrow[scale=1.5,>=stealth]{>}}},postaction={decorate},
line width=.2mm] (-13,-5) -- (13,-5);

\node at (12, -4){$\RR$};

\draw  [decoration={markings,mark=at position 1.0 with
{\arrow[scale=1.5,>=stealth]{>}}},postaction={decorate},
line width=.2mm] (-11,4) -- (-11,-2);

\node at (-13, -1){$p=p_\zeta$};

\draw  [decoration={markings,mark=at position 1.0 with
{\arrow[scale=1.5,>=stealth]{>}}},postaction={decorate},
line width=.2mm] (-14,2) -- (-14,4);

\node at (-13.5, 3){$\zeta$};

\etp
\caption {A pasting diagram from a polygonal subdivision.}\label{fig:pasting-poly}
\end{figure}

\begin{exas}\label{ex:pasting-convex}
(a) The two  sides of Fig. \ref{fig:part-comonad}
are pasting diagrams whose pastings are the compositions of the two paths
in the square in Proposition \ref{prop:PL-ass}(a). The proposition says that these compositions
are equal, which is represented by a single transformation in the middle of the figure.

\vskip .2cm

(b) Let $\zeta\in \CC$, $|\zeta|=1$, be a direction. Then we have the projection $p=p_\zeta:\CC\to\RR$
to the line orthogonal to $\zeta$ which we identify with $\RR$ and orient as in Fig. \ref{fig:pasting-poly}.
 Let $P=\Conv(B)\subset \CC$, $B\subset A$,  $|B|\geq 3$,  be a convex polygon with vertices in $A$.
Assume   that $A$ is in linearly general position including
to $\zeta$-infinity. Then
 $p(x)\neq p(y)$ whenever $[x,y]$ is an edge of  $P$. Therefore  every such edge
$[x,y]$ acquires an orientation from $x$ to $y$, where  $p(x)<p(y)$.
In addition, the boundary $\del P$ becomes decomposed
into the {\em upper boundary} $\del_+(P) = \del_+^\zeta(P)$ and the {\em lower boundary}
$\del_- (P) = \del_-^\zeta(P)$. By definition, $\del_+ (P)$ is ``what is visible
from the  source of the projection'',  see Fig.
\ref{fig:pasting-poly}, and $\del_-(P)$ is the complementary arc (with the same beginning and end).
This makes $P$ into a pasting scheme with one $2$-cell,
 in which the arc
$\del_+(P)$ consists of source edges and $\del_-(P)$ of target edges.

\vskip .2cm

Further, let  $\Pc=\{P_i\}_{i\in I}$ be a
{\em polygonal decomposition} of $P$, i.e., a finite set of convex polygons
$P_i\subset P$, also with vertices in $A$,  such that $P=\bigcup P_i$ and any intersection $P_i\cap P_j$
is a common face (i.e., $\emptyset$, or a vertex or an edge) of both.
By the above,  each $P_i$ becomes a pasting scheme with one $2$-cell.
This makes
  the whole subdivision $\Pc$ into a pasting scheme as well.  See \cite{KV-pasting}.

\end{exas}

\begin{figure}[H]
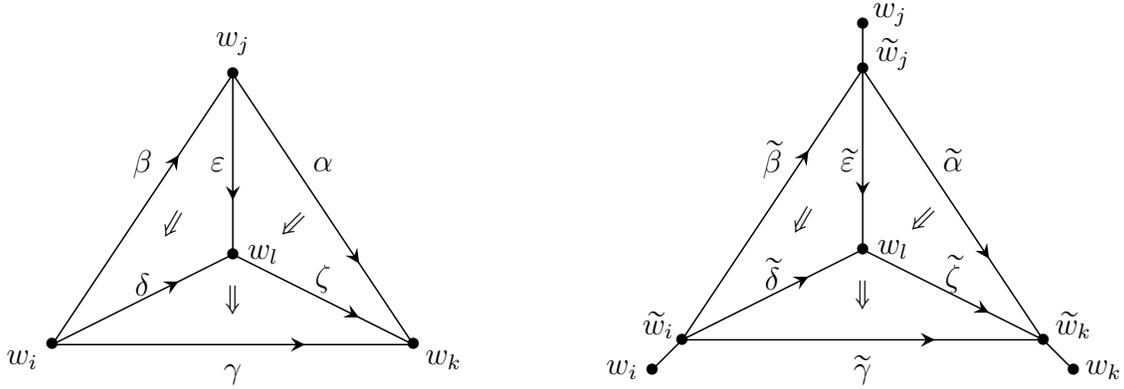

\centering
\btp[scale=0.4]

\node at (0,0){$\bullet$};
\node at (-6,-3){$\bullet$};
\node at (6,-3){$\bullet$};
\node at (0,6){$\bullet$};

\draw  [decoration={markings,mark=at position 0.7 with
{\arrow[scale=1.5,>=stealth]{>}}},postaction={decorate},
line width=.2mm] (-6,-3) -- (6,-3);

 \draw  [decoration={markings,mark=at position 0.7 with
{\arrow[scale=1.5,>=stealth]{>}}},postaction={decorate},
line width=.2mm] (-6,-3) -- (0,6);

\draw  [decoration={markings,mark=at position 0.7 with
{\arrow[scale=1.5,>=stealth]{>}}},postaction={decorate},
line width=.2mm] (0,6) -- (6,-3);

\draw  [decoration={markings,mark=at position 0.7 with
{\arrow[scale=1.5,>=stealth]{>}}},postaction={decorate},
line width=.2mm] (-6,-3) -- (0,0);

 \draw  [decoration={markings,mark=at position 0.7 with
{\arrow[scale=1.5,>=stealth]{>}}},postaction={decorate},
line width=.2mm] (0,6) -- (0,0);

\draw  [decoration={markings,mark=at position 0.7 with
{\arrow[scale=1.5,>=stealth]{>}}},postaction={decorate},
line width=.2mm] (0,0) -- (6,-3);

\node at (-7, -3.5){$w_i$};
\node at (7, -3.5){$w_k$};
\node at (0,7){$w_j$};
\node at (1,0){$w_l$};
\node at (3,3){$\alpha$};
\node at (-3,3){$\beta$};
\node at (0,-4){$\gamma$};
\node at (-3,-1){$\delta$};
\node at (-0.5,3){$\eps$};
\node at (3,-0.8){$\zeta$};

\node at (0, -1.5){$\Downarrow$};
\node[rotate=330] at (-2, 1){$\Downarrow$};
\node[rotate=315] at (2,1) {$\Downarrow$};

\etp
\quad\quad\quad\quad
\btp[scale=0.4]

\node at (0,0){$\bullet$};
\node at (-6,-3){$\bullet$};
\node at (6,-3){$\bullet$};
\node at (0,6){$\bullet$};

\draw  [decoration={markings,mark=at position 0.7 with
{\arrow[scale=1.5,>=stealth]{>}}},postaction={decorate},
line width=.2mm] (-6,-3) -- (6,-3);

 \draw  [decoration={markings,mark=at position 0.7 with
{\arrow[scale=1.5,>=stealth]{>}}},postaction={decorate},
line width=.2mm] (-6,-3) -- (0,6);

\draw  [decoration={markings,mark=at position 0.7 with
{\arrow[scale=1.5,>=stealth]{>}}},postaction={decorate},
line width=.2mm] (0,6) -- (6,-3);

\draw  [decoration={markings,mark=at position 0.7 with
{\arrow[scale=1.5,>=stealth]{>}}},postaction={decorate},
line width=.2mm] (-6,-3) -- (0,0);

 \draw  [decoration={markings,mark=at position 0.7 with
{\arrow[scale=1.5,>=stealth]{>}}},postaction={decorate},
line width=.2mm] (0,6) -- (0,0);

\draw  [decoration={markings,mark=at position 0.7 with
{\arrow[scale=1.5,>=stealth]{>}}},postaction={decorate},
line width=.2mm] (0,0) -- (6,-3);

\node at (-6.8, -2.5){$\wt w_i$};
\node at (7, -2.5){$\wt w_k$};
\node at (1,6.5){$\wt w_j$};
\node at (1,0){$w_l$};
\node at (3,3){$\wt \alpha$};
\node at (-3,3){$\wt \beta$};
\node at (0,-4){$\wt \gamma$};
\node at (-3,-0.8){$\wt \delta$};
\node at (-0.5,3){$\wt \eps$};
\node at (3,-0.8){$\wt \zeta$};

\node at (0, -1.5){$\Downarrow$};
\node[rotate=330] at (-2, 1){$\Downarrow$};
\node[rotate=315] at (2,1) {$\Downarrow$};

\node at (-7, -4){$\bullet$};
\node at (7, -4){$\bullet$};
\node at (0,7.5){$\bullet$};

\draw [line width=0.2mm] (-7, -4) --  (-6,-3);
\draw [line width=0.2mm]  (0,6) -- (0,7.5);
\draw [line width=0.2mm]  (6,-3) -- (7, -4);

\node at (-8, -4){$w_i$};
\node at (8, -4){$w_k$};
\node at (0.9, 7.8){$w_j$};
\etp

\caption{The $3\to 1$ relation for PL arrows and proof by pinching.} \label{fig:PL-3:1}
\end{figure}

\paragraph{ Picard-Lefschetz arrows: the $3\to 1$ move.} Geometrically, associativity and
coassociativity of the Picard-Lefschetz arrows represent the  elementary
$2\to 2$ move of plane triangulations: replacing one of the two possible
triangulations of a $4$-gon into two triangles by the other triangulation.
We have also a companion statement for the $3\to 1$ move: passing
from a triangulation of a triangle into $3$ small triangles, to the ``big'' triangle itself,
see the left part of Fig. \ref{fig:PL-3:1}.

\begin{prop}\label{prop:PL-3:1}
In the situation of Fig. \ref{fig:PL-3:1} , the pasting of the natural transformations
$u_{\zeta, \eps, \alpha}$, $v_{\eps,\beta,\delta}$ and $v_{\zeta,\delta, \gamma}$
is equal to $0$.
\end{prop}

Note that the structure of a schober does not give a non-trivial natural transformation
associated to the ``big'' triangle.

\noindent{\sl Proof:} As in the proof of coassociativity, we pinch the picture at the points
$w_i, w_j, w_k$, see the right part of Fig. \ref{fig:PL-3:1}.
So it suffices to prove that the pasting of the ``inner'' triangle on the right is zero.
In this inner triangle, we have   only one singular point, $w_l$
and  two paths from    $\wt w_i$ to $\wt w_k$,
one above and one below  $w_l$, i.e., we have a Picard-Lefschetz situation.
The first of these is  the composite path
$(\wt\alpha\wt\beta)$  and the second is $\wt\gamma$.
So we have a Picard-Lefschetz triangle, denote it $\Delta$, corresponging
to this Picard-Lefschetz situation.

Note that the pasting of the two  small triangles on the left and on the right
(the lower triangle excluded) is just the Picard-Lefschetz arrow
$u_{\wt\zeta, \wt\delta, (\wt\alpha\wt\beta)}$ involving the composite path.
The the pasting of all three triangles is the composition of
$u_{\wt\zeta, \wt\delta, (\wt\alpha\wt\beta)}$ with $v_{\wt\zeta,\wt\delta, \wt\gamma}$,
which is the composition of two successive morphisms in $\Delta$, so it is zero. \qed


\subsection{Reminder on semi-orthogonal decompositions and gluing functors}
\label{subsec:rem-SOD}

\paragraph{$2$-term SOD's.} We recall the material from \cite{BK-serre}, see also
\cite{huybrechts}, Ch.1.
Let $\Vc$ be a $k$-linear triangulated category. A full triangulated subcategory $\Bc\subset\Vc$
is called {\em strictly full}, if any object of $\Vc$ isomorphic to an object of $\Bc$, is in $\Bc$.
In this section, all subcategories of $\Vc$ will be assumed triangulated and strictly full.

For a class of objects $I\subset \Ob(\Vc)$ we denote by $\< I\>$ the smallest
strictly full triangulated subcategory of $\Vc$ containing $I$.

Given a subcategory $\Bc\subset\Vc$, its {\em left} and {\em right orthogonals}
are the subcategories defined as follows:
\be
\begin{gathered}
^\perp \Bc \,=\, {^\perp}\Bc_\Vc \, := \,\bigl\{ V\in\Vc\, \bigl| \, \Hom(V,B)=0, \,\,\forall \, B\in\Bc\bigr\},
\\
\Bc^\perp  \,=\, \Bc_\Vc^\perp \, := \,\bigl\{ V\in\Vc\, \bigl| \, \Hom(B,V)=0, \,\,\forall \, B\in\Bc\bigr\}.
\end{gathered}
\ee
We can further form the iterated orthogonals $\Bc^{\perp 2} = (\Bc^{\perp})^\perp$,
$\Bc^{\perp (-2)} \,=\,  {^\perp}({^\perp}\Bc)$ etc. General notation: $\Bc^{\perp n}$, $n\in \ZZ$.

\begin{Defi}\label{def:SOD-2}
(a) A pair $(\Bc,\Cc)$ of subcategories in $\Vc$ is called {\em a semi-orthogonal decomposition}
(SOD for short)
of $\Vc$, if $\Hom(\Bc,\Cc)=0$ (i.e., $\Cc\subset\Bc^\perp$) and each object $V\in \Vc$
can be included into an exact triangle
\[
B_V \lra V\lra C_V\lra B_V[1], \quad B_V\in\Bc, \,\, C_V\in \Cc.
\]
(b) A subcategory $\Bc\subset\Vc$ is called {\em left admissible}, if $({^\perp}\Bc, \Bc)$
is an SOD, and {\em right admissible}, if $(\Bc, \Bc^\perp)$ is an SOD.
\end{Defi}

\begin{rems}
(a) The condition of existence of a triangle in part (a) of Definition \ref{def:SOD-2} is equivalent
to $\Vc=\< \Bc, \Cc\>$, i.e., to $\Vc$ being generated by $\Bc$ and $\Cc$ as a triangulated
category. So it is standard to speak about a  ``semiorthogonal decomposition  $\Vc=\< \Bc, \Cc\>$''.

\vskip .2cm

(b) For an SOD  $\Vc=\< \Bc, \Cc\>$, we have $\Cc=\Bc^\perp$ and $\Bc={^\perp}\Cc$.

\vskip .2cm

(c) The triangle in Definition \ref{def:SOD-2}
(a)  in fact depends functorially on $V$, which is reflected in the following fact.
\end{rems}

\begin{prop}
$\Bc\subset \Vc$ is left admissible (resp. right admissible) if and only if the embedding functor
$\eps:\Bc\hra\Vc$ admits a left adjoint $^*\eps$ (resp. right adjoint $\eps^*$). \qed
\end{prop}

For example, for a right admissible $\Bc$ we have $\eps^*(V) = B_V$, the $\Bc$-component
of the canonical triangle for $V$. Indeed, $\Hom_\Vc(\eps(B),V)=\Hom_\Bc(B, B_V)$,
since $\Hom_\Vc(B, C_V)=0$.

\vskip .2cm

We further say that $\Bc$ is {\em $\oo$-admissible}, if all iterated orthogonals
$\Bc^{\perp n}$, $n\in\ZZ$, are left admissible (they are then right admissible as well).
We have the action of the group $\ZZ$ on the set of $\oo$-admissible subcategories in
$\Vc$ generated by the transformation $\Bc\mapsto\Bc^\perp$ with inverse
being $\Bc\mapsto{^\perp}\Bc$.

\begin{prop}
Suppose $\Vc$ has a Serre functor $S$ (in particular, all Hom-spaces in $\Vc$
are finite-dimensional over $\k$). Then for any left and right admissible $\Bc\subset\Vc$
we have $S(\Bc)=\Bc^{\perp\perp}$ and $\Bc$ is $\oo$-admissible.\qed
\end{prop}

For example, if $B\in\Bc$ and $C\in\Bc^\perp$,  then $\Hom(C, S(B))=\Hom(B,C)^*=0$
so $S(B)\subset V^{\perp\perp}$.

\paragraph{$N$-term SOD's.}
Let $\Vc$ be as before.

\begin{Defi}
A sequence $\Bc_\bullet = (\Bc_1, \cdots, \Bc_N)$ of (strictly full, triangulated) subcategories in $\Vc$
is said to form a semi-orthogonal decomposition (SOD) of $\Vc$, if:
\begin{itemize}
\item[(1)] The $\Bc_i$ are semi-orthogonal,  i.e., $\Hom(\Bc_i, \Bc_j) =0$ for $i<j$.

\item[(2)] $\Vc=\< \Bc_1,\cdots, \Bc_N\>$ is generated by the $\Bc_i$.
\end{itemize}
\end{Defi}
Thus, an $N$-term SOD can be seen, in several ways, as a nested sequence of $2$-term SOD's:
\[
\Vc\,=\, \< \cdots \< \Bc_1, \Bc_2\>, \Bc_3\>, \cdots , \Bc_n\>  \,=\,
\< \Bc_1, \< \Bc_2, \cdots, \< \Bc_{N-1}, \Bc_N\> \cdots \>.
\]
Of considerable interest for us will be the action of the braid group $\Br_N$ on the set of
appropriate $N$-term SOD's. Given an SOD $\Vc= \< \Bc_\bullet \> = \< \Bc_1, \cdots, \Bc_n\>$,
and $i=1,\cdots, N-1$, we define its {\em right} and {\em left mutation at the spot}
$i$ to be the sequences of subcategories
\be\label{eq:mut}
\begin{gathered}
\tau_i\Bc_\bullet \,=\, \bigl( \Bc_1,\cdots, \Bc_{i-1}, \Bc_{i+1}, (\Bc_{i+1}^\perp)_{\< \Bc_i, \Bc_{i+1}\> },
\Bc_{i+2}, \cdots, \Bc_N\bigr),
\\
\tau_i^{-1}\Bc_\bullet  \,=\, \bigl( \Bc_1,\cdots, \Bc_{i-1}, ({^\perp}\Bc_i)_{\< \Bc_i, \Bc_{i+1}\> },
 \Bc_{i+1},
\Bc_{i+2}, \cdots, \Bc_N\bigr).
\end{gathered}
\ee
These sequences are semi-orthogonal, but may not generate $\Vc$ so may not be SODs.

\begin{Defi}\label{eq:oo-admis}
(a) A SOD $\Vc=\< \Bc_\bullet \>$ is called $\oo$-{\em admissible}, if all the sequences obtained from
$\Bc_\bullet$ by iterated application of the $\tau_i$ and $\tau_i^{-1}$, $i=1,\cdots, N-1$,
are SOD's.

\vskip .2cm

(b) An {\em $\oo$-admissible filtration} on $\Vc$ is a sequence
\[
F \,=\, \bigl( F_1\Vc \subset F_2\Vc \subset \cdots,\subset F_N\Vc = \Vc\bigr)
\]
of strictly full triangulated subcategories such that
\[
F_i\Vc \,=\, \< \Bc_1,\cdots, \Bc_i\>
\]
for some $\oo$-admissible SOD $\Vc=\<\Bc_1,\cdots, \Bc_N\>$.
\end{Defi}
Note that $\Bc_i$ in the definition is recovered uniquely as $(F_{i-1}\Vc)^{{^\perp}}_{F_i\Vc}$.
Thus considerations of $\oo$-admissible filtrations and $\oo$-admissible SOD's
are equivalent.
The following is well known \cite{BK-serre}.

\begin{prop}
The group $\Br_N$ acts on the set of $\oo$-admissible $N$-term SOD's of $\Vc$
(or, equivalently, on the set of $\oo$-admissible $N$-term filtrations)
so that the
generators $\tau_i$ and their inverse act by \eqref{eq:mut}. \qed
\end{prop}

Further, recall (\cite{kassel-turaev}, Thm. 1.24) that $\Br_N$ has a center which,
for $N\geq 3$,  is freely generated by $\delta^2$, where
$\delta$ is the ``Garside discriminant'' element
\be\label{eq:garside-dis}
\delta \,=\,   \tau_1 (\tau_2\tau_1) (\tau_3\tau_2\tau_1) \cdots (\tau_{N-1}\cdots\tau_1).
\ee
The following is also well known.
\begin{prop}\label{prop:delta-mutation-SOD}
(a) We have
\[
\delta(\Bc_1,\cdots, \Bc_N) \,=\,\bigl( \Bc_N = (\Bc_N^{\perp\perp})_{\Bc_N}, \,
(\Bc_{N-1}^{\perp\perp})_{\<\Bc_{N-1}, \Bc_N\>}, \, (\Bc_{N-2}^{\perp\perp})_{\<\Bc_{N-2}, \Bc_{N-1}, \Bc_N\>},\,
\cdots,  \Bc_1^{\perp\perp}\bigr).
\]
(b) We also have
\[
\delta^2(\Bc_1,\cdots, \Bc_N) \,=\, (\Bc_1^{\perp\perp}, \Bc_2^{\perp\perp},\cdots, \Bc_N^{\perp\perp}).
\]
 (c)
 Assume that $\Vc$ has a Serre functor $S$. Then,
 \[
 \delta^2(\Bc_1,\cdots, \Bc_N) \,=\, (S(\Bc_1), \cdots, S(\Bc_N)).
 \]

\end{prop}

\noindent{\sl Proof:} Part (a) is proved directly by induction on $N$. To see (b), notice that (a)
used twice
implies the statement for the first component:  $(\delta^2\Bc_\bullet)_1 = \Bc_1^{\perp\perp}$.
To deduce it for any other component, say the $i$th, let $\tau_{(1i})=\tau_{i-1}\tau_{i-2}\cdots\tau_1$
be the minimal lift to $\Br_N$ of the transposition $(1i)$. Then, using the fact that $\delta^2$ is
central, we have
\[
(\delta^2\Bc_\bullet)_i = (\tau_{(1i)}\delta^2\Bc_\bullet)_1 = (\delta^2\tau_{(1i)}\Bc_\bullet)_1
=
((\tau_{(1i)}\Bc_\bullet)_1)^{\perp\perp} = \Bc_i^{\perp\perp}.
\]
Finally, (c) follows from (b). Indeed,  for any $\oo$-admissible subcategory $\Bc\subset \Vc$
 we have $\Bc^{\perp\perp} = S(\Bc)$, which is obvious from
\[
\Hom(B, C) \,=\,\Hom(C, S(B))^*, \quad C\in \Bc^\perp. \qed
\]

\begin{conven}
Note that the concept of an SOD can be developed
purely in the framework of classical triangulated categories: no enhancements are needed.
We will use it in the enhanced framework as well: by an SOD in a pre-triangulated category $\Vc$
we will mean an SOD in the triangulated category $H^0(\Vc)$. Since $\Vc$ and $H^0(\Vc)$
have the same objects, this does not cause confusion.
\end{conven}

\paragraph{Gluing $2$-term SODs: gluing functors.}
One can ask: given triangulated categories $\Bc_1,\cdots, \Bc_N$,
what extra data are needed to construct an SOD $\Vc=\< \Bc_1, \cdots, \Bc_N\>$? An obvious invariant
of such an SOD is the collection of functors
\[
h_{ij}: \Bc_i^\op\times\Bc_j\lra\Vect_\k,\quad h_{ij}(B_i, B_j) =\Hom_\Vc(B_i, B_j), \,\,\, i>j.
\]
However, it does not seem possible to use just the data of $h_{ij}$ as triangulated functors to recover (``glue'')
$\Vc$ as a triangulated category.

On the other hand, it is well known how to perform such gluing in the dg-enhanced setting. In this paragraph
and the next subsection
we recall this gluing procedure, starting from the case of $2$-term SODs \cite{kuz-lun}.

\vskip ,2cm

Let $\Bc, \Cc$ be pre-triangulated dg-categories and $f: \Cc \to\Bc$ is a dg-functor. Define a new dg-category
$\Bc\times_f \Cc$ to have, as objects, triples $V=(B, C,u)$, where $B\in\Bc$, $C\in\Cc$ and $u: f(C)\to B$ is
a degree $0$ closed morphism in $\Bc$. The Hom-complex between $V= (B,C,u)$ and
$ V'= (B', C', u')$ is
defined, as a graded vector space to be
\be\label{eq:hom-BFC}
\Hom^\bullet_{\Bc \times_{{\hskip -0.1cm} f}\,  \Cc} (V,V') \,=\, \Hom^\bullet_\Bc(B,B') \oplus \Hom_{\Cc}(C, C')
\oplus \Hom_{\Bc} (f(C), B')[-1],
\ee
with the differential defined in such a way that we have an
exact triangle
\[
\begin{gathered}
\Hom^\bullet_{\Bc \times_{{\hskip -0.1cm} f}\,  \Cc} (V,V') \lra
{
\Hom^\bullet_\Bc (B, B') \oplus
\atop
\oplus \Hom^\bullet_\Cc (C,C')
}
\buildrel \on{dis}\over
\lra \Hom^\bullet_\Bc(f(C), B') \lra  \Hom^\bullet_{\Bc \times_{{\hskip -0.1cm} f}\,  \Cc} (V,V')[1].
\\
\end{gathered}
\]
Here $\on{dis}$  is the obvious ``discrepancy'' map, measuring the failure of a pair of morphisms
$B\to B'$, $C\to C'$ to be compatible with $u$ and $u'$,
see \cite{kuz-lun} \S 4.1 for more details. Informally, putting
$ \Hom^\bullet_{\Bc \times_{{\hskip -0.1cm} f}\,  \Cc} (V,V')$ on the left of $\on{dis}$ in the triangle, is the homotopy
analog of taking its kernel, so we enforce compatibility at the homotopy level.

\begin{prop}[\cite{kuz-lun}]\label{prop:BfC}
(a) $\Bc \times_f \Cc$ is pre-triangulated.

(b) We have a semi-orthogonal decomposition $\Bc\times_f \Cc =\< \Bc, \Cc\>$, with $\Bc\subset\Bc\times_f \Cc$
consisting of  the $(B,0,0)$ and $\Cc$ consisting of the $(0,C,0)$.

(c) The functor
\[
h^\bullet: \Cc^\op\times \Bc \lra \dgVect, \quad h^\bullet(C,B) \,=\, \Hom^\bullet_{\Bc \times_{{\hskip -0.1cm} f}\,  \Cc}(C,B)
\]
is given by
\[
h^\bullet(C,B) \,=\, \Hom^\bullet_\Bc(f(C), B) [-1]. \qed
\]
\end{prop}

\vskip .2cm

Assume now that $f$ has a right adjoint $f^*: \Bc\to\Cc$, so we have the
counit $\eta: f\circ f^*\to\Id_\Bc$. For an object $B\in\Bc$ consider the object
\[
k (B) \,=\,  (B, f^*B, \eta_B: f(f^*B)\to B) \,\in \, \Bc\times_f \Cc.
\]

\begin{prop}\label{prop:mut-f*}
(a)  The correspondence $B\mapsto k(B)$ extends to a dg-functor $k: \Bc\to
{\Bc \times_{{\hskip -0.1cm} f}\,  \Cc}$ which gives quasi-isomorphisms on Hom-complexes.

\vskip .2cm

(b) The category $k(\Bc)$ is identified with $\Cc^\perp=\Bc^{\perp\perp}\subset {\Bc \times_{{\hskip -0.1cm} f}\,  \Cc}$,
and we have an SOD ${\Bc \times_{{\hskip -0.1cm} f}\,  \Cc}= \< \Cc, k(\Bc)\>$.

\vskip .2cm

(c) This extends to an equivalence
\[
{\Cc \times_{{\hskip -0.1cm} f^*}\,  \Bc}\lra {\Bc \times_{{\hskip -0.1cm} f} \,  \Cc}
\]
which takes $\Cc\subset  {\Cc \times_{{\hskip -0.1cm} f^*}\,  \Bc}$ to $\Cc\subset  {\Bc \times_{{\hskip -0.1cm} f}\,  \Cc}$
by the identity functor, and $\Bc\subset  {\Cc \times_{{\hskip -0.1cm} f^*}\,  \Bc}$ to
$k(\Bc) \subset  {\Bc \times_{{\hskip -0.1cm} f}\,  \Cc}$ by the functor $k$.
\end{prop}

\noindent{\sl Proof:} (a)
Notice that for any $B,B'$ we have a canonical embedding of complexes
\[
\Hom^\bullet_\Bc(B,B') \hra \Hom^\bullet_{\Bc \times_{{\hskip -0.1cm} f}\,  \Cc}(k(B), k(B')),
\]
which is, moreover, a quasi-isomorphism. Indeed, by \eqref{eq:hom-BFC},
\[
\Hom^\bullet_{\Bc \times_{{\hskip -0.1cm} f}\,  \Cc}(k(B), k(B'))\,=\,
\Hom_{\Bc}(B,B') \oplus \Hom_\Cc(f^*B, f^* B') \oplus \Hom_\Cc(f^*B, f^*B')[-1])
\]
as a graded vector space. Looking at the differentials, we see that the embedding of the LHS as
the first summand in the RHS is a morphism (embedding) of complexes, and the quotient by the image
is contractible (is the cone of the identity map of the second summand in the RHS).
This defines a dg-functor $k$ with the required property.

\vskip .2cm

(b) Let $C\in\Cc$ and $B\in\Bc$. By definition \eqref{eq:hom-BFC},
\[
\Hom_ {\Bc \times_{{\hskip -0.1cm} f}\,  \Cc}(C, k(B)) \,=\, \Hom_\Cc(C, f^*B) \oplus \Hom_\Bc(f(C), B)[-1]
\]
as a graded vector space. Looking at the differential, we recognize this as the cone of
the adjunction (quasi-) isomorphism) so it is acyclic. This shows an inclusion $k(\Bc)\subset \Cc^\perp$.
To show that this inclusion is an equality and to establish an SOD, it suffices to show that every object
$(B,C, f(C)\buildrel u\over\to B)$ fits into an exact triangle with an object from $\Cc$ on the left and an
object of $\k(\Bc)$ on the right. For this, consider the following closed degree $0$ morphism in
$\Bc\times_f \Cc$, written vertically, with the third component,  as in \eqref{eq:hom-BFC}, equal to $0$:
\[
\xymatrix{
 \bigl(B\ar[d]_{\Id} , & C\ar[d]^{\on{adj}(u)}, & f(C) \ar[r]^u & B\bigr)
\\
\bigl(B, & f^*B, & f(f^*B)\ar[r]_{\hskip 0.5cm \eta_B}& B\bigr)
}
\]
 Its target lies in $k(\Bc)$ and its cone is of the form $(0, C',0)$, so lies in $\Cc$.

\vskip .2cm

(c)
This follows from the identification
\[
\Hom^\bullet_{\Bc \times_{{\hskip -0.1cm} f}\,  \Cc}(k(B), C) \, =  \, \Hom^\bullet_\Cc(f^*(B), C)
\]
which determines the gluing functor to be $f^*$, cf.  Proposition \ref{prop:BfC}(c).
\qed

\begin{cor}
If $f$ admits iterated adjoints of all orders, then the subcategory $\Bc\subset {\Bc \times_{{\hskip -0.1cm} f}\,  \Cc}$
is $\oo$-admissible. \qed
\end{cor}


\subsection{Unitriangular (co)monads and glued  $N$-term SOD's}\label{subsec:unitri-mon}

\paragraph{Unitriangular (co)monads and their (co)algebras.}
We will use a version of the classical categorical concept of a monad \cite{maclane}.

\begin{Defi} \label{def:monad}
Let $\Vc=(\Vc_1, \cdots, \Vc_N)$ be a finite ordered sequence of  pre-triangulated dg-categories.

\vskip .2cm

(a)
A {\em unitriangular (dg-)monad} on $\Vc$  is  a collection $\M=(\M_{ij}: \Vc_i\to\Vc_j)_{i<j}$ of dg-functors
together with closed, degree $0$ natural transformations, called {\em composition maps}
\[
c_{ijk}: \M_{jk} \circ \M_{ij} \lra \M_{ik}, \quad i<j<k,
\]
 such that for any $i<j<k<l$ the diagram below commutes (associativity condition):
\[
\xymatrix{
\M_{kl}\circ \M_{jk}\circ \M_{ij}
\ar[d]_{c_{jkl}\circ \M_{ij}}
\ar[rr]^{\M_{kl}\circ c_{ijk}} && \M_{kl}\circ\ M_{ik}\ar[d]_{c_{ikl}}
\\
\M_{jk}\circ \M_{ij} \ar[rr]_{c_{ijl}}&& M_{il}.
}
\]
We put $\M_{ii}=\Id_{\Vc_i}$.

\vskip .2cm

(b) An $\M$-algebra is a sequence of objects $V=(V_1,\cdots, V_N)$, $V_i\in\Vc_i$, together with closed, degree $0$
morphisms {\em (action maps)}
\[
\alpha_{ij}: \M_{ij}(V_i) \lra V_j, \quad i<j,
\]
such that for any $i<j<k$ the diagram below commutes:
\[
\xymatrix{
\M_{jk}(\M_{ij}(V_i))
\ar[d]_{c_{ijk, V_i}}
\ar[rr]^{\M_{jk}(\alpha_{ij})} && \M_{jk}(V_j) \ar[d]^{\alpha_{jk}}
\\
\M_{ik}(V_i) \ar[rr]_{\alpha_{ik}}&& V_k.
}
\]
\end{Defi}

\begin{rem}
Note that any dg-functor commutes with direct sums.  So, compared to the general concept of
\cite{maclane},  our unitriangular dg-monads $\M$ are of ``linear'' nature and rather resemble associative dg-algebras (in fact,
they can be viewed as ``algebras in the category of functors''). Accordingly,  it would be more natural  to
refer to $\M$-algebras as ``modules'', but we retain the classical terminology of \cite{maclane}.
\end{rem}

Any collection of objects  $U=(U_1,\cdots, U_N)$, $U_i\in V_i$, gives rise to two $\M$-algebras: the
 {\em free $\M$-algebra}
$\M(U)$ and the {\em reduced free $\M$-algebra} $ \M^< (U)$, with
\be
\M(U)_j \,=\,\bigoplus_{i\leq j} \M_{ij}(U_i), \quad  \M^<(U)_j \,=\, \bigoplus_{i<j} \M_{ij}(U_i),
\ee
and the action maps given by the composition maps in $\M$.

\begin{Defi}\label{def:comonad}
(a) A {\em unitriangular (dg-)comonad} $\E$ in $\Vc$ is a unitriangular dg-monad in
$\Vc^\op:=(\Vc_N^\op,\cdots,
\Vc_1^\op)$. Explicitly, $\E$ consists of dg-functors $\E_{ij}: \Vc_i\to\Vc_j$, $i<j$, together with
closed degree $0$ dg-functors ({\em cocomposition maps})
\[
\eps_{ijk}: \E_{ik}\lra \E_{jk}\circ \E_{ij}, \quad i<j<k,
\]
satisfying the coassociativity condition dual to that in Definition \ref{def:monad}(a).
We put $\E_{ii}=\Id_{\Vc_i}$.
To keep notation straight, we will sometimes denote by $\E^\op$ the unitriangular dg-monad in $\Vc^\op$
corresponding to $\E$.

\vskip .2cm

(b) A coalgebra over $\E$ us a collection of objects $V_i\in \Vc_i$ together with closed degree
$0$ morphisms ({\em coaction maps})
\[
\beta_{ij}: V_j \lra \E_{ij}(V_i), \quad i<j,
\]
satisfying the compatibility condition dual to that in Definition \ref{def:monad}(b).
\end{Defi}

Any collection of objects $U = (U_1,\cdots, U_N)$, $U_i\in\Vc_i$, gives rise to  the
{\em free $\E$-coalgebra} $\E(U)$ and the {\em restricted free $\E$-coalgebra $ \E^<(U)$} on $U$
with
\be
\E(U)_j \,=\,\bigoplus_{i\leq j} \E_{ij}(U_i), \quad  \E^<(U)_j \,=\, \bigoplus_{i<j} \E_{ij}(U_i),
\ee
and the coaction maps given by the cocomposition maps in $\E$.

\vskip .2cm

If $\M$ is a unitriangular monad in $(\Vc_1,\cdots, \Vc_N)$, then the collections of the right and
left adjoints $\M^*=(\M_{ij}^*)_{i<j}$ and $^*\M= ({^*\M}_{ij})_{i<j}$ (provided they exist), are
unitriangular comonads in $(\Vc_N,\cdots, \Vc_1)$.  Further, $\M$-algebras are the same as
$\M^*$-coalgebras.

\paragraph{Dg-categories of (co)algebras and SODs.}
Let $\Vc=(\Vc_1,\cdots,\Vc_N)$ be as before. For two sequences of objects
$V=(V_1,\cdots, V_N)$ and $V'=(V'_1,\cdots, V'_N)$, $V_i, V'_i\in\Vc_i$, we define
\[
\Hom^\bullet_\Vc(V, V') \,=\,\bigoplus_{i=1}^N \Hom^\bullet_{\Vc_i}(V_i, V'_i).
\]
In this way we view $\Vc$ as a dg-category (i.e., as the direct sum of the $\Vc_i$).

\vskip .2cm

Let now $\M$ be a unitriangular dg-monad in $\Vc$ and $V, V'$ be $\M$-algebras.
We define the {\em naive complex of morphisms} between $V$ and $V'$ as
\[
\Hom^\bullet_\naive(V, V') \,=\,\bigl\{ \phi\in\Hom^\bullet_\Vc(V, V') \bigl| \, \phi \text{ commutes with the actions}\bigr\}.
\]
This definition is not homotopy invariant and we take the derived functor of it in the following explicit way.

\vskip .2cm
For any $\M$-algebra $V$ we consider the {\em reduced bar-resolution} of $V$ which is the
  complex  of free $\M$-algebras (and differentials being closed, degree $0$ naive morphisms)
\be\label{eq:bar-res}
\xymatrix
 {
\Bar^\bullet(V)\,=\,\biggl\{
\cdots \ar[r] & \M( \M^< (\M^<(V))) \quad
\ar[rrr]^{c_{ \M^<(V)} - \M(c_V) + \M^2(\alpha)} &&&  {\M( \M^<(V))} \ar[rr]^{c_V-\M(\alpha)}&& \M(V) \biggr\}.
}
\ee
The differential  $d$ is the alternating sum of elementary contractions using $c = (c_{ijk})$ or
$\alpha=(\alpha_{ij})$, so $d^2=0$ in virtue of
the axioms.  Note that this complex terminates since $(\M^<)^N(V)=0$. We fix the grading of the complex
so that $\M(V)$ is in degree $0$. Let $B(V)=\Tot(\Bar^\bullet(V))$ be
the total object.

\begin{prop}\label{prop:augm-bar}
The action map $\alpha: \M(V)\to V$ gives rise to a quasi-isomorphisms $B(V)\to V$.
\end{prop}

\noindent {\sl Proof:} We need to show that the augmented complex
\[
\cdots \to \M(\M^<( \M^<(V))) \to \M( \M^<(V)) \to \M(V) \buildrel \alpha\over\lra V,
\]
denote it $B^+$, is exact. This means that for each $j=1,\cdots, N$, the $j$th component $B^+_j$ of $B$
is exact. Explicitly, $B^+_j$ is the complex (in $\Vc_j$)
\[
\cdots \to \bigoplus_{i_1< i_2\leq j} \M_{i_2,j}(\M_{i_1, i_2}(V_{i_1}) \to
\bigoplus_{i\leq j} \M_{ij}(V_i) \to V_j.
\]
Thus,  $B^+_j$ as a graded object of $\Vc_j$, consists of summands
\[
S(i_1, \cdots, i_m) \,=\, \M_{i_m,j} \cdots \M_{i_1, i_2}(V_{i_1})[m],
\]
labelled by sequences $i_1<\cdots < i_m \leq j$. For $p\geq 0$ let $F_p\subset B^+_j$ be the sum of
the $S(i_1,\cdots i_m)$ such that
$\#\{\nu| \,  i_\nu <  j\} \leq p$. It is immediate that $F_p$ is closed under the differential, i.e., is
a subcomplex. So we have a finite increasing filtration $F_0\subset F_1\subset \cdots\subset F_j = B^+_j$
by subcomplexes. Now, each quotient $F_p/F_{p-1}$ is a $2$-term complex with both terms identified
with
\[
\bigoplus_{i_1<\cdots < i_p < j} \M_{i_p, j} \cdots \M_{i_1, i_2}(V_{i_1})
\]
and the differential being the identity. So each $F_p/F_{p-1}$ is exact and therefore $B^+_j$ is exact. \qed

\begin{Defi}
Let $\Vc=(\Vc_1,\cdots, \Vc_N)$ be as before.

\vskip .2cm

(a) For a unitriangular dg-monad $\M$ in $\Vc$ we define the dg-category $\Alg_\M$
to have, as objects, $\M$-algebras and, as Hom-complexes,
\[
\Hom^\bullet_{\Alg_\M}(V,V') \,=\, \Hom^\bullet_\naive (B(V), B(V')).
\]
(b) For a unitriangular dg-comonad $\E$ in $\Vc$ we define the dg-category $\Coalg_\E$ as
$(\Alg_{\E^\op})^\op$,
where $\E^\op$ is the dg-monad in $\Vc^\op$ corresponding to $\E$, see Definition \ref{def:comonad}(a).
\end{Defi}

\begin{rems}
(a) Since $B(V)$ is {\em quasi-free} (i.e., it becomes a free algebra, if we forget differentials in both $B(V)$
and $\M$), the quasi-isomorphism $B(V')\to V'$ induces a quasi-isomorphism
\[
\begin{gathered}
\Hom_{\Alg_\M}^\bullet(V, V') \lra \Hom^\bullet_\naive (B(V), V') \,=\,
\\
=\, \Tot\, \biggl\{ \Hom^\bullet_\Vc(V, V') \buildrel \delta_0\over\lra \Hom^\bullet_\Vc( \M^<(V), V')
\buildrel \delta_1\over\lra \Hom_\Vc(\M^< ( \M^< (V)), V') \buildrel\delta_2\over\lra\cdots\biggr\}.
\end{gathered}
\]
Note that the first differential $\delta_0$ in the last complex is the ``discrepancy map''
\[
\delta_0(\phi) \,=\,\phi\circ\alpha -\alpha'\circ\phi:  \M^<(V)\to V'),
\]
where $\alpha$ and $\alpha'$ are the action maps for $V$ and $V'$. So
$\Ker(\delta_0)=\Hom^\bullet_\naive(V, V')$.

\vskip .2cm

(b) Explicitly, the category $\Coalg_\E$ has, as objects, $\E$-coalgebras, and
\[
\Hom^\bullet_{\Coalg_\E}(V, V')\,=\,\Hom^\bullet_\naive(C(V), C(V')).
\]
Here
\[
C(V)\,\,=\,\,\Tot \, \bigl\{ \E(V) \to\E( \E^<(V))\to\cdots \bigr\}
\]
is the total object of the {\em reduced cobar-resolution} of $V$, defined in the way dual  to
\eqref{eq:bar-res}, and $\Hom^\bullet_\naive$ is the complex of morphisms
commuting with the
coalgebra structures.
\end{rems}

\begin{prop}
(a) Let $\M$ be a quasi-triangular dg-monad in $\Vc=(\Vc_1,\cdots, \Vc_N)$.
The dg-category $\Alg_\M$ is pre-triangulated. It has an SOD $\Alg_\M=\< \Vc_1,\cdots, \Vc_N\>$,
where $\Vc_i$ is identified with the full dg-subcategory in $\Alg_\M$ formed by algebras $V$ with
$V_j=0$ for $j\neq i$.

\vskip .2cm

(b) Dually, let $\E$ be a quasi-triangular dg-comonad in $\Vc$. The category $\Coalg_\E$
is pretriangulated. It has an SOD $\Coalg_\E = \< \Vc_1,\cdots, \Vc_N\>$, where
$\Vc_i$ is identified with the  full dg-category of coalgebras $V$ with $V_j=0$ for $j\neq i$.
\end{prop}

\noindent{\sl Proof:} We prove only (a), since (b) is obtained by duality. We start by showing that
$\Alg_\M$ is pre-triangulated. Let $\Alg_\M^\naive$ be the dg-category with the same objects
as $\Alg_\M$ and Hom-complexes $\Hom^\bullet_\naive(V,V')$. Then $\Alg_\M^\naive$
is pre-triangulated: the total complex $\Tot(U^\bullet)$ of a twisted complex $U^\bullet$
over $\Alg_\M^\naive$ is found componentwise, i.e., separately for each $i=1,\cdots, N$.
The $i$th component of $\Tot(U^\bullet)$ is $\Tot(U^\bullet_i)$, the total object
in $\Vc_i$ of the twisted complex $U^\bullet_i$ over $\Vc_i$.
Now, $\Alg_\M$ is, by definition, the full dg-subcategory in $\Alg^\naive_\M$ on objects
of the form $B(V)$. To prove that it is pre-triangulated, it suffices to show the following.
If, in $U^\bullet$ above, each $U^p$ has the form $B(V^p)$, then $U=\Tot(U^\bullet)$
is {\em homotopy equivalent} to some $B(V)$ (i.e., becomes isomorphic to it
in $H^0(\Alg^\naive_\M)$). For this, we can take $V=U$ and look at the canonical
quasi-isomorphism $\ol\alpha: B(U)\to U$. By our assumptions, $U$ is quasi-free,
and therefore $\ol\alpha$ is a homotopy equivalence.

\vskip .2cm

Next, for $i=1,\cdots, N$, let $\Alg^i_\M$ (resp. $\Alg^{\leq i}_\M$) be
the full dg-subcategory in $\Alg_\M$ formed by algebras $V=(V_1,\cdots, V_N)$ such
that $V_j\neq 0$ only for $j=i$ (resp. for $j\leq i$). Then $\Alg^i_\M$ is isomorphic to $\Vc_i$.
Further, we have a $2$-term SOD $\Alg_\M^{\leq i} = \< \Alg_\M^{\leq i-1}, \Alg_\M^i\>$.
Indeed, for $V\in \Alg_\M^i$ and $V'=\Alg_\M^{\leq i-1}$ we have $\Hom_{\Alg_\M}*V,V')=0$
by the construction. At the same time, for any $V=(V_1, \cdots, V_i,0,\cdots, 0)\in\Alg_\M^{\leq i}$
we have an exact triangle
\[
(V_1, \cdots, V_{i-1},0,\cdots 0) \lra (V_1,\cdots, V_i,0,\cdots, 0) \lra (0,\cdots, 0, V_i, 0,\cdots, 0)
\]
so which shows that we have an SOD. In this way we show by induction that
\[
\Alg_\M\,=\, \< \Alg^1_\M,\cdots, \Alg^N_\M\>\,=\, \< \Vc_1, \cdots, \Vc_N\>.  \qed
\]

\begin{ex}
Let $N=2$. The datum of a uni-triangular monad $\M$ reduces in this case to the datum of a single dg-functor
$f= \M_{12}: \Vc_1\to\Vc_2$ and nothing else. The dg-category $ \Alg_\M $ for such an $\M$ is
precisely the glued dg-category $\Vc_1\times_f \Vc_2$, see Proposition \ref{prop:BfC}.
\end{ex}

\paragraph{Mutations of uni-triangular dg-monads.}
We now describe the analog, for monads, of the concept of mutations of exceptional sequences
and of  their Ext-algebras, see \cite{rudakov, seidel}. Note that the formulas for the mutations below
\eqref{eq:monad-left-mut} are similar to those of Proposition \ref{prop:GMV-mut}. This similarity
is a clear indication of the role of perverse schobers in the theory of mutations.

\vskip .2cm

Let  $\Vc=(\Vc_1,\cdots, \Vc_N)$ be
a sequence of pretriangulated dg-categories. For any $i=1,\cdots, N-1$ we denote
\[
s_i\Vc = (\Vc_1,\cdots, \Vc_{i-1}, \Vc_{i+1},\Vc_i, \Vc_{i+2},\cdots, \Vc_N)
\]
the sequence obtained from $\Vc$  by permuting the $i$th and $(i+1)$st terms.

\vskip .2cm
 Let $\M$ be a uni-triangular
dg-monad in $\Vc$. Assume that each $\M_{ij}$ has both adjoints ${^*\M}_{ij}$ and $\M^*_{ij}$.
We construct two uni-triangular dg-monads $\M'=L_i\M$ and $\M'' = R_i\M$ in $s_i\Vc$
which we call the {\em left} and {\em right mutations} of $\M$.  More precisely,  define the functors
$\M'_{\nu,j}$,
$1\leq \nu < j\leq N$, by
 \be\label{eq:monad-left-mut}
\M'_{\nu, j} =\begin{cases}
\M_{\nu, j},& \text { if } \nu<j\leq i-1;
\\
\Cone\bigl\{ \M_{i,i+1}\circ \M_{\nu,i} \buildrel c_{\nu, i, i+1}\over\lra \M_{\nu, i+1}\bigr\}[-1], &
\text{ if } \nu\leq i-1, \, j=i.
\\
\M_{\nu,i},& \text{ if } \nu\leq i-1, \, j=i+1;
\\
\M_{\nu,j},& \text {if } \nu\leq i-1, \, j\geq i+2;
\\
{^*\M}_{i,i+1}, &\text{ if } \nu=i, \, j=i+1;
\\
\Cone\bigl\{\M_{i+1,j} \buildrel {{^tc}_{i, i+1,j}}\over\lra \M_{ij}\circ ({^*\M}_{i,i+1})\bigr\},
&\text{ if } \nu=i, j\geq i+2;
\\
\M_{i,j},&\text{ if } \nu=i+1,\, j\geq i+2;
\\
\M_{\nu,j}, &\text{ if } i+2\leq\nu < j.
\end{cases}
\ee
Here  ${^tc}_{i, i+1,j}$
corresponds to the composition map $c_{i,i+1,j}: \M_{i+1,j}\circ \M_{i,i+1}\to \M_{i,j}$ of $\M$ by the cyclic rotation isomorphism $\Hom(F\circ G, H) \simeq \Hom (F, H\circ ({^*G}))$, see Corollary \ref{cor:tens-rot}.

\vskip .2cm

We now describe the compositions in $\M'$. Leaving aside the cases when they reduce to the
compositions in $\M$, we consider the following possibilities.

\vskip .2cm

For $\nu < i$ the composition $\M'_{i,i+1}\circ\ M'_{\nu,i} \to \M'_{\nu, i+1}$, which by our definitions, should
be a map
\[
\Cone\bigl\{  {^*\M} _{i,i+1}\circ \M_{i,i+1}\circ \M_{\nu,i} \to {^*\M}_{i, i+1}\circ \M_{\nu, i+1}\bigr \} \lra \M_{\nu,i},
\]
is induced by the following (vertical) map of a (horizontal)  $2$-term complex to a $1$-term complex
\[
\xymatrix{
{^*\M}_{i, i+1}\circ \M_{i,i+1}\circ \M_{\nu,i}\ar[d]_{\eta\circ \M_{\nu,i} }  \ar[r] & {^*\M}_{i, i+1}\circ \M_{\nu, i+1}\ar[d]
\\
\M_{\nu,i} \ar[r]& 0,
}
\]
where $\eta$ is the counit of the adjunction.

\vskip .2cm

For $\nu<i$ and $j>i+1$ the composition $\M'_{i,j}\circ \M_{\nu,i}\to \M'_{\nu,j}$ which, by our definitions,
should be a map
\[
\Cone\{\M_{i+1,j}\circ \M_{i,j}\circ{^*\M}_{i,i+1}\} \circ \Cone \{ \M_{i,i+1}\circ \M_{\nu,i}\to \M_{\nu, i+1}\}[-1]
\lra \M_{\nu,j},
\]
is obtained by considering the double complex
\[
\xymatrix{
\M_{i+1,j}\circ \M_{i,i+1}\circ \M_{\nu,i} \ar[r] \ar[d] & \M_{i+1,j}\circ \M_{\nu, i+1}\ar[d]
\\
\M_{i,j}\circ {^*\M}_{i,i+1}\circ \M_{i,i+1} \circ \M_{\nu,i} \ar[r]& \M_{i,j}\circ {^*\M}_{i,i+1}\circ \M_{\nu, i+1},
}
\]
then mapping it, using the counit of the adjunction and the compositions in $\M$, to the double complex
\[
\xymatrix{
\M_{\nu,j}\ar[r]^\Id\ar[d]_\Id & \M_{\nu,j}\ar[d]^\Id
\\
\M_{\nu,j}\ar[r]_\Id& \M_{\nu,j},
}
\]
and then noticing that the total object of the latter double complex is canonically homotopy equivalent to
 $\M_{\nu,j}$.

 \vskip .2cm

 Dually, we define the functors $\M''_{\nu,j}$ for $1\leq\nu < j\leq N$ by
 \[
 \M''_{\nu,j} = \begin{cases}
 \M_{\nu,j},&\text{ if } \nu<j\leq i-1;
 \\
 \M_{\nu, i+1}, &\text{ if }  j=i,\, \nu \leq i-1;
 \\
 \Cone\bigl\{ \M_{\nu,i} \buildrel c_{i,i+1,j}^t\over\lra \M_{i,i+1}^*\circ \M_{\nu, i+1}\bigr\}[-1], & \text{ if }  \nu\leq i-1, \, j= i+1;
 \\
 \M_{\nu,j},&\text{ if } \nu\leq i-1, \, j\geq i+2;
 \\
 \M^*_{i,i+1},&\text{ if } \nu=i,\, j=i+1;
 \\
 \M_{i+1,j},&\text{ if } \nu=i,\, j\geq i+2;
 \\
 \Cone\bigl\{ \M_{i+1,j}\circ \M_{i,i+1}\buildrel c_{i,i+1,j}\over\lra \M_{ij}\bigr\}, & \text { if } \nu=i+1, j\geq i+2;
 \\
 \M_{\nu,j}, & \text { if }    i+2\leq\nu < j.
 \end{cases}
 \]
 The composition maps for $\M''$ are similar to those for $\M'$ and we leave them to the reader.

 \begin{prop}\label{prop:mut-monads}
 Assuming that the relevant (iterated) adjoints exist, the following hold:

 \vskip .2cm

 (a) Both $L_i\M=\ M'$ and $R_i\M=\M''$ are unitriangular dg-monads in $s_i\Vc$.

 \vskip .2cm

 (b) We have canonical quasi-isomorphisms $L_i R_i \M \simeq R_i L_i \M \simeq \M$.

 \vskip .2cm

 (c) We have canonical quasi-isomorphisms
 \[
 L_i L_j \M \simeq L_j L_i \M, \,\, |i-j|\geq 2, \quad L_i L_{i+1} L_i \M \simeq L_{i+1} L_i L_{i+1}\M, \,\, i=1,\cdots,
 N-2,
 \]
 and similarly for the $R_i$.

 \vskip .2cm

 (d$'$) For any $\M$-algebra $V=(V_1,\cdots, V_N)$ in $\Vc$ the sequence
 \[
 \begin{gathered}
 \lambda_i(V) \,=\, \bigl( V_1,\cdots, V_{i-1}, L_{V_i}(V_{i+1}), V_i, V_{i+2},\cdots, V_N\bigr), \quad  \text{where}
 \\
 L_{V_i}(V_{i+1}) \, := \, \Cone\bigl\{ \M_{i,i+1}(V_i) \lra V_{i+1}\bigr\}[-1] \quad \text{ (``left mutation'')},
 \end{gathered}
 \]
is naturally an $L_i\M$-algebra. This defines
 a quasi-equivalence of pre-triangulated categories $\lambda_i: \Alg_\M\to\Alg_{L_i\M}$. The pullback, under
 $\lambda_i$, of the standard SOD in $\Alg_{L_i\M}$ is the SOD
 \[
 \tau_i \< \Vc_1,\cdots, \Vc_N\> \,=\, \bigl\langle \Vc_1, \cdots, \Vc_{i-1}, {(\Vc^\perp_i)}_{\< \Vc_i, \Vc_{i+1}\>},
 \Vc_i, \Vc_{i+2}, \cdots,\Vc_N\bigr\rangle
 \]
see \eqref{eq:mut}.

\vskip .2cm

(d$''$) Dually, for any $V$ as in (d$'$), the sequence
 \[
 \begin{gathered}
 \rho_i(V) \,=\, \bigl( V_1,\cdots, V_{i-1},  V_{i+1}, R_{V_{i+1}}(V_i),  V_{i+2},\cdots, V_N\bigr),
  \quad  \text{where}
 \\
 R_{V_{i+1}}(V_i) \, := \, \Cone\bigl\{  V_i \lra \M^*_{i, i+1}(V_{i+1})\bigr\}[-1] \quad \text{ (``right mutation'')},
 \end{gathered}
 \]
 is naturally an $R_i\M$-algebra. This defines
 a quasi-equivalence of pre-triangulated categories $\rho_i: \Alg_\M\to\Alg_{R_i\M}$. The pullback, under
 $\rho_i$, of the standard SOD in $\Alg_{R_i\M}$ is the SOD
 \[
 \tau^{-1}_i \< \Vc_1,\cdots, \Vc_N\> \,=\, \bigl\langle \Vc_1, \cdots, \Vc_{i-1}, \Vc_{i+1},
  {({^\perp \Vc}_{i+1})}_{\< \Vc_i, \Vc_{i+1}\>},
\Vc_{i+2}, \cdots,\Vc_N\bigr\rangle
 \]
see again \eqref{eq:mut}.

\end{prop}

\noindent{\sl Proof:} It is completely similar to the (now classical) arguments in the case of exceptional sequences
\cite{rudakov, seidel}, and we leave the details to the reader.  The complicated formulas for $\M'_{\nu, j}$
(resp. $\M''_{\nu, j}$) can be obtained at once by analyzing the Hom-complexes between elements of $\lambda_i(V)$
(resp. $\rho_i(V)$).
\qed

\begin{cor}
Suppose each $\M_{ij}$ has left and right adjoints of all orders. Then:
\vskip .2cm

(a)  Iteration of the $L_i, R_i $ define an action $\sigma\mapsto L_\sigma$ of the braid group $\Br_N$ on the
 category of unitriangular dg-monads, so that $L_{\tau_i}\M = L_i\M$, $L_{\tau_i^{-1}}\M=R_i\M$.
 \vskip .2cm

 (b) The SOD $\< \Vc_1,\cdots,\Vc_N\>$
in $\Alg_\M$ is $\oo$-admissible. \qed
\end{cor}

 For future reference let us mention another, more minor operation on monads and algebras.
 For $p_1,\cdots, p_N\in\ZZ$  and a unitriangular monad $\M$ on $\Vc$ as before,
 we have the {\em shifted monad} $\M[p_1,\cdots, p_N]$ on the same $\Vc$ defined by
 \be\label{eq:shift-monad}
 \M[p_1,\cdots, p_n]_{ij} \,=\, \M_{ij}[p_j-p_i],
 \ee
 with composition maps induced by those in $M$ in an obvious way.
Thus we have an equivalence of dg-categories
\be\label{eq:shift-alg-monad}
\Alg_\M \lra \Alg_{\M[p_1,\cdots, p_n]}, \quad (V_1,\cdots, V_N) \mapsto (V_1[p_1],\cdots, V_N[p_n]).
\ee

\paragraph{Koszul duality for monads.}

 Let $\Vc=(\Vc_1,\cdots, \Vc_N)$ be as before, and $G=(G_{ij}:\Vc_i\to\Vc_j)_{i<j}$
 be any collection of dg-functors. We can form the {\em free (unitriangular dg-) monad} $\FM(G)$
 and the  {\em free (unitriangular dg-) comonad} $\FC(G)$ generated by $G$. As dg-functors,
 \be\label{eq:free-op-coop}
 \begin{gathered}
 \FM(G)_{ij}\,=\,\FC(G)_{ij} \,=\,\, \bigoplus_{p\geq 0} \bigoplus_{i<i_1<\cdots < i_p< j} G_{i,i_1,\cdots, i_p,j},
 \\
  G_{i,i_1,\cdots, i_p,j}\,:= \, G_{i_p,j}\circ G_{i_{p-1}, i_p}\circ\cdots\circ G_{i,i_1}.
 \end{gathered}
 \ee
 The composition in $\FM(G)$ is the formal concatenation. The cocomposition $\eps_{ijk}: \FC(G)_{ik}\to
 \FC(G)_{jk}\circ\FC(G)_{ij}$, $i<j<k$, is given by ``splitting''. This means that the summand
 $G_{i,i_1,\cdots, i_p,k}$ is sent to $0$,  if no $i_\nu$ is equal to $j$ and to
 $G_{i_\nu,\cdots, i_p,k}\circ G_{i,i_1,\cdots, i_\nu}$ (by the identity map), if $i_\nu=j$.

 \vskip .2cm

 Let now $\M$ be a uni-triangular dg-monad in $\Vc$. Consider the free comonad $\FC(M[1])$ with its differential $d$.
 Let $D$ be the differential in  $\FC(M[1])$ whose matrix element $\M_{jk}\circ \M_{ij}\to \M_{ik}$ equals
 $c_{ijk}$
 for each $i<j<k$, and which is extended uniquely to the free coalgebra by compatibility with cocomposition.
 Then $D^2=0$ by associativity of the $c_{ijk}$ and $dD+Dd=0$ by compatibility of the $c_{ijk}$ with the
 differentials in $M$. So we have the uni-triangular dg-comonad
 \[
 \M^! \,=\,\bigl(\FC(\M[1]), d+D\bigr)
 \]
which we call the {\em cobar} (or {\em Koszul}) dual to $\M$. Explicitly, for $i<j$ we have
\[
\M^!_{ij} \,=\,\Tot\biggl\{ \cdots \buildrel D\over\to \bigoplus_{i<i_1<i_2<j} \M_{i_2,j}\circ \M_{i_1,i_2}\circ \M_{i,i_1}
 \buildrel D\over\to \bigoplus_{i<i_1<j} \M_{i_1,j}\circ \M_{i,i_1}  \buildrel D\over\to \M_{ij}\biggr\},
\]
where the (horizontal) grading is such that $\deg(\M_{ij})=-1$. The differential $D$ is the alternating sum
of adjacent composition maps.

\vskip .2cm

Dually, let $\E$ be a uni-triangular dg-comonad on $\Vc$. Consider the free monad $\FM(\E[-1])$ with its differential $d$.
Let $D$ be the differential in $\FM(\E[-1])$ which, on the generators, is given by the cocomposition maps
$\eps_{ijk}: \E_{ik}\to \E_{jk}\circ \E_{ij}$ and is extended to the free monad by the Leibniz rule. As before, $D^2=0$
and $dD+Dd=0$, so we have the uni-triangular dg-monad
\[
\E^\!\,=\,\bigl( \FM(\E[-1]), d_D\bigr),
\]
called the {\em bar-} or {\em Koszul dual} to $\E$. Explicitly,
\[
\E^\!_{ij} \,=\,\Tot\biggl\{ \E_{ij}\buildrel D\over\to \bigoplus_{i<i_1<j} \E_{i_1,j}\circ \E_{i,i_1}\buildrel D\over\to
\bigoplus_{i<i_1<i_2<j} \E_{i_2,j}\circ \E_{i_1, i_2}\circ\ E_{i,i_1} \buildrel D\over\to\cdots\biggr\},
\]
 with the horizontal grading such that $\deg(\E_{ij})=+1$.

\begin{prop}\label{prop:M-!!}
For any uni-triangular dg-monad  $\M$ we have a canonical quasi-isomorphism $\M^{!\!}\to \M$.
Dually, for any uni-triangular dg-comonad $\E$ we have a canonical quasi-isomorphism $\E\to \E^{ \! !}$.
\end{prop}

\noindent{\sl Proof:} This is completely analogous to the proof of the corresponding statement for
the bar-cobar duality for operads \cite{ginzburg-kapranov}, so we give just a bried sketch.
Looking at   $\M^{!\!}_{ij}$ as just a graded functor, we find that it is the direct sum of many shifted copies of
summands of the form $\M_{i,i_1,\cdots, i_p,j}$, see \eqref{eq:free-op-coop}, and there is
just one copy (which is $\M_{ij}$) for $p=0$.
The compositions in $\M$ define,
 for each $i<j$,   the projection, denote it $\pi_{ij}: \M^{!\!}_{ij}\to \M_{ij}$, and these projections
  give  a morphism of dg-monads $\pi: \M^{!\!}\to \M$. To prove that each $\pi_{ij}$
is a quasi-isomorphism, we consider a filtration $F$ in  $\M^{!\!}_{ij}$ with $F_m$ being the sum of
all summands corresponding to the $\M_{i,i_1,\cdots, i_p,j}$ with $p\leq m$. This is compatible with the
differential, and each quotient $F_m/F_{m-1}$ splits into the direct sum of dg-functors of the form
\[
\M_{i,i_1,\cdots, i_m,j} \otimes_\k C(i,i_1,\cdots, i_m,j),
\]
where $C(i,i_1,\cdots, i_m,j)$ is a combinatorial complex of $\k$-vector spaces which is easily found out to be exact
 for $p>0$.
\qed

\vskip .2cm

Next, let $\M$ be a quasi-triangulat dg-monad in $\Vc$ and $V=(V_1,\cdots, V_N)$ be an $\M$-algebra.
Consider the free $\M^!$-coalgebra $\M^!(V)$ with its differential $d$. As before, it has a unique differential $D$,
compatible with the coalgebra structure, whose matrix elements $\M_{ij}(V_i)\to V_j$, $i<j$m are given
by the action maps $\alpha_{ij}$ for $V$. It satisfies $D^2=0$, $dD=Dd=0$, so we have the $M^!$-coalgebra
$V^! = (\M^!(V),d+D)$. Explicitly,
\[
V^!_j \,=\,\Tot \biggl\{ \cdots\lra  \bigoplus_{i_1<i_2<j} \M_{i_1,i_2}(\M_{i_1,j}(V_{i_1})) \lra \bigoplus_{i<j} \M_{ij}(V_i)
 \lra V_j\biggr\}.
\]

\begin{prop}
The functor $V\mapsto V^!$ defines a quasi-equivalence of dg-categories $\Alg_\M\to\Coalg_{\M^!}$.
\end{prop}

\noindent{\sl Proof:} Again, this is completely similar to the corresponding statement for algebras over an operad
\cite{ginzburg-kapranov}. More precisely, we also have the dual construction, which for any uni-triangular
dg-comonad $\E$ and an $\E$-coalgebra $C$ defines an $\E^\!$-coalgebra $C^\! = (\E^\! (C), d+D)$. Then
one proves easily that for any $V\in\Alg_\M$ we have a quasi-isomorphism $V^{!\!}\to V$ and for any
$C\in \Coalg_\E$ we have a quasi-isomorphism $C\to C^{\! !}$. \qed

\vskip .2cm

Note that Koszul duality can be expressed through mutations.

\begin{prop}\label{prop:L-delta-M}
Let $\M$ be a quasi-triangular dg-monad.
The comonad $\M^!$ is identified with $(L_\delta \M)^*$,
 where $\delta\in\Br_N$ is the Garside discriminant
\eqref{eq:garside-dis}
\end{prop}

 This is a direct generalization of a statement about
 mutations of exceptional collections (when instead  of monads we have algebras) which was
 proved   in
   \cite{bondal-polishchuk}.

   \vskip .2cm

\noindent {\sl Proof:} Direct comparison:
 unraveling  the definition of $L_\delta \M$, we find that it
  involves only the left adjoints $^*\M_{ij}$ and after passing to the right adjoints
  we recover $M^!$.
\qed

\begin{cor}
Let $\M$ be a quasi-triangular dg-monad. Assuming the relevant double adjoints exist,
 the dg-monads ${^{**}\M}= ({^{**}\M}_{ij})$
and $\M^{**} = (\M^{**}_{ij})$ are identified with $L_{\delta^2}M$ and $L_{\delta^{-2}} M$
respectively. \qed
\end{cor}

\paragraph{$A_\oo$-monads and their straightening.}
As dg-monads can be seen as associative dg-algebras in the category of functors, one can also consider
the homotopy, or $A_\oo$-version, in which associativity holds only up to a coherent system of
higher homotopies.  This version is similar to the classical concept of $A_\oo$-algebras \cite{keller,
{KS-A-inf}} and we discuss it briefly in our context.

\vskip .2cm

Let $\Vc=(\Vc_1,\cdots, \Vc_N)$ be a sequence of pre-triangulated dg-categories. Denote
$\Vc_i^\sharp$ the graded category obtained from $\Vc_i$ by forgetting the differentials in the Hom-complexes.
We can consider $\Vc_i^\sharp$ as a dg-category with trivial differential.

\begin{Defi}
A uni-triangular $A_\oo$-monad $\M$ in $\Vc$ is a datum of:
\begin{itemize}
\item[(1)] A collection $\M^\sharp$ of dg-functors $(\M_{ij}: \Vc_i^\sharp\to\Vc_j^\sharp)_{i<j}$

\item[(2)] A differential $Q=Q_\M$ in the free comonad $\FC(\M^\sharp [1])$ of degree $+1$, satisfying $Q^2=0$
and compatible with the cocomposition.
\end{itemize}
An $A_\oo$-morphism of unitriangular $A_\oo$-monads $\M\to L$ is, by definition, a morphism of
unitriangular dg-comonads $(\FC(\M^\sharp [1], Q_\M)\to (\FC(L^\sharp[1], Q_L)$.
\end{Defi}

Explicitly, $Q=Q_\M$ is determined by its matrix elements with values in the ``space of cogenerators'' $M^\sharp[1]$.
After moving all the shifts to the right, these matrix elements can be written as:
\[
\begin{gathered}
d_{ij}: \M_{ij}\lra \M_{ij}[-1], \quad i<j;
\\
c_{ijk}: \M_{jk}\circ \M_{ij}\lra \M_{ik}, \quad i<j<k;
\\
c_{ijkl}: \M_{kl}\circ \M_{jk}\circ \M_{ij} \lra \M_{il}[1], \quad i<j<k<l;
\\
\cdots\cdots\cdots\cdots\cdots\cdots\cdots\cdots\cdots\cdots\cdots\cdots
\\
c_{i_1,\cdots, i_p}: \M_{i_{p-1},i_p}\circ\cdots\circ \M_{i_1, i_2}\lra \M_{i_1,i_p}[p-3], \quad i_1 <\cdots < i_p;
\\
\cdots\cdots\cdots\cdots\cdots\cdots\cdots\cdots\cdots\cdots\cdots\cdots\cdots\cdots\cdots
\end{gathered}
\]
The condition $Q^2=0$ unravels, in a standard way, to a series of conditions:
\begin{itemize}
\item[(0)] $d_{ij}^2=0$, so  each $\M_{ij}$ becomes a dg-functor. We denote by $d$
the differentials in all dg-functors obtained from the $\M_{ij}$ by composition.

\item[(1)] $[d, c_{ijk}]=0$, i.e., $c_{ijk}$ is a morphism of dg-functors.

\item[(2)] $[d, c_{ijkl}]$ is the associator for the $c_{ijk}$, i.e., the difference between the two
paths in the square of Definition \ref{def:monad}(a).

\item[($\cdots$)] And so on.
\end{itemize}
Thus the $c_{i_1,\cdots, i_p}$form a system of higher homotopies for associativity of the $c_{ijk}$.

\vskip .2cm

Note that  a uni-triangular  $A_\oo$-monad $\M$ in $\Vc$ gives a uni-triangular
(strictly coassociative) dg-comonad $M^! = (\FC(\M^\sharp[1], Q)$. Therefore
the Koszul dual $\M^{!\!}$ is a strictly associative uni-triangular dg-monad,
which we call the {\em straightening} of $\M$.  Further,  in the same way as in Proposition
\ref{prop:M-!!}, we have a canonical  $A_\oo$-quasi-isomorphism $\M^{!\!}\to \M$
of $A_\oo$-monads.

This procedure of straightening allows us to ignore the $A_\oo$-issues in the sequel, even though
many constructions of monads in the later sections produce, strictly speaking, $A_\oo$-monads.


\subsection{  The Fukaya-Seidel category of a schober on $\CC$}
\label{subsec:monads-glue}

\paragraph{Uni-triangular Barr-Beck theory for spherical functors.}
Let
\be\label{sch-disk-diag2}
\xymatrix{
\Phi_1
\ar[dr]_{a_1}& \cdots & \Phi_N\ar[dl]^{a_N}
\\
& \Psi &
}
\ee
be a diagram of $N$ spherical dg-functors with common target, as in
 \S  \ref{subsec:schob-surf}\ref{par:schob-disk-spid}
 It gives a uni-triangular dg-monad
 \be\label{eq:monad-M-spher}
 \M=\M(a_1,\cdots, a_N) \,=\,\bigl(\M_{ij} \,=\, a^*_j \circ a_i:  \Phi_i\lra\Phi_j\bigr)_{i<j}
 \ee
 on  $\Phi :=(\Phi_1,  \cdots,  \Phi_N)$.
The composition
\[
c_{ijk}: \M_{jk}\circ \M_{ij} \,=\, a_k^*\circ a_j \circ a_j^*\circ a_i \lra a_k^*\circ a_i \,=\, \M_{ik}
\]
is obtained from the counit $\eta_j:  a_j \circ a_j^*\to \Id_{\Phi_j}$ by composing with
$a_k^*$ on the left and $a_i$ on the right.
We will call $\M$ the {\em Fukaya-Seidel monad} associated to the diagram
\eqref{sch-disk-diag2}.

\vskip .2cm

Accordingly, we have a pre-triangulated category $\Alg_\M=\< \Phi_1,\cdots, \Phi_N\>$
with its standard semi-orthogonal decomposition.

\vskip .2cm

For any object $U\in\Psi$ the collection $b(U)=(a_1^*(U),\cdots, a_N^*(U))$ is naturally an $M$-algebra.
The action map
\[
\alpha_{ij}: \M_{ij}(a_i^*(U)) \,=\, a_j^* \, a_i\, a_i^*\, (U) \lra a_j^*(U)
\]
is obtained from the counit $a_i\circ a_i^*\to\Id$. This gives a dg-functor
\be\label{eq:barr=beck}
b: \Psi \lra \Alg_\M
\ee
which we call the {\em Barr-Beck functor}.

\vskip .2cm

Let us find the left adjoint ${^*b}: \Alg_\M\to\Psi$. For an $\M$-algebra $V=(V_1,\cdots, V_N)$, we define
the {\em Beck-Barr complex} $\BBB(V)$ to be the following complex of objects of $\Psi$
(with differentials being closed morphisms of degree $0$):
\[
\cdots\lra  \bigoplus_{i_1<i_2<j} a_j \,a_j^* \,a_{i_2} \,a_{i_2}^*\, a_{i_1}\, (V_{i_1})
\lra
\bigoplus_{i<j} a_j \,a_j^* \,a_i\, (V_i) \lra \bigoplus_j a_j\,(V_j).
\]
The grading is normalized so that the complex terminates in degree $0$. The differential $d$ on the
$(-p)$th term
\[
\bigoplus_{i_1<\cdots < i_p < j} a_j\, a_j^*\, a_{i_p} \, a_{i_p}^* \,\cdots a_{i_2} \, a_{i_2}^* \, a_{i_1}\,
(V_{i_1})
\]
is the alternating sum of the results of applying the counits $a_j\circ a_j^*\to \Id, \cdots,  a_{i_2}\circ a_{i_2}^*
\to\Id$
and the action map $a_{i_2}^* a_{i_1}(V_{i_1})\to V_{i_2}$. It is clear that $d^2=0$.

\begin{prop}\label{prop:*b}
The functor $V\mapsto \Tot(\BBB(V))$ is left adjoint to $b$.
\end{prop}

\noindent{\sl Proof:} By definition of the Hom-complexes in $\Alg_\M$, for any $U\in\Psi$ we have
\[
\begin{gathered}
\Hom^\bullet_{\Alg_\M} (V, b(U)) \,= \Tot\biggl\{
\bigoplus_j \Hom^\bullet_{\Phi_j}(V_j, a_j^*(U)) \to \hfill
\\ \to \bigoplus_{i<j} \Hom^\bullet_{\Phi_j}
(a_j^*\, a_i(V_i), a_j^*(U)) \to
\bigoplus_{i_1<i_2<j} \Hom^\bullet_{\Phi_j}
(a_j^*\, a_{i_2}\, a_{i_2}^* \, a_{i_1}(V_{i_1}), a_j^*(U))\to\cdots\biggr\}
\\
=\Tot\biggl\{ \bigoplus_j \Hom^\bullet_\Psi(a_j(V_j), U) \to \hfill
\\
\to \bigoplus_{i<j} \Hom^\bullet_\Psi
(a_j\, a_j^* \, a_i(V_i), U) \to \bigoplus_{i_1<i_2<j} \Hom^\bullet_\Psi
(a_j \, a_j^* \, a_{i_2}\, a_{i_2}^* a_{i_1}(V_{i_1}), U) \to\cdots \biggr\}
\\
=\,\Hom^\bullet_\Psi (\Tot(\BBB(V)),U). \hfill
\end{gathered}
\]
\qed

\begin{rems}
(a) The complex $\BBB(V)$ has some resemblance with $B^+$, the augmented bar-complex of $V$
(see the proof of Proposition \ref{prop:augm-bar}) or, more precisely, with the direct sum
$\bigoplus_j a_j(B^+_j)$. In fact, the two complexes have the same terms but the differential
in $\BBB(V)$ has more summands than that in $\bigoplus a_j(B^+_j)$. Of course, each
$B^+_j$ and therefore $\bigoplus a_j(B^+_j)$ is exact while $\BBB(V)$ is, in general, not.

\vskip .2cm

(b) The classical Barr-Beck theory \cite{barr-wells} concerns a single adjoint pair
$\xymatrix{
\Cc \ar@<.4ex>[r]^a&\Dc \ar@<.4ex>[l]^{a^*}
}$. In this case $\M= a^*\circ a: \Cc\to\Cc$ is a monad in the classical sense, and we have a functor
$b: \Dc\to\Alg_\M$ sending $U\mapsto a^*(U)$. The Barr-Beck monadicity theorem
\cite{barr-wells} gives a criterion for $b$ to be an equivalence. The complex $\BBB(V)$,
$V\in\Alg_\M$, can always be defined as a simplicial object in $\Dc$. In our uni-triangular context,
our interest is not in $b$ being an equivalence but in a different property.
\end{rems}

\begin{prop}
The functor $b$ is spherical. The spherical {\em cotwist}  $T^{(b)}_\Psi = \Cone\{ \Id_\Psi\to b^*\circ b\}$ is
identified with $T_{N,\Psi}^{-1}\circ\cdots\circ T_{1,\Psi}^{-1}$, where $T_{i,\Psi}=T_\Psi^{(a_i)} =
\Cone\{a_i\circ a_i^*\to \Id_\Psi\}$ is the spherical {\em twist} for $a_i$.
\end{prop}

\noindent{\sl Proof:}
By Corollary \ref{cor:spher-adj} and Remark \ref{rem:spher-adj}
it suffices to show that $^*b$ is spherical and to show that $T_\Psi^{({^*b})}\= T_{1,\Psi}\cdots T_{N,\Psi}$.
So we prove the statements in this latter form.

\vskip .2cm

From the definition of $b$ and from Proposition \ref{prop:*b}, we see that
$T^{({^*}b)}_\Psi$ can be written as $\Tot (C^\bullet)$, where
\[
  C^\bullet\,   =\,\biggl\{ \cdots \lra \bigoplus_{i<j} a_j \, a_j^*\, a_i\, a_i^* \lra \bigoplus_j a_j\, a_j^* \lra \Id\biggr\}
\]
(the grading terminates in degree $0$). This $C^\bullet$ can be obtained as follows.
Take the $2$-term complexes of functors $\{a_i\circ a_i^*\to\Id_\Psi\}$, $i=1,\cdots, N$
(with grading in degrees $-1,0$) and compose them, getting an $N$-fold complex
(commutative $N$-cube)
\[
\{a_1\circ a_1^*\to\Id_\Psi\}\circ\cdots\circ \{a_N\circ a_N^*\to\Id_\Psi\}.
\]
Converting this cube, in a standard way, into a single complex gives $C^\bullet$.
But
\[
\Tot(C^\bullet) \,\=\, T_{1,\Psi} \circ\cdots T_{N,\Psi}.
\]
This proves that $T^{({^*}b)}_\Psi$ is an equivalence and has the claimed form.

\vskip .2cm

Now look at the functor $T^{({^*b})}_{\Alg_\M} = \Cone\{ \Id_{\Alg_\M}\to b\circ{^*b}\}[-1]$.
Let $V=(V_1,\cdots, V_N)\in\Alg_\M$. From  Proposition \ref{prop:*b}, we find that
$T_{\Alg_\M}(V)_p$ is the total object of the double complex
\be\label{eq:T-Alg-M}
\xymatrix{
&&&V_p\ar[d]^e
\\
\cdots \ar[r]&\bigoplus_{i_1<i_2<j} a_p^*\, a_j\, a_j^*\, a_{i_2} a_{i_2}^* a_{i_1}(V_{i_1})
\ar[r]& \bigoplus_{i<j} a_p^* \, a_j a_j^*\, a_i(V_i) \ar[r]& \bigoplus_j a_p^*\, a_j(V_j),
}
\ee
where $V_p$ is in bidegree $(0,0)$ and $e$ takes $V_p$ to $a_p^*\, a_p(V_p)$ by the unit of
th adjunction.

\vskip .2cm

To analyze this, for any $i,j\in\{1,\cdots, N\}$ (not necessarily $i<j$), let
$\M_{ij} = a_j^*\circ a_i: \Phi_i\to\Phi_j$. Then the matrix of functors  $\M^\Box=(\M_{i,j})_{i,j=1}^N$
is a ``matrix monad'', i.e., we have compositions $c_{ijk}: \M_{jk}\circ \M_{ij}\to \M_{ik}$
for any $i,j,k$ and the units $e_i: \Id_{\Phi_i}\to \M_{ii}$, satisfying the usual axioms.

Let $\M^<, \M^\geq$ be the matrices of functors given by
\[
\M^<_{ij}\,=
\begin{cases}
\M_{ij},& \text{ if } i<j;
\\
0,& \text{ otherwise;}
\end{cases}
\quad\quad\quad
\M^\geq_{ij}\,=
\begin{cases}
\M_{ij},& \text{ if } i\geq j;
\\
0,& \text{ otherwise.}
\end{cases}
\]
So formally, $\M^\Box=\M^< + \M^\geq$ is the decomposition of $\M^\Box$ into the strict upper-triangular
and non-strict lower-triangular parts. Our uni-triangular monad $\M$ can be formally written in
this notation as $\M=1+\M^<$.

\vskip .2cm

\noindent{\bf Preliminary step: Symbolic computation.}
To analyze $T_{\Alg_\M}^{({^*b})}$ as a functor, we first do a symbolic computation with matrices,
``at the level of Euler characteristics'', interpreting the total object of  the double complex \eqref{eq:T-Alg-M}
for each $p$   as the alternating
sum of its terms.  In this symbolic way, $T_{\Alg_\M}^{({^*b})}$ is represented by the terminating matrix series
\[
\begin{gathered}
1-\M^\Box + \M^\Box\cdot \M^< - \M^\Box\cdot (\M^<)^2 + \cdots  \hfill
\\
= 1-\M^\Box(1+\M^<)^{-1}\,=\, (1-\M^\geq)(1+\M^<)^{-1}.\hfill
\end{gathered}
\]
Now, $1-\M^\geq$ is lower triangular with diagonal elements corresponding to the functors
$T_{\Phi_i}^{-1} = \Cone\{\Id_{\Phi_i}\to a_i^* \circ a_i\}[-1]$. So, in our  formal computation,
these diagonal elements should be considered invertible and so  $1\-M^\geq$ is invertible, as is
$(1+\M^<)^{-1}$.

\vskip .2cm

\noindent{\bf Upgrade to functors.}
We now upgrade the above symbolic computation to an actual argument with functors. For this
we first define the  functors
\[
T^{\leq}, T^{\geq}: \Alg_\M\lra \Alg_\M.
\]
More precisely, define $T^\leq$ to send $V\in\Alg_\M$
to
\[
\Tot\,\biggl\{ \cdots\lra \M^<(\M^<(V))\lra \M^<(V) \buildrel\alpha\over\lra V\biggr\}.
\]
Here the complex in curly brackets is a subcomplex of the augmented bar-complex of $V$ (see the
proof of Proposition \ref{prop:augm-bar}), but it is not exact. By construction, $T^{\leq}$ is
upper-unitriangular, i.e., it sends each $\Phi_i$ into the subcategory $\< \Phi_1,\cdots, \Phi_i\>$.
In this, it is an analog of the matrix  $(1+\M^<)^{-1}$ in the above preliminary reasoning.
Note further that $T^\leq$ is an equivalence for the same reason a unitriangular matrix is invertible.
That is,   $T^\leq$ preserves the admissible filtration
\[
\<\Phi_1\>\,  \subset \, \<\Phi_1, \Phi_2 \>\,  \subset \, \cdots \,  \subset \, \<\Phi_1, \cdots, \Phi_N\> \,=\,
\Alg_\M
\]
corresponding to our original SOD $(\Phi_1,\cdots, \Phi_N)$,
and induces identity on each quotient $\<\Phi_1,\cdots, \Phi_i\>\bigl/ \<\Phi_1,\cdots, \Phi_{i-1}\> =\Phi_i$.

\vskip .2cm

Further,
  define the functor $\M^\geq: \Alg_\M\to\Alg_\M$ by
\[
(\M^\geq(V))_i \,=\,\bigoplus_{p\geq i} \M_{pi}(V_p).
\]
This is made into an $\M$-algebra as follows. For $i<j$ the map
\[
\M_{ij}\biggl( \bigoplus_{p\geq i}\M_{pi}(V_p)\biggr) \,=\,\bigoplus_{p\geq i} \M_{ij} \M_{pi}(V_p)
\lra \bigoplus_{p\geq j} \M_{pj}(V_p)
\]
sends the summand $\M_{ij}\M_{pi}(V_p)$ to $\M_{pj}(V_p)$ via $c_{ijp}: \M_{ij}\M_{pi}\to \M_{pj}$, if $j\geq p$,
and sends such a summand to $0$, if $j>p$.

We have a natural transformation $\Id\buildrel\eps\over\to \M^\geq$ given by the units of the adjunction
$\Id_{\Phi_i}\to \M_{ii}=a_i^*\circ a_i$. Now define $T^\geq = \Cone(\eps)[-1]$, so
\[
(T^\geq(V))_i \,=\, T_{\Phi_i}^{-1}(V_i)[-1] \,\oplus \,\bigoplus_{p>i} \M_{pi}(V_p)[-1].
\]
By construction, $T^\geq$ preserves the admissible filtration
\[
\, \<\Phi_N\>\,\subset \, \< \Phi_{N-1}, \Phi_N\> \,\subset \, \cdots \,\subset \, \< \Phi_1,\cdots, \Phi_N\> \,=\,
\Alg_\M
\]
corresponding to the mutated SOD $L_\delta(\Phi_1,\cdots, \Phi_N)$, see
Proposition \ref{prop:delta-mutation-SOD},
 and induces on each quotient
\[
\< \Phi_i, \Phi_{i+1}, \cdots, \Phi_N\> \, \bigl/ \, \< \Phi_{i+1}, \cdots, \Phi_N\>\,\simeq \, \Phi_i
\]
an equivalence, namely $T_{\Phi_i}^{-1}$. Therefore $T^\geq$ is an equivalence as well.
It is the functor upgrade of the matrix $1-\M^\geq$ from the above preliminary step.

\vskip .2cm

Finally, we notice that $T_{\Alg_\M}^{({^*b})}=T^\geq\circ T^\leq$ and so it is an equivalence.
Since both $T^{({^*b})}_\Psi$ and $T^{({^*b})}_{\Alg_\M}$ are equivalences, the functor $^*b$ is spherical
and so is $b$.
\qed

\begin{rem}
 The splitting  $T_{\Alg_\M}^{({^*b})} \= T^\geq\circ T^\leq$ is
analogous to the formula $
 T^\Lc_{\on{norm}} \, = \, C^+ T^L (\wt C^-)^{-1}
 $
representing the normalized monodromy
 of the Fourier transform of a perverse sheaf
in terms of the Stokes matrices, see Proposition \ref{prop:mono=2Stokes}.
More precisely, $T^\geq$   is analogous to the product $ C^+ T^L$, the latter being non-strict lower triangular
with isomorphisms (monodromies) on the diagonal.
The functor
 $T^\leq$, whose very shape
is suggestive of a geometric series,  is analogous to $ (\wt C^-)^{-1}$
(strict upper-triangular with identities on the diagonal).

\end{rem}

 \paragraph{The Fukaya-Seidel category.}
 Let now $A=\{w_1,\cdots, w_N\}\subset \CC$ and $\Sen\in\Schob(D,A)$ be a schober on
$\CC$ with singularities in $A$. Topologically, $\CC$ is the same as the standard open disk,
so we use the approach of
\S  \ref{subsec:schob-surf}\ref{par:schob-disk-spid}
That is, we choose a direction at infinity which we think of as a complex number $\zeta$ with
$|\zeta|=1$ and take a ``Vladivostok'' point $\bv=R\zeta$, $R\gg 0$ far away in the direction
$\zeta$. Let $K$ be a $\bv$-spider for $(\CC,A)$. Then $\Sen$
is represented by a diagram \eqref{sch-disk-diag} or \eqref{sch-disk-diag2} in the presense of $K$.
We denote the spherical functor $a_i: \Phi_i\to\Psi$ by $a_{i,K}$ to emphasize its dependence on $K$.
In particular, we have the uni-triangular dg-monad
$\M=\M(K) = \M(a_{1,K},\cdots, a_{N,K})$ as in \eqref{eq:monad-M-spher}.

We define  the {\em Fukaya-Seidel category of}  $\Sen$ {\em relatively to $\bv$ and  $K$ } to be
\be\label{eq:FS=Alg-M}
\FF_K(\CC, \bv; \Sen) \,:= \,  \Alg_{\M(K)}.
\ee
We now study the dependence of this construction on $K$.
Recall that the set $\Kc(\CC,\bv,A)$ of $\bv$-spiders for $(\CC, A)$ is acted upon, simply transitively,
by the group $\Br_N$, see \S \ref{subsec:abs-PL}\ref{par:spider}
By Proposition \ref{prop:schober-K-change}, the effect of this action on the functors $a_i$ is
found, on the generators $\tau_i\in\Br_N$, by
\[
\bigl( a_{1,\tau_i (K)}, \cdots, a_{N,\tau_i (K)}\bigr) \,=\,
\bigl( a_{1,K}, \cdots, a_{i-1, K}, \, T_{i, \Psi} a_{i+1,K}, \,  a_{i, K}, a_{i+2,K}, \cdots, a_{N,K}\bigr).
\]

\begin{prop}\label{prop:mut-monads}
Let $a_i: \Phi_i\to\Psi$, $i=1,\cdots, N$, be spherical dg-functors as in \eqref{sch-disk-diag2}.
For  $i=1,\cdots, N-1$ let
\[
(\ol a_1, \cdots, \ol a_N) \,=\, (a_1, \cdots, a_{i-1}, \, T_{i,\Psi} a_{i+1}, \, a_i, a_{i+2}, \cdots, a_N).
\]
Then we have a natural quasi-isomorphism of dg-monads
\[
\M(\ol a_1, \cdots, \ol a_N) \,\simeq \, (L_i M(a_1, \cdots, a_N))[\delta_{1,i},\cdots, \delta_{N,i}].
\]
 Here $\delta_{\nu,i}$ is the Kronecker symbol and the square brackets  mean the shift of the monad
 \eqref{eq:shift-monad}.

\end{prop}

\noindent {\sl Proof:} This is a direct comparison. Let  $\ol \M= \M(\ol a_1, \cdots, \ol a_N)$.
We express  the $\ol \M_{\nu,j}= \ol a_j^*\,  \ol a_\nu$ through the
$\M_{ij}= a_j^*\, a_\nu$ using the definition of the $\ol a_\nu$
to arrive at the formulas
   \eqref{eq:monad-left-mut} for the left mutation.  The only new $\ol a_\nu$ are
   \[
   \ol a_i \,=\, T_{i,\Psi}\, a_{i+1} \,=\, \Cone\, \{ \Id\to a_i \circ {^*a}_i\}[-1] \,\circ \, a_{i+1}
   \,=\,
   \Cone\, \{ a_{i+1} \to a_i \,\circ \, {^*a}_i \, \circ \, a_{i+1}\} [-1],
   \]
 and $\ol a_{i+1}=a_i$.

 \vskip .2cm

 So $\ol \M_{\nu,j}=\M_{\nu,j}$ for $\nu,j\notin \{i, i+1\}$, which matches the $1$st, $4$th and last lines of
 \eqref{eq:monad-left-mut}.

 \vskip .2cm

 If $\nu\leq i-1$ and $j=i$, then
 \[
 \begin{gathered}
 \ol \M_{\nu,i}\,=\, \ol a_i^8 \, a_\nu \,=\,\Cone \, \{ a_{i+1}^* \circ a_i\circ a_i^* \to a_{i+1}^*\} \,\circ\, a_\nu
 \,=\, \Cone \, \{ a_{i+1}^*\circ a_i\circ a_i^*\circ a_\nu \to a_{i+1}^*\circ a_\nu\}
 \\
 =\, \Cone\, \{ \M_{i, i+1}\circ \M_{\nu,i} \to \M_{\nu, i+1}\}\,=\, \M'_{\nu,i}[1], \hfill
 \end{gathered}
 \]
where $\M'_{\nu,i}$ is given by the second line of  \eqref{eq:monad-left-mut}.

\vskip .2cm

If $\nu\leq i-1$, $j=i+1$, then
\[
\ol \M_{\nu, i+1} \,=\, \ol a_{i+1}^* \ol a_\nu \,=\, a_i^* a_\nu \,=\, \M_{\nu,i} \,=\, \M'_{\nu, i+1},
\]
as given by the third line of  \eqref{eq:monad-left-mut}.

\vskip .2cm

If $\nu=i$, $j=i+1$, then
\[
\begin{gathered}
\ol \M_{i, i+1} \,=\, \ol a_{i+1}^* \ol a_i \,=\, a_i^*\,\circ \, \Cone\, \{ a_{i+1} \to a_i \,\circ \, {^*a}_i
\,\circ \, a_{i+1}\}[-1] \,=
\\
=\, \Cone \, \{ a_i^* \,\circ \, a_{i+1} \to a_i^* \,\circ a_i \,\circ \, {^*a}_i \,\circ\, a_{i+1}\}[-1] \,=\,
\Cone\, \{ a_i^* \to a_i^*\,\circ \, a_i \,\circ\, {^*a}_i\} \,\circ\, a_{i+1} [-1].
\end{gathered}
\]
By Proposition \ref{prop:kuznetsov}(1), the last cone is identified with ${^*a}_i$, so we get
\[
{^*a}_i a_{i+1}[-1] \,=\, {^*\M}_{i, i+1}[-1] \,=\, \M'_{i,i+1}[-1],
\]
see the $5$th line of  \eqref{eq:monad-left-mut}.

\vskip .2cm

In this way we find that $\ol \M_{\nu,j}\,=\, \M'_{\nu,j}[\delta_{i,j}-\delta_{i,\nu}]$ as claimed. \qed

\begin{cor}\label{cor:FS-identified}
The categories $\FF_K(\CC,\zeta;\Sen)$ for different $K\in \Kc(\CC,\bv,A)$, are canonically
(uniquely up to a contractible space of higher homotopies) identified with each other.
\end{cor}

\noindent{\sl Proof:} Follows from Proposition \ref{prop:mut-monads}, and from the
identification \eqref{eq:shift-alg-monad}. \qed

\begin{Defi}\label{def:fuk-sei}
(a) The Fukaya-Seidel category of a schober $\Sen$ on $\CC$ with respect to a direction $\zeta$
at infinity is defined as
\[
\FF(\CC,\zeta ;\Sen) \, :=\, \FF_K(\CC, \bv ; \Sen), \quad \forall \, v=R\zeta,\,  R\gg 0
\text{  and } \, K\in \Kc(\CC, \bv,A).
\]

(b) Similarly, the Fukaya-Seidel category of a schober $\Sen$ on a disk $D$ with respect to
a boundary point $\bv\in\del D$ is defined as
\[
\FF(D, \bv; \Sen) \, := \, \FF_K(D, \bv;  \Sen), \quad\quad \forall \, K\in \Kc(D,\bv, A).
\]
\end{Defi}


\subsection{The Lefschetz schober of a Landau-Ginzburg potential}\label{subsec:lef-schober}

In this section we explain how various constructions related to the ``classical''
Fukaya-Seidel
categories in symplectic geometry \cite{seidel} can be interpreted in our schober language.

\paragraph{The idea of the Lefschetz schober.} We consider the situation of \S
\ref{subsec:PL-class}\ref{par:PL-perv}
That is, $X$ is a Riemann surface and $W: Y\to X$ is an $X$-valued generalized Lefschetz pencil
of relative dimension $n$ (Definition \ref{def:GLP}). Let $A=\{w_1,\cdots, w_N\}\subset X$
be the set of singular values of $W$.
In this case we have the Lefschetz perverse
sheaf
\[
\Lc_W = R^n_\perv W_* \, \, \ul \k_Y \,\in \,\Perv(X,A).
\]
We will construct a natural categorification of $\Lc$
which we call the {\em Lefschetz schober} and denote $\Len_W\in\Schob(X,A)$.
The data  of Definition \ref{def:schober-surf} defining $\Len_W$ will be as follows:

\begin{itemize}
\item The local system  of categories $(\Len_W)_0$ on $X\- A$ will  consist of the Fukaya categories
of the  $W^{-1}(x)$, $x\in X\- A$,  which are considered as
compact real symplectic manifolds with the symplectic form given by a K\"ahler metric.

\item The local system $\bPhi_i(\Len_W)$ on $S^1_{w_i}$ will correspond to the
{\em local Fukaya-Seidel category} $\LFS_i(W)$ with the natural monodromy action, see
below.
\end{itemize}

Thus, in our ``upgrade'' from perverse sheaves to perverse schobers, the Fukaya category categorifies
the middle (co)homology while $\LFS_i(W)$ categorifies, following \cite{seidel},
 the space of vanishing cycles.

 \paragraph{  Reminder on Fukaya categories.} We first recall, following \cite{seidel}, the construction of
 the Fukaya category of a compact K\"ahler manifold. We start with recalling the concept of
 a quadratic differential, necessary to define a dg-category with  $\ZZ$-graded (not just $\ZZ/2$-graded)
 Hom-complexes.

 \vskip .2cm

  Let $V$ be a $\CC$-vector space of dimension $n$. Denote by $G_{\RR \RR}(n, V)$
 the open part of the Grassmannian $G_\RR(n,V)$ of $n$-dimensional real subspaces
 $W\subset V$ formed by $W$ which are totally real, i.e., $W\oplus iW = V$.

  A  {\em quadratic differential} on $V$ is a
 $\CC$-linear map $\eta: (\Lambda^n V)^{\otimes 2} \to \CC$.  A nonzero quadratic differential $\eta$
 on $V$  gives  a map, called the  {\em  phase function}
 \[
 \alpha_\eta: G_{\RR\RR}(n,V) \lra S^1, \quad W \,\,\mapsto \,\,
  \frac{\eta\bigl( (e_1\wedge\ldots \wedge e_n)^{\otimes 2}\bigr) }
{\bigl| \eta \bigl( (e_1\wedge\ldots \wedge e_n)^{\otimes 2}\bigr)  \bigr|},
 \]
where
 $\{e_1, \cdots, e_n\}$  is any $\RR$- basis of $W$.

  \vskip .2cm

  Let now $Z$ be a compact complex manifold of dimension $n$.  A  {\em quadratic differential} on $Z$
 is a holomorphic section $\eta\in\Gamma(Z, (\Omega_Z^n)^{\otimes 2})$.
 Thus, for each $z\in Z$, we have a quadratic differential  $\eta_z$ on the $\CC$-vector
 space $T_zZ$.

 We say that $Z$ is an {\em Enriques manifold}, if it admits a nowhere vanishing quadratic
 differential, i..e, if $(\Omega_Z^n)^{\otimes 2}\simeq \Oc_Z$.  For example, any Calabi-Yau
 manifold ($\Omega^n_Z\simeq \Oc_Z$) is obviously Enriques.

 In the sequel we will assume that $Z$ is Enriques and fix a nowhere vanishing $\eta$.
 Assume further that $Z$ is equipped with a K\"abler metric $h=g+i\omega$, so
 $\omega$ is a symplectic form on $Z$. So we
  can speak about Lagrangan $\RR$-subspaces in the tangent spaces $T_zZ$ and about Lagrangian
  submanifolds in $Z$ (considered as a $2n$-dimensional real manifold).

 Any Lagrangian subspace in any $T_zZ$ is totally real, so the phase
 functions $\alpha_{\eta_z}, z\in Z$, assemble into a map (also called the phase function)
 \[
 \alpha_Z: \on{LG}(TZ) \lra S^1,
 \]
where $\on{LG}(TZ)$ is the Lagrangian Grassmannian bundle associated to the
tangent bundle $TZ$.

\vskip .2cm

The {\em Fukaya category} $\Fuk(Z)$ is a pre-triangulated dg-category defined as the pre-triangulated
envelope of (i.e., the category of twisted complexes over) the $A_\oo$-category $\Fuk'(Z)$ defined as follows.

Objects  of $\Fuk'(Z)$ are
{\it Lagrangian branes}, i.e., data
$
L^\flat = (L, \alpha_L, P_L),
$
where:
\begin{itemize}
\item  $L\subset Z$ is a  compact Lagrangian submanifold;  

\item $\alpha_L: L\lra \RR$ is a smooth map (called a  {\it grading} of $L$) such that
\[
\alpha_L(z) = e^{2\pi i\alpha_X(T_zL)},
\]
 \item $P_L$ is a Pin-structure on $L$.

\end{itemize}
 So, a Lagrangian
brane is a Lagrangian submanifold $L$ decorated by some additional
data.

\vskip .2cm

Let $L_1^\flat, L_2^\flat$ be two Lagrangian branes. Let us suppose
that $L_1$ and $L_2$ intersect transversely. One defines the Hom-complex
 $C^\bullet(L^\flat_1, L^\flat_2) = \Hom^\bullet_{\Fuk'(Z)}(L_1^\flat, L_2^\flat)$ as follows. As a $\k$-vector space
\[
C^\bullet (L^\flat_1, L^\flat_2)\, :=\, \bigoplus_{z\in L_1\cap L_2} \k\cdot z
\]
For each $z\in L_1\cap L_2$ we have two Lagrangian subspaces
$T_z L_1, T_z L_2 \subset T_zZ$, and their {\it index}
$i(T_zL_1, T_z L_2)\in \ZZ$ is defined as in \cite{seidel}  (11.25); to define it
one uses the decorations.
This defines a $\ZZ$-grading on $C(L^\flat_1, L^\flat_2)$, by setting
$\deg(z) = i(T_z L_1, T_z L_2)$.

\vskip .2cm

The Floer differential $d = d_{L_1, L_2}: C(L^\flat_1, L^\flat_2)^i \lra C(L^\flat_1, L^\flat_2)^{i+1}$ is defined by
its matrix elements:
for $z \in L_1\cap L_2$, $\deg(z) = i$
$$
dz = \sum_{y\in (L_1\cap L_2)^{i+1}} n(z,y) y,
$$
$n(z,y)$ being the number, with appropriate signs, of pseudo-holomorphic
(i.e., holomorphic with respect to a generic almost complex deformation of the complex structure)
bigons in $Z$ with vertices $z$ and $y$.

If $L_1$ and $L_2$ do not intersect transversely, one uses
an appropriate Hamiltonian displacement of one of them, to make
them transverse.

\vskip .2cm

 Compositions are defined similarly, using the numbers of pseudo-holomorphic
 polygons, see \cite{seidel}. We note the following:

 \begin {prop}[\cite{PSS}]\label{prop:PSS}
 For any Lagrangian brane $L^\flat\in \Ob(\Fuk(Z))$ we have a
natural isomorphism
\[
 H^*(\Hom_{\Fuk(Z)}(L^\flat,L^\flat)) = H^*(L;\k).
\]
\end{prop}

\paragraph{The local system of Fukaya categories.}
Let now $W: Y\to X$ be a generalized Lefschetz pencil of relative dimension $n$,
with the set of singular values $A\subset X$.
Let $X_0=X\- A$ and $Y_0=W^{-1}(X_0)$, so $W: Y_0\to X_0$  is a  smooth proper map.
We write $Y_x=W^{-1}(x)$ for the fibers of this map.

We assume that both $X$ and $Y$ are Enriques, with
chosen nowhere vanishing quadratic differentials $\eta_X$ and
$\eta_Y$. Then each $Y_x$, $x\in X_0$. is an Enriques manifold, with quadratic differential
$\eta_{Y_x}=(\eta_Y/W^{\eta_X})|_{Y_x}$. We fix a K\"ahler metric
$h$  on $Y$; it induces a metric on each $Y_x$, $x\in X_0$, so the category
$\Fuk(Y_x)$ is defined.

\vskip .2cm

Recall  \eqref{eq:symp-conn}
that the fibration $W: Y_0\to X_0$ carries the  symplectic connection $\nabla_W$.
 In particular, any path $\gamma$ in $X_0$ joining two points $x$ and $x'$ gives rise to a diffeomorphism
 $T_\gamma: Y_x\to Y_{x'}$,  preserving the symplectic structure.
The connection $\nabla_W$ is not flat, i.e., $T_\gamma$, as a diffeomorphism,
can change under small deformations of $\gamma$, but the change is given by  a Hamiltonian
isotopy. Therefore $\nabla_W$ gives rise to a local system
$\Fuk(Y_0/X_0)$ of pre-triangulated dg-categories on $X_0$ whose stalk at any $x\in X_0$ is
$\Fuk(Y_x)$.

\paragraph{The case of a Lefschetz pencil.}\label{par:cas-lef-penc}
 Assume further that $W$ is a Lefschetz pencil,
i.e., all the critical points $y_i\in Y$ are Morse.  Let $D_i$ be a small closed disk around $w_i=W(y_i)$,
choose a point $\eps_i\in\del D_i$
and let $\gamma_i$ a simple path in $D_i$  (``radius'') joining $w_i$ with $\eps_i$.
Let also $Y_i=W^{-1}(w_i)$, a compact Enriques K\"ahler manifold.
As explained in \S \ref{subsec:PL-class}, $\gamma_i$ gives rise to the Lefschetz  thimble $\Th_i\subset W^{-1}(D_i)$
whose projection to $D_i$ is $\gamma_i$. In particular, we have the geometric vanishing cycle $\Delta_i = \Th_i \cap Y_i$,
diffeomorphic to $S^n$, see Fig. \ref{fig:vanish}.
 It is known that $\Delta_i$ is Lagrangian in $Y_i$ and, moreover, has an upgrade to a Lagrangian
 brane $\Delta_i^\flat$,
 which is an object of $\Fuk(Y_i)$, see \cite{seidel} \S 16f. Further, by Proposition \ref{prop:PSS},
 \[
 H^\bullet\, \Hom^\bullet_{\Fuk(Y_i)} (\Delta_i^\flat, \Delta_i^\flat) \,=\, H^\bullet(S^n, \k).
 \]
 whici, together with  the fact that $\Fuk(Y_i)$ is Calabi-Yau, means that
  $\Delta_i^\flat$
 is a spherical object in $\Fuk(Y_i)$,  see Example \ref{ex:spher-ob}.
 This implies that the dg-functor
\[
a_i: \Phi_i := D^b\, \Vect_\k \lra  \Psi_i := \Fuk(Y_i), \quad V\mapsto V\otimes_\k \Delta_i^\flat
\]
 is spherical.

The functor $a_i$ gives rise to a perverse schober $\Len_{W, D_i}$ on $D_i$ with the
only singularity at $w_i$. In particular, on $C_i=\partial D_i$, the schober $\Len_{W, D_i}$ gives a local system of
dg--categories coming from the self-equivalence
\[
T_{\Psi_i} \,\,=\,\,\on{Cone}\bigl\{ a_i\circ a_i^*  \lra  \Id_{\Psi_i}  \bigr\}: \, \Psi_i \lra \Psi_i.
\]
On the other hand, consider the local system $(\Len_W)_0:= \Fuk(Y_0/X_0)$
of dg-categories on $X_0=X\-A$. It has,
  as the stalk at $\eps_i$, the same dg--category
$\Psi_i=\Fuk(Y_i)$, and the monodromy of $\Len_S^\parallel$ around $C_i$ is, by construction, induced by
the symplectomorphism $T_{C_i}: Y_i \to Y_i$, the holonomy of the symplectic connection.
It is known, see  \cite{seidel} \S 16c,  that $T_{C_i}$ is Hamiltonian isotopic to the {\em Dehn twist}
around $\Delta_i$,
which means that the induced self-equivalence on $\Fuk(Y_i)$ is identified with $T_{\Psi_i}$.
Therefore, $\Len_{W, D_i}$ and $\Len_W^\parallel$ agree on $C_i$.

\begin{defi}
We define the Lefschetz schober $\Len_W$ on $X$ to be obtained by gluing the local system
$(\Len_W)_0$ on $X_0$ with
the schobers $\Len_{W, D_i}$ on each $D_i$ along the above identifications on the $C_i$.
\end{defi}

\paragraph{The case of a generalized Lefschetz pencil.} Let now $W: Y\to\CC$ be a
generalized Lefschetz pencil, with $y_i\in Y$ being an isolated singular point with multiplicity
$\mu_i$. Let the disk $D_i \ni w_i=W(y_i)$, as well as $\eps_i\in\del D_i$ and $\gamma_i$
be as in \S \ref{par:cas-lef-penc}

Note that $\mu_i$ is the dimension of  the space of vanishing cycles
of $W$ at $y_i$. The idea of categorifying this space, i.e., upgrading it to a Fukaya-style
pre-triangulated ($A_\oo$-)category $\LFS_i(W)$,   goes back to Kontsevich
\cite{kontsevich-ENS}. We will refer to $\LFS_i(W)$ as the {\em local Fukaya-Seidel category}
of $W$ at $y_i$ and note that there are two natural approaches to its construction.

\vskip .2cm

\noindent {\bf (1) Deformation approach.} Deform $W$ to an actual Lefschetz pencil
$\wt W: \wt Y\to\CC$ so that $y_i$ splits into $\mu_i$ Morse critical points
$y_{i,1},\cdots, y_{i, \mu_i}$, lying close to each other and projecting under $\wt W$
to distinct points $w_{i,1},\cdots, w_{i,\mu_i}\in D_i$. Denoting $\wt Y_i = \wt W^{-1}(D_i)$,
we have a Lefschetz pencil $\wt W_i: \wt Y_i\to D_i$ and so, by the above, the
Lefschetz schober $\Len_{\wt W_i}$ on $D_i$. We then define
\[
\LFS_i(W)_{\eps_i} \,=\, \FF(D_i, \eps_i; \Len_{\wt W_i}),
\]
to be the Fukaya-Seidel category of $D_i$ with coefficients in this schober or, more
directly,  as the Fukaya-Seidel category of $\wt W_i$ in the classical sense \cite{seidel}.

The result should, morally, not depend on the choice of the deformation, since, as far
as local deformations (those in a neighborhood of $y_i$ only) go, different choices
are connected by the braid group action. Such an action translates into a mutation of SOD
on the Fukaya-Seidel category while the category itself is unchanged, see Corollary \ref{cor:FS-identified}.

\vskip .2cm

 The difficulty with this approach is that one can guarantee the existence of
a good deformation only locally (in a neighborhood of $y_i$) and not as a proper map. For
a non-proper map, however, the symplectic connection may be non-integrable,
i.e., the parallel transport of a point can go to infinity.
This makes it difficult to define the local system part of the schober $\Len_{\wt W_i}$.

\vskip .2cm

\noindent{(\bf 2) The direct approach.} As in \cite{kontsevich-ENS},
in Def. 3.7 of  \cite {abouzaid-Morse}
and  in App. A of  \cite{abouzaid-smith} we can take, as objects of $\LFS_i(W)_{\eps_i}$, a
certain clas of (possibly singular) Lagrangian branes which project by $W$ to the radial
path $\eps_i$ and vanish to $y_i$ over $w_i$.

\vskip .2cm

The implementation of this approach involves several subtle issues of symplectic geometry,
for which we refer to \S A.2 of  \cite{abouzaid-smith} and references therein.
Let us assume that we are in  situation where this approach can be implemented. In this case
Proposition A.4 and Theorem A.6 of \cite {abouzaid-smith}  imply that we have a spherical
functor
\[
a_i: \LFS_i(W)_{\eps_i} \lra (\Len^\para_W)_{\eps_i} \,=\, \Fuk(W^{-1}(\eps_i))
\]
whose associated spherical reflection is identified with the effect of  the monodromy around $\del D_i$.
This means that we can construct a schober $\Len_W$ by gluing $\Len_W^\para$ with
the schobers on the $D_i$ associated to the $a_i$, as before.


\section {Algebra of the Infrared associated to a schober on $\CC$,  I}\label{sec:AIR-I}

In this and the next chapter we analyze schobers on $\CC$ in the spirit of the Algebra of the Infrared,
relating curvilinear and rectilinear transport functors. Such analysis involves naturally convex geometry of
the set $A$, but splits naturally into two parts, which use respectively:

\begin{itemize}
\item[(1)] ``Elementary'' convex geometry related to the  data of which elements of $A$ lie
in the convex hulls of which subsets. This data  is formalized by the oriented matroid $\OO(A)$,
see Definition \ref{def:or-datum}.

\item[(2)] ``Advanced'' convex geometry involving the concept of the secondary polytope $\Sigma(A)$.
This polytope depends on more subtle aspects than just the data of $\OO(A)$.

\end{itemize}

In the present Chapter \ref{sec:AIR-I} we concentrate on the first part, with the second part postponed
till the next Chapter  \ref{sec:AIR-II}.

\subsection {Rectilinear approach and the Fourier transform  of a schober.} \label{subsec:schober-rect}

 \paragraph{Rectinilear approach to schobers on $\CC$.}
 We now adopt the point of view of \S \ref{subsec:rec-approach} and emphasize
 the rectilinear and convex-geometrical approach to the (purely topological) concept of
 a schober on $\CC$.

 \vskip .2cm

 We assume that $A$ is in linearly general position and $\Sen\in\Schob(\CC,A)$.
 We distinguish the horizontal direction ($\zeta=1$)  in each circle $S^1_{w_i}$
 and denote $\Phi_i(\Sen)$ the
 stalk of the local system $\bPhi_i(\Sen)$ at $1\in S^1_{w_i}$.
  For any distinct $i,j\in \{1,\cdots, N\}$
 we have the {\em rectilinear transport functor}
 \[
 M_{ij}\, :=\, M_{ij}([w_i, w_j]): \Phi_i(\Sen)\lra \Phi_j(\Sen),
 \]
as in \eqref{eq:rect-trans-def}.

\paragraph{Fourier transform of a schober on $\CC$: topology.}
Let $A=\{w_1,\cdots, w_N\}\subset \CC$. Following the conventions of \S \ref{subsec:FTPS-top}, we
write $\CC=\CC_w$ and denote the dual complex line as $\CC_z$ with coordinate $z$ dual to $w$.
If $\Fc\in\Perv(\CC_w,A)$ is a perverse sheaf with singularities in $A$, then we
can define its Fourier transform $\wc\Fc$ by \eqref {eq:FT-of-F}, using the
 Riemann-Hilbert correspondence. Proposition \ref {prop:wc-M-top-1} says that
 $\wc\Fc\in\Perv(\CC_z,0)$.

\vskip .2cm

Let now  $\Sen\in \Schob(\CC_w,A)$ be a perverse schober with singularities in $A$. We would
like to define its Fourier transform, which should be a new schober $\wc\Sen\in \Schob(\CC_z,0)$.
At present, it is not known how to extend the Riemann-Hilbert correspondence to schobers.
However, the above constructions provide a working  {\em ad hoc}  definition.

\vskip .2cm

First, we define two local systems of triangulated categories $\bPhi(\wc\Sen)$,
$\bPsi(\wc\Sen)$ on the circle $S^1 = \{|\zeta|=1\}$ which we think of as the
circle of directions of $\CC_z$ at $0$. Given $\zeta\in S^1$, let
 $\bv_\zeta = -R\ol\zeta$, $R\gg 0$ be a ``Vladivostok'' point in $\CC_w$
 very far in the direction
$-\ol \zeta$. We put
\[
\bPhi(\wc\Sen)_\zeta \,=\, \Sen_{\bv_\zeta} \,\= \, \bPsi_\oo(\Sen)_{-\ol\zeta}
\]
to be stalk of $\Sen$ at $\bv_\zeta$ which can be thought as the stalk at $-\ol\zeta$
of the local system of nearby cycles of $\Sen$ near $\oo$. Further, we put
\be\label{eq:Psi-Fourier-schob}
\bPsi(\wc\Sen)_\zeta \,=\, \FF(\CC_w, -\ol\zeta; \Sen)
\,=\, \FF_K(\CC,\zeta ;\Sen) \,=\,\Alg_{M_K}
\ee
to be the Fukaya-Seidel category of $\Sen$ in the direction $-\ol\zeta$, see Definition
\ref{def:fuk-sei}(a). Here $K$ is any $\bv_\zeta$-spider for $A$. Further, choosing some
particular value of $\zeta$, we define $\wc\Sen$ to be represented by the
  Barr-Beck
spherical functor
\[
b: \bPhi(\wc\Sen)_\zeta = \Sen_{\bv_\zeta}  \lra
 \Alg_{M_K} = \bPsi(\wc\Sen)_\zeta
\]
of \eqref{eq:barr=beck}.
 Comparison with Proposition \ref {prop:wc-M-top-2} shows that this is indeed a natural
 categorical analog.

 Study of  Fukaya-Seidel categories can be therefore seen as  study of the Fourier transform
 for schobers.

\paragraph{Categorical Stokes structures.}
A categorification of the concept of a Stokes structure (i.e., data defining ``irregular perverse
schobers'' in complex dimension $1$) was introduced by Kuwagaki \cite{kuwagaki}
(and also in an earlier unpublished note
 \cite{Kat5}). Here we recall this categorified
concept in a simple form convenient for us.

\vskip .2cm

Let $\Lambda$ be a local system of finite posets on $S^1$. Thus for each $\zeta\in S^1$ we have
a finite set $\Lambda_\zeta$ with a partial order $\leq_\zeta$. Let us assume the following
two properties of $\Lambda$ (the second one added for simplicity).

\begin{itemize}

\item[($\Lambda1$)] The order
 $\leq_\zeta$ is a total order for any $\zeta$ outside of a finite set $\Delta\subset\Lambda$
 (the anti-Stokes directions).

 \item[($\Lambda2$)] For $\zeta\in\Delta$ the order $\leq_\zeta$ has exactly two incomparable
 elements, whose order switches as we cross $\zeta$. That is, if $\zeta_-,\zeta_+$
 are nearby directions immediately before and after $\zeta$ in the counterclockwise order, then
 the two orders $\leq_{\zeta_+}, \leq_{\zeta_-}$ (which can be considered on the same set
 $\Lambda_\zeta$) have the form (for some $k$)
 \[
 \begin{gathered}
 \lambda_1 <_{\zeta_-} \lambda_2 <_{\zeta_-}\cdots <_{\zeta_-} \lambda_k <_{\zeta_-}\lambda_{k+1}
<_{\zeta_-} \cdots <_{\zeta_-} \lambda_N,
\\
 \lambda_1 <_{\zeta_+} \lambda_2 <_{\zeta_+}\cdots <_{\zeta_+} \lambda_{k+1} <_{\zeta_+}\lambda_{k}
<_{\zeta_+} \cdots <_{\zeta_+} \lambda_N
\end{gathered}
 \]
 \end{itemize}

 \begin{Defi}
 Let $\Lambda$ be a local system of posets on $S^1$ satisfying ($\Lambda1$) and ($\Lambda2$).
 A {\em categorical $\Lambda$-Stokes structure} is a datum of:
 \begin{itemize}
 \item[(1)] A local system $\Len$ of triangulated categories on $S^1$.

 \item[(2)] For each $\zeta\in S^1\-\Delta$, an $\oo$-admissible filtration
 (Definition \ref{eq:oo-admis})  $F^\zeta$ on the triangulated
 category $\Len_\zeta$ indexed by the totally ordered set $(\Lambda_\zeta, \leq_\zeta)$, so that:

 \item[(3)] $F^\zeta$ is covariantly constant on (each component of) $S^1\-\Delta$.

 \item[(4)] If $\zeta\in\Delta$, then in the notation of ($\Lambda2$) we have
 $F^{\zeta_+} = \tau_k F^{\zeta_-}$ (the result of mutation corresponding to the generator
 $\tau_k\in\Br_N$).
 \end{itemize}
 \end{Defi}

 \noindent The definition implies that the quotient categories
 \[
 \gr_\lambda^F \Len_\zeta \,=\, F_\lambda\Len_\zeta \bigl/ F_{\lambda'} \Len_\zeta \,=\, (F_{\lambda'} \Len_\zeta)^\perp_{F\lambda\Len_\zeta}, \quad \lambda' = \max_{\mu <\lambda} \, \mu, \,\,\zeta\in\S^1\-\Delta,
 \]
i.e., the terms of the SOD's associated to the $\oo$-admissible filtrations in
$\Len_\zeta$, $\zeta\in S^1\-\Delta$, unite into a $\Lambda$-graded local system of triangulated
categories $\LL=\gr^F\Len$ on all of $S^1$. Note that in this categorical context, a splitting,
i.e., an identification $\LL\to\Len$ usually does not exist, even locally.

 \paragraph {Fourier transform of a schober: Stokes structure.} We now specialize to the case when
 $\Lambda = \ul A$ is the constant sheaf corresponding to $A=\{w_1,\cdots, w_N\}\subset \CC$
 in strong linearly general position and $\leq_\zeta$ is given by the dominance order of the
 exponentials \eqref{eq:leq-zeta}. Then $\ul A$ satisfies  ($\Lambda1$) and ($\Lambda2$) and we
 will speak about {\em $A$-Stokes structures}.
 The set $\Delta$ of anti-Stokes directions for $\ul A$
 consists of the directions $\ol\zeta_{ij}$, where $\zeta_{ij}$ are the slopes from
 of \eqref{eq:zeta-ij}.

 \vskip .2cm

 Let us construct an $A$-Stokes structure on the local system $\Len = \bPsi(\wc\Sen)$
formed by the Fukaya-Seidel categories \eqref{eq:Psi-Fourier-schob}. Let $\zeta\in S^1$ be
non-anti-Stokes,
  so
we have a picture as in Fig. \ref{fig:horiz-cuts}. Let $\bv_\zeta = -R\ol\zeta$, $R\gg 0$
be very far in the direction $-\ol\zeta$. We can then take the  $\bv_\zeta$-spider for $A$
of the form
\[
K(\bv_\zeta) \,=\, \bigcup_{i=1}^N\,  [w_i, \bv_\zeta],
\]
the union of the straight intervals. As $R\gg 0$, it approximates the
 union $K(\zeta)$ of straight rays
  $K_i(\zeta)$ in the direction $(-\ol\zeta)$
from  Fig. \ref{fig:horiz-cuts}. By definition,
$\bPsi(\wc\Sen)_\zeta$ is the Fukaya-Seidel category
 $\FF_{K(\bv_\zeta)} (\CC_w, \bv_\zeta; \Sen)$, i.e., the category of algebras over
 the uni-triangular monad $M_{K(\bv_\zeta)}$. As such, it comes with a semi-orthogonal
 decomposition or, equivalently, see Definition \ref{eq:oo-admis},
 with an $\oo$-admissible filtration $F^\zeta$.

 \begin{prop}
 The filtrations $F^\zeta$ at the non-anti-Stokes fibers $\bPsi(\wc\Sen)_\zeta$
 form a categorified $A$-Stokes structure.
 \end{prop}

 \noindent{\sl Proof:} Let us look at the effect of crossing an anti-Stokes direction
 $-\ol\zeta_{ij}$. Let $\zeta_-$,   $\zeta_+$ be nearby directions immediately
 before and after $-\ol\zeta_{ij}$ and $K(\bv{\zeta_-})$, $K(\bv_{\zeta_+})$
 be the corresponding rectilinear spiders. Then $K(\bv_{\zeta_+})$
 is isotopic to the spider obtained from $K(\bv_{\zeta_-})$ by an elementary
 braiding of the $i$th and $j$th tentacles (these tentacles will have successive numbers,
 say $k$ and $k+1$ in the numeration according to $<_{\zeta_-}$).  So the
 effect of crossing  $-\ol\zeta_{ij}$ on the filtration is, by Propositions
   \ref{prop:schober-K-change} and \ref{prop:mut-monads}, is the mutation corresponding to
   $\tau_k$. \qed


\subsection{The infrared complex of a schober.} \label{subsec:infra-com-schob}

\paragraph{The setup.}
We now present   categorical analogs of   formulas of \S \ref{subsec:baby-IA}
 expressing  the transport  along a curved path in terms of  a  sum
of rectilinear ones.  This sum will be categorified by a complex of functors
while the proof, involving Picard-Lefschetz identities, will be categorifed by a
diagram of exact triangles which generalize Postnikov systems familiar in homotopy theory
and homological algebra \cite{gelfand-manin}.

\vskip .2cm

\begin{figure}[H]
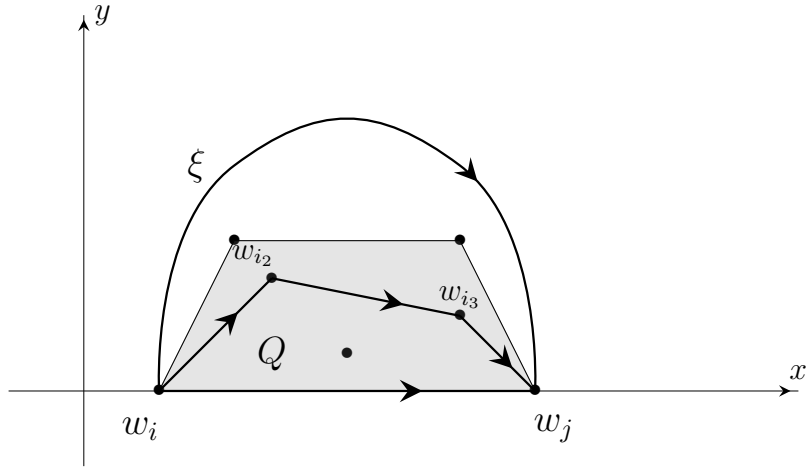

\centering

\btp[scale=0.5]

\node at (0,0){\small$\bullet$};
\node (wi) at (-5,-1){\small$\bullet$};
\node (wi2) at (-2, 2){\small$\bullet$};
\node (wi3) at (3,1){\small$\bullet$};
\node (wj) at (5,-1){\small$\bullet$};

\node at (3,3){\small$\bullet$};
\node at (-3,3){\small$\bullet$};
\node at (-4, -3){};
\node at (3,-4){};

\filldraw [color=gray, opacity=0.2] (5, -1) -- (3,3) --(-3,3) -- (-5, -1);

\draw  (5, -1) -- (3,3) --(-3,3) -- (-5, -1);


\draw  [decoration={markings,mark=at position 0.7 with
{\arrow[scale=2,>=stealth]{>}}},postaction={decorate},
line width=.3mm] (-5, -1) -- (5, -1);

\draw  [decoration={markings,mark=at position 0.7 with
{\arrow[scale=2,>=stealth]{>}}},postaction={decorate},
line width=.3mm] (-5, -1) -- (-2,2);

\draw  [decoration={markings,mark=at position 0.7 with
{\arrow[scale=2,>=stealth]{>}}},postaction={decorate},
line width=.3mm] (-2,2) -- (3,1);

\draw  [decoration={markings,mark=at position 0.7 with
{\arrow[scale=2,>=stealth]{>}}},postaction={decorate},
line width=.3mm] (3,1) -- (5, -1);

\draw  [decoration={markings,mark=at position 0.7 with
{\arrow[scale=2,>=stealth]{>}}},postaction={decorate},
line width=.3mm]   plot  [ smooth, tension=1] coordinates{   (-5, -1)    (-3, 5)  (3, 5)    (5, -1)   } ;

\node at (-2, 0){\large$Q$};
\node at (2, -3){};

\node at (-5.5, -2){\large$w_i$};
\node at (5.5, -2){\large$w_j$};
\node at (-4, 5) {\large$\xi$};

\node at (3.0, 1.5){$w_{i_3}$};
\node at (-2.5, 2.6){$w_{i_2}$};

\draw  [decoration={markings,mark=at position 0.99 with
{\arrow[scale=1.5,>=stealth]{>}}},postaction={decorate},
line width=.1mm] (-9,-1) -- (12,-1);

\draw  [decoration={markings,mark=at position 0.99 with
{\arrow[scale=1.5,>=stealth]{>}}},postaction={decorate},
line width=.1mm] (-7,-3) -- (-7, 9);

\node at (12,-0.5){$x$};
\node at (-6.5,9){$y$};

\etp
\caption{Convex paths as graphs of convex functions. } \label{fig:conv-paths}
\end{figure}

We assume that $A=\{w_1,\cdots, w_N\}\subset \CC$ is in linearly general position,
denote $Q=\Conv(A)$ and  fix a schober $\Sen\in\Schob(\CC,A)$.
  Consider the situation depicted on the right of the Fig. \ref{fig:circum}, reproduced below on the Fig. \ref{fig:conv-paths},
   so we have an edge
$[w_i, w_j]$ of $Q$ and  a path $\xi$ from $w_i$ to
$w_j$ which circumnavigates $Q$.  We would like to express the transport functor
$M(\xi): \Phi_i\to \Phi_j$ in terms of compositions of rectilinear transports associated to convex paths.

Accordingly, we will categorify the sum formula of
  Proposition \ref{prop:baby-IA}.

\paragraph{ Paths as graphs of functions.} We further assume that there is
  an affine coordinate system $(x,y)$ on $\CC=\RR^2$
such that:
\begin{itemize}
\item $w_i$ and $w_j$ lie on the horizontal axis $y=0$.

\item The projection of $Q$ to the horizontal axis along the vertical one  is the interval $[w_i, w_j]$ and not a bigger set.
Further, no edges of $Q$ are vertical, see Fig. \ref{fig:conv-paths}.
\end{itemize}

\noindent This assumption will be automatically satisfied in our main application. In general it can  always be ensured
by peforming a real projective transformation of $\RR^2$ (preserving straight  lines and therefore preserving convex geometry).

Let $\Pc$ be the (infinite-dimensional) convex cone formed by all convex functions $f: [w_i,w_j]\to\RR_{\geq 0}$
such that $f(w_i)=f(w_j)=0$.
Each such $f$ gives a path $\gamma_f$, the graph of $f$, joining $w_i$ and $w_j$.
We refer to paths $\gamma$ of the form $\gamma=\gamma_f$ as {\em convex paths}.
In particular, we can assume our
circumnavigating path $\xi$ to be convex.

\vskip .2cm

Each point $w_\nu\in A$, $\nu\neq i,j$, gives an affine hyperplane $H_\nu$ in $\Pc$ formed by all convex paths
$\gamma$ passing through $w_\nu$.
If $w_\nu$ has coordinates $(x_\nu, y_\nu)$ in our system, then $H_\nu$ consists of
 all convex functions $f$ such that $f(x_\nu)=y_\nu$. This  condition  is given
 by vanishing of the affine-linear functional
 \be
 l_\nu: \Pc\lra \RR, \quad l_\nu(f) = f(x_\nu)-y_\nu.
 \ee
We get an affine hyperplane arrangement $\Hc=\{H_\nu\}$ in $\Pc$.
In particular, we can speak about {\em flats} of $\Hc$, i.e., intersections $F_I:=\bigcap_{\nu\in I} H_\nu$
for vanious subsets $I\subset \{1,\cdots, N\}\- \{i,j\}$.  This includes $I=\emptyset$ which corresponds to
$F_\emptyset=\Pc$.  Each flat $F$ has well-defined {\em codimenion} $\codim(F)$, defined  as the
number of the $H_\nu$ containing $F$. This definition  is meaningful because of the following.

\begin{prop}\label{prop:H-trans}
The hyperplanes from $\Hc$ intersect transversely.
\end{prop}

\noindent{\sl Proof:} The statement means that for any  $I\subset J\subset \{1,\cdots, N\}\- \{i,j\}$ such that  $I\neq J$
and $F_I\neq\emptyset$,
the inclusion $F_J \subset F_I$ is proper. But this is clear, as among convex paths passing through all the
$w_\nu, \nu\in I$ (which exist by the  assumption $F_I\neq\emptyset$) we can always find a path not passing through any other
$w_\mu, \mu\notin I$. \qed

\vskip .2cm

Further, the arrangement $\Hc$ subdivides $\Pc$ into {\em faces} which are maximal subsets
on which each equation $l_\nu$ has the same sign $+,-,0$.
Faces are locally closed convex sets.
A face $\Gamma$ can be considered
as an isotopy class of convex paths $\gamma$ joining $w_i$ with $w_j$, with isotopies understood
relatively to the set $A$.
  The set of faces will be denoted by $P$. For a face $\Gamma\in P$ we denote by $\Aff(\Gamma)$
  its affine envelope which is the minimal flat of $\Hc$ containing $\Gamma$. Thus $\Gamma$
  is open on $\Aff(\Gamma)$, and we denote $\codim(\Gamma) := \codim (\Aff(\Gamma))$.

  In particular, open faces (i.e., faces of codimension $0$) correspond to
  {\em empty paths}, i.e., paths not containing any $w_\nu$, $\nu\neq i,j$.

  \paragraph{The infrared Postnikov system.}
 Each face $\Gamma\in P$ gives a transport functor $M(\Gamma): \Phi_i\to\Phi_j$.
  Explicitly, we choose any path $\gamma\in \Gamma$ and denote
  $w_{k_1} = w_i, w_{k_2},  \cdots, w_{k_p}=w_j$ all the elements of $A$ lying on $\gamma$,
  in order from $w_i$ to $w_j$.
  Denote $\gamma_\nu$ be the segment of $\gamma$ between $w_{k_{\nu}}$ and $w_{k_{\nu+1}}$.
  Then
  \[
  M(\Gamma) \, := \, M_{k_{p-1}, k_p}(\gamma_{p-1}) \circ \cdots \circ M_{k_1, k_2}(\gamma_1).
  \]
  Following Convention \ref{conv:clock}, the composition of the transport functors above
  is defined by   identifying  the stalks of the intermediate local systems $\bPhi$
  using clockwise monodromies, i.e., along the arcs going inside $Q$.

\begin{figure}[H]
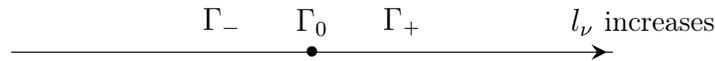

  \centering
  \btp[scale=0.4]
 \draw  [decoration={markings,mark=at position 0.99 with
{\arrow[scale=2,>=stealth]{>}}},postaction={decorate},
line width=.1mm] (-10,0) -- (10,0);

\node at (0,0){\small$\bullet$};
\node at (-3,1){$\Gamma_-$};
\node at (3,1){$\Gamma_+$};
\node at (0,1){$\Gamma_0$};
\node at (11,1){$l_\nu$ \small{increases}};
  \etp
  \caption{A wall-crossing triple.}\label{fig:wct}
  \end{figure}

  \begin{Defi}
  Three faces $\Gamma_-, \Gamma_0, \Gamma_+\in P$ are said to form a {\em wall-crossing triple},
  if:
  \begin{itemize}
  \item[(1)] $\Gamma_-$ and $\Gamma_+$ are open in the same flat $F= \Aff(\Gamma_-) = \Aff(\Gamma_+)$.
  In particular, they have the same codimension, say $r$.

  \item[(2)]  $\Gamma_0$ has codimension $r+1$ and separates $\Gamma_-$ and $\Gamma_+$ in $F$.
  That is, $\Gamma_0$ is contained in the closures of both $\Gamma_-$ and
  $\Gamma_+$, and $\Aff(\Gamma_0)$ is the intersection of $F$ and one more
  (unique, by Proposition \ref{prop:H-trans}) hyperplane $H_\nu$.

  \item[(3)]
  The equation $l_\nu$ of $H_\nu$ is negative on $\Gamma_-$, zero on $\Gamma_0$ and positive on $\Gamma_+$. See Fig. \ref{fig:wct}.
  \end{itemize}
  \end{Defi}

  Thus, crossing the ``wall'' $\Gamma_0$ in the direction   from $\Gamma_-$ to $\Gamma_+$ correspond, in terms
  of a path $\gamma$, to moving $\gamma$ across one point $w_\nu$,  in the direction from below to above.
  Crossing just this one wall and no others means that  no other points on $A$ leave  $\gamma$
  or appear on  $\gamma$. Therefore we have a Picard-Lefschetz situation and  the corresponding
  Picard-Lefschetz triangle
  gives rise to
   the {\em wall-crossing triangle} of functors
 \be\label{eq:WC-triangle}
 M(\Gamma_0) \buildrel v_{\Gamma_0, \Gamma_+}\over\lra M(\Gamma_+)
 \buildrel z_{\Gamma_+, \Gamma_-}\over
 \lra M(\Gamma_-) \buildrel u_{\Gamma_-, \Gamma_0}\over\lra M(\Gamma_0)[1].
 \ee

 \begin{Defi}\label{def:IPS}
 The diagram (in the homotopy category of dg-functors $\Phi_i\to\Phi_j$)  formed by the functors $M(\Gamma), \Gamma\in P$ and
 the natural transformations $ v_{\Gamma_0, \Gamma_+}$, $z_{\Gamma_+, \Gamma_-}$
 and $u_{\Gamma_+, \Gamma_0}$ for all wall-crossing triples of faces $(\Gamma_-, \Gamma_0, \Gamma_+)$,
 will be called the {\em infrared Postnikov system} and denoted $\Pen$.
 \end{Defi}

 Note that any arrow of $\Pen$ is contained
 in a unique exact triangle.

 \paragraph{Examples of infrared Postnikov systems.}

 \begin{ex}[(The maximally non-convex position)]\label{ex:max-nonconv}
 Let $A$ be as in Fig. \ref{fig:max-nonconv} on the left.
 That is, $Q$ is a triangle with one side $[w_1, w_N]$, and all the other points
 $w_2,\cdots, w_{N-1}\in A$ are positioned very close to the perpendicular to that side
 (but not exactly on the perpendicular, to satisfy the assumption that $A$ is in linearly general position).

 The face structure of the hyperplane arrangement $\Pc$ is depicted on the right of
  Fig. \ref{fig:max-nonconv} in the form of a $1$-dimensional transversal slice.
 Open faces correspond to   empty paths, i.e., paths containing only $w_1$ and $w_N$;
 there are $N-1$ such faces, corresponding to the paths $\xi_1=[w_1,w_N]$,
 $\xi_2$ and so on up to the circumnavigating path $\xi=\xi_{N-1}$.
 Faces of codimension $1$ correspond to polygonal paths $[12N] =[w_1, w_2]\cup [w_2,w_N]$,
 $[13N]=[w_1,w_3]\cup [w_3, w_N]$ and so on up to $[1, N-1, N]$.

 \begin{figure}[H]
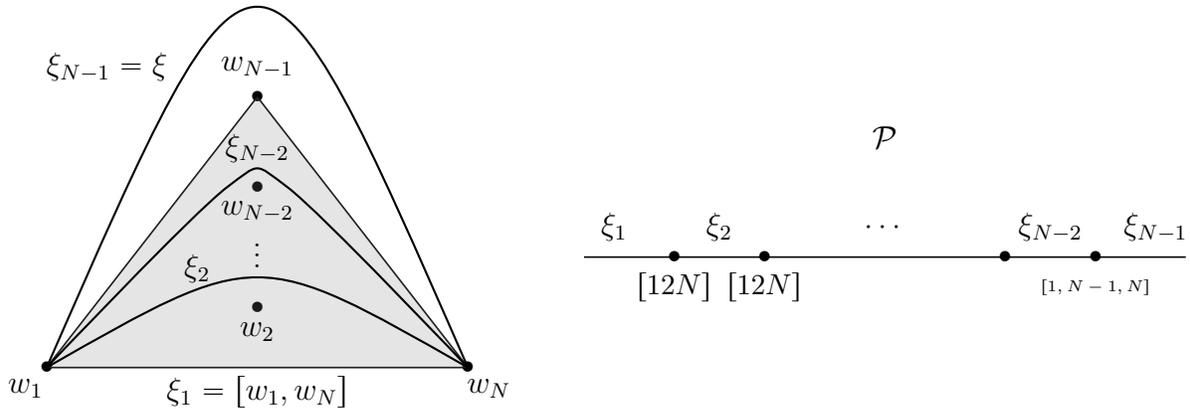

 \centering
 \btp[scale=0.4]
 \node at (-7,0){\small$\bullet$};
  \node at (7,0){\small$\bullet$};
\node at (0,2){\small$\bullet$};
   \node at (0,4){$\vdots$};
   \node at (0,6){\small$\bullet$};
    \node at (0,9){\small$\bullet$};
 \draw [line width=0.6] (-7,0) -- (0,9) -- (7,0) --(-7,0);
 \filldraw [color=gray, opacity=0.2] (-7,0) -- (0,9) -- (7,0) --(-7,0);

 \node at (-7.7, -0.7){$w_1$};
  \node at (7.7, -0.7){$w_N$};

  \node at (0,1.2){$w_2$};
   \node at (0,5.2){$w_{N-2}$};
  \node at (0,10){$w_{N-1}$};

  \draw[line width=0.8] (-7,0) ..controls (0,16) .. (7,0);
   \draw[line width=0.8] (-7,0) ..controls (0,4) .. (7,0);

  \draw[ line width = 0.8] plot [smooth,tension=0.5] coordinates{
 (-7,0) (-1,6) (1,6) (7,0)
};

 \node at  (-2,3.3){$\xi_2$};
 \node at (0,7.3){$\xi_{N-2}$};
\node at (-5,10){$\xi_{N-1}=\xi$};
\node at (0, -0.8){$\xi_1=[w_1, w_N]$};
 \etp
 \quad\quad
 \btp[scale=0.4]
 \draw[line width=0.6] (-10,0) -- (10,0);
 \node (A) at (0,-5){};
   \node at (-7,0){\small$\bullet$};
  \node at (-4,0){\small$\bullet$};
    \node at (4,0){\small$\bullet$};
  \node at (7,0){\small$\bullet$};

  \node at (0,4){$\Pc$};
 \node at (-9,1){$\xi_1$};
 \node at (-7,-1){\small$[12N]$};
 \node at (-5.5,1){$\xi_2$};
 \node at (-4,-1){\small$[12N]$};

 \node at (0,1) {$\cdots$};

 \node at (7,-1){\tiny$[1,N-1,N]$};
 \node at (9,1) {$\xi_{N-1}$};
 \node at (5.5, 1){$\xi_{N-2}$};

 \etp
\caption{The maximally nonconvex position.} \label{fig:max-nonconv}
\end{figure}
The diagram $\Pen$ in this case is a Postnikov system in the usual sense, i.e.,
(see  \cite{gelfand-manin} or \cite{KS-arr} \S 1A) a sequence
of interlocking exact triangles of the form:

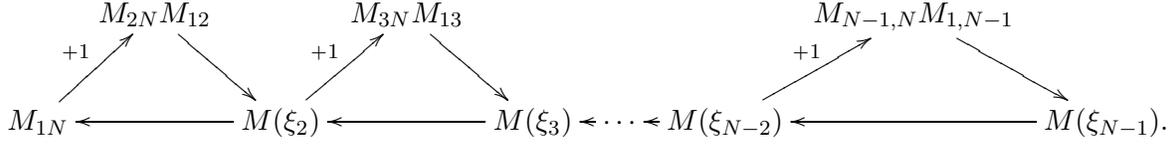
\[
\small{
 \xymatrix@C=0.5em{
 &M_{2N}M_{12} \ar[dr]
 &&M_{3N}M_{13}\ar[dr]
 &&&&M_{N-1,N}M_{1,N-1}\ar[dr]&
 \\
 M_{1N} \ar[ur]^{+1} &&
 \ar[ll] M(\xi_2) \ar[ur]^{+1} &&
 \ar[ll] M(\xi_3) & \cdots \ar[l] &
 \ar[l] M(\xi_{N-2}) \ar[ur]^{+1} &&
 \ar[ll] M(\xi_{N-1}).
}
}
\]

As with any  Postnikov system, $\Pen$ can be viewed in two different ways, which emphasize two
different parts of the diagram.

\vskip .2cm

\noindent {\bf  $\Pen$ as an analog of a cofiltration. } The part of $\Pen$ corresponding to empty paths (open cells)
\[
M(\xi)=M(\xi_{N-1}) \lra \cdots \lra M(\xi_2) \lra M(\xi_1) = M_{1N}
\]
can be seen as an analog of a cofiltration in $M(\xi)$. The polygonal transfers, i.e., $M_{1N}$ and
  $M_{1i}M_{iN}$,
play the role of kernels of the cofiltration. So $\Pen$ represents $M(\xi)$ as ``assembled'' out of
these kernels.

\vskip .2cm

\noindent{\bf $\Pen$ as a complex.} Composing the adjacent non-horizontal arrows in $\Pc$, we get
a  diagram in the derived category involving only polygonal transfers:
\[
 R^\bullet \,=\,\bigl\{ M_{1N} \lra M_{2N}M_{12}[1] \lra M_{3N}M_{13}[2] \lra \cdots \lra M_{N-1, N}M_{1,N-1}[N-2]\bigr\}.
\]

It is a cochain complex: the composition of any two consecutive arrows is zero. So $\Pen$
represents $M(\xi)$ as a ``total object'' of the complex $I^\bullet$.
  \end{ex}

  \begin{ex}[(The convex position)] Now assume that $A=\{w_1,\cdots, w_N\}$ consists of vertices
  of a convex $N$-gon, numbered clockwise, see
   Fig. \ref{fig:max-conv} for $N=4$.

   \begin{figure}[H]
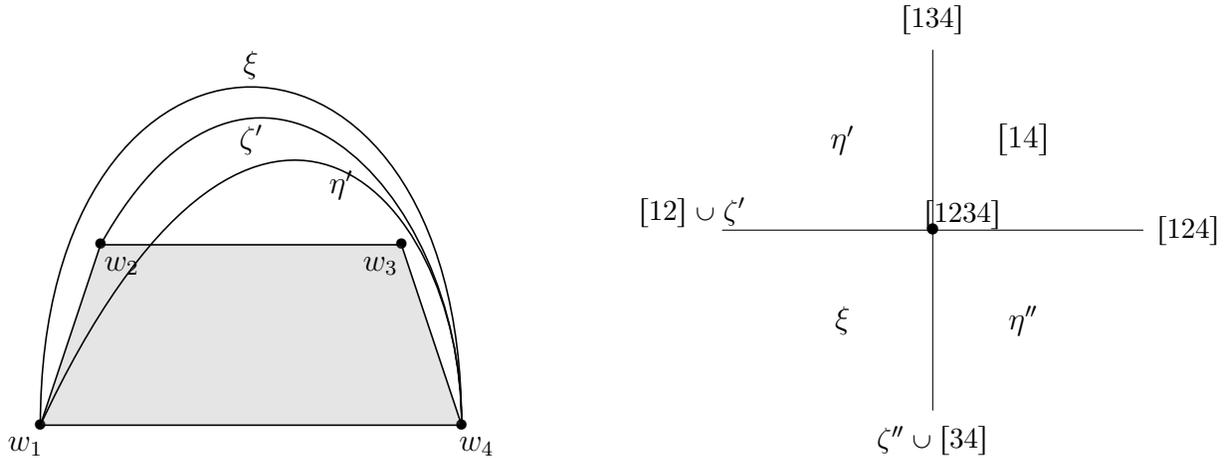

   \centering
   \btp[scale=0.4]

  \node at (-7,0){\small$\bullet$};
   \node at (-5,6){\small$\bullet$};
 \node at (5,6){\small$\bullet$};
  \node at (7,0){\small$\bullet$};

   \filldraw [color=gray, opacity=0.2] (-7,0) -- (-5,6) -- (5,6) -- (7,0) -- (-7,0);
   \draw [line width=0.6] (-7,0) -- (-5,6) -- (5,6) -- (7,0) -- (-7,0);

   \node at (7.5, -0.7){$w_4$};
      \node at (-7.5, -0.7){$w_1$};
 \node at (4.3, 5.3){$w_3$};
  \node at (-4.3, 5.3){$w_2$};
\draw[line width =0.6] (-7,0) .. controls (-7,15) and (7,15) .. (7,0);

 \draw[line width =0.6] (-5,6)  .. controls  (0,15) and (7,8)     .. (7,0);

  \draw[line width =0.6] (-7,0)  .. controls  (0,15) and (7,8) .. (7,0);
 \node at (3,8){$\eta'$};
  \node at (0,9.5){$\zeta'$};
  \node at (0,12){$\xi$};
  \etp
  \quad\quad\quad\quad
  \btp[scale=0.4]
  \node at (0,0){$\bullet$};
  \draw (0,-6) -- (0,6);
  \draw (-7,0) -- (7,0);

  \node at (1,0.5){\small$[1234]$};

  \node at (-3,-3){$\xi$};
  \node at (-3,3){$\eta'$};
  \node at (3,3){$[14]$};
  \node at (3,-3){$\eta''$};
 \node at (8.5,0){\small$[124]$};
 \node at (0,7){\small$[134]$};
 \node at (-8, 0.6){\small$[12]\cup \zeta'$};
 \node at (0,-7){\small$\zeta''\cup[34]$};

  \etp
   \caption{Convex position: $N=4$.}\label{fig:max-conv}
   \end{figure}

    In this case $\Hc$ consists of $N-2$ hyperplanes  intersecting transversely along a flat
  of codimension $N-2$, so  we have $3^{N-2}$ faces. For $N=4$ they are depicted on the right of
  Fig. \ref{fig:max-conv}. Here, along with the  paths $\xi$ (the circumnavigating path),
 $\eta'$ and $\zeta'$ depicted on the left figure,
  we use the paths $\eta''$ and $\zeta''$ which are symmetric to $\eta'$
and $\zeta'$ with respect to the vertical axis of symmetry of $Q$. For example, $\zeta''$ joins
$w_1$ with $w_3$, going over $w_2$. We also write $[ij]$ for $[w_i, w_j]$, as well as $[ijk]$ for
$[w_i, w_j]\cup [w_j, w_k]$ etc.

\vskip .2cm

  The infrared Postnikov system
  $\Pen$ is a commutative  $3\times 3\times\cdots\times 3$-diagram
  consisting of exact triangles (written linearly)  in each of the $N-2$ directions. For example, for $N=4$ we
  get a $3\times 3$ diagram with all rows and columns being exact triangles:
  \be\label{eq:3x3}
  \xymatrix{
  M(\xi)\ar[d]\ar[r]& M(\eta') \ar[d] \ar[r]& M(\zeta') M_{12} \ar[d]
  \\
  M(\eta'') \ar[r]\ar[d] & M_{14} \ar[r]\ar[d] & M_{24}M_{12}\ar[d]
  \\
  M_{34}M(\zeta'') \ar[r]& M_{34}M_{13}\ar[r]& M_{34}M_{23}M_{12}.
  }
  \ee
  This is not a Postnikov system in the  standard sense but it can be analyzed along the same lines as in
  Example \ref{ex:max-nonconv}:

  \vskip .2cm

  \noindent{\bf $\Pen$ as an analog of a filtration.} The part of $\Pen$ formed by empty paths, corresponding to
  the upper left $2\times 2$ square in \eqref{eq:3x3}, can be seen as an analog of a cofiltration in $M(\xi)$.
  This ``cofiltration'', however, is labelled not by a totally ordered set but by a poset:
  the  poset of all subsets in $\{2,\cdots, N-1\}$. Note that we are working in the derived category, and
  the naive abelian-categorical
   definition \eqref{eq:gr-alpha-F} of the kernels of a cofiltration labelled by a poset,
 does not  admit an immediate generalization to this case.

 \vskip .2cm

 \noindent{\bf $\Pen$ as a complex.} The part of the diagram $\Pen$  ``opposite'' to the above,
 i.e., corresponding to the lower  right $2\times 2$ square in  \eqref{eq:3x3},
 is a commutative $(N-2)$-dimensional
 hypercube formed by polygonal transfers. This hypercube can be converted
into a cochain complex $R^\bullet$ in a standard way.  The entire diagram can be seen as expressing
$M(\xi)$ as a total object of $R^\bullet$ in the sense that we will make precise later.   \end{ex}

\paragraph{The infrared complex.}\label{par:infra-com}
We now develop the cochain complex point of view on the   infrared Postnikov system $\Pen$
in the general context of Definition \ref{def:IPS}.

As in \eqref{eq:polyg-path}, let $\Lambda = \Lambda (i,j)\subset P$ be the set of (isotopy classes of) polygonal convex paths $\gamma$ with vertices in $A$ joining $w_i$ and $w_j$.
 Note that each such isotopy class contains precisely one polygonal path, so we can think of $L$
consisting of polygonal paths themselves.
For such a $\gamma$ we put $\bl(\gamma) = A\cap\gamma \-\{w_i,w_j\}$ to be the set of
the intermediate vertices of $\gamma$ and $l(\gamma)=|\bl(\gamma)|$ the number of such points, cf.
  Proposition  \ref{prop:baby-IA-counter}.

  \vskip .2cm

Consider also the {\em height function}
\[
h: \Lambda \lra \ZZ_+, \quad h(\gamma) \,=\, |\bh(\gamma)|, \quad \bh(\gamma)
 \,:=\, A\cap \Conv(\gamma)\- \{w_i, w_j\},
 \]
counting the number of ``new'' points (i.e., points other than $w_i, w_j$) on or under $\gamma$.
Thus $h(\gamma)-l(\gamma)$ is the number of ``new'' points strictly under $\gamma$.
Let also $\OR(\gamma) = \bigwedge^\max( \k^{\bh(\gamma)})$ be the $1$-dimensional ``orientation space''
of $\Gamma$.

\vskip .2cm

Let  $\Lambda _m\subset \Lambda $ be  the set of paths with height $m$.

 \vskip .2cm

We now define a cochain complex of functors
\be\label{eq:complex-R}
\begin{gathered}
R^\bullet \,=\, \bigl\{0 \to R^0  \buildrel q\over \to R^{1}
 \buildrel q\over  \to R^2 \buildrel q\over \to \cdots \big\}, \\
 R^{m}\,:=\, \bigoplus_{\gamma\in \Lambda_m} M(\gamma)\otimes_\k \OR(\gamma)
\,\in\,  \dgFun(\Phi_i, \Phi_j).
\end{gathered}
\ee
The differential $q: R^{m} \to R^{m+1}$ is defined as follows.

\vskip .2cm

\begin{figure}[H]
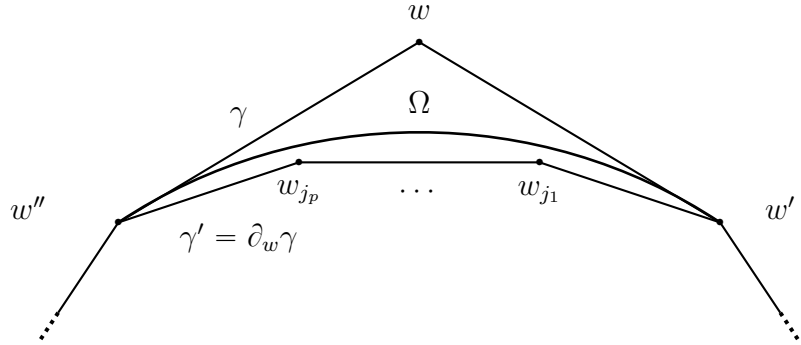

\centering
\btp[scale=0.4]

\node at (-10,3){\tiny$\bullet$};
\node at (10,3){\tiny$\bullet$};
\node at (4,5){\tiny$\bullet$};
\node at (-4,5){\tiny$\bullet$};
\node at (0,9){\tiny$\bullet$};
\draw[line width=0.8] (-12,0) -- (-10,3) -- (0,9) -- (10,3) -- (12,0);
\draw[line width=0.8] (-10,3) -- (-4,5) -- (4,5)  -- (10,3);
\draw[dotted, line width = 1.5] (-12,0) -- (-12.66, -1);
\draw[dotted, line width = 1.5] (12,0) -- (12.66, -1);

\node at (12, 3.5) {$w'$};
\node at (-13, 3.5) {$w''$};
\node at (0,10){$w$};
\node at (-4,4){$w_{j_p}$};
\node at (4,4) {$w_{j_1}$};
\node at (0,4){$\cdots$};

\draw[line width = 1] (-10,3) .. controls (-4,7) and (4, 7)  .. (10,3);

\node at (-6,6.5){$\gamma$};

\node at (-6,2.5){$\gamma' = \del_w\gamma$};

\node at (0,7){$ \Omega$};

\etp

\caption{Reduction $\gamma'=\del_w\gamma$ of a path $\gamma$ at a vertex $w$.}\label{fig:neighborly}
\end{figure}

Let $\gamma\in \Lambda_{m+1}$ and
$w\in\bl(\gamma)$ be an intermediate vertex. The {\em reduction} of $\gamma$ at $w$,
denoted $\del_w\gamma$, is a new convex path obtained by removing the vertex $w$ from
$\gamma$ and taking the convex hull of the remaining vertices, see Fig. \ref{fig:neighborly}.

\vskip .2cm Thus, a $2$-segment part $\ol \gamma = [w',w,w'']$
of $\gamma$ is replaced by a $p+1$-segment path
\[
\ol\gamma' \,=\, [w', w_{j_1}, \cdots, w_{j_p}, w'']
 \quad p\geq 0,
\]
of $\gamma'=\del_w\gamma$
so that there are no elements of $A$ inside  the nonconvex polygon bounded by $\ol\gamma$ and
$\ol\gamma'$.   By construction,
\[
\bh(\del_w\gamma) = \bh(\gamma)\-\{w\}, \quad h(\del_w\gamma)=h(\gamma)-1,
\quad l(\del_w\gamma) = l(\gamma)+p-1.
\]

Denote by $\ol\Omega$ the (isotopy class of the) curved path in the above polygon
that joins $w'$ and $w''$ avoiding all elements of $A$. Let
$\Omega$ be the path joining $w_i$ and $w_j$ which extends $\ol\Omega$ by coinciding with
 $\gamma$ and $\gamma'$
before $w'$ and after $w''$.

\vskip .2cm

The curvilinear triangle formed by $\ol\gamma$ and $\ol\Omega$ gives, by
 \eqref{eq:WC-triangle},
the Picard-Lefschetz ``comultiplication'' arrow
$v_{\ol\Omega, \ol\gamma}: M(\ol\Omega) \to M(\ol\gamma)[1]$.
Note that we also have a curvilinear polygon formed by $\ol\Omega$ and the polygonal path $\ol\gamma'$.
Because of the coassociativity of the ``multiplication'' Picard-Lefschetz $u$-arrows
(Proposition \ref{prop:PL-ass}), this latter triangle gives the iterated multiplication arrow
$u_{\ol\gamma', \ol\Omega}: M(\ol\gamma') \to M(\ol\Omega)$.
Let
\[
M(\gamma') \buildrel v_{\gamma', \Omega} \over\lra M(\Omega)
\buildrel u_{\Omega, \gamma}\over\lra M(\gamma)[1]
\]
be the natural transformations  between the  functors
$\Phi_i\to\Phi_j$ obtained by pre- and post-composing
$v_{\ol\gamma,\ol\Omega}$ and $u_{\ol\Omega, \ol\gamma'}$ by the identity natural transformations
of the transport functors corresponding to the common parts of $\gamma, \gamma'$ and
$\Omega$,
(unchanged under our move).

\vskip .2cm

Let now $\gamma\in\Lambda_{m+1}$ and $\gamma'\in\Lambda_{m}$. We define
the matrix element
\[
q_{\gamma', \gamma}: M(\gamma')\otimes\OR(\gamma') \lra M(\gamma)\otimes \OR(\gamma)[1]
\]
of $q$ to be $0$ unless $\gamma'=\del_w\gamma$ for some $w\in\bl(\gamma)$, in which
case we put

\[
q_{\gamma', \gamma}\,=\,  (u_{\Omega, \gamma}\circ  v_{\gamma', \Omega}) \otimes_\k \wedge w.
 \]
Here
\[
\wedge w: \OR(\gamma) = \bigwedge\nolimits^\max(\k^{\bh(\gamma')}) \lra \OR(\gamma) = \bigwedge\nolimits^{\max} (\k^{\bh(\gamma)})
\]
is the Grassman multiplication by an element of $(\k^{\bh(\Gamma)})$ corresponding to $w$.

\begin{prop}
The differential $q$ with the matrix elements $q_{\gamma', \gamma}$ satisfies $q^2=0$
and so $(R^\bullet, q)$ is a complex over
$H^0\, \dgFun(\Phi_i, \Phi_j)$, the  homotopy category of dg-functors $\Phi_i\to\Phi_j$.

\end{prop}

\noindent{\sl Proof:} By definition, each matrix element of $q^2$ has the form
\[
(q^2)_{\gamma'', \gamma} \,=\,\sum_{\gamma'\in \Lambda_{m+1}} q_{\gamma', \gamma}
 q_{\gamma'',\gamma'}, \quad
\gamma\in \Lambda_{m+2}, \,\,\gamma'\in \Lambda_{m}.
\]
For a summand in the RHS to be nonzero we must have $\gamma'=\del_w\gamma$ for some $w\in\bl(\gamma)$
and $\gamma''=\del_y \gamma'$ for some $y\in \bl(\gamma')$.
There are three possibilities:
\begin{itemize}
\item[(1)] $y\in \{w_{j_1}, \cdots, w_{j_p}\}$   is one of the ``new''  vertices which become exposed after we remove
$w$.

\item[(2)] $y$ is one of the ``old'' vertices of $\gamma$
which remain as vertices of $\gamma'$. Here one can subdivide further:
\begin{itemize}
\item[(2a)] $r\in\{w,w'\}$ is at one of the ends of the new segment $\ol\gamma'$.

\item[(2b)]   $r\notin\{w,w'\}$  is separated from the new segment $\ol\gamma'$.
\end{itemize}
\end{itemize}

\noindent In the case (1) there is only one possible path $\gamma'$, so
 $(q^2)_{\gamma'', \gamma} = q_{\gamma', \gamma} q_{\gamma'', \gamma'}$.
 Since each of the two factors here is itself a composition of a $u$-arrow and a $v$-arrow, the
 vanishing of their composition follows from the face that the composition of a $v$-arrow and a $u$-arrow
 in the same Picard-Lefschetz triangle is zero.

 \vskip .2cm

 In the case (2a) or (2b) there are two possible paths $\gamma'$, namely $\del_{w}(\gamma)$ and
 $\del_y(\gamma)$. In this case $(q^2)_{\gamma'', \gamma} $ will be the sum of two
equal summands with opposite signs.  \qed


\subsection{ The  infrared monad  and the Fukaya-Seidel monad }\label{subsec:infra-FS-mon}

\paragraph{The Fukaya-Seidel monad as circumnavigation.} The ``circumnavigating'' transport functors
studied in \S \ref{subsec:infra-com-schob} can be viewed as the components of the Fukaya-Seidel monad
$M$ from \eqref{eq:monad-M-spher} defining the Fukaya-Seidel category of a schober $\Sen\in\Schob(\CC,A)$.

\vskip .2cm

\begin{figure}[H]
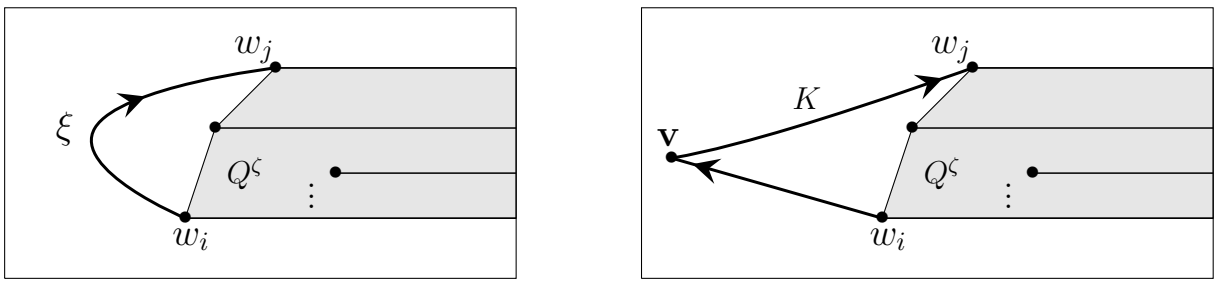

\centering

\btp[scale=.4]

\draw (-8, -4) -- (-8,5) -- (9,5) -- (9, -4) -- (-8,-4);
\node at (1,3){$\bullet$};
\node at (-1,1) {$\bullet$};
\node at (3,-0.5) {$\bullet$};
\node at (-2,-2) {$\bullet$};

\draw
[ line width=.2mm] (9,3) -- (1,3);

\draw [ line width=.2mm] (-1,1) -- (9,1);
\draw [line width=.2mm] (3,-0.5) -- (9,-0.5);
\draw [line width=.2mm] (9,-2) -- (-2,-2);

\node at (0.3, 3.6){\large$w_j$};
 \node at (2.2,  -1){\large$\vdots$};
\node at (-1.8,  -2.7){\large$w_i$};

\filldraw [opacity = 0.1] (9,3) -- (1,3) -- (-1,1) -- (-2,-2) -- (9,-2) -- (9,3);
\draw   (9,3) -- (1,3) -- (-1,1) -- (-2,-2) -- (9,-2) -- (9,3);

\draw
[decoration={markings,mark=at position 0.6 with
{\arrow[scale=2,>=stealth]{>}}},postaction={decorate},
line width=.4mm]
plot  [ smooth, tension=1] coordinates{ (-2, -2) ( -5, 1) (1,3)} ;

\node at (-6,1){\large $\xi$};

\node at (0,-0.5){$Q^{\zeta}$};

\etp
\quad\quad\quad\quad
\btp[scale=.4]

\draw (-10, -4) -- (-10,5) -- (9,5) -- (9, -4) -- (-10,-4);
\node at (1,3){$\bullet$};
\node at (-1,1) {$\bullet$};
\node at (3,-0.5) {$\bullet$};
\node at (-2,-2) {$\bullet$};

\draw
[ line width=.2mm] (9,3) -- (1,3);

\draw [ line width=.2mm] (9,1) -- (-1,1);
\draw [line width=.2mm] (3,-0.5) -- (9,-0.5);
\draw [line width=.2mm] (-2,-2) -- (9,-2);

\node at (0.3, 3.6){\large$w_j$};
 \node at (2.2,  -1){\large$\vdots$};
\node at (-1.8,  -2.7){\large$w_i$};

\draw
[decoration={markings,mark=at position 0.9 with
{\arrow[scale=2,>=stealth]{>}}},postaction={decorate},
line width=.4mm]
plot  [ smooth, tension=1] coordinates{   (-9,0)   ( -5, 1) (1,3)} ;

\draw
[decoration={markings,mark=at position 0.9 with
{\arrow[scale=2,>=stealth]{>}}},postaction={decorate},
line width=.4mm]
plot  [ smooth, tension=1] coordinates{(-2, -2)  
 (-9,0)} ;

\filldraw [opacity = 0.1] (9,3) -- (1,3) -- (-1,1) -- (-2,-2) -- (9,-2) -- (9,3);
\draw   (9,3) -- (1,3) -- (-1,1) -- (-2,-2) -- (9,-2) -- (9,3);

\node at (0,-0.5){$Q^{\zeta}$};

\node at (-9,0){$\bullet$};
\node at (-9.1, 0.7){\large$\bv$};

\node at (-4.5,2){$K$};
\etp

\caption{Circumnavigation = Fukaya-Seidel monad. } \label{fig:lasso}
\end{figure}

More precisely, fix a direction $\zeta$, $|\zeta|=1$ and suppose $A=\{w_1,\cdots, w_N\}$ is in linearly
general position including $\zeta$-infinity, so that the $\zeta$-direction rays
$w_i+\zeta\RR_+$ do not intersect, see Definition \ref{def:conv-pos}.
Let us number the $w_i$ so that $\Im(\zeta^{-1}w_1) < \cdots < \Im(\zeta^{-1} w_N)$
and consider the Fukaya-Seidel category $\FF(\CC,-\zeta;\Sen)$ corresponding to the opposite
direction $(-\zeta)$. This category is obtained by taking  the ``Vladivostok'' point
$\bv = \zeta R$, $R\gg 0$ very far in the direction $(-\zeta)$ and forming the $\bv$-spider $K$
consisting of the straight intervals $[\bv, w_i]$, almost parallet to each other, see the right of
Fig. \ref{fig:lasso}. This spider defines spherical functors $a_i:  \Phi_i =  \Phi_i(\Sen\to\Sen_\bv$
and we have the Fukaya-Seidel monad $\M= M(K) = (\M_{ij})_{i<j}$, where
$\M_{ij}= a_j^*\circ a_i: \Phi_i\to\Phi_j$. The category $\FF(\CC, -\zeta;\Sen)$ is then the
category of $\M(K)$-algebras.

\vskip .2cm

The nature of $K$ implies that $\M_{ij}$ is nothing but the circumnavigating transport  functor
from $w_i$ to $w_i$, see the left of Fig. \ref{fig:lasso}.

\paragraph{The infrared monad.}  The iterated rectilinear transforms studied
in  \S \ref{subsec:infra-com-schob}, form a monad as well. This is in fact a replacement,
in the general schober situation,  of the
$A_\oo$-algebra  $R_\oo$ defined in \cite{GMW1, GMW2}, see also
  \cite{KaKoSo} \S 10.

\vskip .2cm

More precisely, let $Q^\zeta = \Conv^\zeta(A)$ (the shaded region in Fig. \ref{fig:lasso}).
 We adjust slightly the definitions of   \S  \ref{subsec:infra-com-schob}\ref{par:infra-com}
 That is,  for any $i<j$  we let $\Lambda(i,j)$
 be the set of $\zeta$-convex polygonal paths $\gamma$ joining $w_i$ and $w_j$.
 For $\gamma\in\Lambda(i,j)$
 we have the iterated transport $M(\gamma):\Phi_i\to\Phi_j$.  Recall (Convention \ref{conv:clock})
 that in the definition of $M(\gamma)$ we identify the stalks of the local systems $\bPhi$ at
 the intermediate points using the clockwise monodromy: along the arcs going inside $Q^\zeta$.

 We also
  denote by  $l(\gamma)$ the number of intermediate vertices of $\gamma$ and
 $\bh(\gamma)$ the set  $A\cap\Conv^\zeta(\gamma\-\{w_i,w_j\}$ and $l(\gamma)=|\bl(\gamma)|$.
 Define the dg-functor $\R_{ij}:\Phi_i\to\Phi_j$ by
 \be\label{eq:R-ij}
 \begin{gathered}
 \R_{ij}\,=\,\bigoplus_{\gamma\in \Lambda(i,j)} R_\gamma,
 \\
 \R_\gamma\,=\,
 M(\gamma)\otimes \OR(\gamma),\quad
 \OR(\gamma) := \bigwedge\nolimits^\max (\k^{\bh(\gamma)}).
 \end{gathered}
 \ee
This is similar to the second line of \eqref{eq:complex-R}, except we do not consider the differential $q$ yet.

\vskip .2cm

Let $i<j<k$ and $\gamma \in\Lambda (j,k)$, $\gamma''\in \Lambda(i,j)$. If $\gamma \cup \gamma'$ is
$\zeta$-convex, then it belongs to $\Lambda (i,k)$, and we have a natural
 map
\[
c_{\gamma, \gamma' }: \R_\gamma \otimes\ R_{\gamma '}\lra \R_{\gamma \cup \gamma'}.
\]
It is  defined as the  tensor product of the  composition map
  \[
  M(\gamma')\circ M(\gamma) \lra M(\gamma'\cup\gamma),
  \]
  of the isomorphism
  \[
   \OR (\gamma')\otimes \OR(\gamma) \lra \OR(\gamma'\cup\gamma)
  \]
 obtained by contracting the two copies of the basis vector of $\k^A$.

 Let us define the map $c_{ijk}: R_{jk}\otimes R_{ij}\to R_{ik}$ by
\[
c_{ijk}|_{\R_\gamma \otimes \R_{\gamma'}} =  \begin{cases}
c_{\gamma,\gamma'}, &\text{ if } \gamma \cup \gamma' \text{ is convex};
\\
0,& \text { otherwise.}
\end{cases}
\]

 \begin{prop}\label{prop:infra-mon}
 The maps $c_{ijk}$ make $\R=\R^\zeta=(\R_{ij})_{i<j}$ into a uni-triangular monad.
 \qed
 \end{prop}

 \noindent{\sl Proof:} Since the $c_{ijk}$ are given by composition of functors, they are associative. \qed

 \vskip .2cm

 We call $\R$ the {\em infrared monad} associated to $\Sen$ and the direction $\zeta$.

\begin{ex}[(Maximally concave position)]
\label{ex:max-conc}
Suppose that $A$ lies on a single $(-\zeta)$-convex polygonal path. Then from the
point of view of $(+\zeta)$-convexity, $A$ is ``maximally concave'', see Fig. \ref{fig:max-conc}.

\begin{figure}[H]
\centering
\btp[scale=0.5]
\node at (4,3){$\bullet$};
\node at (6,1){$\bullet$};
\node at (5, -2){$\bullet$};
\node at (1, -4){$\bullet$};
\draw (4,3) -- (12,3);
\draw (6,1) -- (12,1);
\draw (5, -2) -- (12, -2);
\draw (1,-4) -- (12, -4);

\draw [line width= 1] (4,3) -- (6,1) -- (5,-2) -- (1, -4);

\node at (4, 4){\large$w_j$};
\node at (1.3, -5){\large$w_i$};

\filldraw [color=gray, opacity=0.2] (12,3) -- (4,3) -- (6,1) -- (5, -2) -- (1, -4)-- (12, -4);

\draw
[decoration={markings,mark=at position 0.9 with
{\arrow[scale=2,>=stealth]{>}}},postaction={decorate},
line width=.4mm]
plot  [ smooth, tension=1] coordinates{(1,-4)  
(4,3)} ;

\node at (1, -0.7){\large$\R_{ij}$};

\draw
[decoration={markings,mark=at position 0.9 with
{\arrow[scale=2,>=stealth]{>}}},postaction={decorate},
line width=.4mm]
plot  [ smooth, tension=1] coordinates{(1,-4)  ( -5, -0.7)
(4,3)} ;

\node at (-6.5, -0.7){\large$\M_{ij}$};

\etp

\caption{Maximally concave position: $\M_{ij}=\R_{ij} = M_{ij}$.}\label{fig:max-conc}
\end{figure}

\noindent  In this case, any $\zeta$-convex polygonal path consists of a single segment
$[w_i, w_j]$.
This segment is isotopic to the circumnavigating path passing through any faraway point
in the direction $(-\zeta)$. Therefore in this case the components of the
 monads $M$ and $R$ coincide:  $\M_{ij}=\R_{ij} =M_{ij}$ is the rectilinear transport.
  But the composition maps in $\R$ are zero,
 while in $\M$ they do not have to be zero.

\end{ex}

 \paragraph{Two layers of infrared analysis.} By the  {\em infrared analysis} of a schober $\Sen$
we will mean, somewhat simplistically,  the study of the relations between the infrared monad $R$ and the Fukaya-Seidel monad $M$ or of ways to obtain (the components of) $M$ from
(those of) $\R$.
 In fact, there can be two layers of such analysis:
\begin{itemize}
\item[(1)] Deforming the differential in each component $\R_{ij}$ to get a dg-functor
quasi-isomorphic to $\M_{ij}$.

\item[(2)] Deforming the compositions $c_{ijk}$ in the monad $\R$ to get an $A_\oo$-monad
quasi-isomorphic to $\M$.
\end{itemize}

\noindent These two steps involve two levels of convex geometry mentioned in the beginning
of this chapter. Here we consider the first step which is more elementary.

\vskip .2cm

Similarly to \eqref{eq:complex-R}, we have a differential $q$ in each $\R_{ij}$, thus upgrading
$\R_{ij}$ to a cochain complex
\be\label{eq:complex-R2}
\begin{gathered}
\R_{ij}^\bullet \,=\, \bigl\{ 0 \to \R_{ij}^0  \buildrel q\over \to \R_{ij}^1
 \buildrel q\over  \to \R_{ij}^2 \buildrel q\over \to \cdots\big\}, \\
 \R_{ij}^{m}\,:=\, \bigoplus_{\gamma\in \Lambda_m(i,j)} M(\gamma)\otimes_\k \OR(\gamma)
\,\in\,  \dgFun(\Phi_i, \Phi_j).
\end{gathered}
\ee
We would like to say that this complex ``represents'' the dg-functor $M_{ij}$. However,
 $\R^\bullet_{ij}$, as we defined it, is a complex over the triangulated category
 $H^0 \dgFun(\Phi_i, \Phi_j)$ and there is, a priori,  no canonical way of associating to
 it a single object of this category. Here, a general discussion may be useful.

\paragraph{Totalization of  complexes and Maurer-Cartan elements.}\label{par:tot-compl}

Let $\Cc$ be a pre-triangulated dg-category and $(C^\bullet, q)$ be a bounded complex over
the triangulated category $H^0 \Cc$. We would like to associate to it a ``total object''
$\Tot(C^\bullet)\in \Cc$ which would categorify the formal alternating sum $\sum (-1)^i C^i$.

\begin{rem}
In this section we use slightly different notation - the differential $q$ in a complex now has degree zero.
\end{rem}

 The reason that it is not automatic, is as follows.
\vskip .2cm

By our assumptions,  each $C^i$ is an object of $\Cc$ and
each $q: C^i\to C^{i+1}$ is an element of $H^0\Hom^\bullet_\Cc (C^i, C^{i+1})$
so that $q^2=0$ in $H^0\Hom^\bullet_\Cc (C^i, C^{i+2})$.
This means that if we lift   $q$ to a system $q_1$ of actual closed degree $0$ morphisms
$ q_{1,i}:  C^i \to C^{i+1}$, then $q_1^2 $  does not have to be zero but
has the form $[d, q_2]$ where $q_2$ is a system
of degree $(-1)$-morphisms $q_{2,i}: C^i\to C^{i+2}$. To continue in a meaningful
way, we need  to extend
$q_1, q_2$ to a full  coherence data. Such data can
 be formulated in terms of Maurer-Cartan elements.

\vskip .2cm

Let $(\gen^\bullet,d)$ be a dg-Lie algebra over $\k$. In particular, $\gen^0$, the degree $0$
component, is a Lie algebra in the usual sense.
We recall that a {\em Maurer-Cartan element}
in $\gen^\bullet$ is an element $\alpha\in \gen^1$ such that $d\alpha + {1\over 2} [\alpha, \alpha]=0$.
The set of Maurer-Cartan elements in $\gen$ is denoted $\MC(\gen^\bullet)$.

\vskip .2cm

In particular, let $(E^\bullet, d)$ be an associative dg-algebra. Then $E^\bullet$ can be considered
as a dg-Lie algebra via the graded commutator. In particular, we can speak about Maurer-Cartan
elements in $E^\bullet$. Explicitly, they are elements $\alpha\in E^1$ such that $d\alpha+ \alpha^2=0$.
The degree $0$ part $E^0$ is an associative algebra in the usual sense, so we have the
group $E^{0,\times}$ of invertible elements in $E^0$. This group acts on $\MC(E)$ by conjugation
 Two elements of $\MC(E^\bullet)$ lying in the same orbit of this action, are called
{\em gauge equivalent}.

\vskip .2cm

Returning to  our situation of $(C^\bullet, q)$,
consider the associative dg-algebra
\[
\Ec^\bullet \,=\,   \bigoplus_{i\leq j} \Hom^\bullet_\Cc
(C^i, C^j)[j-i-1],
\]
with the differential $[d,-]$. A Maurer-Cartan element $\wt q\in\MC(\Ec^\bullet)$
 is a sum $\wt q=\sum_{j\geq 1} q_j$ with
$q_j\in\bigoplus_i \Hom^{1-j}(C^i, C^{i+j})$ and the condition $[d,\wt q] + \wt q^2=0$
unravels into a system of equations of which the first two are $[d, q_1]=0$ and
 $q_1^2 + [d,q_2]=0$.

\begin{Defi}
An {\em MC-lifting} of the complex $(C^\bullet, q)$  is a Maurer-Cartan element
$\wt q \in\Ec^1$ such that the cohomology class of $q_1$ is $q$.
\end{Defi}

\begin{ex}
Let $\Cc = D^b(\Vect_\k)$ be the standard dg-enhancement of the bounded derived
category of $\k$-vector spaces. That is, objects of $\Cc$ are bounded complexes
$E^\bullet$ of vector spaces and $\Hom^\bullet(E^\bullet_1, E^\bullet_2)$
is the standard Hom-complex (with differential given by graded commutation with
the differentials in $E^\bullet_1$ and $E^\bullet_2$. A morphism in $H^0(\Cc)$ is a homotopy
class of morphisms of complexes. Thus, a complex over $H^0(\Cc)$ is a
bigraded vector space  $C^{\bullet\bullet}$ together with an actual differential $q_0: C^{i,j}\to C^{i,j+1}$,
of degree $(0,1)$, satisfying $q_0^2=0$ (so that each $C^i = (C^{i,\bullet}, q_0)$ is a complex) and a system of
homotopy classes of morphisms of complexes
 $q: C^{i,\bullet}\to C^{i+1,\bullet}$, squaring to $0$ up to homotopy.

An MC-lifting  of this data is a differential $\wt q= \sum_{i\geq 0} q_i$ in $C^{\bullet\bullet}$
with $q_n$ of  degree $(n, 1-n)$,  satisfying $\wt q^2=0$.
Since $\wt q$ is homogeneous of degree $1$ with respect to the total grading, we
 have a new complex $(T^\bullet, \wt q)$, where $T^m=\bigoplus_{i+j=m} C^{ij}$.

 \end{ex}

In the general case,
An MC-lifting $\wt q$ can be seen as a structure of a {\em twisted complex}
\cite{BK, seidel}  in $\Cc$  with the same
terms as $C^\bullet$. We will denote this twisted complex by $\wt C^\bullet$.
The pre-triangulated structure on $\Cc$ gives then the totalization $\Tot(\wt C^\bullet)$.
We have the following general statement.

\begin{prop}\label{prop:MC=PS}
Let $C^\bullet = \{ C^{-m} \buildrel q\over\to\cdots \buildrel q\over\to C^0\}$
be situated in degrees $[-m,0]$ (this does not restrict the generality). Then,
we have a bijection between the  sets consisting of:
\begin{itemize}
\item[(i)] Gauge equivalence
classes of MC-liftings of $C^\bullet$.

\item[(ii)]  Isomorphism classes of
Postnikov systems in $H^0\Cc$
\[
\small{
 \xymatrix@C=1em{
 &C^{-m} [-m+1]\ar[dr]^{h_m}
 &&&&C^{-2}[-1] \ar[dr]^{h_2} \ar@{-->}[rr]^q
 &&C^{-1} \ar[dr]^{h_1=q}&
 \\
 D^{[-m,0]}
 \ar[ur]^{+1}&& \ar[ll] D^{[-m+1,0]} &\cdots \ar[l] & \ar[l] D^{[-2,0]} \ar[ur]^{+1}
 && \ar[ll] D^{[-1,0]}  \ar[ur]^{+1}&& \ar[ll] D^0=C^0
}
}
\]
with the property that each  degree $1$ map  $C^{-i}[-i+1] \to C^{-i-1}[-i]$, $i=2, \cdots, m$,
obtained as the $\searrow \nearrow$ composition  (as indicated by the dotted horizontal arrow)
 is equal to $q: C^{-i}\to C^{-i+1}$
 \end{itemize}
\end{prop}

\noindent{\sl Proof (sketch):} Suppose that we have a MC-lifting  $\wt q =\sum_{j\geq 1} q_j$.
Let $C^{[-i,0]}$ be the truncation of $C^\bullet$ obtained by retaining only the
graded components $C^l$ for $l\in[-i,0]$. Then the components of $\wt q$ that
preserve $C^{[-i,0]}$, give a MC-lifting for it,  which we denote  $\wt q^{[-i,0]}$.
This gives rise to a twisted
complex $\wt C^{[-i,0]}$ with the same terms as $C^{[-i,0]}$. Let $D^{[-i,0]}=\Tot(\wt C^{[-i,0]})$.
Now, we have a filtration of twisted complexes
\[
C^0 \subset \wt C^{[-1,0]} \subset \wt C^{[-2,0]}\subset \cdots \subset \wt C^{[-m,0]}
\]
with quotients being $C^{-i}[-i]$. This filtration gives the desired Postnikov system
for the total objects.

\vskip .2cm

Conversely, suppose we have a Postnikov system as in (ii). Let $(q_1)_{-1,0}\in\Hom^0_\Cc(C^{-1}, C^0)$
be a closed element
representing the component $(q)_{-1,0}: C^{-1}\to C^0$ of $q$. Then we can identify
$D^{[-1,0]}$ with the canonical cone of $(q_{1})_{-1,0}$ given by the pre-triangulated structure on $\Cc$.
Then, the complex $\Hom^\bullet_\Cc(C^{-2}[-1], D^{[-1,0]})$ can be identified with the cone
of the morphism of complexes
\[
\Hom^\bullet_\Cc(C^{-2}, C^{-1}) \lra\Hom^\bullet_\Cc(C^{-2}, C^0)
\]
induced by $(q_1)_{-1,0}$. So the morphism $h_2$ can be thought as the degree $0$ cohomology class
of this latter cone. A $0$-cocycle $\wt h_2$
 representing $h_2$ consists of two components, which we denote
$(q_1)_{-2,-1}\in\Hom^0_\Cc(C^{-2}, C^{-1})$ and $(q_2)_{-2,0}\in\Hom^{-1}_\Cc (C^{-2}, C^0)$.
Further, we can think of $D^{[-2,0]}$ as the canonical cone of $\wt h_2$, so
$\Hom^\bullet_\Cc(C^{-3}[-2], D^{[-2,0]})$ is the cone of the moprhism
\[
\Hom^\bullet_\Cc(C^{-3}[-2], C^{-2}[-1]) \lra\Hom^\bullet_\Cc (C^{-3}[-2], D^{[-1,0]})
\]
and $h^3$ is a degree $0$ cohomology class of this cone. Choosing a $0$-cocycle $\wt h_3$
representing $h_3$, we decompose it into, first, the component denoted $(q_1)_{-3,-2}\in\Hom^0_\Cc (C^{-3}, C^{-2})$ and second, the component in $\Hom^0_\Cc(C^{-3}[-2], D^{[-1,0]})$
which, using the representation of $D^{[-1,0]}$ as a cone, is further split into two components
$(q_2)_{-3,-1}\in \Hom^{-1}_\Cc(C^{-3}, C^{-1})$ and $(q_3)_{-3,0}\in \Hom^{-2}_\Cc (C^{-3}, C^0)$.
Continuing like this, we get components
\[
(q_j)_{i, i+j}\,\in\, \Hom^{1-j}_\Cc (C^i, C^{i+j}), \quad j\geq 1.
\]
Denote $q_j\in \Hom^{1-j}_\Cc (C^\bullet, C^\bullet)$  the morphism with components
$(q_j)_{i,i+j}$ and  put $\wt q=\sum_{j\geq 1} q_j$.
The conditions that each $\wt h_j$ is a   closed (of degree $0$) element of the
corresponding Hom-complex,
imply together, in the standard way,  that
$\wt q$ is a Maurer-Cartan element. \qed

\paragraph{From the infrared Postnikov system to an MC-lifing?}
We are interested in applying the above to $C^\bullet = (\R^\bullet_{ij},q)$ being the
$(i,j)$th infrared complex \eqref{eq:complex-R2}. We would like to construct an MC-lifting
of the differential $q$.

\vskip .2cm

 Note that we have the infrared Postnikov system
$\Pen$ of Definition \ref{def:IPS}.  But Proposition \ref{prop:MC=PS} does not apply directly
since $\Pen$
 is a more sophisticated diagram than the
``linear'' Postnikov systems appearing there.
Only in the maximally nonconvex situation of Example \ref{ex:max-nonconv} we get a linear diagram.
In general, $\Pen$ does contain linear Postnikov systems (typically more than one)  and
we can apply  Proposition \ref{prop:MC=PS}  to those.
However, the  MC-liftings $\wt q$  thus obtained need not not compatible with the grading on $R_{ij}^\bullet$, i.e., may
contain components which  act within the  individual $R_{ij}^m$. Nevertheless, we think that the following is
true.

\begin{conj}\label{conj:MC-simple}
For each schober $\Sen\in\Schob(\CC,A)$,   there is an MC-lifting $\wt q$ of the differential $q$ in
 $\R^\bullet_{ij}$, canonical
  up to gauge equivalence, with the following property:  The total object of $(\R^\bullet_{ij},\wt q)$ is quasi-isomorphic to the component $M_{ij}$ of
 the Fukaya-Seidel  monad.
\end{conj}

 This should follow from a more systematic definition of schobers which will include  various
 coherence data from the very outset. We leave this question for future study.


\section{Algebra of the Infrared associated to a schober on $\CC$,  II}\label{sec:AIR-II}

In this chapter we consider the more advanced aspects of the ``infrared'' approach
to schobers, which use the notion of secondary polytopes. This corresponds to the
approach via webs in \cite{GMW1, GMW2}, see the discussion in \cite{KaKoSo}
as well as above, in Introduction
\S\ref{par:web-vs-sec}

\subsection{Background on secondary polytopes}\label{subsec:back-sec}

\paragraph{Basic definitions.}
The concept of secondary polytopes makes sense in the Euclidean space $\RR^n$ of any number of
dimensions \cite{GKZ}. In this paper we will be only interested in the case $n=2$ when $\RR^2$ appears as
the complex plane $\CC$. To fix the terminology, we give a brief reminder, under
this and  some other additional assumptions, see
\cite{GKZ} for a systematic exposition.

\vskip .2cm

Let $A=\{w_1,\cdots, w_N\}\subset \CC$, $N\geq 3$,  be a set in linearly general position. We denote
$Q=\Conv(A)$, a convex polygon.
Denote $\Vert(Q)$ the set of vertices of $Q$.
We will refer to the pair $(Q,A)$ as a {\em marked polygon}.
Thus a marking of a given polygon $Q$ consists in a choice of a set $A\subset Q$ in linearly
general position which contains $\Vert(Q)$.

\vskip .2cm

A {\em subpolygon} of $(Q,A)$ is a polygon $Q'$ such that $\Vert(Q')\subset A$. We write
$Q'\subset (Q,A)$.
A {\em marked subpolygon} of $(Q,A)$ is a marked polygon $(Q', A')$ such that $A'\subset A$.
We write $(Q',A')\subset (Q,A)$.
A marked subpolygon $(Q', A')\subset (Q,A)$ will be called {\em geometric}, if $A'=Q'\cap A$.
A {\em marked triangle} is a marked polygon $(Q,A)$ such that $Q$ is a triangle and
$A=\Vert(Q)$.

\vskip .2cm

A {\em polygonal decomposition} of a marked polygon $(Q,A)$ is a finite collection
$\Pc=\{(Q_\nu, A_\nu)\}$ of marked subpolygons in $(Q,A)$ such that
$Q=\bigcup_\nu Q_\nu$ and each intersection $Q_\nu\cap Q_\mu$, $\mu\neq\nu$,
is a common face of both (i.e.,  either empty, or a common vertex, or a common edge).
A {\em triangulation} of $(Q,A)$ is a polygonal decomposition $\Tc=\{(Q_\nu, A_\nu)\}$
such that each $(Q_\nu, A_\nu)$ is a marked triangle.
The set of polygonal decompositions is partially ordered by the relation $\preceq$
of {\em refinement}: we have $\Pc\preceq \Qc$, if $\Pc$ induces a polygonal decomposition
of any marked polygon from $\Qc$. In such a case $\Qc$ gives a (necessarily) regular decomposition
$\Qc_\nu$ of each $(Q_\nu, A_\nu)$.

\begin{Defi}
Let $\Pc=\{(Q_\nu, A_\nu)\}$ be a polygonal decomposition of $(Q,A)$.
A function $\psi: A\to\RR$ is called:

\begin{itemize}

\item [(a)]  {\em $\Pc$-piecewise affine}, if there is a (necessarily unique) continuous
function $f: Q\to\RR$ which is affine-linear on any $Q_\nu$ and such that $f(w)=\psi(w)$
for any $\nu$ and any $w\in A_\nu$.

\item [(b)] {\em $\Pc$-normal}, if
\begin{itemize}

\item  [(b1)] $\psi$ is  $\Pc$-piecewise-affine.

\item [(b2)] The function $f$ in (b) not piecewise affine
on any subset of $Q$ which strictly contains one of the $Q_\nu$.
In other words, $f$ does break along any common edge of any two polygons in $\Pc$.

\item[(b3)] The function $f$ is convex.

\item[(b4)] For $w\in A$ not lying in any of the $A_\nu$, we have $f(w) < \psi (w)$.
\end{itemize}
\end{itemize}
\end{Defi}

We denote by $\RR^A$ the space of all functions $\psi: A\to\RR$.  For a polygonal decomposition $\Pc$
we denote by $D_\Pc\subset \RR^A$ the  linear space formed  by $\Pc$-piecewise affine functions
and $C_\Pc\subset D_\Pc$ the set of $\Pc$-normal functions. It is clear that $C_\Pc$ is an open
convex cone in $D_\Pc$. A decomposition $\Pc$ is called {\em regular}, if $C_\Pc\neq\emptyset$.
The following is by now standard.

\begin{prop}\label{prop:R^A=UC_P}
The space $\RR^A$ is the disjoint union of the cones $C_\Pc$ over all (regular)  decompositions $\Pc$
of $(Q,A)$.
\end{prop}

\noindent{\sl Proof:} Given $\psi: A\to\RR$, we construct, first of all, a piecewise-linear convex function
$f: Q\to \RR$ by looking at the unbounded polyhedron
\[
G_\psi \,=\,\Conv\bigl\{ (w,t)\bigl| w\in A,\,  t\geq \psi(w)\bigr\} \,\subset \, \CC\times \RR = \RR^3
\]
and taking the finite (i.e., the bottom) part of its boundary as the graph of $f$.
The decomposition $\Pc=\Pc_\psi$
is then given by the maximal regions $Q_\nu$  of piecewise affine behavior of $f$, and the marking
of $Q_\nu$ is read off (b4).  Cf. \cite{KaKoSo} Fig.2. \qed

\vskip .2cm

Let $V$ be an $\RR$-vector space. By a {\em fan} in $V$ we will mean
a   decomposition of  $V$  into a  disjoint union of locally closed convex cones.

\begin{exas}
(a) Proposition \ref{prop:R^A=UC_P} means that we have a fan in $\RR^A$ formed by the
cones $C_\Pc$.  It is  called the {\em secondary fan}
of $A$ and denoted $\Sec(A)$.

\vskip .2cm

(b) Let $P\subset V$ be a (bounded) convex polytope. It gives rise to the {\em normal fan} $\Nc(P)$ in $V^*$
whose cones $N_F(P)$ are labelled by the faces $F\subset P$ of all dimensions. By definition, $N_F P$,
called the {\em normal cone} of $P$ at $F$,
consists of linear functionals $l: V\to\RR$ such that $l|_P$ achieves   minimum precisely on the face $F$.
\end{exas}

The concept of secondary polytopes connects the examples of type (a) and (b) for $V=\RR^A$. In this
case we identify $V^*$ with $\RR^A$ as well, using the delta-functions at $w\in A$ as an orthonormal basis.

\begin{propdef}
The secondary fan $\Sec(A)$ is realized as the normal fan of a convex polytope $\Sigma(A)\subset\RR^A$
called the {\em secondary polytope} of $A$. More precisely:

(a) $\Sigma(A)$ is a convex polytope of dimension $|A|-3$.

(b) Faces $F_\Pc$ of $\Sigma(A)$ are labelled by regular subdivisions $\Pc$ of $(Q,A)$.
The inclusion of faces $F_\Pc\subset F_\Qc$ correspons to refinement $\Pc\preceq\Qc$.

(c) In particular, vertices $F_\Tc$ of $\Sigma(A)$ correspond to regular triangulations $\Tc$ of $(Q,A)$.

(d) The dimension of $F_\Pc$ is equal to the codimension of the cone $C_\Pc$ or, what is the same,
of the subspace $D_\Pc$, in $\RR^A$.
\end{propdef}

The proof proceeds via an explicit construction of the vertices $F_\Tc\in\RR^A$, see
\cite{GKZ}.

\paragraph{ Faces of small dimension and codimension.}

\begin{exas}\label{ex:circuit}
Suppose that $A$ is a {\em circuit}, i.e.,
$|A|=4$.  Then $\Sigma(A)$ has dimension $1$, so it is an interval. There are precisely $2$
triangulations of $(Q,A)$, both regular.  These triangulations are coarse subdivisions of $(Q,A)$.
There are two possibilities:
\begin{itemize}
\item[(a)] $A$ is in convex position,
so $Q$ is a $4$-gon. The two triangulations  of $(Q,A)$ are just the two says of splitting $Q$ into two triangles.

\item[(b)] One point of $A = \{w_1, \cdots, w_4\}$, say, $w_4$, lies in the convex hull of the $3$ others. In this case $Q$ is a triangle (but $(Q,A)$ is not a marked triangle).
One triangulation of $(Q,A)$ splits $Q$ into $3$ smaller  triangles , and the other consists of
a single marked triangle $(Q, \{w_1, w_2, w_3\})$.

\end{itemize}

  \end{exas}

\begin{figure}[H]
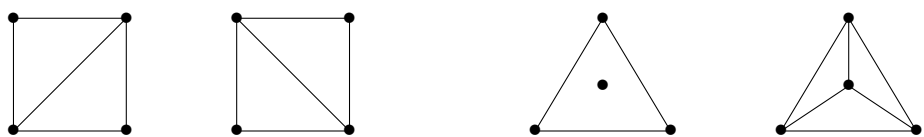

\centering

\btp[scale=.3, baseline=(current  bounding  box.center)]

\node (1) at (0,0){\small$\bullet$};
\node (2) at (5,0) {\small$\bullet$};
\node (3) at (5,5) {\small$\bullet$};
\node (4) at (0,5) {\small$\bullet$};
\draw (0,0) -- (0,5) -- (5,5) -- (5,0) -- (0,0);
\draw (0,0) -- (5,5);

\etp
\hskip 1cm
\btp[scale=.3, baseline=(current  bounding  box.center)]

\node (1) at (0,0) {\small$\bullet$};
\node (2) at (5,0) {\small$\bullet$};
\node (3) at (5,5) {\small$\bullet$};
\node (4) at (0,5) {\small$\bullet$};
\draw (0,0) -- (0,5) -- (5,5) -- (5,0) -- (0,0);
\draw (0,5) -- (5,0);

\etp
\hskip 2cm
\btp[scale=.3, baseline=(current  bounding  box.center)]

\node (1) at (-3,0) {\small$\bullet$};
\node (2) at (3,0) {\small$\bullet$};
\node (3) at (0,5){\small$\bullet$};
\node (4) at (0,2){\small$\bullet$};
\draw (-3,0) --(3,0) -- (0,5) -- (-3,0);
\etp
\hskip 1cm
\btp[scale=.3, baseline=(current  bounding  box.center)]

\node (1) at (-3,0){\small$\bullet$};
\node (2) at (3,0){\small$\bullet$};
\node (3) at (0,5){\small$\bullet$};
\node (4) at (0,2){\small$\bullet$};
\draw (-3,0) --(3,0) -- (0,5) -- (-3,0);
\draw (-3,0) -- (0,2) -- (0,5);
\draw (0,2) -- (3,0);
\etp

\caption{Circuits and flips.}\label{fig:circ-flip}
\end{figure}

\noindent So for $A$ a being a circuit, the interval $\Sigma(A)$   represents one of the standard  $2\leftrightarrow  2$ and $3\leftrightarrow 1$ flips of plane triangulations,
see Fig. \ref{fig:circ-flip}.  Compare also with the diagrams  in Figs. \ref{fig:part-comonad} and
\ref{fig:PL-3:1}.
For a general  $A$ it is known that  the edges of $\Sigma(A)$ correspond to flips of regular triangulations
along such circuits.

\begin{Defi}\label{def:coarse}
A regular decomposition $\Pc$ of $(Q,A)$ is called {\em coarse}, if $F_\Pc\subset\Sigma(A)$ is a face
of codimension $1$.
\end{Defi}

\begin{rem}
$\Pc$ being coarse means that  there is just one, up to rescaling and adding affine-linear functions,
$\Pc$-piecewise-affine continuous function $f: Q\to\RR$. So a simple way to construct a coarse
subdivision would be to cut $Q$ by a straight line into two subpolygons, so that $f$ breaks only along this line.
The triangulations in  Example  \ref{ex:circuit}(a) give an instance of this construction. However,  \ref{ex:circuit}(b) gives
a coarse subdivision into $3$ triangles. This shows that coarse subdivisions can  involve many polygons.

\end{rem}


\subsection{Variation of secondary polytopes and exceptional walls }\label{subsec:W-II-kind}

\paragraph{The phenomenon.}
Let $A\in\CC^N_\genr$ be in linearly general position.  The partially ordered set of all polygonal decompositions of
$(Q,A)$ (regular on not) depends only on the convex geometry of the set $A$, i.e., on the oriented matroid $\OO(A)$,
see Definition \ref{def:or-datum}. In particular, it remains unchanged inside any configuration chamber (connected component) $C\subset\CC^N_\genr$.

\vskip .2cm

However, the property of a given decomposition to be regular and therefore to correspond to a face of $\Sigma(A)$
is more subtle and   may not hold or fail uniformly on any configuration chamber $C$. This means that,
as $A$ varies in $C$, the structure of $\Sigma(A)$ can  change. This phenomenon, known in the
theory of secondary polytopes, was rediscovered in the dual language in \cite{GMW1} and called
{\em exceptional wall-crossing}.

\vskip .2cm

This means that we have a further chamber structure on the chamber $C$ given by a closed subset $D\subset C$
(union of ``exceptional walls'') so that $C\- D$ is dense and the structure of  $\Sigma(A)$
is unchanged on each connected component of $C\- D$.
In this section we quantify this picture in our case of two dimensions.

\paragraph{ The deformation complex of a subdivision.}
For any subset $I\subset\CC$ let $\Aff(I)\subset\RR^I$ be the  subspace formed by the restrictions to $I$ of
$\RR$-affine-linear functions $\CC\to\RR$. Thus $\dim_\RR \Aff(I)$ can be $1,2$ or $3$.

\vskip .2cm

Let $\Pc$ be a polygonal subdivision of $(Q,A)$, regular or not.
Denote by  $\Pc_0, \Pc_1, \Pc_2$  the
sets of all vertices, resp. edges, resp. polygons from $\Pc$ which do not lie on the boundary $\del Q$.
For any $i=0,1,2$ and any $\sigma\in\Pc_i$
we denote $\OR(\sigma) = H^i(\sigma, \del\sigma;  \RR)$ the $1$-dimensional orientation $\RR$-vector space of $\sigma$.
Thus, for $\sigma$ being a polygon ($i=2$), the space $\OR(\sigma)$ is canonically trivialized by the standard orientation of
$\CC$. Also, for $i=0$ ($\sigma$ a point), the space  $\OR(\sigma)$ is trivialized as well.   We now construct a
$3$-term cochain  complex
\be
\Def^\bullet(\Pc) \,=\,\biggl\{ \bigoplus_{\sigma\in \Pc_2} \Aff(\sigma)\otimes\OR(\sigma)
\buildrel d\over
\lra \bigoplus_ {\sigma\in\Pc_1} \Aff(\sigma)\otimes \OR(\sigma) \buildrel d\over
\lra \bigoplus_{\sigma\in\Pc_0 } \Aff(\sigma)\otimes\OR(\sigma) \biggr\},
\ee
with grading normalized so that the leftmost term is in degree $0$.
More precisely, the differential $d$ is defined to send
$\Aff(\sigma)\otimes \OR(\sigma)$ to the sum of $\Aff(\tau)\otimes\OR(\tau)$
where $\tau$ runs over codimension $1$ faces of $\sigma$ (i.e., over the edges of $\sigma$, if $\sigma$
is a polygon and over the two vertices of $\sigma$, if $\sigma$ is an edge). We define the matrix element
\[
d_{\sigma,\tau}: \Aff(\sigma)\otimes \OR(\sigma)\lra \Aff(\tau)\otimes\OR(\tau)
\]
to be equal to $r_{\sigma,\tau}\otimes \eps_{\sigma,\tau}$, where:
\begin{itemize}
\item  $r_{\sigma, \tau}: \Aff(\sigma)\to\Aff(\tau)$ is the
restriction  of    functions  from  $\sigma$ to $\tau$, and

\item  $\eps_{\sigma,\tau}: \OR(\sigma)\to\OR(\tau)$ is  the isomorphism given by canonical co-orientation
of $\tau$ in $\sigma$.
\end{itemize}

\begin{prop}
(a) We have $d^2=0$, i.e., $\Def^\bullet(\Pc)$ is indeed a cochain complex.

\vskip .2cm

(b) The $0$th cohomology space $H^0 \, \Def^\bullet(\Pc)$ is identified with $D_\Pc$,
the space of $\Pc$-piecewise-affine functions.

\vskip .2cm

(c)  We have $H^2\, \Def^\bullet(\Pc)=0$.
\qed
\end{prop}

\noindent{\sl Proof:}
(a) Let   $\sigma\in\Pc_2$ be  a polygon and let  $(\sigma, f)$
be an element of the $\sigma$th summand of $\Def^0(\Pc)$, so that $f\in\Aff(\sigma)$. Then
Then $d^2(\sigma,f)$ is the sum over vertices $w$ of $\sigma$, and the summand corresponding
to each such $w$ is, in its turn, the sum over the two edges of $\sigma$ bounded by $w$, of the
identical summands $f(w)$ but with the opposite signs coming from the orientations. So $d^2(\sigma,f)=0$.

\vskip .2cm

(b)
By definition, $H^0\,  \Def^\bullet(\Pc)$ consists of collections $(f_\sigma)$ of affine-linear
functions on the polygons $\sigma$ of $\Pc$  which match over the internal edges, i.e., precisely of $\Pc$-affine-linear functions.

\vskip .2cm

(c) For any $i$ and any $\sigma\in\Pc_i$ we have the embedding $\RR\subset \Aff(\sigma)$ of constant functions.
This gives a subcomplex
\be\label{eq:complex-C}
C^\bullet \,=\, \biggl\{ \bigoplus_{\sigma\in \Pc_2}  \OR(\sigma)
\lra \bigoplus_ {\sigma\in\Pc_1} \OR(\sigma)
\lra \bigoplus_{\sigma\in\Pc_0 } \OR(\sigma) \biggr\}\,\,\subset \,\, \Def^\bullet(\Pc).
\ee
This $C^\bullet$ is, up to shift, the relative cellular chain complex $C_\bullet(Q,\del Q; \RR)$ for the
cell decomposition given by $\Pc$.  So
\[
H^i  (C^\bullet) \,=\,H_{2-i}(Q,\del Q; \RR) \,=\,\begin{cases}
\RR,& \text{ if i=0};
\\
0, & \text{ otherwise.}
\end{cases}
\]
Note that $C^2=\Def^2(\Pc)$, so $\Def^\bullet(\Pc)/C^\bullet$ is situated in degrees $0$ and $1$, while
$H^2(C^\bullet)=0$ by the above. Therefore $H^2 \Def^\bullet(\Pc)=0$ as well. \qed

\vskip .2cm

Let us  denote by $\exc(\Pc)$ and call
the {\em exceptionality} of $\Pc$,  the dimension of $H^1\, \Def^\bullet(\Pc)$.
We say that $\Pc$ is {\em exceptional}, if $\exc(\Pc)\neq 0$.

\begin{prop}\label{prop:dim-FP-exc}
Suppose $\Pc$ is regular and $F_\Pc$ be the corresponding face of $\Sigma(A)$. Then
\[
\codim \, F_\Pc \,=\, 3|\Pc_2| - 2|\Pc_1| + |\Pc_0| -3 + \exc(\Pc).
\]

\end{prop}

\noindent{\sl Proof:} If $\Pc$ is regular, then $\codim\, F_\Pc = \dim D_\Pc-3$
is the dimension of the space of $\Pc$-piecewise-affine functions modulo globally affine functions.
To get our statement, it remains to
equate the alternating sum of dimensions of the terms of $ \Def^\bullet(\Pc)$ with the
alternating sum of the dimensions of the cohomology. \qed

\vskip .2cm

So  in the case  $\exc(\Pc)\neq 0$  we have  an ``unexpected jump'' in the codimension of $F_\Pc$.

\paragraph{The deformation complex and deformations of $A$.}
For a subset $I\subset\CC$ let $\Lin(I):=\Aff(I)/\RR$ be the quotient of $\Aff(I)$ by the subspace of
constant functions. Thus $\dim_\RR\Lin(I)$ can be $2,1$ or $0$. Denote $L=\Lin(\CC) = \Hom_\RR(\CC,\RR)$.

\vskip .2cm
Let $\Pc$ be any polygonal subdivision of $(Q,A)$.
Let $\ol\Def=\ol\Def^\bullet (\Pc):=\Def^\bullet (\Pc)/C^\bullet$ be the quotient of $\Def(\Pc)$ by the
subcomplex $C^\bullet$ of constant functions \eqref{eq:complex-C}. Explicitly,
\[
\ol\Def\,=\,\biggl\{ \bigoplus_{\sigma\in\Pc_2} \Lin(\sigma)\otimes\OR(\sigma) \lra
\bigoplus_{\sigma\in\Pc_1} \Lin(\sigma)\otimes\OR(\sigma)\biggr\},
\]
with the leftmost term in degree $0$. Note that for $\sigma\in\Pc_2$  we have $\Lin(\sigma)=L$,
while for $\sigma\in \Pc_1$ we have a short exact sequence
\[
0\to N^*_{\sigma/\CC} \lra L\lra \Lin(\sigma)\to 0,
\]
where $N^*_{\sigma/\CC}$ is the conormal space (i.e., the dual to the normal space $N_{\sigma/\CC}$)
of  $\sigma$ inside $\CC$.
Therefore we have a short exact sequence (with maps written vertically) of complexes (with differentials
written horizontally):
\be\label{eq:SES-Def}
\xymatrix{
0 \ar[d] &0\ar[d] &0\ar[d] &0\ar[d]
\\
\Kc =
\ar[d] & \biggl\{ \hskip 1.5cm 0 \ar[d]
\hskip 1.5cm \ar[r]& \bigoplus_{\sigma\in\Pc_1} N^*_{\sigma/C}\otimes\OR(\sigma) \ar[d]
\ar[r]^d& \bigoplus_{\sigma\in\Pc_0} L\otimes\OR(\sigma) \ar[d] \biggr\}
\\
C^\bullet\otimes L=
\ar[d] & \biggl\{ \bigoplus_{\sigma\in\Pc_2} L\otimes\OR(\sigma) \ar[r] \ar[d] &
\bigoplus_{\sigma\in\Pc_1} L\otimes\OR(\sigma)\ar[r] \ar[d] &  \bigoplus_{\sigma\in\Pc_0} L\otimes\OR(\sigma)\ar[d]
\biggr\}
\\
\ol\Def =
\ar[d] & \biggl\{ \bigoplus_{\sigma\in\Pc_2} L\otimes\OR(\sigma) \ar[r]. \ar[d] &  \bigoplus_{\sigma\in\Pc_1}
\Lin(\sigma) \otimes\OR(\sigma)\ar[d] \ar[r]  & \hskip 1.5cm 0\ar[d]\hskip 1.5cm  \biggr\}
\\
0&0&0&0
}
\ee
where in the middle we have the ``constant'' complex $C^\bullet\otimes L$, and $\Kc$ is defined as the
kernel of the projection.  Let us analyze this sequence.

\vskip .2cm

First of all,  $C^\bullet\otimes L$ has $H^0=L$ and no other cohomology. So  the long exact
sequence of cohomology corresponding to \eqref{eq:SES-Def} gives a short exact sequence
\be\label{SES-break}
0\to L\lra H^0(\ol\Def) \buildrel\delta\over\lra H^1(\Kc)\to 0
\ee
and an isomorphism
\[
H^1(\ol\Def) \buildrel\simeq\over\lra H^2(\Kc).
\]
Now, $H^0(\ol\Def) = D_\Pc/\RR$ is the space of $\Pc$-piecewise affine functions modulo constants.
So $H^0(\Def)/L = D_\Pc/\Aff(\CC)$ is the quotient of $D_\Pc$ by the subspace of global affine-linear functions.

\vskip .2cm
Further, $\Kc^1 = \bigoplus_{\sigma\in\Pc_1} N^*_{\sigma/\CC}\otimes\OR(\sigma)$ is the space of
``normal derivative data''  of linear functions with respect to various $\sigma\in\Pc_1$. So the coboundary
\[
\delta:  H^0(\ol\Def) \lra H^1(\Kc)\subset\Kc^1
\]
is the ``break data'' map. It sends a piecewise affine function $f=(f_\sigma)_{\sigma\in\Pc_2}$
(modulo global affine-linear functions) to the collection of breaks of $f$ at all the intermediate edges.

\vskip .2cm

Note that for $\sigma\in \Pc_1$ the orientation line $\OR(\sigma)$  is the same as the tangent line to $\sigma$,
and,  since $\OR(\sigma)^{\otimes 2}\simeq\RR$, it is also identified with the cotangent line,
i.e., with $\Lin(\sigma)$. Therefore $N^*_{\sigma/\CC}\otimes\OR(\sigma)$ is identified with $\RR$,
in particular, we can speak about positive and negative elements of it.

\begin{prop}\label{prop:reg-K+}
Let $\Pc$ be a polygonal subdivision of $(Q,A)$.
Let $\Kc^1_{>0}\subset \Kc^1$ be the octant formed by vectors
$(g_\sigma) \in\bigoplus_{\sigma\in\Pc_1} N^*_{\sigma/\CC}\otimes \OR(\sigma)$ with all $g_\sigma$
positive.
Then $C_P/\Aff(\CC) \simeq  H^1(\Kc)\cap \Kc^1_+$. In particular, $\Pc$ is regular if and only if
$H^1(\Kc)\cap \Kc^1_{>0}\neq\emptyset$.
\end{prop}

\noindent{\sl Proof:}
For any $\Pc$-piecewise affine function $f\in D_\Pc$
let $\ol f\in H^0(\ol\Def)$ be the image of $f$ (the class of $f$ modulo global affine linear functions).
The above analysis implies that
$C_\Pc$ consists of
$f\in D_\Pc $ such that $\delta(\ol f)$ has all components positive. Our claim follows from  this
via the short exact sequence \eqref{SES-break}. \qed

\vskip .2cm

Let us now consider the above picture from the point of view of deforming $A$.
When we move $A$
inside a given connected
component $C\subset \CC^N_\neq$, we can uniquely move the subdivision $\Pc$  as well,
as the convex geometry of $A$ does not change.
The space of infinitesimal deformations of $A$ is $T_A\CC^N=\CC^N$.
Let us call a {\em restricted deformation} a deformation that moves only points of $\Pc_0$,
leaving all the other elements of $A$ (that is, the vertices of $Q$ and the points inside $Q$
not in $\Pc_0$) unchanged.  The corresponding space of infinitesimal deformation is
$T_A^\res = \CC^{\Pc_0}$ which we view as the leftmost (degree $(-2)$) term of the dual complex
\[
\Kc^* \,=\,\biggl\{ \bigoplus_{w\in\Pc_0} L^*\otimes \OR(w) \buildrel \del\over\lra
\bigoplus_{\sigma\in\Pc_1} N_{\sigma/\CC}\otimes\OR(\sigma)\biggr\}
\]
Indeed, $L^*=\CC$ and for a point $w$ we have $\OR(w)=\RR$.

\begin{prop}\label{prop:del=rotation}
$H^{-2}(\Kc^*)=\Ker(\del)$, viewed as a subspace of $T_A^\res\CC^N$, consists of infinitesimal
deformations of $\Pc_0$ such that the corresponding infinitesimal deformation of $\Pc$ has all the edges
moving in a parallel way.
\end{prop}

We call restricted infinitesimal deformations of $A$ appearing in the proposition,  {\em parallel deformations}.

\vskip .2cm

\noindent{\sl Proof:}
Let $\xi = (\xi_w)_{w\in\Pc_0}$ be an infinitesimal deformation of $A$, considered as an element of $(\Kc^*)^{-2}$. Let us define $\xi_w=0$ for $w\in A\- \Pc_0$.

Let $\sigma =[w,w']\in \Pc_1$ so that either one or both   vertices $w,w'$ lie in $\Pc_0$.
Choose an orientation of $\sigma$, say from $w$ to $w'$.
Denote $p_\sigma: L^*=\CC\to N_{\sigma/\CC}$ be the projection to the normal space.
With these preliminaries,
$\sigma$-component of $\del(\xi)$ is, by definition $p_\sigma(\xi_w)-p_\sigma(\xi_{w'})\in N_{\sigma/\CC}$.
The vanishing of this difference means that, under the infinitesimal move of $[w,w']$ given by
moving $w$ according to $\xi_w$ and $w'$ according to $\xi_{w'}$, the direction of the interval
does not change. So vanishing of all the components of $\del(\xi)$ means that the deformation
is parallel. \qed

\begin{cor}
A subdivision $\Pc$ of $(Q,A)$ is exceptional, if and only if it admits parallel  deformations.
\end{cor}

\noindent{\sl Proof:} By definition, $\Pc$ is exceptional, if  the space
\[
H^1(\Def^\bullet(\Pc))\,\simeq \, H^1(\ol\Def) \,\simeq \, H^2(\Kc)
\]
is nonzero, in other words if $H^{-2}(\Kc^*)\neq 0$. \qed

\begin{figure}[H]
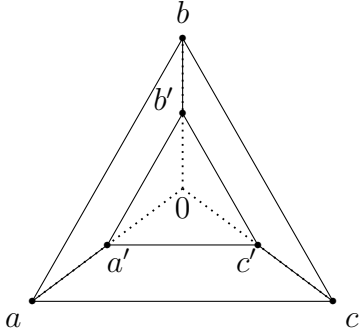

 \centering
 \btp[scale=0.5]
 \node at (0,4){\tiny$\bullet$};
  \node at (4,-3){\tiny$\bullet$};
 \node at (-4, -3){\tiny$\bullet$};

  \node at (2, -3/2){\tiny$\bullet$};
    \node at (-2, -3/2){\tiny$\bullet$};
  \node at (0,2){\tiny$\bullet$};

  \draw (0,4) -- (4,-3) -- (-4,-3) -- (0,4);
\draw (2, -3/2) -- (-2, -3/2) -- (0,2) -- (2, -3/2);

\draw (0,2) -- (0,4);
\draw (2, -3/2) -- (4, -3);
\draw (-2, -3/2) -- (-4, -3);

\draw[dotted, line width=0.7] (0,0) -- (0,4);
\draw[dotted, line width=0.7] (0,0) -- (4,-3);
\draw[dotted, line width=0.7] (0,0) -- (-4,-3);

\node at (-4.5, -3.5){$a$};
\node at (0,4.7){$b$};
\node at (4.5, -3.5){$c$};

\node at (-1.7, -1.9){$a'$};
\node at (1.7, -1.9){$c'$};

\node at (-.5, 2.4){$b'$};

\node at (0, -0.5){$0$};
 \etp
 \caption{An exceptional subdivision with an obvious parallel deformation. }\label{fig:parallel-def}
 \end{figure}

\begin{ex}
Consider two ``concentric''  triangles in $\CC$, i.e., a triangle $[a,b, c]$ contaning $0$ inside and
a contracted triangle $[a'=\lambda a, b'=\lambda b, c'=\lambda c]$, $\lambda <1$. Let
$A=\{a,b,c,a', b', c' \}$, so $Q=[a,b,c]$ and let $\Pc$ be the subdivision of $(Q,A)$ in  Fig.
\ref{fig:parallel-def}.
Then $\Pc$ is exceptional.  Its parallel restricted deformation is obtained by varying $\lambda$.
\end{ex}

\paragraph{ Exceptionally general position.}
We would like to specify regions in $\CC^N_\genr$ where ``the structure of the secondary polytope
does not change''. A precise way of expressing this is as follows.

\vskip .2cm

For any convex polytope $P$ let $\Fen(P)$ be the poset of
nonempty faces of $P$ (including $P$ itself), ordered by inclusion.
 Two polytopes $P_1, P_2$
 are called {\em combinatorially equivalent}, if  there is an isomorphism of posets $\Fen(P_1)\to \Fen(P_2)$.
 For example, all plane $n$-gons with given $n$ are combinatorially equivalent.

\begin{Defi} We  say that $A\in \CC^N_\genr$ is in {\em exceptionally general position},
if,  for  any
regular subdivision $\Pc$ of $(Q, A)$,   we have $ \exc(\Pc)=0$.
Denote by $\CC^N_\egen\subset\CC^N_\genr$ the set of $A$ in exceptionally general position
and put $D=\CC^N_\genr\-\CC^N_\egen$.
\end{Defi}

\begin{prop}\label{prop:Egen-dense}
(a) $\CC^N_\egen$ is an $\RR$-Zariski open    subset in $\CC^N_\genr$, i.e., $D$ is an
$\RR$-Zariski closed subset (zero locus of a system of real polynomials).

\vskip .2cm

(b) $\CC^N_\egen$ is dense in $\CC^N_\genr$ with respect to the usual transcendental topology.
In other words, $\codim(D)\geq 1$.

\vskip .2cm

(c) For any connected component $E\subset \CC^N_\egen$, the combinatorial equivalence class of the
polytopes $\Sigma(A)$, $A\in E$, is the same.
\end{prop}

\begin{rem}
The set $D$ can be seen as the union of  ``exceptional walls'' of \cite{GMW1} \S 8.3. However, we cannot exclude
the  possibility that
some irreducible components of $D$ have codimension $>1$, in which case they would not represent walls in
the strict sense.
\end{rem}

\noindent{\sl Proof of Proposition \ref{prop:Egen-dense}:}  (a) Let $A$ vary inside a component $C$ of $\CC^N_\genr$.
As $A$ moves, we can move $\Pc$ with $A$. This means that the deformation complex
$\Def^\bullet(\Pc)$ also moves, i.e., we get a complex of  $\RR$-algebraic vector bundles on $C$,
considered as an $\RR$-algebraic variety.
Therefore the subset in $C$ over which some particular cohomology group of this complex,
in our case $H^1$, is non-zero, is $\RR$-Zariski closed, given by vanishing of certain
minors of the matrices of the differentials. So $\CC^N_\egen\cap C$ is  $\RR$-Zariski open in $C$.

\vskip .2cm

(b) To show that $\CC^N_\egen$ is dense, it is enough to show that for any $A\in\CC^N_\genr$ and
any exceptional decomposition $\Pc$ of $(Q,A)$, there is a nearby deformation (actual, not infinitesimal)
$A'$ of $A$
so that the corresponding decomposition $\Pc'$ of $(Q',A')$ is no longer regular.

\begin{figure}[H]
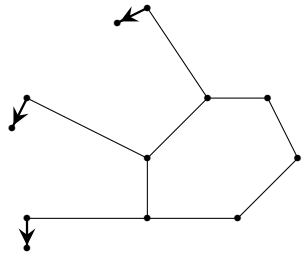

\centering
\btp[scale=0.4]

\node at (0,0){\tiny$\bullet$};
\node at (2,2){\tiny$\bullet$};
\node at (4,2){\tiny$\bullet$};
\node at (5,0){\tiny$\bullet$};
\node at (3,-2){\tiny$\bullet$};
\node at (0,-2){\tiny$\bullet$};

\draw (0,0) -- (2,2) -- (4,2) -- (5,0) -- (3, -2) -- (0, -2) -- (0,0) ;

\node at (0,5){\tiny$\bullet$};
\node at (-4,2){\tiny$\bullet$};
\node at (-4,-2){\tiny$\bullet$};


\draw (-4,-2) -- (0,-2);
\draw (-4,2) -- (0,0);
\draw (0,5) -- (2,2);

 \draw  [  decoration={markings,mark=at position 0.9 with
{\arrow[scale=1.5,>=stealth]{>}}},postaction={decorate},
line width=.3mm] (0,5) -- (-1,4.5);

\node at (-1,4.5){\tiny$\bullet$};

 \draw  [  decoration={markings,mark=at position 0.9 with
{\arrow[scale=1.5,>=stealth]{>}}},postaction={decorate},
line width=.3mm] (-4,2) -- (-4.5, 1);

\node at (-4.5,1){\tiny$\bullet$};

 \draw  [  decoration={markings,mark=at position 0.9 with
{\arrow[scale=1.5,>=stealth]{>}}},postaction={decorate},
line width=.3mm] (-4,-2) -- (-4, -3);

\node at (-4,-3){\tiny$\bullet$};
\etp
\caption{The deformation $A'$ of $A$ obtained by moving the vertices of $Q$ counterclockwise.}
\label{fig:def-A'}
\end{figure}

So we suppose that $\Pc$ is exceptional, i.e., admits a parallel   infinitesimal deformation
$\xi= (\xi_w)_{w\in\Pc_0}$ and exhibit  an $A'$ as above.
In searching for such an $A'$, we are allowed to move the vertices of $Q$ as well.
Let $A'$ be obtained by moving all the vertices of $Q$ a little in the counterclockwise direction,
while keeping all the other points, in particular, those from $\Pc_0$, unchanged, see Fig. \ref{fig:def-A'}.

Let $\Pc'$ be  the subdivision of $(Q', A')$ obtained by deforming $\Pc$. Note that $\Pc'_0=\Pc_0$.
We claim that  $\Pc'$
is not regular. Indeed, the restricted deformation of $\Pc$ given by
the same collection of vectors $\xi$, will now  not be parallel,  and instead we have the following property:
\begin{itemize}
\item [(*)] If an edge $[w,w']$ of $\Pc'$ has both vertices in $\Pc'_0=\Pc_0$, then the direction
of $[w,w']$ does not change, as before. If, however, $[w,w']$ has only one vertex in $\Pc_0$
and the other one on $\del Q'$, then the effect of $\xi$ will rotate the edge in the counterclockwise
direction, see Fig. \ref{fig:def-A'}.

\end{itemize}

\noindent We claim that (*) implies that $\Pc'$ is not regular. Let us use
Proposition \ref{prop:reg-K+} and denote $\Kc(\Pc')$ the complex $\Kc$ for $\Pc'$.
Let $d$ be the differential in $\Kc(\Pc')$, as in \eqref{eq:SES-Def} and $\del$ be the dual
differential in  $\Kc^*(\Pc')$, as above.
By  \ref{prop:reg-K+}, regularity of $\Pc'$ means that the subspace $\Ker(d)=H^1\Kc({\Pc'})$
intersects the  strictly positive octant $\Kc^1(\Pc')_{ >0}$, so that the corresponding
intersection is essentially the cone $C_\Pc$. In dual terms, this means that:
\begin{itemize}
\item[(**)]  $\Im(\del)$ intersects the non-strictlly positive octant $\Kc^{1 }(\Pc)^*_{ \geq 0}$ only at $0$
or along a single half-line.
\end{itemize}
Indeed, the dual cone to $\Ker(d)\cap \Kc^1(\Pc')_{ \geq 0}$ is the image, denote it $I$,
of the non-strict octant $\Kc^{1}(\Pc')^*_{ \geq 0}$
in $\Coker(\del) = \Kc^{1*}(\Pc')/\Im(\del)$. So nonemptiness of  $\Ker(d)\cap \Kc^1(\Pc')_{ >0}$,
means that $I$ does not contain a linear space, which  is equivalent to (**).

\vskip .2cm

But, as shown in the proof of Proposition \ref{prop:del=rotation},  $\del$ is given by the infinitesimal rotation of
the edges, so (*) means that $\Im(\del)$
contains a vector which have some components $0$ and other (at least 3) components
strictly positive. This means that $\Pc'$ cannot be regular.

\vskip .2cm

This proves part (b) of Proposition  \ref{prop:Egen-dense}. Let us show part (c).
Let $A$ move along a path in $E$. Then all  subdivisions $\Pc$ of $(Q,A)$
move as well, so we write $\Pc(A)$ for the subdivision for a given $A$. Now, since all regular subdivisions
$\Pc(A)$
remain non-exceptional, the dimension of the corresonding face $F_{\Pc(A)}\subset \Sigma(A)$ remains the same by
 Proposition \ref{prop:dim-FP-exc}. This means that the polytope $\Sigma(A)$ cannot undergo a combinatorial
 change.  \qed

\paragraph{Factorization property of secondary polytopes.}
.
\begin{prop}
If $A\in\CC^N_\egen$, then for  every regular subdivision $\Pc=(Q_\nu, A_\nu)$, the
face $F_\Pc$ of $\Sigma(A)$ is combinatorially equivalent to
the product $\prod_\nu \Sigma (A_\nu)$.
\end{prop}

This proposition corrects the treatment in \cite{KaKoSo}, where the condition of
exceptionally general position was overlooked.

\vskip .2cm

\noindent{\sl Proof:}   The faces   of $F_\Pc$ are the $F_\Qc$ where $\Qc\preceq \Pc$
is a regular subdivision   refining $\Pc$.
Such $\Qc$ consists of a system of subdivisions  $\Pc_\nu$ of $(Q_\nu, A_\nu)$
which are necessarily regular, since $\Pc$ is regular. This gives an embedding  of posets
\[
a: \Fen( F_\Pc) \,\,\hra\,\,  \prod_\nu \Fen(\Sigma(A_\nu))\,=\, \Fen\biggl(\prod_\nu \Sigma(A_\nu)\biggr).
\]
We claim that $a$ is a bijection, hence an isomorphism of posets.
Since $A\in\CC^N_\egen$, the subdivision $\Pc$ is non-exceptional, and so
Proposition \ref{prop:dim-FP-exc}  calculates $\dim F_\Pc$ in combinatorial terms.
This calculation implies that $\dim F_\Pc$
equal the sum
of the dimensions of the $\Sigma(A_\nu)$, i.e., the dimension of $\prod_\nu \Sigma(A_\nu)$
(we leave the elementary computation to the reader).  So bijectivity of $a$
is a particular case of the following fact.

\begin{lem}
Let $P_1, P_2$ be two convex polytopes of the same dimension, and $b: \Fen(P_1)\hra \Fen(P_2)$
be an embedding of posets. Then $b$ is bijective.
\end{lem}

\noindent{\sl Proof of the Lemma:}
Let $\delta = \dim P_1=\dim P_2$. The statement is clearly true for $\delta=0,1$,
so we argue by induction, assuming the lemma true for the case when the common
dimension of $P_1$ and $P_2$ is  less than $\delta$.

\vskip .2cm

For a poset $(S,\leq)$ we call the {\em height} of $S$ and
denote by $h(S)$ the maximal $d\geq 0$ such that there is a chain $s_0<\cdots < s_s$
of strictly increasing elements of $S$ (it may be $h(S)=\oo$).

\vskip .2cm

For any convex polytope $P$ we have $\dim P=h(\Fen(P))$. In particular, if $F$ is a face of $P_1$,
then $\dim(F)$ is the height of the interval $\Fen(F)=\Fen(P_1)_{\leq F}$. Since $b$ is an embedding,
this means that $\dim b(F)\geq \dim F$. Since $\dim(P_1)=\dim(P_2)$, we conclude that
$\dim b(F)=\dim F$ for any face $F$ of $P_1$.

\vskip .2cm

In particular, $b$ takes the maximal face $P_1\in\Fen(P_1)$ to $P_2\in\Fen(P_2)$, so $P_2\in \Im(b)$.
Let us prove that all  proper (non-maximal) faces of $P_2$ lie in $\Im(b)$ as well.

\vskip .2cm

By the inductive assumption, for any proper   face $F\in \Fen(P_1)$ the embedding
$b$ induces an isomorphism
\[
b_F:   \Fen(P_1)_{\leq F}= \Fen(F)  \lra  \Fen(b(F))= \Fen(P_2)_{\leq b(F)}.
\]
Now, the geometric realization of the poset of proper faces of $P_i$, $i=1,2$,  is homeomorphic to the sphere
$S^{\delta-1}$. So the isomorphisms $b_F$ for proper $F\in\Fen(P_1)$ give, after passing to
geometric realizations,  a continuous embedding
$S^{\delta-1}\hra S^{\delta-1}$. By topological reasons, such an embedding must be a homeomorphism.
Therefore $b$ is surjective on proper faces as well. \qed


\subsection{The  $\Lie_\oo$-algebra  associated to a schober on $\CC$}\label{subsec:L-oo-andR_oo}

The main feature of the Infrared Algebra formalism \cite{GMW1,  GMW2}
is a $\Lie_\oo$-algebra $\gen$ associated to a ``theory''. As shown in \cite{KaKoSo},
$\gen$ can be defined using secondary polytopes. In this section we show how such a datum can be
constructed from any schober $\Sen\in\Schob(\CC,A)$.

We consider the situation of \S \ref{subsec:schober-rect}.
Let   $A=\{w_1,\cdots, w_N\}\subset \CC$ be in linearly general position
and $Q=\Conv(A)$.

\paragraph{Framed subpolygons, their content and pasting.}
We start by  adding some complements to the discussion of pasting diagrams in
\S \ref{subsec:PL-tri}\ref{par:pasting-NT}

 By a {\em subpolygon} in $Q$ we mean a set of the form $P=\Conv(B)$,
where  $B\subset A$ has cardinality at least $3$.
Let  $P\subset Q$ be a subpolygon.  By a {\em framing} of $P$ we mean a datum $f$ of
two vertices $\alpha,\omega$ of $P$ (possibly equal) called the {\em $0$-dimensional source}
and {\em target} of $P$ (with respect to $f$). These vertices subdivide the boundary $\del P$ into
two complementary polygonal arcs
\be\label{eq:del(P,f)}
\begin{gathered}
\del_+(P,f) \,=\, [\alpha = w_{i_0}, w_{i_1},\cdots, w_{i_r}=\omega],
\\
\del_-(P,f) \,=\, [\alpha = w_{j_0}, w_{j_1}, \cdots, w_{j_s}=\omega], \quad i_0=j_0, \, i_r=j_s,
\end{gathered}
\ee
with common beginning and end. More precisely, $\del_+(P,f)$ is that of the two arcs
which goes in the counterclockwise direction around $P$, and $\del_-(P,f)$ is
the arc going in the  clockwise direction. These two arcs are referred to as the
{\em $1$-dimensional source} and {\em target} of $P$ with respect to $f$.
In this way a choice of a framing makes $P$ into a pasting diagram with one $2$-dimensional
cell.

\begin{ex}\label{ex:f-zeta}
As in Example \ref{ex:pasting-convex}(b), let $\zeta$ be a direction so that $A$ is
in linearly general position including $\zeta$-infinity and $p_\zeta: \CC\to \RR$
be the projection to a line orthogonal to $\zeta$. This datum defines a framing $f_\zeta$
on $P$, with $\alpha$ and $\omega$ being the vertices where $p_\zeta$ achieves the
minimum and maximum, and $\del_+(P, f_\zeta) = \del_+^\zeta(P)$ resp.
$\del_-(P,f_\zeta)=\del_-^\zeta(P)$ is the upper, resp. lower part of the boundary
if looked at from $\zeta R$, $R\gg 0$, see Fig. \ref{fig:pasting-poly}.
\end{ex}

\begin{Defi}
The {\em content} of a framed polygon $(P,f)$ is defined as
\[
\cont(P,f) \,=\, | A\cap (P\- \del_+(P,f))|,
\]
i.e., as the number of elements of $A$ lying either in the interior of $P$ or
strictly inside $\del_-(P,f)$.
\end{Defi}

Let now  $\zeta$ be generic as above and $\Pc = (P_i)_{i\in I}$ be a polygonal decomposition of $P$ into subpolygons,
as in Example  \ref{ex:pasting-convex}(b). By Example \ref{ex:f-zeta},  we get framings
(all denoted $f_\zeta$) on $P$ and on each $P_i$. The following is straightforward.

\begin{prop}
In the situation described we have $\cont(P,f_\zeta)=\sum_i \cont(P_i,f_\zeta)$. \qed
\end{prop}

\begin{figure}[H]
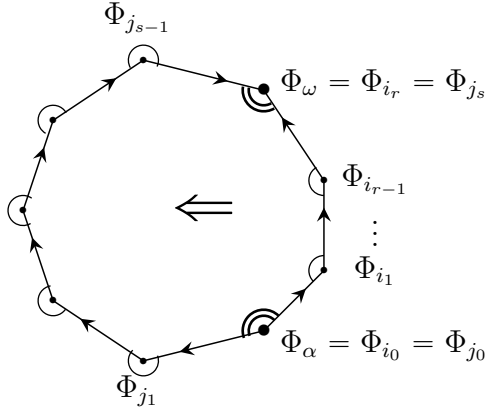

\centering

 \btp[scale=0.4]

\node at (4,4){$\bullet$};

\centerarc[line width=1](4,4)(170:300:0.5)
\centerarc[line width=1](4,4)(170:300:0.7)

\node at (4,-4){$\bullet$};

\centerarc[line width=1](4,-4)(45:200:0.5)
\centerarc[line width=1](4,-4)(45:200:0.7)

\node at (6,1){\tiny$\bullet$};
\node at (6,-2){\tiny$\bullet$};
\node at (0,-5){\tiny$\bullet$};
\node at (0,5){\tiny$\bullet$};
\node at (-3,3){\tiny$\bullet$};
\node at (-3,-3){\tiny$\bullet$};
\node at (-4,0){\tiny$\bullet$};

\drarr(6,1) -- (4,4);
\drarr (6, -2)  --  (6,1) ;
\drarr(4, -4) --  (6, -2);
\drarr (4, -4)  -- (0,-5) ;
\drarr (0, -5)  -- (-3, -3) ;
\drarr (-3,-3)  -- (-4,0) ;
\drarr  (0,5)  -- (4,4);
\drarr  (-3,3) -- (0,5);
\drarr (-4,0) --  (-3,3) ;

\node at (2,0){\huge$\Leftarrow$};

\node at (8, -4.5) {$\Phi_\alpha=\Phi_{i_0}=\Phi_{j_0}$};

\node at (8, 4.2){$\Phi_\omega= \Phi_{i_r}=\Phi_{j_s}$};
\node at (7.7, 1){$\Phi_{i_{r-1}}$};
\node at (7.7, -2){$\Phi_{i_{1}}$};
\node at (7.7, -0.5){$\vdots$};
\node at (-0.2,-6) {$\Phi_{j_{1}}$};
\node at (-0.2,6.4){$\Phi_{j_{s-1}}$};

\centerarc[line width=0.5](6,1)(120:270:0.5)

\centerarc[line width=0.5](6,-2)(90:225:0.5)

 \centerarc[line width=0.5](-3,3)(30:240:0.5)

 \centerarc[line width=0.5](0,5)(-25:210:0.5)

\centerarc[line width=0.5](-4,0)(55:300:0.5)

\centerarc[line width=0.5](-3, -3)(100:325:0.5)

 \centerarc[line width=0.5](0,-5)(135:380:0.5)

\etp
\caption{The intertwiner dg-space and the identifications of the $\Phi$-categories. }\label{fig:intertw-real}
\end{figure}

 \paragraph{The space of intertwiners associated to a polygon.}
 Let now  $\Sen\in\Schob(\CC,A)$ be a schober with singularities in $A$.
 Let $P\subset Q$ be a subpolygon and $f$ be a framing of $P$. Keeping the notation of
 \eqref{eq:del(P,f)}, we define the iterated transport functors
 \[
  M\bigl(\del_+(P,f)\bigr), \,   M\bigl(\del_-(P,f)\bigr): \Phi_\alpha(\Sen)\lra\Phi_\omega(\Sen)
 \]
 associated to the positive and negative boundary paths given by $f$. To avoid confusion, let
 us emphasize the following:
 \begin{itemize}
 \item In the definition of $M\bigl(\del_\pm(P,f)\bigr)$ we follow Convention \ref{conv:clock},
 i.e., we identify the stalks of the local systems $\bPhi$ at the intermediate points by using
 clockwise monodromy from the incoming to the outgoing direction. This means that
 for $M\bigl(\del_+(P,f)\bigr)$ we follow the arc going outside of the polygon $P$
 and for $M\bigl(\del_-(P,f)\bigr)$ we follow the arc going inside $P$, as shown in Fig.
 \ref{fig:intertw-real}.

 \item To identify the sources of the two functors $M\bigl(\del_\pm(P,f)\bigr)$ with each
 other (since a priori the sources are the stalks of the local system $\bPhi_\alpha(\Sen)$
 at the directions corresponding to the two edges of $P$ coming out of $\alpha$),
 we use, similarly to \eqref{eq:3-isotop}, the monodromy along the arc going inside $P$.
 In the similar way (inside arc) we identify the targets of these two functors. This is
 shown by two thick double arcs in Fig.
 \ref{fig:intertw-real}.

 \end{itemize}

\noindent   We define
 the {\em (dg-)space of intertwiners} associated to $(P,f)$ to be
  \be \label{eq:def-space-inter}
 \begin{gathered}
\II(P,f) \,=\, \Hom^\bullet_{\dgFun(\Phi_\alpha, \Phi_\omega)}
\biggl(  M\bigl(\del_+(P,f)\bigr), \, M\bigl(\del_-(P,f)\bigr)\biggr)
[\cont(P,f)].
\end{gathered}
 \ee

 \begin{prop}\label{prop:I(P)-move}
 The spaces $\II(P,f)$ for different framings of the same polygon $P$ are canonically identified.
 \end{prop}

We will therefore use the notation $\II(P)$ for $\II(P,f)$ with any framing $f$.

\vskip .2cm

 \noindent{\sl Proof:}
 Let us pass from a given framing $f=(\alpha,\omega)$ to a new framing $f= (\alpha, \omega')$
 with the same $\alpha$ and $\omega=w_{i_{r-1}}$ moved from $\omega = w_{i_r}=w_{j_s}$
 one step clockwise. This move shortens $\del_+$ and enlarges $\del_-$ by one,
 so $\cont(P,f')=\cont(P,f)+1$. Let
  \[
 F= M_{i_{r-1}, i_r}, \,\,\, G= M(\gamma), \,\,\, H=M\bigl(\del_-(P,f)\bigr),
 \]
 where $\gamma$ is the part of $\del_+(P,f)$ going from $\alpha$ to
 $w_{i_{r-1}}$. Then
 \[
G=M(\del_+(P,f')),\,\,\, F^* = M_{i_{r-1}, i_r}^* \=  M_{i_r, i_{r-1}} T_{\Phi_r}^{-1}[-1],
 \]
 the last identification following from Proposition \ref{prop:adj-trans}.
 Recall that  the monodromy of the local system $\bPhi_{i_r}$ is
 $\TT:= T_{\Phi_{i_r}}[2]$,  see \S \ref{sec:schob-surf}\ref{par:schober-D}
 Using this notation,
 \[
 F^* H \,=\, M_{i_r, r_{r-1}} \TT^{-1} M(\del_-(P,f))[1] \,=\, M(\del_-(P, f')) [1].
 \]
 Therefore  the cyclic rotation isomorphism
\[
 \Hom^\bullet(FG,H)\, \=\,  \Hom^\bullet(G, F^* H)
 \]
  given by (the dg-version of) Corollary \ref{cor:tens-rot}, gives
 an identification $\II(P,f)=\II(P,f')$.

 This shows the invariance of $\II(P,f)$ under
  our elementary move and thus, inductively,
 under any move that changes $\omega$ without changing $\alpha$. In a similar way
 we prove an identification under any move that changes $\alpha$ without changing $\omega$.
 \qed

\paragraph{Pasting of the intertwiners.} We now consider the situation of
 Example \ref{ex:pasting-convex}(b). That is,
let  $P\subset Q$ be a subpolygon and
$\Pc=\{P_\nu\}_{\nu\in N}$
be its polygonal subdivision into several other subpolygons $P_\nu$.  Denote for short
\be\label{eq:I(Pc)}
\II(\Pc) \,=\,  \bigotimes_{\nu\in N} \II(P_\nu).
\ee
In this situation
we have the canonical {\em pasting map}
\be\label{eq:pasting-Pi}
\Pi_\Pc:  \II(\Pc)   \lra \II(P).
\ee
Its definition is as follows.
As in Example
\ref{ex:pasting-convex}(b), we choose a generic direction $\zeta$ and the corresponding
 linear projection $p_\zeta: \CC\to\RR$. This  defines a framing $f_\zeta$ on each $P_\nu$ and
 on $P$, so $\Pc$ is made into a pasting scheme.
 Therefore,  whenever we assign a category  $\Cc_i$  to each vertex $w_i$, a functor $F_{ij}:\Cc_i\to\Cc_j$ to each oriented edge
$[w_i, w_j]$ and a natural
transformation $U_\nu$
to each $P_\nu$, we can speak about the pasted natural transformation associated to $P$ itself.

\vskip .2cm

We apply  this to the  situation where:
\begin{itemize}
\item   $\Cc_i$  is an appropriate stalk of  the local system of categories  $\bPhi_i$;

\item  $F_{ij}$  is  the transport functor $M_{ij}$;

\item $U_\nu\in \II(P_\nu)$, $\nu\in N$, is a given familly of elements of the intertwiner spaces. More precisely,
each $U_\nu$ is viewed, via Proposition \ref{prop:I(P)-move},
as  a natural transformation acting from $M\bigl(\del_+(P_\nu, f_\zeta)\bigr)$ to
$M\bigl(\del_-(P_\nu, f_\zeta)\bigr)$.
\end{itemize}

 \noindent We define $\Pi_\Pc(\bigotimes_\nu U_\nu)$ as the result of the
pasting.
A different choice of the direction $\zeta$ amounts
to applying the tensor rotation (changing the source and target path) of each $U_\nu$ and thus
the result will differ by tensor rotation only, i.e., will be the same, according to our identifications.

\paragraph{Reminder on $\Lie_\oo$-algebras.}
Let $\hen$ be a graded $\k$-vector space.
Recall  that a $\Lie_\oo$-structure on $\hen$ can be defined
  an algebra differential
$D$ in the completed symmetric algebra $\wh S^{\geq 1} (\hen^*[-1])$,
continuous with respect to the natural topology\footnote{This is imposed to make sure
that $D$ is obtained from from an endomorphism  $\Delta$ of  $S^{\geq 1}(\gen[1])$ by dualization.
We can  of course work directly with $\Delta$ but then we have to  use the {\em coalgebra}, not algebra structure on
$S^{\geq 1}(\gen[1])$, similarly to the use of comonads in \S \ref{subsec:unitri-mon}.
}
and satisfying $D^2=0$.
Explicitly, $D$ is uniquely defined, by the Leibniz rule, by its value on the
space of generators $\hen^*[-1]$, i.e., by the data of linear maps
$D_n: \hen^*[-1]\to S^n(\hen^*[-1])$ given for $n\geq 1$. The dual map to $D_n$ can
be considered as  linear map
\[
l_n: \hen^{\otimes n} \lra \hen,   \quad \deg(l_n) = 1-n
\]
 The symmetry conditions on the $D_n$ and the condition $D^2=0$ imply a series of
well known relations on the $l_n$, see \cite{stasheff, getzler, KS-def},  of which we recall the following:
\begin{itemize}
\item[(1)] $l_1: \hen\to \hen$ satisfies $l_1^2=0$ so it makes $\hen$ into a cochain complex.
Further, each $l_n$ is a morphism of complexes with respect to the differentials
given by $l_1$.

\item[(2)] $l_2: \hen^{\otimes 2} \to \hen$ is antisymmetric and makes the cohomology space
$H^\bullet_{l_1}(V)$ into a graded Lie algebra.

\item[(3)] $l_3: \hen^{\otimes 3}\to \hen$
provides an explicit chain homotopy for the Jacobi identity
in (2), and so on.
\end{itemize}

\noindent A $\Lie_\oo$-algebra with $l_n=0$ for $n>2$ is the same as a dg-Lie algebra.

\vskip .2cm

Sometimes it is convenient to  start with a dg-vector space $(\hen,d)$. In this case by
``a $\Lie_\oo$-algebra structure on   $(\hen,d)$'' one means a $\Lie_\oo$-structure on
$\hen$ (considered as just a graded vector space) such that $l_1=d$.

\vskip .2cm

It is also convenient to write the value of $l_n$ on decomposable tensors as an
$n$-ary multilinear operation (bracket)
\[
l_n(v_1\otimes\cdots \otimes v_n) \,=\, [v_1,\cdots, v_n].
\]
By a {\em bracket monomial} in variables $x_1,\cdots,x_m$ one means any operation
obtained by iteration of brackets, for example $[[x_1, x_2], [x_3,x_4,x_5]]$.
Such a monomial defines a map  (called its evaluation) $\hen^m\to \hen$
 for any $\Lie_\oo$-algebra $\hen$.
 We say that $\hen$ is  {\em nilpotent} \cite{getzler}, if there is $m\geq 2$
(the degree of nilpotency) such that any bracket monomial in $m'\geq m$ variables
 has evaluation on $\hen$ identically equal to $0$.

\paragraph{The $\Lie_\oo$-algebra $\gen$.}
We now adapt the constructions of \cite{KaKoSo} \S 10 to our situation of a schober $\Sen\in\Schob(\CC,A)$.
Assume   that $A$ is in exceptionally general position.
For any subset $B\subset A$ with $|B|\geq 3$ we define $P=P_B=\Conv(B)$
and
$\II(B) = \II(P_B)$.

\vskip .2cm

Recall that the secondary polytope $\Sigma(B)$ has dimension $|B|-3$ and its faces $F_\Pc$
are labelled by regular polygonal decompositions $\Pc=(P_\nu, B_\nu)$ of $(P,B)$.
In particular, we have the dg-vector space $\II(\Pc)$ defined in  \eqref{eq:I(Pc)}.
Further, an inclusion of faces $F_{\Pc'}\subset\ol  F_\Pc$ means that
we have a refimenent relation $\Pc' \prec \Pc$ and so the pasting maps \eqref{eq:pasting-Pi}
give rise to a map
 \[
\Pi_{\Pc', \Pc}: \II(\Pc') \lra \II(\Pc).
\]
The ``2-dimensional associativity'' of pasting means that the
maps $\Pi_{\Pc',\Pc}$  are transitive for triple refinement
$\Pc''\prec\Pc'\prec \Pc$ and so define a cellular complex of sheaves $\Nc_{B}$
on $\Sigma(B)$. That is, the stalk of $\Nc_{B}$ on $F_\Pc$ is $\II(\Pc)$ and the generalization map
from $F_{\Pc'}$ to $F_\Pc$ is $\Pi_{\Pc', \Pc}$.

\begin{prop}
The complex $\Nc_{B}$ is factorizing, i.e., for any regular polygonal decomposition
 $\Pc=(P_\nu, B_\nu)$ as above,
the restriction of $\Nc_{B}$ to $F_\Pc = \prod \Sigma(B_\nu)$ is identified with $\boxtimes_\nu  \Nc_{B_\nu}$.
\end{prop}

\noindent{\sl Proof:} Follows directly from the definition of the the stalks and generalization maps. \qed

\vskip .2cm

Let us now define
\be\label{eq:gen-Q}
E_{B} \,=\, \II(B)\otimes \OR(\Sigma(B))[-\dim\Sigma(B)-1], \quad E_{P} := E_{P\cap A}, \,\, P\subset (Q,A).
\ee

\begin{thm}\label{thm:g-l-oo}
The differentials in the cellular cochain complexes of the $\Nc_{B}$ unite to make the dg-vector space
$\bgen = \bigoplus_{(P.B)\subset (Q,A)} E_{B}$ into a $\Lie_\oo$-algebra.
This $\Lie_\oo$-algebra is nilpotent.
The subspace
\[
\gen\, =\, \gen(\Sen)\,=\, \bigoplus_{P\subset (Q,A)} E_{P}
\]
is a $\Lie_\oo$-subalgebra in $\bgen$.
\qed
\end{thm}

\noindent{\sl Proof:} The construction and argument are similar to those of \cite{KaKoSo}.
  Let us specialize  the general discussion of $\Lie_\oo$-structures above to $\hen=\bgen$
and  denote the product in the symmetric algebra by
$\odot$.
 For any   subdivision $\Pc=(P_\nu, B_\nu)$ of a marked subpolytope $(P,B)$
we have the subspace
\[
\begin{gathered}
\II(\Pc)^*\otimes \OR(F_\Pc) [\dim\, F_\Pc]\,=\, \bigotimes_\nu \II(B_\nu)^*
\otimes \OR(\Sigma(B_\nu)) [\dim\, \Sigma(B_\nu)] \, \=\,
\\
\=\,
\bigodot_\nu \bigl( E_{B_\nu}^*[-1]\bigr) \,\subset \, \wh S^\bullet(\bgen{}^*[-1])
\end{gathered}
\]
Recall (Definition \ref{def:coarse}) the concept of a coarse subdivision
of a marked subpolygon $(P,B)$. Such subdivisions $\Pc$ (corresponding to taut webs of
\cite{GMW1}) have the property that the corresponding face $F_\Pc$ of $\Sigma(B)$ has
codimension $1$. In this case we have the identification (coorientation)
$\kappa_\Pc: \OR(F_\Pc) \to \OR(\Sigma(B))$ given by the outer normal to $F_\Pc$.
Therefore we have a degree $1$ map
\[
D_{P,\Pc}: E_B^*[-1] \lra \bigodot_\nu  \bigl( E_{B_\nu}^*[-1]\bigr)
\]
given by the tensor product of
$\Pi_{\Pc, P}^t$, the map dual to the pasting \eqref {eq:pasting-Pi},  and
 of the coorientation $\kappa_\Pc$.

\vskip .2cm

We now define $D$ on each summand $E_B^*[-1]$ to be the sum over all coarse
subdivisions $\Pc$ of $(P,B)$, of the maps $D_{P, \Pc}$. The Leibniz rule extends
this uniquely to an algebra differential of $\wh S^\bullet({\bgen}{}^* [-1])$,
which commutes with the differential induced from that in $\bgen$ by construction
(all pasting maps commute with differentials). The property $D^2=0$ follows from
the similar, obvious,  property  $\delta^2=0$ for the differentials $\delta$ in  the cellular
 cochain complexes of the factorizing complexes of sheaves
$\Nc_B$ on the $\Sigma(B)$.

\vskip .2cm

More precisely,  it is enough to prove that $D^2=0$
on the generators.  Let  $(D^2)_{B}^{B_1,\cdots, B_n}$ be the  matrix element of $D^2$  acting from the summand
$E_B^*[-1]$ of the space of generators to the  summand
$\bigotimes_\nu (E_{B_\nu}^*[-1])$ in $\wh S^\bullet({\bgen}{}^* [-1])$.
In order for $(D^2)_{B}^{B_1,\cdots, B_n}$ to not be $0$ from the outset,
we must have that  $\Pc = \{(P_\nu, B_\nu)\}$ is a regular polygonal subdivision of
$(P,B)$ and the corresponding face $F_\Pc\subset\Sigma(B)$ has codimension $2$.
Now, if this is true, then $(D^2)_{B}^{B_1,\cdots, B_n}$ is identical with the
  matrix element of $\delta^2$ corresponding to the maximal cell of
$\Sigma(B)$ and a codimension $2$ cell  $F_\Pc$. So it is zero and therefore
$D^2=0$ identically.

\vskip .2cm

This  makes $\bgen$ into a $\Lie_\oo$-algebra.
The claims that it is nilpotent and that
$\gen$ is a $\Lie_\oo$-subalgebra
follow in the same way as in \cite{KaKoSo} by analyzing the brackets, see
  the next Remark \ref{rem:gen-brackets}. \qed

\begin{rem}\label{rem:gen-brackets}
It follows from the construction that the $\Lie_\oo$-operations in $\bgen$   are given by ``pasting in coarse decompositions''. That is, the restriction of
$l_n: (\bgen)^{\otimes n}\to\bgen$ to a summand $E_{B_1}\otimes\cdots \otimes E_{B_n}$
is zero unless $\Pc=\{(P_\nu, B_\nu)\}_{\nu=1}^n$, $P_\nu := \Conv(B_\nu)$,
is a coarse regular polygonal decomposition of $(P,B)$ where $P = \bigcup P_\nu$
and $B=\bigcup B_\nu$
(in particular, unless $P$ is convex).  In this is the case, then
$l_n$ is the product of the pasting $\Pi_\Pc$ and the coorientation $\kappa$.
This implies:
\begin{itemize}
\item[(a)]  The maximal degree of a nonvanishing bracket monomial in $\bgen$ is
 bounded by the maximal number of triangles   in a triangulation of $(Q,A)$. Therefore
$\bgen$ is nilpotent.

\item[(b)] If  each $B_\nu= P_\nu\cap A$ represents a summand in $\gen\subset \bgen$,
then the value of $l_n$ on the tensor product of the $E_{B_\nu}$ also
lies in a summand from $\gen$,  as $B=P\cap A$ in this case. Therefore
$\gen\subset\bgen$ is a $\Lie_\oo$-subalgebra.

\end{itemize}
\end{rem}


\subsection{Deforming the infrared monad $\R$} \label{subsec:def-R}

\paragraph{The   deformation complex of $\R$.}
We fix a direction $\zeta$ such that $A$ is in linearly general position including $\zeta$-infinity.
Let $\Sen\in\Schob(\CC,A)$ and
 $\R=\R^\zeta$ be the corresponding infrared monad, see Proposition \ref{prop:infra-mon}.

  We  define the {\em deformation complex}, or {\em ordered Hochschild complex} of the monad $\R$
as
\[
\begin{gathered}
\vC^\bullet =
\vC^\bullet (\R,\R)\,=\, \biggl\{\bigoplus_{i<j} \Hom^\bullet(\R_{ij}, \R_{ij}) \to
 \bigoplus_{i<j<k} \Hom^\bullet(\R_{jk}\R_{ij}, \R_{ik}) \to
 \\
\hskip 8cm \to \bigoplus_{i<j<k<l} \Hom^\bullet (\R_{kl}\R_{jk}\R_{ij}, \R_{il})\to\cdots\biggr\},
\end{gathered}
\]
with the (horizontal) grading starting in degree $0$. Compare with \cite{KaKoSo} \S 12.
Elements of the component
\[
\vC^p \,=\bigoplus_{i_0<\cdots < i_{p+1}} \Hom^\bullet\bigl(\R_{i_p, i_{p+1}} \cdots \R_{i_0, i_1}, \, \R_{i_0, i_{p+1}}\bigr)
\]
can be seen as  $(p+1)$-ary ``operations'' having as inputs, the $p+1$ functors $\R_{i_\nu, i_{n+1}}$
and as output, the functor $\R_{i_0, i_{p+1}}$.
Similarly to the case of the Hochschild complex of an associative (dg-)algebra, $\vC^\bullet$ is a dg-Lie algebra
with respect to the {\em Gerstenhaber bracket}  \cite {gerstenhaber}
\[
[f,g] \,=\, \sum_{\nu=1}^{p+1} f\circ_\nu g \,-(-1)^{\deg(f)\deg(g)} \sum_{\nu=1}^{q+1} g\circ_\nu f, \quad
f\in \vC^p, \, g\in\vC^q.
\]
Here $f\circ_\nu g$ is the ``operadic composition'' cf.  \cite{markl}, obtained by substituting the output of the operation $g$
into the $\nu$th input of the operation $f$ (whenever this makes sense, and defined to be zero otherwise).

\begin{ex}\label{ex:R-beta}
As in the case with associative (dg-)algebras, a Maurer-Cartan element $\alpha\in\vC^1$ gives a unitriangular $A_\oo$-deformation
$\R(\alpha)$ of $\R$. More precisely, $\alpha$, being of total degree $1$,  consists of components
\[
\begin{gathered}
\alpha^{(1)} = (\alpha^{(1)}_{ij}) \in \bigoplus_{i<j} \,\, \Hom^1(\R_{ij}, \R_{ij}), \quad
\alpha^{(2)} = (\alpha^{(2)}_{ijk}) \in \bigoplus_{i<j<k} \Hom^0(\R_{jk}\R_{ij}, \R_{ik}),
\\
\alpha^{(3)} = (\alpha^{(3)}_{ijkl}) \in\bigoplus_{i<j<k<l}
\Hom^{-1}(\R_{kl}\R_{jk}\R_{ij}, \R_{il}), \cdots
\end{gathered}
\]
The deformed monad $\R(\alpha)$ has:
\begin{itemize}
\item[(0)]  The same functors $\R(\alpha)_{ij}= \R_{ij}$. but with differential $d+\alpha^{(1)}_{ij}$.

\item[(1)]  The composition maps
$c_{ijk}(\alpha): \R_{jk}\R_{ij}\to \R_{ik}$
defined by $c_{ijk}(\alpha) = c_{ijk} + \alpha^{(2)}_{ijk}$.

\item[(2)] The $c_{ijk}(\alpha)$ may not be strictly associative but the
$\alpha^{(3)}_{ijkl}$ define the homotopy for the associativity  and so on to give a full $A_\oo$-structure.

\end{itemize}
\end{ex}

\paragraph {Reminder on $\Lie_\oo$-morphisms.} Let $\hen, \ken$ be two $\Lie_\oo$-algebras given by
differentials $D^\hen$ and $D^\ken$ in $\wh S^{\geq 1}(\hen^*[-1])$ and $\wh S^{\geq 1}(\ken^*[-1])$
respectively. Recall that a $\Lie_\oo$-{\em moprhism} $\phi: \hen\to \ken$ can be
seen as a morphism of commutative dg-algebras
\[
\phi^\sharp: \bigl(\wh S^{\geq 1}(\ken^*[-1]),
D^\ken\bigr) \lra \bigl(\wh S^{\geq 1}(\hen^*[-1]), D^\hen\bigr)
\]
continuous with respect to the natural topologies on the source and target.
 As above,  we can view the $\Lie_\oo$-structures on $U$ and $V$ in terms of
the multilinear operations $l_n^\hen: \hen^{\otimes n}\to \hen$ and
$l_n^\ken: \ken^{\otimes n}\to \ken$, $n\geq 1$.
 We also denote both differentials $l_1^\hen$ and $l_1^\ken$ by $d$
and write the higher $l_n$ as brackets.

The morphism $\phi^\sharp$, considered as just a morphism of free commutative algebras,
is determined by its values on the generators which dualize into a sequence of
degree $0$ maps
\be\label{eq:wt-phi}
\wt \phi_n: S^n(\hen[1])[-1] \lra \ken
\ee
or, equivalently,
 maps
\[
\phi_n: \hen^{\otimes n} \lra \ken, \,\,\, n\geq 1, \quad \deg(\phi_n)= n-1
\]
 satisfying the symmetry/antisymmetry conditions derived from the above.
The condition that $\phi^\sharp$ commutes with the differentials, gives a series of
identities on the $\phi_n$, on which the first two are:
\begin{itemize}
\item[(1)] $\phi_1: \hen\to \ken$ commutes with the differentials.

\item[(2)] For $u_1, u_2\in \hen$ we have
\[
\phi_1([u_1, u_2]) -  [\phi_1(u_1), \phi_1(u_2)] \,=\, d(\phi_2(u_1, u_2))
-\bigl( \phi_2(du_1, u_2) + (-1)^{\deg(u_1)}\phi_2(u_1, du_2)\bigr).
\]
In particular, the map $H^\bullet_d(\hen)\to H^\bullet_d(\ken)$ defined by $\phi_1$,
is a morphism of graded Lie algebras.
\end{itemize}

\noindent See \cite{KS-def} \S 3.2.6 or \cite{lada-markl}, Def. 5.2 for the full list of these identities, the second reference treating the case
(which alone will be relevant for us)
when $\ken$ is a dg-Lie algebra, i.e., $l_n^\ken=0$ for $n>2$.

\vskip .2cm

A $\Lie_\oo$-morphism $\phi$ is called {\em of finite order}, if $\phi_n=0$ for $n\gg 0$.

\paragraph{ From  $\gen$ to the deformation complex of $\R$.}

Next, adapting the construction of \cite{GMW1, GMW2} as interpreted in \cite{KaKoSo} \S 10-11, we construct
a $\Lie_\oo$-morphism
$\eta: \gen \to \vC^\bullet$.

\vskip .2cm

In order to do this, for each $n\geq 1$ we need to define a map
$\eta_n: \gen^{\otimes n}\to \vC^\bullet$. Recall that the source and target of this desired $\eta_n$
are direct sums of summands labelled by certain combinatorial data:

\begin{itemize}
\item  We have  $\gen=\bigoplus_P E_P$, see
\eqref{eq:gen-Q}.  Therefore
   \[
  \gen^{\otimes n} = \bigoplus_{P_1,\cdots, P_n} E_{P_1}\otimes\cdots \otimes E_{p_n},
  \]
 the summands  being  labelled by  sequences of  subpolygons $P_1,\cdots, P_n \subset Q$.

\item  We have   $R_{ij} = \bigoplus_{\gamma \in \Lambda (i,j)} R_\gamma$, the summands
being labelled by
$\zeta$-convex polygonal paths $\gamma$ joining $w_i$ and $w_j$, see \eqref{eq:R-ij}.
Therefore
\be\label{eq:vc-bullet}
\vC^\bullet = \bigoplus_{\gamma_1,\cdots, \gamma_m, \gamma}
\Hom(\R_{\gamma_1}\circ\cdots \R_{\gamma_n}, \R_\gamma),
\ee
the summands being labelled by sequences of $\zeta$-convex paths
$\gamma_1,\cdots, \gamma_m,\gamma$ such that
 $\gamma_1,\cdots, \gamma_m$  form a  composable sequence
 and $\gamma$  has the same beginning and end as the composition
$\gamma_1\circ\cdots
\circ\gamma_n$.
\end{itemize}

\noindent
Therefore $\eta_n$ should be given by  matrix elements
\be\label{eq:eta-matrix}
\eta_{P_1,\cdots, P_n}^{\gamma_1,\cdots, \gamma_m|\gamma}: \,\,
\bigotimes_{\nu=1}^n E_{P_\nu} \lra \Hom(\R_{\gamma_1}\circ\cdots \R_{\gamma_n}, \R_\gamma),
\ee
where $P_1,\cdots P_n$ and $\gamma_1,\cdots, \gamma_m,\gamma$ are as above.
As in \cite{KaKoSo}, we will define them by analyzing coarse subdivisions of
unbounded polygons.

\vskip .2cm

More precisely, let us choose a Vladivostok point $\bv= r\zeta, \, r\gg 0$ far away in the direction
$\zeta$ and let $\ol A=A\cup\{\bv\}$ and $\ol Q=\Conv(\ol A)$. Assume that $\ol A$ is in exceptionally general position.
A subpolygon $G\subset \ol Q$ will be called {\em bounded}, if $\bv\notin G$ and {\em unbounded},
if $\bv\in G$.

\vskip .2cm

A bounded subpolygon is the same as a subpolygon in $A$.
An unbounded subpolygon $G$ can be  thought of as approximately $\zeta$-convex,
in particular, its  {\em unbounded edges} (those containing $\bv$) can be depicted
as half-lines going in direction $\zeta$, see Fig. \ref{fig:1-fin}. Further, such $G$
is uniquely determined by a $\zeta$-convex path $\gamma$. Explicitly,
$\gamma$ is the bounded part of the boundary $\del G$.  Alternatively,
$\gamma = \del_-^\zeta(G)$ is the ``bottom part'' of $\del G$ if   $\zeta\oo$ is positioned
vertically.  In what follows let us write $\del_\pm$ for $\del_\pm^\zeta$.

\vskip .2cm

 In the other direction,  $G=\Conv(\gamma\cup\{\bv\})$, or $G=\Conv^\zeta(\gamma)$ in the
 parallel representation of Fig.   \ref{fig:1-fin}.
 Therefore instead of  convex paths $\gamma_1,\cdots, \gamma_m, \gamma$
 in \eqref{eq:vc-bullet}
 we can use the  corresponding unbounded polygons $G_1,\cdots, G_m, G$.

 \vskip .2cm

 Let now $G$ be an unbounded polygon
 and
 $\ol\Pc$ be a polygonal decomposition of $G$.
 It consists of several bounded polygons which we denote $P_1,\cdots, P_n$
 (the precise numeration is not important)
 and several unbounded polygons which we denote $G_1,\cdots, G_m$
 and number consecutively clockwise around the boundary of $G$, see
 Fig. \ref{fig:1-fin}. Let
 $\gamma_j = \del_-G_j$ be the $\zeta$-convex path corresponding to $G_j$,
 oriented counterclockwise around the boundary of $G_j$. Then $\gamma_1,\cdots, \gamma_m$
 form a composable sequence.
 Let $\gamma = \del_- G$. Then $\gamma$ has the same beginning and end as
 $\gamma_1\circ \cdots \circ  \gamma_m$.
 Note that $\gamma$ may have nontrivial overlap with $\gamma_1$ along
 the  ``left handle'' $\lambda$ and with $\gamma_m$ along the ``right handle'' $\rho$, see
  Fig. \ref{fig:1-fin}.

  \vskip .2cm

   Further, our choice of $\zeta$ gives a framing $f_\zeta$ of each $P_\nu$, $\nu=1,\cdots, n$.
 In particular we have the positive and negative boundaries $\del_\pm P_\nu$
   and  the pasting map
\be\label{eq:pasting-G}
\Pi_{\ol\Pc}: \bigotimes_{\nu=1}^n \Hom^\bullet \bigl(M(\del_+P_\nu), M(\del_-P_\nu)\bigr)
\lra \Hom^\bullet (M_{\gamma_1}\circ\cdots\circ M_{\gamma_m}, M_\gamma).
\ee

\vskip .2cm

Now recal \eqref{eq:R-ij} that
$
\R_\gamma\,=\,
 M(\gamma)\otimes \OR(\gamma).
 $
 Further, we have a canonical identification
 \be\label{eq:or-gamma-sigma}
 \OR(\gamma)=\bigwedge\nolimits^\max \k^{G\cap A} \lra \OR(\Sigma(G\cap \ol A))
 \ee
 coming from two steps:
 \begin{itemize}
 \item[(1)] The point $\bv\in\ol A$ is distinguished, so wedge product with it on the left
 gives an identitification $\bigwedge^\max (\k^{G\cap A})\to \bigwedge^\max(\k^{G\cap \ol A})$.

 \item[(2)] By construction \cite{GKZ}, the polytope $\Sigma(G\cap \ol A)$ is
 embedded into $\k^{G\cap \ol A}$ and its affine span there is a vector subspace of codimension $3$
 the quotient by which is $\Aff(\ol A)$, the space of restriction to $\ol A$ of affine $\RR$-linear
 functions $\CC\to\RR$. This quotient subspace has a canonical orientation
 coming from the standard orientation of $\CC$. So the orientation space  of
 $\Sigma(G\cap \ol A)$, being the same as the orientation space of its affine span,
 is identified with the orientation space of $\k^{G\cap \ol A}$.
  \end{itemize}

\begin{Defi}\label{def:associated}
Let us call a sequence of bounded subpolygons $P_1,\cdots, P_n$ and
a sequence of $\zeta$-convex paths $\gamma_1,\cdots, \gamma_m,\gamma$
{\em associated}, if there exists a  (necessarily unique)
 coarse subdivision $\ol\Pc$ of an unbounded subpolygon $G$
 into bounded subpolygons $P_\nu$ (in some order) and unbounded polyygons
 $G_1,\cdots, G_m$ (in the counterclockwise order as above) so that $\gamma_j=\del_-G_j$
 and $\gamma = \del_-G$.
 \end{Defi}

 Taking into account the identification \eqref{eq:or-gamma-sigma}
 we see, just as in the definition of $\gen$, that
  in the situation of Definition \ref{def:associated} we have a canonical coorientation
 \[
 \kappa_{\ol \Pc}: \bigotimes_{\nu=1}^n \OR(\Sigma(P_\nu\cap\ol A)) \,\otimes
 \,
 \bigotimes_{j=1}^m  \OR(\gamma_j) \,\lra \, \OR(\gamma).
 \]
 It is just the coorientation of the codimension $1$ face $F_{\ol\Pc}$ of the polytope
 $\Sigma(G\cap\ol A)$
 given by the outer normal to this face.
 Note also that
  \be\label{eq:EP-explicitly}
E_{P_\nu} = \Hom^\bullet \bigl(M(\del_+P_\nu), M(\del_-P_\nu) \bigr)\otimes \OR(\Sigma(A\cap P_\nu))[2-|A\cap \del_+P|].
\ee

\bef[H]
\centering
\btp[scale=0.3]
\node (1) at (0,0){};
\node (2) at (4,6){};
\node (3) at (6,2){};
\node (4) at (6,-1){};
\node (5) at (4,-6){};

\fill(1) circle (0.15);
\fill(2) circle (0.15);
\fill(3) circle (0.15);
\fill(4) circle (0.15);
\fill(5) circle (0.15);

\node (7) at (7,-7){};
\fill (7) circle (0.15);

\node (8) at (6,8){};
\fill (8) circle (0.15);

\node (9) at (9, 9){};
\fill (9) circle (0.15);

\draw (15,9) -- (9,9);
\draw[line width=0.8mm, color=gray] (9,9) -- (6,8) -- (4,6);

\draw (0,0) -- (4,6) -- (6,2) -- (6,-1) -- (4,-6) -- (0,0); 
\draw  (15,-1) -- (6,-1);
\draw (6,2) -- (15,2);
\draw (0,0) -- (6,2);

\draw (7,-7) -- (15,-7);
\draw[line width=0.8mm, color=gray] (7,-7) -- (4,-6);
\draw[line width=0.8mm] (4,-6) -- (6,-1) -- (6,2) -- (4,6);

\node at (3.5,3){$P_1$};
\node at (3.5, 1){$\vdots$};
\node at (3.5,-1.5){$P_n$};
\node at (6.3,-4.5){$\gamma_m$};
\node at (7, 1 ){$  \gamma_2 $};
\node at (7.3, 0){$\cdots$};
\node at (6,5.5){$\gamma_1$};
\node at (10,-4){$G_m$};
\node at (10, 1){$ G_2 $};
 \node at (10,5){$G_1$};

\node[text width=3cm, left] (LH) at (-4,-9){Left handle $\lambda
$
};

\node[text width=3cm, left] (RH) at (-4,12){Right handle
 $\rho
 $
 };
\node (rh) at (6,8.5){};
\draw[->, line width=0.8mm, color=gray] (RH) -- (rh);

\node (lh) at (5.5,-7){};
\draw[->, line width=0.8mm, color=gray] (LH) -- (lh);

\node at (14,10){$G$};

\draw  [decoration={markings,mark=at position 1.0 with
{\arrow[scale=1.5,>=stealth]{>}}},postaction={decorate},
line width=.2mm]
(12,-10) -- (15,-10);

\node at (13.5,-11){\large$\zeta$};

\etp
\caption{A  subdivision of an unbounded polygon $G$.}
\label{fig:1-fin}
\enf

We now define the matrix element
$\eta_{P_1,\cdots, P_n}^{\gamma_1,\cdots, \gamma_m|\gamma}$
in \eqref {eq:eta-matrix} to be $0$ unless  $P_1,\cdots, P_n$ and $\gamma_1,\cdots,\gamma_m,\gamma$ are associated.  If they are associated via a coarse subdivision
 $\ol \Pc$,
  we define the matrix element
  to be given by  the pasting map \eqref{eq:pasting-G}, tensored with the coorientation
  $\kappa_{\ol\Pc}$. Comparison of the shifts in \eqref{eq:EP-explicitly} and with the
  horizontal grading in $\vC^\bullet$ shows that in this way we get a degree $0$ map of graded
  vector spaces.

\begin{thm}\label{thm:Lie-oo-map}
The matrix elements $\eta_{P_1,\cdots, P_n}^{\gamma_1,\cdots, \gamma_m|\gamma}$ thus defined
form a finite order $\Lie_\oo$-morphism $\eta: \gen\to \vC^\bullet$.
\end{thm}

\paragraph{Proof of Theorem \ref{thm:Lie-oo-map}.}  We start with a general discussion.
Let $\hen$ be a $\Lie_\oo$-algebra, $U$ be an associative dg-algebra and
\[
C^\bullet = C^\bullet(U,U) \,=\, \bigl\{ \Hom^\bullet_\k(U,U) \to \Hom^\bullet_\k(U^{\otimes 2}, U) \to \cdots\bigr\}
\]
be the standard Hochschild complex of $U$, with the horizontal grading starting in degree $0$.
Similarly to the above (and more classically), $C^\bullet$ is a dg-Lie algebra,
so we can speak about $\Lie_\oo$-morphisms $\phi: \hen\to C^\bullet$.   Let us recall a description of such
morphisms in terms closer to those used earlier.
\vskip .2cm

The dg-algebra structure on $U$ is encoded by a differential $D^U$ in the completed tensor
algebra $\wh T^{\geq 1}(U^*[-1])$, which has the form $D^U=D^U_1+D^U_2$,
where  $D_1^U$ is  the differential induced by that in $U$
while $D^U_2$ is the graded derivation given on the space of generators by the map
$U^*\to U^*\otimes U^*$ dual to the multiplication in $U$.

As before, let us encode the $\Lie_\oo$-structure on $\hen$ by a differential $D^\hen$ in
$\wh S^{\geq 1}(\hen^*[-1])$.

\begin{prop}\label{prop:l-oo-mor-D}
In the above notation, we have a bijection between:
\begin{itemize}
\item[(i)]  $\Lie_\oo$-morphisms $\phi: \hen\to C^\bullet(U,U)$.

\item[(ii)]   Algebra differentials  $D$ in
the completed tensor product
$\wh S^{\geq 1}(\hen^*[-1])\wh\otimes \, \wh T^{\geq 1}(U^*[-1])$
which are continuous in the natural topology, satisfy $D^2=0$ and
whose values on the generators
have the form:
\[
\begin{gathered}
D|_{\hen^*[-1]\otimes 1} = D^\hen|_{\hen^*[-1]}\otimes 1,
\quad D|_{1\otimes U^*[-1]} =  1\otimes D^U|_{ U^*[-1]} + \sum_{m,n\geq 1} D^\phi_{n,m},
\\
D^\phi_{n,m}: 1\otimes U^*[-1] \lra S^n(\hen^*[-1]) \otimes T^m(U^*[-1]).
\end{gathered}
\]
\end{itemize}
  \end{prop}

  \noindent{\sl Proof:} Given data in (ii), we dualize $D^\phi_{n,m}$ and add a shift by $-1$
   to get a  degree $0$ map
  \[
 \wt \phi_{n,m}: S^n(\hen[1])[-1] \lra \Hom^\bullet \bigl(T^m(U[1]), U([1]) \bigl)\,=\, C^m(U,U) [-m],
  \]
  the target being the $m$th term of the Hochschild complex. These maps unite into a collection
  of  degree $0$ maps $\wt \phi_n: S^n(\hen[1])[-1] \to C^\bullet(U,U)$. The condition that the $\wt \phi_n$
  form a $\Lie_\oo$-morphism, see \eqref{eq:wt-phi},  is then expressed by $D^2=0$. \qed

  \vskip .2cm

  We now adapt this approach to our situation which, on the face of it, has two differences:
  \begin{itemize}
  \item[(1)] Instead of a dg-algebra $U$ we have a dg-monad, i.e., a dg-endofunctor $R: \Vc\to \Vc$ with
  a  (closed, degree $0$)    associative natural transformation $R R\to R$. Here $\Vc = \bigoplus_{i=1}^N\Phi_i$.

  \item[(2)] We consider not the full but  the ordered Hochschild complex
  with respect to the ordering on the summands $\Phi_i$ and the   uni-triangular form $R=(R_{ij})_{i<j}$ of $R$.
  \end{itemize}

 \noindent To account for (1), we associate to the dg-monad  $R$ an associative dg-agebra $U(R)$ as follows.
 Let $F: \Cc \to \Dc$ be a dg-functor of dg-categories. Consider the dg-vector space
 \[
 U(F) \,=\, \bigoplus_{c\in \Cc, d\in \Dc} \Hom^\bullet_\Dc(d, F(c) ).
 \]
 Given two composable dg-functors $\Cc\buildrel F\over\to\Dc\buildrel G\over\to\Ec$, we have the morphism
 of dg-vector spaces $U(G)\otimes U(F) \lra U(GF)$ which consists of the composition maps
 \[
 \Hom^\bullet_\Ec(e, G(d)) \otimes\Hom^\bullet_\Dc(d, F(c)) \lra \Hom^\bullet_\Ec(GF(c),e),\quad
 g\otimes f \,\mapsto \,  G(f)\circ g,
 \]
 and is zero one any other tensor products of summands. Further, a closed, degree $0$
 natural transformation $t: F_1\to F_2$ between dg-functors $F_1, F_2: \Cc \to \Dc$ gives rise
 to a morphism of dg-vector spaces $U(t): U(F_1)\to U(F_2)$.

  \vskip .2cm

 Therefore the (infrared) monad structure $\R\R\to \R$ on
 the endofunctor $\R:\Vc\to\Vc$ gives an associative dg-algebra structure $U(\R)\otimes U(\R)\to U(\R)$.
 To account for the difference (2) above, note that the ordered Hochschild complex $\vC^\bullet$ of $\R$
 is in fact a $\Lie_\oo$-subalgebra in $C^\bullet (U(\R), U(\R))$.
 Therefore
 our morphisms  of dg-vector spaces $\eta_n: \gen^{\otimes n} \to \vC^\bullet$ give
morphisms of dg-vector spaces $U(\eta_n): \gen^{\otimes n} \to C^\bullet(U(\R), U(\R))$ and
it is enough to show that the $U(\eta_n)$ form a $\Lie_\oo$-morphism.

\vskip .2cm

For this, we can use Proposition \ref{prop:l-oo-mor-D} and argue similarly to the proof of Theorem
\ref{thm:g-l-oo}. That is, we construct a differential $D$ in $S^{\geq 1}(\gen^*[-1])\otimes T^{\geq 1}(U(\R)^*[-1])$
(no completions needed in our particular case) incorporating
the $U(\eta_n)$,  by looking at the faces of the secondary polytopes.

\vskip .2cm

More precisely, $\eta_n$ is defined in terms of  coarse (corresponding to codimension $1$ faces of
the secondary polytope) regular decompositions $\ol\Pc$
 of  all possible unbounded polygons
$G$, consisting   of $n$ bounded polygons and some number $m$ of unbounded ones, see
Fig. \ref{fig:1-fin}.  The bounded polygons do not have any natural order, but the unbounded ones
are ordered counterclockwise as in Fig. \ref{fig:1-fin}. Therefore the
formula for $\eta_n$  defines, for each $m\geq 1$,  a morphism
\[
D^{U(\eta)}_{n,m}: U^*[-1]  \lra S^n(\gen^*[-1]) \otimes T^m(U^*[-1]),
\]
 (we use the order on the unbounded polygons to lift the value into the tensor algebra).
 This, together with other components, defines a degree $1$ graded derivation $D$ of $S^{\geq 1}(\gen^*[-1])\otimes
 T^{\geq 1}(U^*[-1])$, and all we need to show is $D^2=0$.

 \vskip .2cm

 This follows just like in the proof of Theorem \ref{thm:g-l-oo},  from comparison of matrix elements of $D^2$
 with those of the
   $\delta^2$
 where $\delta$ runs over the differentials in the cellular chain complexes of  appropriate cellular sheaves
 on the secondary polytopes $\Sigma(G\cap \ol A)$. \qed


\subsection { Maurer-Cartan elements in the infrared $\Lie_\oo$-algebra}\label{subsec:MC}

\paragraph{Maurer-Cartan elements in  $\Lie_\oo$-algebras.} Let $\hen$ be a nilpotent
$\Lie_\oo$-algebra. An element $\beta\in \hen^1$ (of degree $1$) is called a
{\em Maurer-Cartan element}, if
\[
\sum_{n=1}^\infty{1\over n!}  l_n(\beta \otimes\cdots\otimes \beta) \,=\, d\beta + {1\over 2!} [\beta,\beta]
+ {1\over 3!} [\beta,\beta,\beta] + \cdots \,=\,0.
\]
This extends the more familiar concept for dg-Lie algebras, discussed in \S
\ref{subsec:infra-FS-mon}\ref{par:tot-compl}
We denote $\MC(\hen)$ the set of Maurer-Cartan elements in $\hen$.
We also refer to \cite{getzler} for the concept of {\em gauge equivalence} of Maurer-Cartan elements.

\vskip .2cm

Let $\phi: \hen\to \ken$ be a finite order $\Lie_\oo$-morphism of nilpotent $\Lie_\oo$-algebras,
with components $\phi_n: \hen^{\otimes n}\to \ken$, $n\geq 1$. Then for
a Maurer-Cartan element $\beta \in \MC(\hen)$ the element
\[
\phi_*(\beta) \,=\, \sum_{n=1}^\infty {1\over n!} \, \phi_n(\beta \otimes\cdots\otimes \beta)
\]
is a Maurer-Cartan element in $\ken$, cf.  \cite{KS-def} \S 3.2.6.
  This  gives  a map $\phi_*: \MC(\hen)\to\MC(\ken)$.

  \paragraph{Maurer-Cartan elements in $\gen$ explicitly.}
  We specialize the above to the $\Lie_\oo$-morphism $\eta:\gen\to \vC^\bullet(\R,\R)$.
  Example \ref{ex:R-beta} implies that any Maurer-Cartan element $\beta$ of $\gen$
  gives an $A_\oo$-deformation $\R(\beta) :=\R(\eta_*(\beta))$ of
  the infrared monad $\R$.
   For this reason let us analyze the data encoded in a datum of $\beta\in\MC(\gen)$
  more explicitly.

  \vskip .2cm

  We have $\gen=\bigoplus_P E_P$, where $P$ runs over subpolygons $P\subset Q$
  (with vertices in $A$). Therefore $\beta = (\beta_P)$ is a system of
  elements $\beta_P\in E_P$, one for each subpolygon $P$.
  Let us write $\Sigma(P)=\Sigma(A\cap P)$ for the secondary polytope of $P$.
  Let us also choose a generic direction $\zeta$ which gives a framing $f_\zeta$ on each $P$,
  in particular, it gives the positive and negative parts $\del_\pm P$ of the boundary of $P$.

  \vskip .2cm

  Recalling \eqref{eq:gen-Q} the definition of $E_P$, taking into account the shifts  as in
  \eqref{eq:EP-explicitly} and
  setting aside the orientation factors (amounting to $\pm 1$), we can say that $\beta_P$
  is a natural transformation
  \[
  \beta_P \,\in\, \Hom^\bullet\bigl(M(\del_+P), M(\del_-P)\bigr), \quad
  \deg (\beta_P) \,=\, 3-|A\cap \del_+P|.
  \]
  So the range of possible degrees of the components  $\beta_P$  is $1,0,-1,\cdots$.
  It is convenient to analyze the Maurer-Cartan conditions
  on them recursively, starting from some ``basic'' components of  degrees $1$ and $0$ and
  viewing the rest as a  ``system of higher coherence data'' on these basic components.
 A subpolygon $P\subset (Q,A)$ will be called {\em empty},
  if it contains no elements of $A$ other than the vertices.

  \vskip .2cm

\begin{figure}[H]
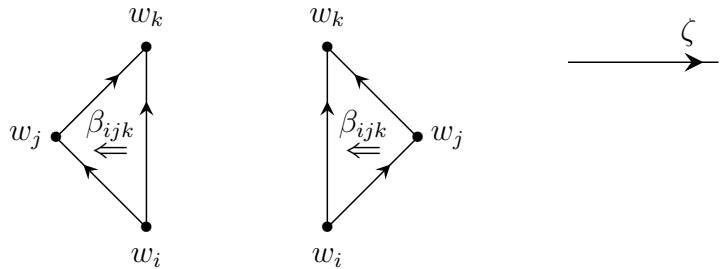

  \centering
  \btp[scale=0.4]

  \node at (0,0){\small$\bullet$};
   \node at (3,3){\small$\bullet$};
    \node at (3,-3){\small$\bullet$};

  \drarr (3,-3) -- (0,0);
 \drarr (3,-3) -- (3,3);
 \drarr (0,0) -- (3,3) ;

 \node at (1.8, -0.5){\large$\Leftarrow$};
 \node at (1.8, +0.4){$\beta_{ijk}$};

 \node at (3, -4){$w_i$};
  \node at (3, 4){$w_k$};
   \node at (-1,0){$w_j$};
   \etp
   \quad\quad\quad\quad
    \btp[scale=0.4]

  \node at (0,0){\small$\bullet$};
   \node at (-3,3){\small$\bullet$};
    \node at (-3,-3){\small$\bullet$};

  \drarr (-3,-3) -- (0,0);
 \drarr (-3,-3) -- (-3,3);
 \drarr (0,0) -- (-3,3) ;

 \node at (-1.8, -0.5){\large$\Leftarrow$};
 \node at (-1.8, +0.4){$\beta_{ijk}$};

 \node at (-3, -4){$w_i$};
  \node at (-3, 4){$w_k$};
   \node at (1,0){$w_j$};

   \draw
[decoration={markings,mark=at position 0.9 with
{\arrow[scale=2,>=stealth]{>}}},postaction={decorate},
line width=.2mm] (5,2.5) -- (10,2.5);
\node at (9,3.5){$\zeta$};
   \etp
   \caption{Empty triangles give comultiplications and multiplications.}\label{fig:empty-tri}
   \end{figure}

  \noindent{\bf Basic level: multiplication and comultiplication.}
  The ``basic level'' of $\beta$ consists of the $\beta_P = \beta_{ijk}$ where $P=\Conv\{w_i, w_j, w_k\}$
  is an empty triangle. We number them so that $w_i <_\zeta w_j <_\zeta w_k$, see Fig.
  \ref {fig:empty-tri}.
  There are two possibilities:
  \begin{itemize}
  \item[(1)] $\del_+P = [w_i,w_k]$ consists of a single edge, see the left of Fig. \ref {fig:empty-tri}.
  In this case $\beta_{ijk}: M_{ik}\to M_{jk} M_{ij}$ is a natural transformation of degree $1$,
  which we view as a comultiplication. The Maurer-Cartan condition on $\beta_{ijk}$ amounts to
  $d(\beta_{ijk})=0$.

  \item[(0)] $\del_+P = [w_i,w_j]\cup w_j,w_k]$ consists of two edges, see the right of Fig. \ref {fig:empty-tri}.
  In this case $\beta_{ijk}: M_{jk} M_{ij}\to M_{ik}$ is a natural transformation of degree $0$,
  which we view as a  multiplication. The Maurer-Cartan condition on $\beta_{ijk}$ amounts to
  $d(\beta_{ijk})=0$.

  \end{itemize}

     \vskip .2cm

  \noindent {\bf Extremal empty polygons:   $A_\oo$-structures.} We call an empty $n$-gon $P$
  {\em extremal of the first kind},  if  $\del_+P$   consists of a single edge.   $\beta_P$ has degree $1$.
  We say that $P$ is  {\em extremal of the second kind}, if  $\del_-P$ consists of a single edge.
 In this case $\beta_P$ has degree $3-n$.

  \vskip .2cm

  The Maurer-Cartan conditions of the   $\beta_P$ corresponding to extremal  empty $P$ of the first kind
  amount to saying that these $\beta_P$ form an $A_\oo$ (homotopy co-associative) structure for the
  comultiplications $\beta_{ijk}$ corresponding to empty triangles of the first kind.

  \begin{figure}[H]
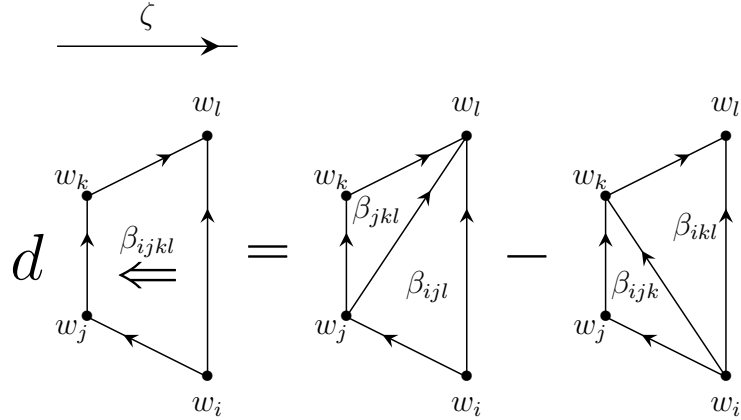

  \centering

     \btp[scale=0.4]

     \node at (2,-4){\small$\bullet$};
   \node at (2,4){\small$\bullet$};
    \node at (-2,2){\small$\bullet$};
     \node at (-2,-2){\small$\bullet$};

 \drarr (2,-4) -- (2,4);
  \drarr (2,-4) -- (-2,-2);
  \drarr (-2,-2) -- (-2,2);
 \drarr (-2,2) -- (2,4);

\node at (2,-5){$w_i$};
\node at (2,5){$w_l$};
 \node at (-2.5,-2.5){$w_j$};
\node at (-2.5,2.55){$w_k$};

\node at (0, -0.7){\huge$\Leftarrow$};
\node at (0,0.5){$\beta_{ijkl}$};

\node at (-4,0){\huge $d$};
\node at (4,0){\huge$=$};

\draw
[decoration={markings,mark=at position 0.9 with
{\arrow[scale=2,>=stealth]{>}}},postaction={decorate},
line width=.2mm] (-3,7) -- (3,7);
\node at (0,8){$\zeta$};
     \etp
  \btp[scale=0.4]

     \node at (2,-4){\small$\bullet$};
   \node at (2,4){\small$\bullet$};
    \node at (-2,2){\small$\bullet$};
     \node at (-2,-2){\small$\bullet$};

 \drarr (2,-4) -- (2,4);
  \drarr (2,-4) -- (-2,-2);
  \drarr (-2,-2) -- (-2,2);
 \drarr (-2,2) -- (2,4);

\node at (2,-5){$w_i$};
\node at (2,5){$w_l$};
 \node at (-2.5,-2.5){$w_j$};
\node at (-2.5,2.55){$w_k$};

\node at (4,0){\huge$-$};
\drarr (-2,-2) -- (2,4);

\node at (0.7,-1){$\beta_{ijl}$};
\node at (-1,1.5){$\beta_{jkl}$};

\etp
  \btp[scale=0.4]

     \node at (2,-4){\small$\bullet$};
   \node at (2,4){\small$\bullet$};
    \node at (-2,2){\small$\bullet$};
     \node at (-2,-2){\small$\bullet$};

 \drarr (2,-4) -- (2,4);
  \drarr (2,-4) -- (-2,-2);
  \drarr (-2,-2) -- (-2,2);
 \drarr (-2,2) -- (2,4);

\node at (2,-5){$w_i$};
\node at (2,5){$w_l$};
 \node at (-2.5,-2.5){$w_j$};
\node at (-2.5,2.55){$w_k$};

 \drarr (2,-4) -- (-2,2);

 \node at (-1,-1){$\beta_{ijk}$};
 \node at (1,1){$\beta_{ikl}$};

\etp

  \caption { An extremal empty $4$-gon of the first kind gives a homotopy for coassociativity. } \label{fig:comult}
  \end{figure}

  \vskip .2cm

  For example, let
  $P$ be  an empty $4$-gon of the first kind. Labeling the vertices $w_i, w_j, w_k, w_l$
  as in Fig. \ref{fig:comult}, we get a degree $0$ natural transformation
  $\beta_P=\beta_{ijkl}: M_{il}\to M_{kl} M_{jk} M_{ij}$. The Maurer-Cartan condition on $\beta_{ijkl}$
  is that it ensures homotopy coassociativity of the elementary (triangular) comultiplications:
  $d(\beta_{ijkl})$ is the difference between two possible bracketing of four factors which correspond
  to two possible trianguilations of the $4$-gon $P$, see Fig. \ref{fig:comult}.

  \vskip .2cm

\begin{figure}[H]
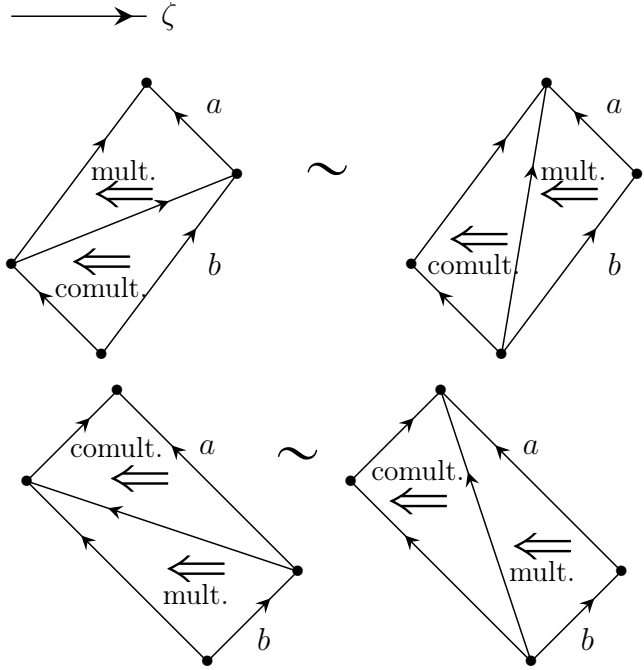

  \centering
  \btp[scale=0.3]

  \node at (0,0){\small$\bullet$};
    \node at (4,-4){\small$\bullet$};
  \node at (10,4){\small$\bullet$};
   \node at (6,8){\small$\bullet$};

  \drarr (4,-4)--(0,0);
   \drarr (4,-4)--(10,4);
   \drarr (0,0)--(6,8);
    \drarr (10,4)--(6,8);
    \drarr (0,0)--(10,4);

   \node at (4,0){\huge$\Leftarrow$};
  \node at (4,-1){\small{comult.}};

  \node at (5,3){\huge$\Leftarrow$};
  \node at (5,4.2){\small{mult.}};

  \draw
[decoration={markings,mark=at position 0.9 with
{\arrow[scale=2,>=stealth]{>}}},postaction={decorate},
line width=.2mm] (0,11) -- (6,11);
\node at (7,11){$\zeta$};
\node at (14,4){\huge$\sim$};
\node at (9 ,7){$a$};

\node at (9,0){$b$};
  \etp \quad
   \btp[scale=0.3]

  \node at (0,0){\small$\bullet$};
    \node at (4,-4){\small$\bullet$};
  \node at (10,4){\small$\bullet$};
   \node at (6,8){\small$\bullet$};

  \drarr (4,-4)--(0,0);
   \drarr (4,-4)--(10,4);
   \drarr (0,0)--(6,8);
    \drarr (10,4)--(6,8);
    \drarr (4,-4)--(6,8);

   \node at (3,1){\huge$\Leftarrow$};
  \node at (2.8,0){\small{comult.}};

  \node at (7,3){\huge$\Leftarrow$};
  \node at (7.2,4.2){\small{mult.}};

\node at (9 ,7){$a$};

\node at (9,0){$b$};

  \etp
  \quad\quad
  \btp[scale=0.3]

  \node at (0,0){\small$\bullet$};
   \node at (-4,-4){\small$\bullet$};
 \node at (-8,8){\small$\bullet$};
  \node at (-12,4){\small$\bullet$};

  \drarr (-4,-4) -- (0,0);
\drarr (-4,-4) -- (-12,4);
\drarr (-12,4) -- (-8,8);
\drarr (0,0) -- (-8,8);
\drarr (0,0) -- (-12,4);

\node at (-4.5, 0){\huge$\Leftarrow$};
\node at (-4.5, -1){\small{mult.}};

\node at (-7, 4){\huge$\Leftarrow$};
\node at (-8,5.5){\small{comult.}};

\node at (0,5){\huge$\sim$};

\node at (-4,5.5){$a$};
\node at (-1.5,-3){$b$};

  \etp
   \btp[scale=0.3]

  \node at (0,0){\small$\bullet$};
   \node at (-4,-4){\small$\bullet$};
 \node at (-8,8){\small$\bullet$};
  \node at (-12,4){\small$\bullet$};

  \drarr (-4,-4) -- (0,0);
\drarr (-4,-4) -- (-12,4);
\drarr (-12,4) -- (-8,8);
\drarr (0,0) -- (-8,8);
\drarr (-4,-4) -- (-8,8);

\node at (-3.5, 1){\huge$\Leftarrow$};
\node at (-3.5, -0){\small{mult.}};

\node at (-9, 3){\huge$\Leftarrow$};
\node at (-9,4.5){\small{comult.}};

 \node at (-4,5.5){$a$};
\node at (-1.5,-3){$b$};

  \etp

  \caption{Homotopy relations between multiplication and comultiplication. } \label{fig:mu-Delta}
  \end{figure}

  \noindent{\bf  The Frobenius condition on the muiltiplication and comiltiplication.}
  Consider now empty $4$-gons $P$ with $\del_+P$ consisting of two segments, as in Fig.
  \ref{fig:mu-Delta}. The corresponding $\beta_P$ have degree $0$ and can be seen as
  homotopies establishing certain relations between the multiplication and comultiplication
  arrows obtained by triangulating $P$ into two empty triangles.
  There can be two different types of position depicted on the top and bottom of Fig.   \ref{fig:mu-Delta}.
   The relations we obtain are  analogous to the so-called
 {\em Frobenius conditions} in the theory of Frobenius algebras \cite{kock}.
 More precisely, suppose $H$ is a vector space with a multiplication $H\otimes H\to H$, written $a\otimes b\mapsto ab$
 and a comultiplication $\Delta: H\to H\otimes H$, written in Sweedler's notation as
 $\Delta(a) = a^{(1)}\otimes a^{(2)}$, the latter being a shorthand for $\sum_i a^{(1)}_i\otimes a^{(2)}_i$.
 Then the top and bottom part of Fig. \ref{fig:mu-Delta}, read right to left, can be written as
 \[
 \Delta(ab) \sim ab^{(1)}\otimes b^{(2)}, \quad \Delta(ab) \sim a^{(1)}\otimes a^{(2)} b
 \]
 (with $\sim$ meaning being homotopic).
 Here the elements $a$ and $b$ are symbolically depicted on the $\del_+$-arrows in Fig. \ref{fig:mu-Delta}. Compare with  \cite{kock}, Lemma 2.3.19.

\vskip .2cm

\noindent{\bf A triangle with one point inside:  $3\to 1$ move.} Let now $P = \Conv\{w_i, w_j, w_k\}$
 be a triangle with
one element $w_l$ of $A$ inside.   The Maurer-Cartan condition on $\beta_P$ has then the form
that $d(\beta_P)$ equals to the pasting of the diagram of multiplications and comultiplications
formed by the three triangles with vertex $w_l$ partitining $P$. In other words, this pasting
is  {\em homotopically zero}.
In the case when
  $\del_+P$ consists of two segments,
$\beta_P$ has degree $0$ and the pasting diagram we obtain is identical to that depicted in Fig. \ref{fig:PL-3:1}.
So this $\beta_P$ gives an analog of the $3\to 1$ move identity of Proposition \ref {prop:PL-3:1}.

\paragraph{The infrared Maurer-Cartan element (conjectural).}
From the above analysys we see that
the conditions on the basic components of a Maurer-Cartan element $\beta\in\gen^1$ are completely
analogous to those satisfied by the Picard-Lefschetz arrows in \S \ref{subsec:PL-tri}.
The difference is that the Picard-Lefschetz arrows as we defined them, live in the homotopy category of
dg-functors (are cohomology classes of closed natural transformations of degrees $1$ or $0$),
 while the components of $\beta$ are natural transformations themselves.
 This makes natural the following conjecture.

 \begin{conj}\label{conj:MC}
 For any schober $\Sen\in\Schob(\CC,A)$ there exists a Maurer-Cartan element
 $\beta = \beta_\Sen \in \gen_A(\Sen)^1$, defined canonically up to gauge equivalence
 so that for any generic direction $\zeta$ we have  the following properties:
 \begin{itemize}
 \item[(1a)] For an empty triangle $P = \Conv\{w_i, w_j, w_k\}$
 with $\del_+P$ consisting of one edge, as in the left  of Fig. \ref{fig:empty-tri},
 the cohomology class of $\beta_P$ in $H^1\, \Hom^\bullet (M_{ik}, M_{jk} M_{ij})$
 is equal to the  Picard-Lefschetz arrow $u_{ijk} = u_{[w_j, w_k], [w_i, w_j], [w_i,w_k]}$,
 see \eqref {eq:PL-arrows}.

 \item[(1b)] For an empty triangle $P = \Conv\{w_i, w_j, w_k\}$
 with $\del_+P$ consisting of two  edges, as in  the right of  Fig. \ref{fig:empty-tri},
 the cohomology class of $\beta_P$ in $H^0\, \Hom^\bullet ( M_{jk} M_{ij}, M_{ik})$
 is equal to the  Picard-Lefschetz arrow $v_{ijk} = v_{[w_j, w_k], [w_i, w_j], [w_i,w_k]}$,
 see \eqref {eq:PL-arrows}.

 \item[(2)] The deformation $\R^\zeta(\beta)$ of the infrared monad $\R^\zeta$ corresponding to
 $\beta$, is quasi-isomorphic to the Fukaya-Seidel monad $\M(\bv)$ corresponding to
 the Vladivostok point $\bv=-r\zeta$, $r\gg 0$  situated far away in the direction $(-\zeta)$
 and the $\bv$-spider $K_\bv$ being the union of straight segments $[\bv, w_i]$
 see Fig. \ref{fig:lasso}.
 \end{itemize}
 \end{conj}

 This conjecture should follow from a more systematic definition of perverse schobers where
 all the coherence data will be given at the outset as a part of the structure.

 \begin{rem}
 Note that Conjecture \ref{conj:MC-simple} about an MC-lifting of the differential in the infrared complex
 is a particular case of Conjecture \ref {conj:MC} obtained by focussing on the effect
 of $\beta$ just on the differentials in the complexes $\R_{ij}$ (level (0) of Example
 \ref {ex:R-beta}). Indeed, Fig. \ref{fig:neighborly},  describing
  the construction of the differential $q$ in the infrared complex, can be seen as  a particular case of
 Fig. \ref{fig:1-fin} corresponding to $m=1$.
 \end{rem}

 \begin{ex}
 Consider the ``maximally concave'' situation of Example \ref{ex:max-conc}. That is, assume that $A$
 and $\zeta$ are positioned so that every element of $A$ is a vertex of $\Conv^{-\zeta}(A)$,
 see Fig. \ref{fig:max-conc}.
 Then every subpolygon $P\subset (Q,A)$ is, in the terminology of this section,
  empty and extremal of the second kind.
 As we saw, in this case $\R_{ij}=M_{ij}$ is the rectilinear transport functor.  Moreover, it is identified
 with $\M_{ij}= a_j^* \, a_i$ by deforming the interval $[w_i, w_j]$ to the Vladivostok path
 $[w_i,\bv]\cup [\bv, w_j]$, as in Fig. \ref {fig:lasso}. However,  all the compositions
 $c_{ijk}^\R: \R_{jk}\circ \R_{ij}\to \R_{ik}$, $i<j<k$,  are zero, while
 $c_{ijk}^\M: \M_{jk}\circ \M_{ij}\to\M_{ik}$ is given by the counit for the adjoint pair $(a_j, a_j^*)$.
 Note that under our system of identifications $\M_{ij}\= M_{ij}$, the cohomology class of
 the closed degree $0$ natural transformation
 $c_{ijk}^\M$ corresponds to  the Picard-Lefschetz arrow $v_{ijk}$. Moreover, the $c_{ijk}^\M$
 are associative by construction. Therefore, putting
 $\beta_P = c_{ijk}^M$ for $P=\Conv\{w_i, w_j, w_k\}$ being a triangle and
 $\beta_P=0$ for all $n$-gons $P$ with $n\geq 4$, we get a Maurer-Cartan element $\beta$
 and $R(\beta) \= \M$, thus establishing the conjecture in this case.

 \end{ex}

 A natural strategy of proving Conjecture \ref {conj:MC} would be to start from  $A=A(0)$ in a maximally
 concave position above and then deform $A$ via a path $A(t)= (w_1(t),\cdots, w_N(t))$, $t\in [0,1]$,
 to a new position $A'=A(1)$. We  can then
  deform $\Sen$ accordingly, i.e.,
 ``isomonodromically''  to a new  schober $\Sen(t)\in\Schob(\CC, A(t))$ for each $t$.
 If the path $A(t)$ is generic,  the nature of $\gen$, $\R$ and $\M$ will be unchanged except
 for finitely many ``wall-crossing'' points where we cross either a wall of collinearity/horizontality
 (\S \ref {subsec:iso-recti}) or an exceptional wall (\S \ref{subsec:W-II-kind}), and the task
 is to describe explicitly the modification of the Maurer-Cartan element $\beta(t)$ under such
 wall-crossing. We plan to return to this in a future work.


\appendix

\renewcommand{\theparagraph}{\arabic{paragraph}.}
\numberwithin{equation}{section}

\section {Conventions on (enhanced) triangulated categories}\label{app:enh}

\paragraph{DG-enhancements.}
Let $\k$ be a field of characteristic $0$.
We will work with $\k$-linear triangulated categories (all Hom's are $\k$-vector spaces).
The $n$-fold shift functor of a triangulated
category $\VV$ will be denoted by $x\mapsto x[n]$.  An exact triangle in $\VV$ will be written
as
\[
x\lra y\lra z\lra x[1].
\]
In many cases we need an additional structure on a triangulated category called an
{\em enhancement}. In this paper we are mostly interested in differential graded (dg-)enhancements
\cite{BK} so we briefly recall the terminology and conventions.

\vskip .2cm

Recall that a ($\k$-linear)  {\em dg-category}  is a category $\Vc$ in which
each $\Hom$-set is made into a dg-vector space (=cochain complex) $\Hom^\bullet_\Vc(x,y)$
over $\k$, so that the differentials satisfy the Leibniz rules with respect to composition.
We refer to \cite { tabuada, toen-morita, toen-dgcat} for a systematic backround on dg-categories.

A dg-category $\Vc$ give a $\k$-linear category $H^0(\Vc)$ by passing to the $0$th
cohomology of the Hom-complexes:
\[
\Ob(H^0(\Vc))\,=\, \Ob(\Vc), \quad \Hom_{H^0(\Vc)}(x,y) \,=\, H^0 \Hom^\bullet_\Vc(x,y),
\]
and a similar category $H^\bullet(\Vc)$ obtained by passing to the full cohomology spaces,
not just $H^0$.

For example, all dg-vector spaces over $\k$ form a dg-category $\dgVect$.

\vskip .2cm

A {\em dg-functor} $F: \Vc\to\Wc$ between dg-categories is a functor preserving
the structures (vector space, differential, grading) on the Hom-complexes. All dg-functors
$\Vc\to\Wc$ form a new dg-categoru $\dgFun(\Vc, \Wc)$. Any dg-category $\Vc$ has the
{\em Yoneda embedding}
\[
\Upsilon: \Vc\lra \dgFun(\Vc^\op, \dgVect), \quad x\mapsto \Hom^\bullet_\Vc(-,x).
\]
A {\em dg-functor} $F: \Vc\to\Wc$ is called a {\em quasi-equivalence}, if
$H^\bullet(F): H^\bullet(\Vc)\to H^\bullet(\Wc)$ is an equivalence of categories.

\vskip .2cm

Any dg-category $\Vc$ has the {\em pre-triangulated envelope} $\Pretr(\Vc)$, see \cite{BK, tabuada}.
It is a new dg-category whose objects are {\em twisted complexes} over $\Vc$,
i.e., data $V^\bullet = \bigl( (V^i)_{i\in \ZZ}, (q_{ij})_{i<j}\bigr)$, where $V^i$ are objects of
$\Vc$ and $q_{ij}\in \Hom^{j-i}_\Vc(V^i, V^j)$ are such that
\[
dq_{ij} + \sum_{i<k<j} q_{jk} q_{ij} \,= \, 0.
\]
The category $\Tr(\Vc):= H^0(\Pretr(\Vc))$ is always triangulated. A morphism $\wt f\in\Hom^0_{\Pretr(\Vc)}(x,y)$ which is closed (i.e., $d\wt f=0$), has the {\em canonical cone}
$\Cone(\wt f)$ which fits into an exact triangle in $\Tr(\Vc)$
\[
x\buildrel f\over \lra y\lra \Cone(\wt f)\lra x[1].
\]
Here $f$ is the morphism in $\Tr(\Vc)$ represented by the cohomology class of $\wt f$.

\vskip .2cm

A dg-category $\Vc$ is called {\em pre-triangulated}, if the canonical embedding
$ \kappa: \Vc\to\Pretr(\Vc)$ is a quasi-equivalence. In this case $H^0(\Vc)$ is triangulated.
Moreover,  for a pre-triangulated $\Vc$ we have  a quasi-inverse dg-functor to $\kappa$ which
is called the  {\em total object}  functor $\Tot: \Pretr(\Vc)\to\Vc$.
Note that a closed, degree $0$
morphism $\wt f: x\to y$ is a particular type of a twisted complex: in this case
$V^\bullet$ has $V^{-1}=x$, $V^0=y$ and $q_{-1,0}=\wt f$. The corresponding total
object $\Tot(V^\bullet)$ is simply the cone of $\wt f$.
This means that we can speak about canonical cones of closed, degree $0$ morphisms of
$\Vc$.

\vskip .2cm

if $\Vc, \Wc$ are dg-categories and $\Wc$ is pre-triangulated, then so is $\dgFun(\Vc, \Wc)$.
In particular, we can speak about cones of natural transformations of dg-functors.
Note that $\Pretr(\Vc)$ can be seen as the full dg-subcategory in $\dgFun(\Vc^\op,\dgVect)$ formed
by taking iterated cones, and shifts, starting from objects of $\Upsilon(\Vc)$.

\vskip .2cm

A {\em dg-enhancement} of a ($\k$-linear) triangulated category $\VV$ is an  identification
of triangulated categories
$\VV\simeq H^0(\Vc)$ for a pre-triangulated category $\Vc$.

\begin{conven}
In the main body of the paper, we will often not distinguish notationally between $\VV$ and $\Vc$,
saying that $\Vc$ is an ``enhanced triangulated category'' (meaning ``$\Vc$ is an enhancement
of $\VV$'') or speaking about ``exact triangles in $\Vc$'' (which belong, strictly speaking, to  $\VV$)
etc.

\end{conven}

\begin{ex}\label{ex:bimodules}
Let $R$ be a dg-algebra (=dg-category with one object). The category $\mod_R=\dgFun(R^\op,\dgVect)$ of right dg-modules over $R$ is pre-triangulated. A dg-module $P$ is called
{\em perfect}, if it can be obtained from free modules (copies of $R$) by iterated application
of shifts, taking of cones and of homotopy direct summands, see \cite{toen-morita, toen-vaquie}.
Let $\Perf_R\subset \mod_R$ be the full dg-subcategory of perfect dg-modules.

Given two dg-algebras $R_1, R_2$ and an $(R_1, R_2)$-dg-module $M$, perfect over $R_1$
and $R_2$, we have the dg-functor
\[
F_M: \Perf_{R_1} \lra \Perf_{R_2}, \quad P\mapsto P\otimes^L_{R_1} M
\]
(the derived tensor product). See \cite{toen-morita} for a more general treatment, where
dg-algebras are replaced by dg-categories.

In most of considerations of the present paper, one can assume the enhanced triangulated
categories to be of the form $\Perf_R$ and the functors to be of the form $F_M$,
thus simplifying the foundational issues.
\end{ex}

\paragraph{$\oo$-categorical enhancements.}\label{par:oo-cat-enh}
 The theory of dg-categories can be embedded into
the theory of $\oo$-categories (or weak Kan complexes) \cite{HTT} which provides a more general
approach to ``enhancements''. While we do not use this theory in the present paper,
it will be used in the more systematic treatment of schobers in \cite{DKSS}, so let us
say a few words about the relation of the two approaches.

An $\oo$-category $\Vs$ gives an ordinary category $\ho(\Vs)$ known as the {\em homotopy
category} of $\Vs$. From the point of view on $\oo$-categories as simplicially enriched categories
( \cite{HTT} \S 1.1.5), an $\oo$-categorical enhancement of a category $\VV$ represents
each $\Hom_\VV(x,y)$ as $\pi_0$ of a simplicial set (instead of $H^0$ of a cochain complex
as is the case with dg-enhancements).

To a dg-category $\Vc$ one can associate an $\oo$-category $\Ndg\Vc$ known as the
{\em dg-nerve} of $\Vc$. Viewed as a weak Kan complex, $\Ndg\Vc$ consists of
Sugawara (weakly commutative) simplices of closed, degree $0$ morphisms of $\Vc$,
see \cite{HS, HA}, and $\ho(\Ndg \Vc)= H^0(\Vc)$.

\vskip .2cm

The role of pre-triangulated dg-categories is played, in the $\oo$-categorical context,
by {\em stable $\oo$-categories} \cite{HA}. If $\Vc$ is a pre-triangulated dg-category,
then $\Ndg\Vc$ is stable \cite{faonte},  so the two approaches are compatible.

The use of $\oo$-categorical enhancements can clarify the treatment of several issues related to
perverse schobers.
For example, small dg-categories do not form a dg-category, but small $\oo$-categories do form an
$\oo$-category. So the construction of $\Ndg\Vc$ leads to a definition of an $\oo$-category
$\dgCat$ formed by small dg-categories.  This means that perverse schobers
should, in a more systematic theory, form an $\oo$-category.

\paragraph{Adjunctions for ordinary categories.}\label{par:adj}
Let $\Vc, \Wc$ be categories.
Recall \cite{maclane} that a pair $(F,G)$ of functors in the opposite directions
$
\xymatrix{
\Vc \ar@<.4ex>[r]^{F}&\Wc \ar@<.4ex>[l] ^{G}
},
$
is called an {\em adjoint pair}, if it is equipped with natural isomorphisms
\be\label{eq:adjunction}
\alpha_{x,y}: \Hom_\Wc(F(x), y) \buildrel\sim\over\lra \Hom_\Vc(x, G(y)), \quad x\in\Vc, \, y\in\Wc.
\ee
One says that $G$ is {\em right adjoint} to $F$ and writes $G=F^*$, and also that
$F$ is {\em left adjoint} to $G$ and writes $F={^*}G$. Equivalently, the data of
$(\alpha_{x,y})$ are encoded in the {\em unit} and {\em counit morphisms of functors}
\be
e=e_F: \Id_\Vc\lra GF, \quad \eta=\eta_F: FG\lra \Id_\Wc,
\ee
satisfying the standard axioms \cite{maclane}. For example, for $x\in\Vc$ the
morphism $e_x: x\to GF(x)$ is $\alpha_{x, F(x)}(\Id_{F(x)})$.

Note that the adjoints, if they exist, are defined canonically (i.e., uniquely up to
a unique isomorphism). For example if $F$ is given, then \eqref{eq:adjunction}
means that $G(y)$ represents the functor $x\mapsto \Hom(F(x),y)$, and
the representing object of a functor, if  it exists, is defined canonically.

Note further that passing to the adjoints  reverses the order of compositions:
\be
(F_1\circ F_2)^* \,=\, F_2^*\circ F_1^*, \quad {^*}(G_1\circ G_2) \,=\, {^*}G_2 \circ {^*} G_1.
\ee
We can also speak about iterated adjoints $F^{**} = (F^*)^*$, $^{**}G = {^*}({^*}G)$ etc.
provided they exist.

\begin{rem}
The concept of adjoint functors between categories is a particular case of the
concept of adjoint $1$-morphisms between objects of a $2$-category $\Cen$
(corresponding to $\Cen=\Cc at$ being the $2$-category of categories). Most
of the constructions related to adjoints generalize to that $2$-categorical
setting. On the other hand, a $2$-category $\Cen$ with one object $\pt$
is the same as a monoidal category $(\Vc,\otimes)$ whose objects
are $1$-morphisms $\pt\to\pt$ and $\otimes$ is the composition of $1$-morphisms.
The $2$-categorical concept of left and right adjoints specializes then to the
concept of {\em left} and {\em right duals} $^*V, V^*$ of an object $V\in\Vc$,
see, e.g.,  \cite{etingof} \S 2.10. The features below are perhaps more familiar
in such  ``tensor context''.
\end{rem}

\vskip .2cm

Given two functors $F_1, F_2: \Vc \to\Wc$ and a natural transformation
$k: F_1\to F_2$, we have the {\em right} and {\em left transposes}
\be\label{eq:transposed-functors}
k^t: F_2^* \lra F_1^*, \quad {^t}k: {^*F}_2 \lra {^*}F_1.
\ee
For example, $k^t$ is the composition
\[
\xymatrix{
F_2^* \ar[rr]^{\hskip -0.5cm e_{F_1}\circ F_2^*}&& F_1^* F_1 F_2^*
\ar[rr]^{F_1^*\circ k\circ F_2^*} && F_1^* F_2 F_2^*
\ar[rr]^{\hskip 0.3cm F_1^* \circ \eta_{F_2}} && F_1^*,
}
\]
and $^tk$ is defined similarly.

\begin{prop}\label{prop:FGH-iso}
(a) Passing the the   left or right transposes gives natural identifications
\[
\Hom(^*F_2, ^*F_1) \,\simeq \,
\Hom(F_1, F_2) \,\simeq \, \Hom(F_2^*, F_1^*).
\]

\vskip .2cm

(b) The identifications in (a) are compatible with composition. That is,
for natural transformations $F_1\buildrel k_1\over\to F_2\buildrel k_2\over \to F_3$
we have
\[
(k_2k_1)^t \,=\, (k_1^t) (k_2^t), \quad ^t(k_2k_1) \,=\, (^tk_2) (^tk_1)
\]
(provided the adjoints exist).

\vskip .2cm

(c) Let $\Cc, \Dc, \Ec$ be categories and
$\Cc\buildrel G\over\to \Dc\buildrel F\over\to \Ec$ and $\Cc\buildrel H\over\to\Ec$ be functors.
Then we have natural idenfitications (partial transposes), provided the adjoints exist:
\[
\Hom(^*H FG, \Id_\Dc)\,\simeq \,
\Hom(FG, H) \,\simeq \, \Hom (G,  F^* H).
\]
\end{prop}
\noindent{\sl Proof:} Similar to the tensor case, see \cite{etingof} (Ex. 2.10.7 and Prop. 2.10.8).
For example, the last identification in (c) is the composition
\[
\Hom(FG,H)\buildrel F^*\circ(-)\over\lra \Hom(F^*FG, F^*H)
\buildrel \rfloor e_F\over\lra \Hom(G, F^*H),
\]
where the first arrow is the application of $F$ and the second is the substitution of
$e_F: \Id\to F^*F$.
\qed

\begin{cor}\label{cor:tens-rot}
Assuming that the relevant  adjoints exist,
in the situation of Proposition \ref{prop:FGH-iso} we have natural identifications
(``tensor rotation'')
\[
\begin{gathered}
\Hom(FG,H) = \hfill
\\
=   \Hom(G, F^* H) = \Hom(GH^*, F^*) = \Hom(H^*, G^* F^*) = \cdots
= \Hom(F^{**}G^{**}, H^{**}) =\cdots \hfill
\\
= \Hom(F, H(^*G)) = \Hom(^*H F, (^*G)) = \Hom(^*H, ^*G(^*F))=\cdots = \Hom(^{**}F (^{**}G), (^{**}H))
=\cdots
\end{gathered}
\]
\qed
\end{cor}

\paragraph{Adjunctions at the dg-level.} \label{par:app-dgadj}

We will need the formalism of adjoints in the setting of dg-categories and dg-functors.
In this case there are extra subtleties: the $\alpha_{x,y}$ need to be not isomorphisms but
only quasi-isomorphisms
of Hom-complexes, the naturality should hold not ``on the nose''  but only up to a coherent
system of homotopies and so on. A practical way of handling these issues was proposed
in \cite{AL-adj} for dg-functors coming from bimodules, cf. Example \ref {ex:bimodules}. In this
paper we will mostly ignore these coherence issues, assuming implicitly the ansatz of
\cite{AL-adj}.

A more satisfactory general treatment of adjunctions can be provided in the framework
of $\oo$-categories \cite{HTT} \S 5.2 and will be used in \cite{DKSS}.

\vskip .2cm

We will use the following statement which is clear in the bimodule setting of  \cite{AL-adj}.

\begin{prop}\label{prop:cone-kt}
Suppose $\Vc, \Wc$ are pre-triangulated dg-categories and $F_1, F_2: \Vc\to \Wc$
are dg-functors with right adjoints $F_1^*, F_2^*$. Let $k: F_1\to F_2$ be a closed
natural transformation of degree $0$, so that we have a triangle of dg-functors
\[
\Cone(k)[-1] \lra F_1 \buildrel k\over\lra F_2\lra \Cone(k). \leqno (*)
\]
Then $\Cone(k)^*$, the right adjoint to $\Cone(k)$, exists and is identified with
$\Cone(k^t)[1]$, so that the triangle
\[
\Cone(k)^* \lra F_2^* \buildrel k^t \over\lra F_1^*\lra \Cone(k)^*[1]
\]
dual to (*), is identified with the triangle
\[
\Cone(k^t)[-1] \lra F_2^* \buildrel k^t\over\lra  F_1^* \lra \Cone(k^t).
\]
\qed
\end{prop}

\paragraph{Serre functors. Calabi-Yau categories.}
Let $\Vc$ be a $k$-linear triangulated category with all $\Hom_\Vc(x,y)$ finite-dimensional.
Recall \cite{BK-serre} that a {\em Serre functor} for $\Vc$ is a self-equivalence $S:\Vc\to\Vc$
equipped with natural isomorphisms (the {\em Serre structure})
\be\label{eq:def-serre}
\rho_{x,y}: \Hom_\Vc(x,y)^* \lra\Hom_\Vc(y, S(x)).
\ee
The basic example is $\Vc=D^b\Coh(X)$, the coherent derived category of a smooth projective
variety $X/\k$ of dimension $n$. Then $S(\Fc)=\Fc\otimes\Omega^n_X[n]$ os
a Serre functor, with \eqref{eq:def-serre} being the classical Serre duality.

Note that $S$, if it exists, is defined canonically (uniquely up to a unique isomorphism of functors).
This is because  \eqref{eq:def-serre}  means that $S(x)$ is the representing object for
the functor $y\mapsto \Hom(x,y)^*$. However, when we say that a given functor
$F:\Vc\to\Vc$  ``is'' a Serre functor, we mean that it {\em is made into} a Serre functor by
equipping it with an appropriate Serre structure. This Serre structure does not have to be
unique, it is unique only up to an automorphism of $F$ as a functor.

For example, we say that $\Vc$ is a {\em Calabi-Yau category of dimension} $n\in\ZZ$,
if the shift functor $x\mapsto x[n]$ is (made into) a Serre functor, i.e., if $\Vc$
is equipped with natural isomorphisms ({\em Calabi-Yau structure})
\be
\kappa_{x,y}: \Hom_\Vc(x,y)^* \lra \Hom_\Vc(y, x[n]).
\ee
For instance, in the example $\Vc=D^b\Coh(X)$, the functor $x\mapsto x[n]$, $n=\dim(X)$,
can be made into a Serre functor, if $X$ admits a global nowhere vanishing volume form
$\omega\in\Gamma(X,\Omega^n)$ (so it is a Calabi-Yau manifold in the  usual sense).
A choice of such a form (unique only up to a scalar)
gives a Calabi-Yau structure.

\vskip .2cm

We will also use the dg-enhanced versions of these concepts. That is, let $\Vc$ be
a pre-triangulated dg-category such that each complex $\Hom^\bullet_\Vc(x,y)$
has the total cohomology space finite-dimensional. Then by a Serre functor for
$\Vc$ we mean a dg-functor $S:\Vc\to\Vc$ which is a quasi-equivalence and is
equipped with natural quasi-isomorphisms of Hom-compleces $\rho_{x,y}$,
like in \eqref{eq:def-serre} but with $\Hom^\bullet$ instead of $\Hom$.
As this amounts to (homotopy) representing dg-functors of the form
$y\mapsto \Hom^\bullet(x,y)^*$, the object $S(x)$ is  defined homotopy canonically
(uniquely up to a contractible space of choices) and so one can use the $\oo$-categorical
approach to address various issues of coherence etc.

Similarly, we will speak about enhanced versions of Calabi-Yau structures.

\paragraph{ Local systems of (pre-triangulated) dg-categories.}
Let $X$ be a topological space. The natural categorification of the concept of a sheaf (of sets,  vector spaces etc.)
on $X$ is the concept of a {\em stack of categories}. Let us explain it in a convenient language analogous to
that of sheaves, see, e.g., \cite{breen} for a more systematic treatment.

\begin{Defi}\label{def:prestack}
A {\em pre-stack} (or a {\em pseudo-functor on open sets}) on  $X$ is a datum $\Fen$ of:

\begin{itemize}
\item[(0)] For any open set $U\subset X$, a category $\Fc(U)$.

\item[(1)] For any any inclusion $U_1\subset U$, a functor
(called the {\em restriction})  $r_{U,U_1}: \Fen(U)\to \Fen(U_1)$, so that $r_{UU}=\Id$.

\item[(2)] For any triple inclusion $U_2\subset U_1\subset U$, an isomorphism of functors ({\em transitivity of restrictions})
$\rho_{U,U_1,U_2}: r_{U_1,U_2} \circ r_{U,U_1}\to r_{U,U_2}$, so that:

\item[(3)] For any quadruple inclusion $U_3\subset U_2 \subset U_1\subset U_0=U $ we have a compatibility
(cocycle) condition on the $4$ isomorphisms $\rho_{U_i, U_j, U_k}$, $0\leq i< j<k\leq 3$.
\end{itemize}

\end{Defi}

 \noindent A {\em stack} is a pre-stack satisfying the descent conditions for any open covering
 $U=\bigcup_a U_a$ which say that objects and morphisms in $\Fc(U)$ can be glued out of
  families of such on the $U_a$ equipped with compatibility data on the intersections.

  \vskip .2cm

  A standard example of a stack is the correspondence $U\mapsto \on{Sh}(U)$, the category of all
  sheaves on $U$. The descent conditions in this case mean that a sheaf on $U$ can be glued
  out of sheaves on the $U_a$ with compatible identifications on the intersections.

    \vskip .2cm

    Let $X$ be an $n$-dimensional $C^\oo$-manifold, possible with boundary.
    A {\em local system
    of categories} i(or a {\em locally constant stack}) on $X$ is a stack $\Fc$ with the following property:
    whenever $U_1\subset U$
    are two {\em balls} in $X$, the functor $r_{U, U_1}$ is an equivalence of categories.
   Here by
      a   ball in $X$
    we mean an open set homeomorphic either to the $n$-dimensional open ball
    $B_n=\{\|x\|< 1\}\subset
    \RR^n$, or to the ``ball at the boundary'' $B_n\cap \{x_1\geq 0\}$.

    \vskip .2cm

    Given a local system $\Fen$ on $X$, for each $x\in X$ we have the category $\Fen_x$, the
    {\em stalk} of $\Fen$ at $x$, defined as $\Fen(U)$ for sufficiently small balls containing $U$.
    Given
      a path $\gamma$ joining two points $x,y\in X$,
    we have  the {\em monodromy  } (or {\em categorical parallel transport}) functor
    $T(\gamma): \Fen_x\to\Fen_y$ by covering $\gamma$
    by a sequence of small balls in a standard way.
    These objects are defined {\em canonically} in the following sense.
    The category  $\Fen_x$ is defined uniquely up to
    an equivalence of categories  which, in its turn, is defined uniquely up to a unique isomorphism of functors,
    and  the functor $T(\gamma)$ is defined uniquely up to a unique isomorphism.

    \vskip .2cm

    Alternatively, a local system can be defined as a pseudo-functor $U\mapsto \Fen(U)$ as above
    but defined only on disks and such that any $r_{U,U_1}$ is an equivalence.

    \vskip .2cm

It is  straightforward to generalize to the case  of dg-categories. As mentioned above,
all dg-categories
are united into an $\oo$-category $\dgCat$, see
\S \ref{app:enh}\ref{par:oo-cat-enh}
A   {\em (homotopy) stack} of dg-categories  is, formally, a sheaf $\Fen$ on $X$ with values in $\dgCat$
in the sense of \cite{HTT}.

\vskip .2cm

Explicitly, this means that $\Fen$ is, similarly to Definition \ref{def:prestack}
 a datum consistinng of dg-categories $\Fc(U)$ at level (0),  dg-functors
$r_{U,U_1}$ at level (1), but at level (2) we have {\em closed degree $0$} natural transformations $\rho_{U,U_1, U_2}$, and level (3) is replaced by an infinite chain of higher coherence conditions.
Similarly, the descent conditions all involve infinite systems of higher coherences.

\vskip .2cm

A {\em local system} of dg-categories on a manifold $X$  is a stack $\Fc$ such that
$r_{U,U_1}: \Fc (U)\to \Fc(U_1)$ is a {\em quasi-equivalence} of dg-categories whenever $U,U_1$
are disks. As before, we can speak about stalks $\Fen_x$ and monodromies $T(\gamma)$
which are dg-categories and dg-functors defined {\em homotopy canonically}, i.e.,
up to a contractible space of possible choices.

\vskip .2cm

If $X$ is a complex manifold, a local system of pre-triangulated dg-categories should be
certainly considered as an example of a perverse schober. However, it does not seem possible
to realize more general perverse schobers
 as stacks of dg-categories in the sense we described.


\vskip 2cm

M.K.: Kavli IPMU, 5-1-5 Kashiwanoha, Kashiwa, Chiba, 277-8583 Japan. \hfill\break
 Email: {\tt mikhail.kapranov@protonmail.com}

 \vskip .2cm

 Y.S.: Department of Mathematics, Kansas State University, Manhattan KS 66508 USA.
 \hfill\break Email: {\tt soibel@math.ksu.edu}

 \vskip .2cm

 L.S.: Department of Mathematics, Higher School of Economics, 109028
 Moscow, Russia. \hfill\break
 Email:   {\tt amligart@yandex.ru}

\end{document}